\documentclass[oneside,noinfoline,leqno]{article}

\usepackage{imsart}

\usepackage[nottoc]{tocbibind}
\usepackage{latexsym}
\usepackage{amssymb,exscale}
\usepackage{amsmath}
\usepackage{amsthm}
\usepackage[all]{xy}
\usepackage{graphicx}
\usepackage{multicol}
\usepackage{makeidx}
\usepackage{footnpag}

\usepackage[hypertex,linktocpage,pagebackref]{hyperref}

\numberwithin{equation}{subsection}

\renewcommand{\theequation}{\thesection\alph{subsection}\arabic{equation}}

\swapnumbers
\theoremstyle{definition}
\newtheorem{theorem}[equation]{Theorem}
\newtheorem{lemma}[equation]{Lemma}
\newtheorem{proposition}[equation]{Proposition}
\newtheorem{corollary}[equation]{Corollary}
\newtheorem{definition}[equation]{Definition}
\newtheorem{example}[equation]{Example}
\newtheorem{remark}[equation]{Remark}
\newtheorem{question}[equation]{Question}

\newcommand{\mysection}[1]{\section[#1]{\raggedright#1}}
\newcommand{\mysubsection}[1]{\subsection[#1]{\raggedright#1}}

\newcommand{\beginproof}{\begin{proof}[Proof \textup{(sketch)}]}
\newcommand{\proofend}{\end{proof}}
\newcommand{\proofendnoqed}{\renewcommand{\qed}{}\end{proof}}

\newcommand{\I}{{\rm i}}
\newcommand{\D}{\mathrm{d}}
\newcommand{\E}{\mathrm{e}}

\newcommand{\ti}{\tilde}

\renewcommand{\(}{\bigl(}
\renewcommand{\)}{\bigr)\vphantom{)}}

\newcommand{\ip}[2]{\langle#1,#2\rangle}

\newcommand{\imply}{\;\;\;\Longrightarrow\;\;\;}
\newcommand{\imp}{$ \Longrightarrow $ }

\newcommand{\pd}{\partial}

\newcommand{\splt}{\text{\textup{split}}}
\newcommand{\stick}{\text{\textup{stick}}}
\newcommand{\stable}{\text{\textup{stable}}}
\newcommand{\white}{\text{\textup{white}}}
\newcommand{\cls}{\text{\textup{cls}}}
\newcommand{\co}{\text{\textup{c}}}

\newcommand{\Comp}{\operatorname{Comp}}
\newcommand{\Lip}{\operatorname{Lip}}
\newcommand{\trace}{\operatorname{trace}}
\newcommand{\Homeo}{\operatorname{Homeo}}
\newcommand{\Diff}{\operatorname{Diff}}
\newcommand{\Exp}{\operatorname{Exp}}
\newcommand{\Log}{\operatorname{Log}}
\newcommand{\Aut}{\operatorname{Aut}}

\newcommand{\Var}{\operatorname{Var}}
\newcommand{\Cov}{\operatorname{Cov}}

\newcommand{\N}{\operatorname{N}}
\newcommand{\U}{{\operatorname{U}}}
\newcommand{\V}{\operatorname{V}}
\newcommand{\Poisson}{{\operatorname{Poisson}}}

\renewcommand{\P}{\operatorname{P}}
\renewcommand{\Im}{\operatorname{Im}}

\newcommand{\modO}{{\operatorname{mod}\,0}}

\newcommand{\One}{\mathbf1}

\renewcommand{\phi}{\varphi}
\newcommand{\eps}{\varepsilon}
\newcommand{\si}{\sigma}
\newcommand{\ga}{\gamma}
\newcommand{\om}{\omega}
\newcommand{\Om}{\Omega}
\newcommand{\de}{\delta}
\newcommand{\De}{\Delta}
\newcommand{\al}{\alpha}
\newcommand{\be}{\beta}
\newcommand{\AAA}{\mathsf A}
\newcommand{\Pc}{\mathcal P}
\newcommand{\cH}{\mathcal H}
\newcommand{\X}{\mathcal X}
\newcommand{\cT}{\mathcal T}

\newcommand{\cM}{\mathcal M}

\newcommand{\Ec}{\mathcal E}
\newcommand{\F}{\mathcal F}
\newcommand{\A}{\mathcal A}
\newcommand{\B}{\mathcal B}
\newcommand{\la}{\lambda}
\newcommand{\La}{\Lambda}

\newcommand{\dist}{\operatorname{dist}}
\newcommand{\const}{{\mathrm{const}}}
\newcommand{\Ex}{\mathbb E\,}
\renewcommand{\Pr}[1]{\mathbb{P}\mskip1.5mu\(\mskip1.5mu#1\mskip1.5mu\)}
\newcommand{\M}{{\operatorname{M}}}
\newcommand{\Prob}{\mathbb P}
\newcommand{\R}{\mathbb R}
\newcommand{\C}{\mathbb C}

\newcommand{\Z}{\mathbb Z}
\newcommand{\T}{\mathbb T}
\newcommand{\cE}[2]{\mathbb{E}\mskip1.5mu\(\mskip1.5mu#1\mskip1.5mu
 \big|\mskip1.5mu#2\mskip1.5mu\)}
\newcommand{\CE}[2]{\mathbb{E}\mskip1.5mu\bigg(\mskip1.5mu#1\mskip1.5mu
 \bigg|\mskip1.5mu#2\mskip1.5mu\bigg)}
\newcommand{\cP}[2]{\mathbb{P}\mskip1.5mu\(\mskip1.5mu#1\mskip1.5mu
 \big|\mskip1.5mu#2\mskip1.5mu\)}

\newcommand{\Tl}{\rlap{$\scriptstyle\mskip-1mu
 \overleftarrow{\raise.25ex\hbox{$\textstyle\phantom{T}$}}
$}T}
\newcommand{\Tr}{\rlap{$\scriptstyle\mskip-0.5mu
 \overrightarrow{\raise.25ex\hbox{$\textstyle\phantom{T}$}}
$}T}

\newcommand{\sif}{$\sigma$\nobreakdash-field}

\newcommand{\sifinite}{$\sigma$\nobreakdash-\hspace{0pt}finite}

\newcommand{\almost}[1]{$#1$\nobreakdash-\hspace{0pt}almost}
\newcommand{\particle}[1]{$#1$\nobreakdash-\hspace{0pt}particle}
\newcommand{\valued}[1]{$#1$\nobreakdash-\hspace{0pt}valued}
\newcommand{\measurable}[1]{$#1$\nobreakdash-\hspace{0pt}measurable}
\newcommand{\flow}[1]{$#1$\nobreakdash-\hspace{0pt}flow}
\newcommand{\dimensional}[1]{$#1$\nobreakdash-\hspace{0pt}dimensional}
\newcommand{\Smap}{S\nobreakdash-map}
\newcommand{\Skernel}{S\nobreakdash-kernel}

\makeatletter
  {\end{multicols}}%
\makeatother

\makeatletter
\renewcommand*\l@section[2]{%
  \ifnum \c@tocdepth >\z@
    \addpenalty\@secpenalty
    \addvspace{0.25em \@plus\p@}
    \setlength\@tempdima{2.5em}
    \begingroup
      \parindent \z@ \rightskip \@pnumwidth
      \parfillskip -\@pnumwidth
      \leavevmode \bfseries
      \advance\leftskip\@tempdima
      \hskip -\leftskip
      #1\nobreak\hfil \nobreak\hb@xt@\@pnumwidth{\hss #2}\par
    \endgroup
  \fi}
\renewcommand*\numberline[1]{\hb@xt@\@tempdima{\hfil#1\hskip1em}}
\makeatother

\begin{document}

\begin{frontmatter}

\title{Nonclassical stochastic flows\\
 and continuous products}
\runtitle{Nonclassical flows and products}
\author{\fnms{Boris}
\snm{Tsirelson}
\ead[label=e1]{tsirel@post.tau.ac.il}}
\address{School of Mathematics, Tel Aviv University,\\
69978 Tel Aviv, Israel\\
\printead{e1}\\
url:
\href{{http://www.tau.ac.il/~tsirel/}}{http://www.tau.ac.il/$\sim$tsirel/}}

\runauthor{B.~Tsirelson}

\begin{abstract}
Contrary to the classical wisdom, processes with independent values
(defined properly) are much more diverse than white noises combined
with Poisson point processes, and product systems are much more
diverse than Fock spaces.

This text is a survey of recent progress in constructing and
investigating nonclassical stochastic flows and continuous products of
probability spaces and Hilbert spaces.
\end{abstract}



\end{frontmatter}

\setcounter{tocdepth}{2}

\tableofcontents

{
\renewcommand{\theequation}{0.\arabic{equation}}
\section*{Introduction}
\addcontentsline{toc}{section}{Introduction}
The famous Brownian motion $ B = (B_t)_{t\in[0,\infty)} $ is
an especially remarkable bizarre random function. Various fine
properties of Brownian sample paths have been investigated, but are beyond the
scope of this survey. Stochastic differential equations are a different (and
maybe more important) way of using $ B $; the stochastic differentials $
\D B_t $ are more relevant here than the Brownian path. The latter is rather
an infinitely divisible reservoir of independent random variables. Of course,
the path is not a differentiable function, but may be differentiated as a
generalized function (Schwarz distribution), giving the white noise $
\frac{\D}{\D t} B_t $. On the other hand, the Brownian motion may be thought
of as the scaling limit (as $ n \to \infty $) of the random walk,
\begin{equation}\label{1.3}
B_{k/n}^{(n)} = \frac{ \tau_0 + \dots + \tau_{k-1} }{ \sqrt n } \, ,
\end{equation}
where $ \tau_0, \tau_1, \dots $ are random signs (that is, i.i.d.\ random
variables taking on two values $ \pm 1 $ with probabilities $ 1/2, 1/2
$). Accordingly, the white noise may be thought of as the scaling limit of the
random locally integrable function (treated as a random Schwarz distribution)
$ W^{(n)} $,
\[
W^{(n)} (t) = \sqrt n \, \tau_k \qquad \text{for } t \in \(\tfrac k n,
\tfrac{k+1}n \) \, ;
\]
indeed, $ W^{(n)} $ converges in distribution (as $ n \to \infty $) to the
white noise.  This is a classical wisdom: the scaling limit of random signs
\emph{is} the white noise. Similarly, the scaling limit of a two-dimensional
array of random signs \emph{is} the white noise over the plane $ \R^2 $. These
are examples of processes with independent values. Conceptually, nothing is
simpler than independent values; but technically, they cannot be treated as
random functions.

A spectacular achievement of percolation theory (S.~Smirnov, 2001) is
existence (and conformal invariance) of the scaling limit of critical
site percolation on the triangular lattice (see for instance \cite{SW}
and references therein). The model is based on a two-dimensional array
of random signs (colors of vertices). Does it mean that the scaling
limit is driven by the white noise over the plane? No, it does
not. The percolation model uses the random signs in a nonclassical
way. However, the profound theory of percolation is only touched on in the
last section, \ref{11b}, of this survey. In order to understand the phenomenon
of nonclassicality we examine it on three levels of complexity: toy models
(Sect.~\ref{sec:1}), concentrated singularity (Sect.~\ref{sec:4}), distributed
singularity (Sect.~\ref{sec:7}). Some of the examined models may be treated as
scaling limits of \emph{oriented} percolation.

Contrary to the classical wisdom, processes with independent values
(defined properly) are much more diverse than white noises, Poisson
point processes and their combinations, time derivatives of L\'evy
processes. The L\'evy-It\^o theorem does not lie; processes with
independent increments are indeed exhausted by Brownian motions,
Poisson processes and their combinations. The scope of the classical
theory is limited by its treatment of independent values via
independent increments belonging to $ \R $ or another linear space (or
commutative group).

\begin{sloppypar}
Nowadays, independent increments are investigated also in
non-commutative groups and semigroups, consisting of homeomorphisms or
more general maps (say, $ \R \to \R $), kernels (that is, maps from
points to measures), bounded linear operators in a Hilbert space,
etc. These are relevant to stochastic flows. Some flows, being smooth
enough, are strong solutions of stochastic differential equations;
these flows are classical. Other flows contain some singularities
(turbulence, coalescence, stickiness, splitting etc.); these flows
tend to be nonclassical. They still are scaling limits of discrete
models driven by random signs, but these signs are used in a
nonclassical way.
\end{sloppypar}

The intuitive idea of a process with independent values appeared to be
deeper than its classical treatment. A general formalization of the
idea is given in Sect.~\ref{sec:2}. The distinction between classical and
nonclassical is formalized in Sect.~\ref{sec:5} by means of the concept
of stability/sensitivity, framed for the discrete case by computer scientists
and (in parallel) for the continuous case by probabilists. A model is
classical if and only if all random variables are stable. See also
Sect.~\ref{sec:6} for an equivalent definition in terms of L\'evy processes.

In discrete time, independence corresponds to the product of
probability spaces. In continuous time, a process with independent
values corresponds to a continuous product of probability spaces. The
corresponding Hilbert spaces $ L_2 $ of square integrable random
variables form a continuous tensor product of Hilbert spaces. Such
products are a notion well-known in analysis (the theory of operator
algebras) and relevant to noncommutative probability (and quantum
theory). A part of the (noncommutative) theory of such products has
recently been in close contact with the (commutative) probability theory
and is surveyed here (Sections \ref{sec:3}, \ref{sec:6}(d--g), \ref{9d},
\ref{sec:10}). The classical case is well-known as Fock spaces, or
type $ I $ Arveson systems. The nonclassical case consists of Arveson
systems of types $ II $ and $ III $. Their existence was revealed (in
different terms) in 1987 by R.~Powers \cite{Po87} (see also
\cite{Po99} and Chapter 13 of recent monograph \cite{Ar} by
W.~Arveson) while the corresponding probabilistic theory was in
latency; the idea was discussed repeatedly by A.~Vershik and
J.~Feldman, but the only publication \cite{Fe} was far from revealing
existence of nonclassical systems. The first nonclassical continuous
product of probability spaces was published in 1998 \cite{TV} by
A.~Vershik and the author. (I was introduced to the topic by
A.~Vershik in 1994.)

The classical part is well understood in both setups (commutative and
noncommutative), see Sect.~\ref{sec:6}. Some results on classical systems may
be found there; however, this survey cannot be used as an introduction to
classical systems. L\'evy-Khinchin-It\^o theory, Fock spaces and their
interrelations are covered elsewhere. About Arveson systems and their
relevance to quantum theory see the recent book by Arveson
\cite{Ar}. Sometimes the classical part is trivial, which is called `black
noise' in the commutative setup and `type $ III $' in the noncommutative
setup. Nonsmooth stochastic flows, homogeneous in space and time, can lead to
black noises, see Sect.~\ref{sec:7}. The interaction between noises theory and
turbulence theory enriches both sides.

Stochastic flows of unitary operators in a Hilbert space belong to
commutative probability, but are closely connected with noncommutative
probability. The connection is used in Sect.~\ref{sec:8} for proving
that such flows are classical.

The `classical/nonclassical' dichotomy is the starting point of a
more detailed classification, see Sect.~\ref{sec:9}.

Rich sources of examples for the noncommutative theory (type $ II $ and type $
III $) are found by the commutative theory, see Sect.~\ref{sec:10}.

Hopefully this survey contributes to the interaction between the
commutative and the noncommutative.

}

\mysection{Singularity concentrated in time (toy models)}
\label{sec:1}
\mysubsection{Preliminaries: flows}
\label{1a}

Independent increments are more interesting for us than values of a
process $ (X_t) $ with independent increments. Thus we denote
\[
X_{s,t} = X_t - X_s \quad \text{for } s < t \, ,
\]
use two properties
\begin{gather}
X_{r,t} = X_{r,s} + X_{s,t} \quad \text{for } r<s<t \, ,
 \label{1a1} \\
X_{t_1,t_2}, X_{t_2,t_3}, \dots, X_{t_{n-1},t_n} \text{ are
 independent for } t_1 < t_2 < \dots < t_n \, , \label{1a2}
\end{gather}
and discard $ X_t $. A two-parameter family $ (X_{s,t})_{s<t} $ of
random variables $ X_{s,t} : \Om \to \R $ satisfying \eqref{1a1},
\eqref{1a2} will be called an \flow{\R}%
\index{flow}\index{stochastic flow}
(or, more formally, a
stochastic \flow{(\R,+)}; here $ (\R,+) $ is the additive group of
real numbers). Why not just `a flow in $ \R $'? Since the latter is
widely used for a family of random diffeomorphisms of $ \R $ (or more
general maps, kernels etc). On the other hand, $ \R $ acts on itself
by shifts ($ y \mapsto x+y $), which justifies calling $
(X_{s,t})_{s<t} $ a stochastic flow. Similarly we may consider \flow{G}s for
some other groups $ G $ (such as $ \R^n $, $ \Z $, $ \Z_m $, $ \T $, see
below). Often, indices $ s,t $ of $ X_{s,t} $ run
over $ [0,\infty) $ each, but any other linearly ordered set may be specified,
if needed. Full generality is postponed to Sect.~\ref{sec:2}.

\mysubsection{Two examples}
\label{1b}

Tracing nonclassical behavior to bare bones we get very simple models
shown here. They may be discrete or continuous, this is a matter of
taste.

A stationary (not just `with stationary increments') random walk (or
Brownian motion) is impossible on such
groups as $ \Z $ or $ \R $, but possible on compact groups such as the
finite cyclic group $ \Z_m = \Z / m\Z $ or the circle $ \T = \R / \Z $.

\subsubsection*{Discrete example}

We choose $ m \in \{ 2,3,\dots \} $, take the group $ \Z_m $ and introduce a
\flow{\Z_m} $ X = (X_{s,t})_{s<t;s,t\in T} $ over the time set $ T = \{
0,1,2,\dots \} \cup \{\infty\} $ (with the natural linear order) as follows:
\begin{itemize}
\item[]
$ (X_{s,\infty})_{s=0,1,2,\dots} $ is a random walk on $ \Z_m $ started
according to the uniform distribution;
\item[]
$ X_{s,t} = X_{t,\infty} - X_{s,\infty} $ are its increments; $ X_{s,s+1} $ is
$ 0 $ or $ 1 $ with probabilities $ 1/2, 1/2 $. 
\end{itemize}
I claim that the flow $ X $ is both usual and strange.

The flow $ X $ is usual in the sense that it is the limit of a sequence of
very simple and natural \flow{\Z_m}s $ X^{(n)} = (X_{s,t}^{(n)} )_{s<t;s,t\in
T} $ constructed as follows:
\begin{gather*}
\Pr{ X_{s,s+1}^{(n)} = 0 } = \frac12 = \Pr{ X_{s,s+1}^{(n)} = 1 }
 \quad \text{for } s=0,\dots,n-1 \, ; \\
\Pr{ X_{s,\infty}^{(n)} = 0 } = 1 \quad \text{for } s = n, n+1, \dots
\end{gather*}
These conditions (together with \eqref{1a1}, \eqref{1a2}) determine
uniquely the joint distribution of all $
X_{s,t} $ for $ s,t \in T $, $ s<t $. Namely, $ X_{s,t}^{(n)} $ may be
thought of as increments of a random walk (in $ \Z_m $) stopped at the
instant $ n $. (Till now, $ \Z $ could be used instead of $ \Z_m $,
but the next claim would be violated; indeed, $ X_{0,\infty}^{(n)} $
would not be tight.) Random processes $ X^{(n)} $ converge in
distribution (for $ n \to \infty $) to a random process $ X $. (It
means weak convergence of finite-dimensional distributions, or
equivalently, probability measures on the compact space $
(\Z_m)^{ \{(s,t)\in T\times T : s<t \} } $, a product of countably
many finite topological spaces.) The limiting process $ X =
(X_{s,t})_{s<t;s,t\in T} $ is again a \flow{\Z_m}.

The flow $ X $ is strange in several respects (in contrast to $ X^{(n)} $). It
is not stable, as explained
in Sect.~\ref{1c} below. It leads to unusual states of an infinite chain of
quantum bits (spins), see Sect.~\ref{1d}. It generates nonclassical
continuous products, see Sections \ref{sec:2}--\ref{sec:6}. Here is a rather
informal discussion of the flow $ X $.

The random variable $ X_{0,\infty} $, distributed uniformly on $ \Z_m
$, is independent of the whole finite-time part of the process, $
(X_{s,t})_{s<t<\infty} $. (The same holds for each $ X_{t,\infty} $
separately, but surely not for $ X_{0,\infty} - X_{1,\infty} = X_{0,1}
$.) Therefore $
X_{0,\infty} $ is not a function of the i.i.d.\ sequence $
(X_{t,t+1})_{t<\infty} $. A paradox! You may guess that an additional
random variable $ X_{\infty-,\infty} $, independent of all $ X_{t,t+1}
$, squeezes somehow through a gap between finite numbers and
infinity. However, such an explanation does not work. It cannot happen
that $ X_{s,\infty} = f_s ( X_{s,s+1}, X_{s+1,s+2},\dots; \,
X_{\infty-,\infty} ) $ for all $ s $. Here is a proof (sketch).
Assume that it happens. The conditional distribution of $ X_{0,\infty}
$ given $ X_{0,1}, \dots, X_{s-1,s} $ and $ X_{\infty-,\infty} $ is
uniform on $ \Z_m $, since
\[
X_{0,\infty} = X_{0,s} + X_{s,t} + f_t ( X_{t,t+1}, X_{t+1,t+2},\dots;
\, X_{\infty-,\infty} ) \, ,
\]
and the conditional distribution of $ X_{s,t} $ is nearly uniform for
large $ t $. Thus, $ X_{0,\infty} $ is independent of $ X_{0,1},
X_{1,2}, \dots $ and $ X_{\infty-,\infty} $; a contradiction.

We see that some flows cannot be locally parameterized by independent
random variables. The group $ \Z_m $ is essential; every \flow{\Z}
(or \flow{\R}, or \flow{\R^n}) $ (X_{s,t})_{s<t;s,t\in T} $ can be
locally parameterized by $ X_{s,s+1} $ ($ s=0,1,\dots $) and $
X_{\infty-,\infty} \linebreak[0]
 = X_{0,\infty} - \lim_{k\to\infty} (X_{0,k}-c_k) $
for appropriate centering constants $ c_1, c_2, \dots $ (see also
Corollary \ref{6a1a}). It is also
essential that the time set $ T = \{ 0,1,2,\dots \} \cup \{\infty\} $
contains a limit point ($ \infty $). Otherwise, say, for $ T = \Z $,
every flow $ (X_{s,t})_{s<t} $ reduces to independent random variables
$ X_{s,s+1} $. Time reversal does not matter; the same phenomenon
manifests itself for $ T = \{ -k: k=0,1,\dots \} \cup \{-\infty\} $,
as well as $ T = \{ 2^{-k} : k=1,2,\dots \} \cup \{0\} $, or just $ t
= [0,\infty) $, the latter being used below.

\subsubsection*{Continuous example}

We take the group $ \T = \R / \Z $ (the circle). For every $ \eps
> 0 $ we construct a \flow{\T} $ Y^{(\eps)} = (Y_{s,t}^{(\eps)}
)_{s<t;s,t\in [0,\infty)} $ as follows:
\begin{gather*}
Y_{s,t}^{(\eps)} = B \Big( \ln \frac t \eps \Big) - B \Big( \ln \frac
 s \eps \Big) \pmod{1} \qquad \text{for } \eps \le s < t <
 \infty \, , \\
Y_{0,t}^{(\eps)} = 0 \quad \text{for } t \in [0,\eps] \, ;
\end{gather*}
here $ \(B(t)\)_{t\in[0,\infty)} $ is the usual Brownian motion.
These conditions determine uniquely the joint distribution of all $
Y_{s,t}^{(\eps)} $ for $ s,t \in [0,\infty) $, $ s<t $. Namely, $
Y_{s,t}^{(\eps)} $ may be thought of as increments of a process,
Brownian (on the circle) in logarithmic time after the instant $ \eps
$, but constant before $ \eps $. Random processes $ Y^{(\eps)} $
converge in distribution (as $ \eps \to 0 $) to a random process $ Y
$ (weak convergence of finite-dimensional distributions is meant), and
the limiting process $ Y = (Y_{s,t})_{s<t;s,t\in[0,\infty)} $ is again
a \flow{\T}. The random variable $ Y_{0,1} $, distributed uniformly on
$ \T $, is independent of all $ Y_{s,t} $ for $ 0 < s < t $. We cannot
locally parametrize $ Y $ by increments of a Brownian motion (and
possibly an additional random variable $ Y_{0,0+} $ independent of the
Brownian motion).

The one-parameter random process $ (Y_{0,\E^t})_{t\in\R} $ is a
stationary Brownian motion in $ \T $. The complex-valued random
process $ (Z_t)_{t\in[0,\infty)} $,
\[
Z_t = \sqrt t \, \E^{2\pi \I Y_{0,t}} \quad \text{for } t \in [0,\infty)
\]
is a continuous martingale, and satisfies the stochastic differential
equation
\[
\D Z_t = \frac{\I}{\sqrt t} Z_t \, \D B_t \, , \quad Z_0 = 0 \, ,
\]
where $ (B_t)_{t\in[0,\infty)} $ is the usual Brownian
motion. However, the random variable $ Z_1 $ is independent of the
whole Brownian motion $ (B_t) $. The weak solution of the stochastic
differential equation is not a strong solution. See also \cite{Yor92},
\cite{ES}, and \cite[Sect.~1a]{Ts03}.

\mysubsection{Stability and sensitivity}
\label{1c}

Stability and sensitivity of Boolean functions of many Boolean
variables were introduced in 1999 by Benjamini, Kalai and Schramm
\cite{BKS} and applied to percolation, random graphs etc. They
introduce errors (perturbation) into a given Boolean array by flipping
each Boolean variable with a small probability (independently of
others), and observe the effect of these errors by comparing the new
(perturbed) value of a given Boolean function with its original
(unperturbed) value. They prove that percolation is sensitive!
Surprisingly, their `stability' is basically the same as our
`classicality'. See also \cite{ST}.

The \flow{\Z_m} $ X $ of Sect.~\ref{1b} (denote it here by $
X^{\text{\ref{1b}}} $) contains i.i.d.\ random variables $ X_{s,s+1} $ that
are Boolean in the sense that each one takes on two values $ 0 $ and $
1 $, with probabilities $1/2, 1/2 $. However, $ X_{0,\infty} $ is not
a function of these (Boolean) variables. Here is a proper
formalization of the idea. A pair of two correlated \flow{G}s (one
`unperturbed', the other `perturbed') is a \flow{(G \times G)}; here $
G \times G $ is the direct product, that is, the set of all pairs $
(g_1,g_2) $ for $ g_1,g_2 \in G $ with the group operation $ (g_1,g_2)
(g_3,g_4) = (g_1 g_3, g_2 g_4) $. (For $ G = \Z_m $ we prefer additive
notation: $ (g_1,g_2) + (g_3,g_4) = (g_1 + g_3, g_2 + g_4) \in \Z_m
\oplus \Z_m $.) Let $ X = (X_{s,t})_{s<t;s,t\in T} $ be a \flow{(G
\times G)}; we have $ X_{s,t} = ( X'_{s,t}, X''_{s,t} ) $, that is, $ X
= (X',X'') $, where $ X', X'' $ are \flow{G}s on the same
probability space $ (\Om,\F,P) $. Note that \eqref{1a2} constrains the pair,
not just $ X' $ and $ X'' $ separately. Let $ \F', \F'' $ be sub-\sif s of $
\F $ generated by $ X', X'' $ respectively. We introduce the maximal
correlation\index{maximal correlation}
\[
\rho_{\max} (X) = \rho_{\max} (X',X'') = \rho_{\max} (\F',\F'') = \sup
| \Ex(fg) | \, ,
\]
where the supremum is taken over all $ f \in L_2(\Om,\F',P) $, $ g \in
L_2(\Om,\F'',P) $ such that $ \Ex f = 0 $, $ \Ex g = 0 $, $ \Var f \le
1 $, $ \Var g \le 1 $. The idea of a (non-degenerate) perturbation of
a flow may be formalized by the condition $ \rho_{\max} (X) < 1 $.

\begin{proposition}\label{1c1}
Let $ X = (X',X'') $ be a \flow{(\Z_m \oplus \Z_m)} such that $ X',X''
$ both are distributed like $ X^{\text{\ref{1b}}} $ (the discrete
\flow{\Z_m} of Sect.~\ref{1b}). If $ { \rho_{\max}(X) < 1 } $ then
random variables $ X'_{0,\infty} $ and $ X''_{0,\infty} $ are
independent.
\end{proposition}

A small perturbation has a dramatic effect on the random variable $
X^{\text{\ref{1b}}}_{0,\infty} $; this is instability (and moreover,
sensitivity). All flows in Proposition \ref{1c1} use the
time set $ T = \{0,1,2,\dots\} \cup \{\infty\} $. Nothing like that
happens on $ T = \{0,1,2,\dots\} $ or $ T=\Z $. Also, the group $ \Z_m
$ is essential; nothing like that happens for \flow{\Z}s (or
\flow{\R}s, or \flow{\R^n}s; see also Corollary \ref{6a1a}).
\emph{Sketch of the proof of Proposition \textup{\ref{1c1}} for the
special case $ m=2 $:}
\begin{multline*}
| \Ex (-1)^{X'_{0,\infty}} (-1)^{X''_{0,\infty}} | = \\
= | \Ex (-1)^{X'_{t,\infty}} (-1)^{X''_{t,\infty}} | \cdot
\prod_{s=0}^{t-1} | \Ex (-1)^{X'_{s,s+1}} (-1)^{X''_{s,s+1}} | \le 1
\cdot ( \rho_{\max}(X) )^t \xrightarrow[t\to\infty]{} 0 \, .
\end{multline*}

The same can be said about the other (continuous) example $
Y^{\text{\ref{1b}}} $ of Sect.~\ref{1b}. The random variable $
Y^{\text{\ref{1b}}}_{0,1} $ is sensitive.

See also Sect.~\ref{sec:5}, especially \ref{5d}.

\mysubsection{Hilbert spaces, quantum bits (spins)}
\label{1d}

Stochastic flows belong to (and are of interest to) commutative probability
theory; in addition, they can help to noncommutative probability theory by
providing new models with unusual properties, which will be discussed in
detail in Sections \ref{sec:3} and \ref{sec:10}. As a simple illustration of
the idea, this approach is applied below to our first toy model.

We return to $ X = X^{\text{\ref{1b}}} $ (the discrete \flow{\Z_m} of
Sect.~\ref{1b}) and consider the Hilbert space $ H $ of all square
integrable complex-valued measurable functions of random variables $
X_{s,t} $, with the norm
\[
\| f ( X_{0,1}, X_{1,2}, \dots; X_{0,\infty} ) \|^2 = \Ex | f (
X_{0,1}, X_{1,2}, \dots; X_{0,\infty} ) |^2
\]
(other $ X_{s,t} $ are redundant). Equivalently, $ H = L_2 (\Om,\F,P)
$ provided that $ \F $ is the \sif\ generated by $ X $.

We may split $ X $ at the instant $ 1 $ in two independent components:
the past, --- just a single random variable $ X_{0,1} $; and the
future, --- all $ X_{s,t} $ for $ 1 \le s < t $, $ s,t \in T
$. Accordingly, $ H $ splits into the tensor product, $ H = H_{0,1}
\otimes H_{1,\infty} $. The Hilbert space $ H_{0,1} $ is
two-dimensional (since $ X_{0,1} $ takes on two values, $ 0 $ and $ 1
$), spanned by two orthonormal vectors $ 2^{-1/2} (1-X_{0,1}) $ and $
2^{-1/2} X_{0,1} $. Using this basis, we may treat the famous Pauli
spin matrices
\[
\si_1 = \begin{pmatrix} 0 & 1 \\ 1 & 0 \end{pmatrix} \, , \quad
\si_2 = \begin{pmatrix} 0 & -\I \\ \I & 0 \end{pmatrix} \, , \quad
\si_3 = \begin{pmatrix} 1 & 0 \\ 0 & -1 \end{pmatrix}
\]
as operators on $ H_{0,1} $, and also on $ H $ (identifying $ \si_k $
with $ \si_k \otimes \One $). Thus, the $ 2 \times 2 $ matrix algebra
$ \M_2(\C) $ acts on $ H $.

\begin{sloppypar}
Similarly,
\[
H = H_{0,1} \otimes \dots \otimes H_{t-1,t} \otimes H_{t,\infty}
\]
for any $ t = 1,2,\dots $; each $ H_{s-1,s} $ is two-dimensional, and
$ \{ 2^{-1/2} (1-X_{s-1,s}), \linebreak[0] 2^{-1/2} X_{s-1,s} \} $ is
its orthonormal basis. We get commuting copies of $ \M_2(\C) $;
\begin{gather*}
\si_k^{t-1,t} : H \to H \quad \text{for } k=1,2,3 \text{ and }
 t=1,2,\dots \, ; \\
[ \si_k^{s-1,s}, \si_k^{t-1,t} ] = 0 \quad \text{for } s \ne t \, .
\end{gather*}
The random variable $ \exp \( \frac{2\pi \I}m X_{0,\infty} \) $ is a
factorizing vector,
\[
\exp \Big( \frac{2\pi\I}m X_{0,\infty} \Big) = \exp \Big(
\frac{2\pi\I}m X_{0,1} \Big) \otimes \dots \otimes \exp \Big(
\frac{2\pi\I}m X_{t-1,t} \Big) \otimes \exp \Big( \frac{2\pi\I}m
X_{t,\infty} \Big) \, .
\]
The first factor is not a basis vector but a linear combination,
\[
\exp \Big( \frac{2\pi\I}m X_{0,1} \Big) = \frac1{\sqrt2}
\begin{pmatrix} 1 \\ \exp\frac{2\pi\I}m \end{pmatrix} \, ;
\]
each factor $ \exp \( \frac{2\pi\I}m X_{t-1,t} \) $ is a copy of this
vector. The quantum state on a local algebra,
\[
M \mapsto \ip{M\psi}{\psi} \, , \quad M \in M_{2^t} (\C) =
\underbrace{ M_2(\C) \otimes \dots \otimes M_2(\C) }_{t}
\]
corresponding to the vector $ \psi = \exp \( \frac{2\pi\I}m
X_{0,\infty} \) \in H $ is equal to the quantum state corresponding to
the vector
\[
\begin{pmatrix} 2^{-1/2} \\ 2^{-1/2} \exp\frac{2\pi\I}m \end{pmatrix}
\otimes \dots \otimes \begin{pmatrix} 2^{-1/2} \\ 2^{-1/2}
\exp\frac{2\pi\I}m \end{pmatrix} =
\begin{pmatrix} 2^{-1/2} \\ 2^{-1/2} \exp\frac{2\pi\I}m
\end{pmatrix}^{\otimes t} \in \C^{2^t} \, . 
\]
The vector $ \exp \( \frac{2\pi\I}m X_{0,\infty} \) \in H $ may be
interpreted as a factorizing state $ \( \begin{smallmatrix} 2^{-1/2}
\\ 2^{-1/2} \exp\frac{2\pi\I}m \end{smallmatrix} \)^{\otimes\infty} $
of an infinite collection of spins. Similarly, for each $ k =
0,1,\dots,m-1 $ the vector $ \exp \( \frac{2\pi\I k}m X_{0,\infty} \)
\in H $ may be  interpreted as $ \( \begin{smallmatrix} 2^{-1/2} \\
2^{-1/2} \exp\frac{2\pi\I k}m \end{smallmatrix} \)^{\otimes\infty} $.
\end{sloppypar}

Vectors of $ H $ of the form
\[
f(X_{0,1}, X_{1,2}, \dots) \exp \Big( \frac{2\pi\I k}m X_{0,\infty}
\Big)
\]
are a subspace $ H_k \subset H $ invariant under all local operators;
and $ H = H_0 \oplus \dots \oplus H_{m-1} $. Each $ H_k $ is
irreducible in the sense that it contains no nontrivial subspace
invariant under all local operators.

The operator $ R $ defined by
\[
R f(X_{0,1}, X_{1,2}, \dots; X_{0,\infty}) = f(X_{0,1}, X_{1,2},
\dots; X_{0,\infty}+1) 
\]
(where $ X_{0,\infty}+1 $ is treated $ \bmod m $) commutes with all
local operators. Every operator commuting with all local operators is
a function of $ R $. The subspaces $ H_k $ are
eigenspaces of $ R $, their eigenvalues being $ \exp \( \frac{2\pi\I
k}m \) $. In a more physical language, the group $ \{ R^k :
k=0,1,\dots,m-1 \} $ is the gauge group, and $ H_k $ are
superselection sectors. Each sector has its own asymptotic behavior of
remote spins, according to $ \exp \( \frac{2\pi\I}m X_{0,\infty} \) = \(
\begin{smallmatrix} 2^{-1/2} \\ 2^{-1/2} \exp\frac{2\pi\I}m \end{smallmatrix}
\)^{\otimes\infty} $.

In contrast, the classical model based on functions of i.i.d.\ random
variables $ X_{s,s+1} $, may be identified with just one sector (the zero
sector, $ k=0 $) of our model.

See also Sect.~\ref{sec:3}, \cite[Appendix]{Ts99},
\cite[Sect.~8.4]{Li}.

\mysection{Singularity concentrated in space (examples)}
\label{sec:4}
\mysubsection{Preliminaries: convolution semigroups, stationary flows, noises}
\label{4aa}

A weakly continuous (one-parameter)
convolution semigroup in $ \R $\index{convolution semigroup}
is a family $
(\mu_t)_{t\in(0,\infty)} $ of probability measures $ \mu_t $ on $ \R $
such that $ \mu_s * \mu_t = \mu_{s+t} $ for all $ s,t \in (0,\infty)
$, and $ \lim_{t\to0} \mu_t \( (-\eps,\eps) \) = 1 $ for all $ \eps >
0 $. Two basic cases are normal distributions $ \N(0,t) $ and Poisson
distributions $ \P(t) $. They correspond to the Brownian motion and the
Poisson process, respectively. Every convolution semigroup decomposes
into a combination of these two basic cases, and corresponds to a
process with independent increments; the process decomposes into
Brownian and Poisson processes. That is the classical theory
(L\'evy-Khinchin-It\^o).

Topological groups $ G $ other than $ \R $ may be used as well. Simple
examples are $ \R^n $, $ \Z $, $ \Z_m $, $ \T $ mentioned in
Sect.~\ref{sec:1}. More advanced, noncommutative topological groups consisting
of diffeomorphisms, unitary operators etc.\ will be used in
Sect.~\ref{sec:8}. Topological semigroups $ G $ are also useful; for instance,
the multiplicative semigroup $ (\C,\cdot) $ of complex numbers (including $ 0
$) will appear in Sect.~\ref{6b}. More advanced examples, consisting of
conformal endomorphisms or (not just invertible) linear operators, will appear
in Sect.~\ref{sec:8}.

In the present section, elements of $ G $ are some quite simple maps $ \R \to
\R $ or $ [0,\infty) \to [0,\infty) $, given by explicit formulas with (at
most) three parameters. The binary operation $ G \times G \to G $ is the
composition of maps; $ G $ is a noncommutative semigroup. The maps are
non-invertible; $ G $ is not a group. The maps are discontinuous; the binary
operation $ G \times G \to G $ is also discontinuous, thus, $ G $ is not a
topological semigroup. It is a semigroup and a (finite-dimensional)
topological space, but still, not a topological semigroup! We axiomatize such
objects as follows.

\begin{definition}\label{2c6}
A \emph{topo-semigroup}\index{topo-semigroup}
is a semigroup $ G $ equipped with a topology such that

(a) $ G $ is a separable metrizable topological space, (complete or) Borel
measurable in its completion;\footnote{%
 The choice of a metric does not matter (even though it can change the
 completion).}

(b) the binary operation is a Borel map $ G \times G \to G $;

(c) the semigroup $ G $ contains a unit $ 1 $, and
\[
\text{if } x_n \to 1 \text{ and } z_n \to 1 \text{ then } x_n y z_n \to y
\]
for all $ x_1, x_2, \dots; y; z_1, z_2, \dots \in G $.
\end{definition}

The convolution $ \mu*\nu $ of two probability measures $ \mu, \nu $
on a topo-semigroup $ G $ is defined evidently (as the image of $ \mu
\times \nu $ under $ (x,y) \mapsto xy $). A (one-parameter)
convolution semigroup $ (\mu_t)_{t\in(0,\infty)} $ in $ G $ is defined
accordingly. We call it \emph{weakly continuous,} if
\begin{equation}\label{4aa*}
\mu_t (U) \to 1 \quad \text{as } t \to 0
\end{equation}
for every neighborhood $ U $ of the unit $ 1 $ of $ G $. It follows
easily that $ \int f \, d\mu_t $ is continuous in $ t \in (0,\infty) $
(and tends to $ f(1) $ as $ t \to 0 $) for every bounded continuous
function $ f : G \to \R $.

When $ G $ is a group, we have two equivalent languages: processes with
independent increments, and flows. However, neither the formula $ X_{s,t} =
X_t - X_s $ nor its noncommutative counterpart, $ X_{s,t} = X_s^{-1} X_t $,
can be used in a semigroup. The only appropriate language is, stochastic
flows.

\begin{proposition}\label{2c7}
Let $ G $ be a topo-semigroup and $ (\mu_t)_{t>0} $ a weakly
continuous convolution semigroup in $ G $. Then there exists a family $
(X_{s,t})_{-\infty<s<t<\infty} $ of \valued{G} random variables $ X_{s,t} $
satisfying the four conditions:

(a) $ X_{s,t} $ is distributed $ \mu_{t-s} $, whenever $ -\infty<s<t<\infty $;

(b) $ X_{t_1,t_2}, X_{t_2,t_3}, \dots, X_{t_{n-1},t_n} $ are independent for $
t_1 < t_2 < \dots < t_n $;

(c) $ X_{r,t} = X_{r,s} X_{s,t} $ a.s.\ whenever $ r < s < t $;

(d) $ X_{s_n,t_n} \to X_{s,t} $ in probability whenever $ s < t $, $ s_n
\downarrow s $ and $ t_n \uparrow t $.
\end{proposition}

The proof is postponed to Sect.~\ref{sec:2}.

The family $ (X_{s,t})_{s<t} $ (determined uniquely in distribution) may be
called the (stationary) \flow{G} corresponding to $ (\mu_t)_{t>0} $.

Denoting by $ \F_{s,t} $ the sub-\sif\ generated by random variables $ X_{u,v}
$ for all $ u,v $ such that $ s \le u < v \le t $, we get
\[
\F_{r,s} \otimes \F_{s,t} = \F_{r,t} \quad \text{whenever } r<s<t \,
.
\]
(That is, $ \F_{r,s} $ and $ \F_{s,t} $ are independent and generate $
\F_{r,t} $.) In addition, the time shift by $ h $ sends $ \F_{s,t} $ to $
\F_{s+h,t+h} $ (stationarity); see Sect.~\ref{sec:2} for details. Such a
family of sub-\sif s (and time shifts) is called a noise (see Def.~\ref{2c1}
later).

A probabilist might feel that noises are too abstract; \sif s do not
catch distributions. (Similarly a geometer might complain that
topological invariants do not catch volumes.) However, they do! The
delusion is suggested by the discrete-time counterpart. Indeed, the
product of countably many copies of a probability space does not
distinguish any specific random variable (or distribution). Continuous
time is quite different. Consider for example the white noise $
(\F_{s,t}) $, $ (T_h) $, corresponding to the \flow{\R} $ X_{s,t} =
B_t - B_s $ of Brownian increments. At first sight, $ X_{s,t} $ cannot
be reconstructed from $ (\F_{s,t}) $ and $ (T_h) $, but in fact they
can! The conditions
\begin{gather*}
X_{s,t} \text{ is \measurable{\F_{s,t}}} \, , \\
X_{r,s} + X_{s,t} = X_{r,t} \, , \\
X_{s,t} \circ T_h = X_{s+h,t+h} \, , \\
\Ex X_{s,t} = 0 \, , \quad \Ex X^2_{s,t} = t-s
\end{gather*}
determine them uniquely up to a sign; $ X_{s,t} = \pm (B_t - B_s)
$. For the Poisson noise the situation is similar. However, for a
L\'evy process with different jump sizes, only their rates are encoded
in $ (\F_{s,t}) $, $ (T_h) $; the sizes are lost.

Every weakly continuous convolution semigroup in a topo-semigroup leads to a
stationarity flow and further to a noise. (See Sect.~\ref{sec:2} for details.)

\mysubsection{Coalescence: another way to the white noise}
\label{4a}

A model described here is itself of little interest, but helps to
understand more interesting models introduced afterwards.

Functions $ [0,\infty) \to [0,\infty) $ of the form $ f_{a,b} $,
\[
\begin{gathered}
f_{a,b} (x) = a + \max (x,b) \, ,
\end{gathered}
\qquad
\begin{gathered}\includegraphics{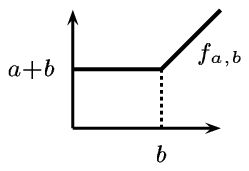}\end{gathered}
\]
for $ a,b \in \R $, $ b \ge 0 $, $ a+b \ge 0 $, form a semigroup $ G
$. That is, the composition $ fg = g \circ f : x \mapsto g(f(x)) $ of
two such functions is such a function, again:
\[
\begin{gathered}
f_{a_1,b_1} f_{a_2,b_2} = f_{a,b} \, , \quad
\begin{aligned}
 a &= a_1 + a_2 \, , \\
 b &= \max ( b_1, b_2-a_1 ) \, .
\end{aligned}
\end{gathered}
\]
Equipped with the evident topology, $ G $ is a two-dimensional
topological semigroup. The following probability distributions are a
weakly continuous convolution semigroup $ (\mu_t)_{t>0} $ in $ G $:
\begin{equation*}\begin{split}
\frac{ \mu_t (\D a \D b) }{ \D a \D b } = \frac{ 2 (a+2b) }{
\sqrt{2\pi} \, t^{3/2} } \exp \bigg( \! - \frac{ (a+2b)^2 }{ 2t }
\bigg) \, .
\end{split}\end{equation*}
It leads to a stationary \flow{G} $ (X_{s,t})_{s<t} $; $ X_{s,t} =
f_{a_{s,t},b_{s,t}} $.

The map $ (a,b) \mapsto a $ is a homomorphism $ G \to (\R,+) $. It
sends $ \mu_t $ to the normal distribution $ \N(0,t) $, which means
that $ a_{s,t} $ is nothing but the increment of the standard Brownian
motion $ (a_{0,t})_t $ in $ \R $. It appears that
\[
\begin{aligned}
b_{r,t} &= - \min_{s\in[r,t]} a_{r,s} \, , \\
a_{r,t} + b_{r,t} &= \max_{s\in[r,t]} a_{s,t} \, . \\
\end{aligned} \qquad
\begin{gathered}\includegraphics{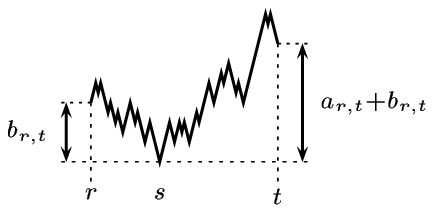}\end{gathered}
\]
The `two-dimensional nature' of the flow is a delusion; the second
dimension $ b $ reduces to the first dimension $ a $. The noise
generated by this \flow{G} is (isomorphic to) the white noise.

The \flow{G} $ (X_{s,t})_{s<t} $ may be treated as the scaling limit
of a discrete-time \flow{G} formed by (compositions of) two
functions $ f_+, f_- : \Z_+ \to \Z_+ $ (chosen equiprobably),
\[
\begin{gathered}\includegraphics{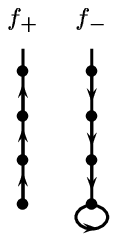}\end{gathered}
\qquad
\begin{gathered}
f_+ (x) = x+1 \, , \quad f_- (x) = \max ( 0, x-1 ) \, .
\end{gathered}
\]
The semigroup spanned by $ f_-, f_+ $ may also be treated as the
semigroup (with unit, non-commutative) defined by two generators $
f_-, f_+ $ and a single relation $ f_+ f_- = 1 $. (The second
relation $ f_- f_+ = 1 $ would turn the semigroup into $ \Z $, giving
in the scaling limit the homomorphism $ G \to \R $ mentioned above.)

Our $ G $ is not just a semigroup, but a semigroup of maps; it acts on $
[0,\infty) $. Thus, any \flow{G} $ (X_{s,t})_{s<t} $ leads to the so-called
one-point motion, the random process $ (X_{0,t}(x))_{t>0} $, provided that a
starting point $ x \in [0,\infty) $ is chosen. Similarly, the two-point motion
is the two-dimensional random process $ \( X_{0,t}(x_1), X_{0,t}(x_2) \) $;
and so on. For our specific \flow{G}, the one-point motion is
(distributed like) the reflecting Brownian motion (starting at $ x $). Two
particles starting at $ x_1, x_2 $ ($ x_1 < x_2 $) keep their distance
($ X_{0,t}(x_2) - X_{0,t}(x_1) = x_2 - x_1 $) as long as the boundary
is not hit ($ X_{0,t}(x_1) > 0 $). In general, the distance decreases
in time. At some instant $ s $ (when $ b_{0,s} $ reaches $ x_2 $) the
two particles coalesce at the boundary point ($ X_{0,s}(x_1) =
X_{0,s}(x_2) = 0 $) and never diverge afterwards ($ X_{0,t}(x_1) =
X_{0,t}(x_2) $ for all $ t \in [s,\infty) $).

\mysubsection{Splitting: a nonclassical noise}
\label{4b}

Functions $ \R \to \R $ of two forms, $ f^-_{a,b} $ and $ f^+_{a,b} $,
\begin{gather*}
 \begin{gathered}
  f^-_{a,b} (x) = \begin{cases}
   x-a &\text{for $ x \in (-\infty,-b) $},\\
   -a-b &\text{for $ x \in [-b,b] $},\\
   x+a &\text{for $ x \in (b,\infty) $};
 \end{cases}
 \end{gathered}
 \qquad
 \begin{gathered}\includegraphics{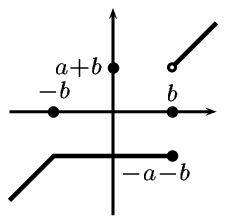}\end{gathered}
 \displaybreak[0]\\
 \begin{gathered}
  f^+_{a,b} (x) = \begin{cases}
   x-a &\text{for $ x \in (-\infty,-b) $},\\
   a+b &\text{for $ x \in [-b,b] $},\\
   x+a &\text{for $ x \in (b,\infty) $}
 \end{cases}
 \end{gathered}
 \qquad
 \begin{gathered}\includegraphics{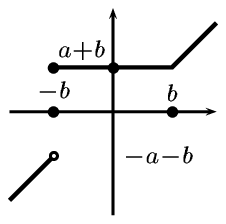}\end{gathered}
\end{gather*}
for $ a,b \in \R $, $ b \ge 0 $, $ a+b \ge 0 $, form a semigroup $ G
$. It is not a topological semigroup, since the composition is not
continuous, but it is a topo-semigroup (as defined by \ref{2c6}). The
map $ f^-_{a,b} \mapsto f^{\text{\ref{4a}}}_{a,b} $, $ f^+_{a,b} \mapsto
f^{\text{\ref{4a}}}_{a,b} $ is a homomorphism $ G \to G^{\text{\ref{4a}}} $; here
$ G^{\text{\ref{4a}}} $ stands for the semigroup
denoted by $ G $ in Sect.~\ref{4a}. We define a measure $ \mu_t $ on $
G $ by two conditions: first, the homomorphism $ G \to G^{\text{\ref{4a}}} $
sends $ \mu_t $ to $ \mu_t^{\text{\ref{4a}}} $, and second, $ \mu_t $ is
invariant under the map $ f^-_{a,b} \mapsto f^+_{a,b} $, $ f^+_{a,b}
\mapsto f^-_{a,b} $. In other words, $ a $ and $ b $ are distributed
as in Sect.~\ref{4a}, while the third parameter is `$ - $' or `$ + $'
with probabilities $ 1/2 $, $ 1/2 $, independently of $ a,b $. These
distributions are a convolution semigroup. Proposition \ref{2c7} gives
us a stationary \flow{G} $ (X_{s,t})_{s<t} $ and a noise, ---
\emph{the noise of splitting.}\index{noise!of splitting}
It is a nonclassical noise! (See Sect.~\ref{5d}.)

The \flow{G} $ (X_{s,t})_{s<t} $ may be treated as the scaling limit
of a discrete-time \flow{G} formed by (compositions of) two
functions $ f_+, f_- : \Z + \frac12 \to \Z + \frac12 $ (chosen
equiprobably),
\[
\begin{gathered}
\begin{gathered} f_- (x) = x-1 \, , \\ f_+ (x) = x+1 \end{gathered}
 \quad \text{for } x \in \( \Z + \tfrac12 \) \cap (0,\infty) \, , \\ 
 f_- (-x) = - f_- (x) \, , \quad f_+ (-x) = - f_+ (x) \, .
\end{gathered}
\qquad
\begin{gathered}\includegraphics{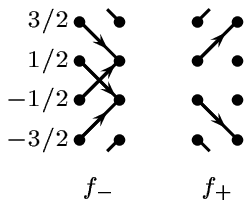}\end{gathered}
\]
They satisfy the relation $ f_+ f_- = 1 $ and generate the same
(discrete) semigroup as in Sect.~\ref{4a}, but the scaling limit is
different, since here (in contrast to \ref{4a}) the product $ (f_-)^b
(f_+)^{a+b} $ for large $ a,b $ is sensitive to $ (-1)^b $.

The \flow{G} $ (X_{s,t})_{s<t} $ is intertwined with the
\flow{G^{\text{\ref{4a}}}} $ (X^{\text{\ref{4a}}}_{s,t})_{s<t} $ by the map $ \R
\to [0,\infty) $, $ x \mapsto |x| $. Indeed, $ | f^\pm_{a,b} (x) | =
f^{\text{\ref{4a}}}_{a,b} (|x|) $. The radial part $ |X_{0,t}(x)| $ is
(distributed like) the coalescing flow of Sect.~\ref{4a}. The sign of 
$ X_{0,t}(x) $, being independent of the radial motion, is chosen anew
each time when the radial motion starts an excursion. The one-point
motion is just the standard Brownian motion in $ \R $.

Two particles starting at $ x_1, x_2 $ ($ |x_1| < |x_2| $) keep the value
$ |X_{0,t}(x_2)| - |X_{0,t}(x_1)| = |x_2| - |x_1| $ as long as
$ |X_{0,t}(x_1)| > 0 $. In general, the value decreases in time. At some
instant the two particles coalesce at the origin and never diverge afterwards.
Before the coalescence the second particle does not hit the origin, while the
first particle chooses the sign anew each time when it starts an excursion. It
does so after the coalescence, too, but now --- together with the second
particle.

Similarly we may take the space set as the union $ \{ z \in \C : z^3
\in [0,\infty) \} $ of three (or more) rays on the complex plane and
define a splitting flow such that its radial part is the coalescing
flow, and the argument (the angular part) is chosen anew (with
probabilities $ 1/3, 1/3, 1/3 $) each time when starting an
excursion. Then the one-point motion is a complex-valued martingale
known as the spider martingale, see \cite[Sect.~2]{BEKSY}.

The noise of splitting was introduced and investigated by J.~Warren
\cite{Wa1}. See also \cite{Wat0}, \cite[Example 1d1]{Ts03}, and
Sections \ref{4d}, \ref{5d} of this survey.

\mysubsection{Stickiness: a time-asymmetric noise}
\label{4c}

Functions $ [0,\infty) \to [0,\infty) $ of the form $ f_{a,b,c} $,
\begin{equation*}
\begin{gathered}
f_{a,b,c} (x) = \begin{cases}
 c &\text{for $ x \in [0,b] $}, \\
 x+a &\text{for $ x \in (b,\infty) $}
\end{cases}
\end{gathered}
\qquad
\begin{gathered}\includegraphics{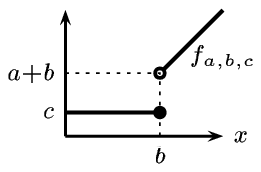}\end{gathered}
\end{equation*}
for $ a,b,c \in \R $, $ b \ge 0 $, $ 0 \le c \le a+b $, form a
semigroup $ G $ (a topo-semigroup, not topological). The map $
f_{a,b,c} \mapsto f^{\text{\ref{4a}}}_{a,b} $ is a homomorphism $ G \to
G^{\text{\ref{4a}}} $. In fact, $ G^{\text{\ref{4a}}} = \{ f_{a,b,a+b} \} $ is a
sub-semigroup of $ G $, therefore the convolution semigroup $
(\mu_t^{\text{\ref{4a}}})_{t>0} $ in $ G^{\text{\ref{4a}}} $ is also a convolution
semigroup in $ G $; it is a degenerate case ($ \la=0 $) of a family of
convolution semigroups $ (\mu_t^{(\la)})_{t>0} $ on $ G $; the
parameter $ \la $ runs over $ (0,\infty) $. Namely,
\begin{equation*}
\begin{gathered}
\mu_t^{(\la)} \text{ is the joint distribution of $ a,b $ and $
 c=\max(0,a+b-\la\eta) $}, \\ 
\text{where the pair $ (a,b) $ is distributed $ \mu_t^{\text{\ref{4a}}} $},
 \\
 \text{while } \eta \text{ is independent of $ (a,b) $ and distributed
 $ \Exp(1) $} \, ;
\end{gathered}
\end{equation*}
that is, $ \Pr{ \eta > c } = \E^{-c} $ for $ c \in [0,\infty) $.
It is indeed a convolution semigroup, due to a property of the
composition in $ G $ ($ c = a_2 + c_1 $ if $ c_1 > b_2 $, otherwise $
c_2 $; about $ a,b $ see Sect.~\ref{4a}): for every $ a_1, b_1, a_2,
b_2 $, if $ c_1 \sim \max(0,a_1+b_1-\la\eta) $ and $ c_2 \sim
\max(0,a_2+b_2-\la\eta) $ are independent then $ c \sim
\max(0,a+b-\la\eta) $.
\[
\begin{gathered}\includegraphics{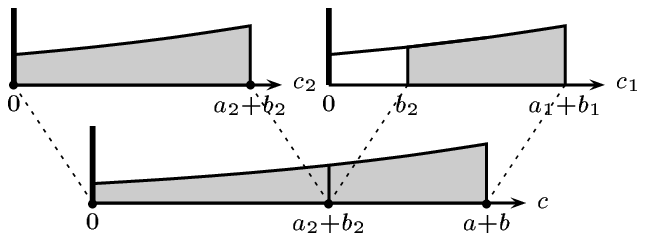}\end{gathered} \quad
 \parbox{5cm}{\small Verifying the $ c $-component of the equality $
  \mu_s^{(\la)} * \mu_t^{(\la)} = \mu_{s+t}^{(\la)} $.\\
  (The case $ a_1 + b_1 > b_2 $ is shown; the other case is trivial.)}
\]
Note that the measure $ \mu_t^{(\la)} $ has an absolutely continuous
part (its three-dimensio\-nal density can be written explicitly, using
the two-dimensional density of $ \mu_t^{\text{\ref{4a}}} $ and the
one-dimensional exponential density of $ \eta $) and a singular part
concentrated on the plane $ c=0 $; the singular part has a
two-dimensional density (it can also be written
explicitly). Proposition \ref{2c7} gives us a stationary \flow{G} $
(X_{s,t})_{s<t} $ and a noise, --- 
\emph{the noise of stickiness.}\index{noise!of stickiness}
It is a nonclassical noise. Moreover, the noise is time-asymmetric!
(See Sect.~\ref{4e}.)

The \flow{G} $ (X_{s,t})_{s<t} $ may be treated as the scaling limit
of a discrete-time \flow{G} formed by (compositions of) three
functions $ f_+, f_-, f_* : \Z_+ \to \Z_+ $:
\begin{equation*}
\begin{gathered}
 \includegraphics{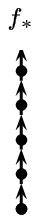}\hspace{-2mm}
 \includegraphics{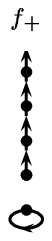}\hspace{-2mm}
 \includegraphics{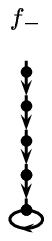}
\end{gathered}
\quad
\begin{gathered}
f_* (x) = x+1 \, , \quad f_- (x) = \max ( 0, x-1 ) \, , \\
f_+ (x) = \begin{cases}
 x+1 &\text{ for $ x>0 $},\\
 0 &\text{for $ x=0 $}.
\end{cases}
\end{gathered}
\end{equation*}
The functions are chosen with probabilities $ \Pr{f_-} = \frac12 $, $
\Pr{f_*} = \frac1{2\la} \sqrt{\De t} $, $ \Pr{f_+} = \frac12 -
\frac1{2\la} \sqrt{\De t} ) $, where $ \De t $ is the time pitch
(tending to $ 0 $ in the scaling limit); the space pitch is equal to $
\sqrt{\De t} $. The semigroup spanned by $ f_-, f_+, f_* $ may also be
treated as the semigroup defined by three generators $ f_-, f_+, f_* $
and three relations $ f_+ f_- = 1 $, $ f_* f_- = 1 $, $ f_* f_+ = f_*
f_* $. 

The one-point motion $ (X_{0,t}(x))_{t>0} $ of our \flow{G} is
(distributed like) the sticky Brownian motion (starting at $ x $). A
particle spends a positive time at the origin, but never sits there
during a time interval. Two particles keep a constant distance until
one of them reaches the origin. Generally, the distance is
non-monotone. But ultimately the two particles coalesce.

The noise of stickiness was introduced and investigated by J.~Warren
\cite{Wa}. See also \cite[Sect.~4]{Ts03}, and Sections \ref{4e},
\ref{5d} of this survey.

\mysubsection{Warren's noise of splitting}
\label{4d}

The noise of splitting consists of \sif s generated by random
variables $ a_{s,t} $, $ b_{s,t} $ and $ \tau_{s,t} $ according to the
parameters $ a,b,\tau $ of an element $ f_{a,b}^\tau $ of the
semigroup $ G $ ($ = G^{\text{\ref{4b}}} $); $ b \ge 0 $, $ a+b \ge 0 $, and
$ \tau = \pm 1 $. We may drop $ b_{s,t} $ but not $ \tau_{s,t} $. The
binary operation in $ G $ is such that (assuming $ r<s<t $) $
\tau_{r,t} $ is either $ \tau_{r,s} $ or $ \tau_{s,t} $ depending on
whether the minimum of the Brownian motion $ B_u = a_{0,u} $ on $
[r,t] $ is reached on $ [r,s] $ or $ [s,t] $. It means that the random
sign $ \tau_{r,t} $ may be assigned to the minimizer $ s \in [r,t] $
of the Brownian motion on $ [r,t] $. ``This is a noise richer than
white noise: in addition to the increments of a Brownian motion $ B $
it carries a countable collection of independent Bernoulli random
variables which are attached to the local minima of $ B $'' \cite[the
last phrase]{Wa1}.

It may seem that these Bernoulli random
variables appear suddenly, having no precursors in the past (like
jumps of the Poisson process). However, this is a delusion.

\begin{definition}
A noise (or, more generally, a continuous product of probability spaces, see
\ref{2b1}) $ (\F_{s,t})_{s<t} $ is
\emph{predictable,}%
\index{predictable (noise; cont.\ prod.\ of prob.\ spaces)}
if the filtration $ (\F_{-\infty,t})_{t\in\R} $
admits of no discontinuous martingales.
\end{definition}

Equivalently: for every stopping time $ T $ (w.r.t.\ the filtration $
(\F_{-\infty,t})_{t\in\R} $) there exist stopping times $ T_n $ such
that $ T_n < T $ and $ T_n \to T $ a.s.

The white noise is predictable; the Poisson noise is not.

\emph{The noise of splitting is predictable.}

What is wrong in saying `each one of these Bernoulli random variables
appears suddenly at the corresponding instant'? The very beginning
`each one of these' is misleading. We cannot number them in real
time. Rather, we can consider (say) $ \tau_{0,1} $, the Bernoulli
random variable attached to the minimizer of $ B $ on $ [0,1] $. Its
conditional expectation, given $ \F_{-\infty,t} $ ($ 0<t<1 $), does
not jump, since we do not know (at $ t $) whether the minimum was
already reached or not; the corresponding probability is continuous in
$ t $.

Is there anything special in local minima of the Brownian path? Any
other random dense countable set could be used equally well, if it
satisfies two conditions, locality and stationarity, formalized
below. However, what should we mean by a `random dense countable set'?
The set of all dense countable subsets of $ \R $ does not carry a
natural structure of a standard measurable space. (Could you imagine a
function of the set of all Brownian local minimizers that gives a
non-degenerate random variable?) They form a \emph{singular space} in
the sense of Kechris \cite[\S2]{Ke99}: a `bad' quotient space of a
`good' space by a `good' equivalence relation. Several possible
interpretations of `good' and `bad' are discussed in \cite{Ke99}, but
we restrict ourselves to few noise-related examples.

\emph{Please consult Sect.~\textup{\ref{2aa}} for some general notions and
notations used below.}

\begin{sloppypar}
The space $ \R^\infty $ of all infinite sequences $ (t_1,t_2,\dots) $
of real numbers is naturally a standard measurable space. The
group $ S_\infty $ of all bijective maps $ \{1,2,\dots\} \to \{1,2,\dots\} $
acts on $ \R^\infty $ by permutations: $ (t_1,t_2,\dots) \mapsto
(t_{n_1},t_{n_2},\dots) $. The Borel subset $ \R^\infty_{\ne} \subset
\R^\infty $ of all sequences of \emph{pairwise different} numbers $
t_1,t_2,\dots $ is $ S_\infty $-invariant, and the set of orbits $
\R^\infty_{\ne} / S_\infty $ may be identified with the set of all
countable subsets of $ \R $. The same for $ (a,b)^\infty_{\ne} /
S_\infty $ and countable subsets of a given interval $ (a,b) \subset
\R $.
\end{sloppypar}

The group $ L_0(\Om \to S_\infty) $ of \emph{random permutations} acts on
the space $ L_0(\Om \to \R^\infty) $ of \emph{random sequences}. The
subset $ L_0(\Om \to \R^\infty_{\ne}) $ is invariant under random
permutations. We treat the quotient space $ L_0(\Om \to \R^\infty_{\ne}) /
L_0(\Om \to S_\infty) $ as a well-defined substitute of the ill-defined $
L_0 ( \Om \to \R^\infty_{\ne} / S_\infty) $. A random countable set is
treated as a random sequence up to a random permutation.

Local minimizers of a Brownian path are such a random set; that is,
they admit a measurable enumeration. Here is a simple construction for
$ (0,1) $. First, $ t_1(\om) $ is the minimizer on the whole $ (0,1) $
(unique almost sure). Second, if $ t_1(\om) \in (0,1/2) $ then $
t_2(\om) $ is the minimizer on $ (1/2,1) $, otherwise --- on $ (0,1/2)
$. Third, $ t_3(\om) $ is the minimizer on the first of the four
intervals $ (0,1/4) $, $ (1/4,1/2) $, $ (1/2,3/4) $ and $ (3/4,1) $
that contains neither $ t_1(\om) $ nor $ t_2(\om) $. And so on.

Random sets $ M_{s,t} $ of Brownian minimizers on intervals $ (s,t)
\subset \R $ satisfy $ M_{r,s} \cup M_{s,t} = M_{r,t} $ for $
r<s<t $ (almost sure, $ s $ is not a local minimizer), and $ M_{s,t} $
depends only on the increments of $ B $ on $ (s,t) $,
which may be said in terms of $ (\F_{s,t})_{s<t} $. Namely, $
M_{s,t} $ is defined via \measurable{\F_{s,t}} random sequences
modulo \measurable{\F_{s,t}} random permutations,
\[
M_{s,t} \in L_0 ( \F_{s,t} \to (s,t)^\infty_{\ne} ) / L_0 ( \F_{s,t}
\to S_\infty ) \, .
\]

\begin{definition}
A \emph{stationary local random dense countable set}%
\index{random dense countable set|see{set\dots}}%
\index{set, random dense}
(over a given noise) is a family $ (N_{s,t})_{s<t} $ of random sets
\[
N_{s,t} \in L_0 ( \F_{s,t} \to (s,t)^\infty_{\ne} ) / L_0 ( \F_{s,t}
\to S_\infty )
\]
satisfying
\begin{gather*}
N_{r,s} \cup N_{s,t} = N_{r,t} \quad \text{a.s.,} \\
N_{s,t} \circ T_h = N_{s+h,t+h} \quad \text{a.s.}
\end{gather*}
whenever $ r<s<t $, $ h \in \R $. (Here $ T_h $ are time shifts, see
\ref{2c1}.)
\end{definition}

Brownian minimizers are an example of a stationarity local random
dense countable set over the white noise. Brownian maximizers are
another example. Their union is the third example.

\begin{question}\label{4d4}\index{question}
Do these three examples exhaust \emph{all} stationarity local random
dense countable sets over the white noise?
\end{question}

New examples could lead to new noises.

\mysubsection{Warren's noise of stickiness, made by a Poisson snake}
\label{4e}

The noise of stickiness consists of \sif s generated by random
variables $ a_{s,t} $, $ b_{s,t} $ and $ c_{s,t} $ according to the
parameters $ a,b,c $ of an element $ f_{a,b,c} $ of the semigroup $ G
$ ($ = G^{\text{\ref{4c}}} $); $ b \ge 0 $, $ 0 \le c \le a+b $. We may drop
$ b_{s,t} $ but not $ c_{s,t} $. 

Consider the (random) set $ C_t = \{ c_{s,t} : s \in (-\infty,t) \}
\setminus \{ 0 \} $; its points will be called `spots'. For a small $
\De t $ usually (with probability $ 1 - O(\sqrt{\De t}) $) $
c_{t,t+\De t} = 0 $ (since $ a+b-\la\eta < 0 $, recall Sect.~\ref{4c}),
therefore $ C_{t+\De t} = ( C_t + a_{t,t+\De t} ) \cap ( a_{t,t+\De t}
+ b_{t,t+\De t}, \infty ) $. We see that the spots move up and down,
driven by Brownian increments. The boundary annihilates the spots
that hit it. However, sometimes the boundary creates new spots. It
happens (with probability $ \sim \const \cdot \sqrt{\De t} $) when $
c_{t,t+\De t} > 0 $.

An observer that moves according to the Brownian increments sees a set
$ C_t - a_{0,t} $ of fixed spots on the changing ray $
(-a_{0,t},\infty) $. The spotted ray may be called a \emph{Poisson
snake.}\index{snake}\index{Poisson snake}
The movement of its endpoint ($ -a_{0,t} $) is Brownian. When
the snake shortens, some spots disappear on the moving boundary. When
the snake lengthens, new spots appear on the moving boundary. It
happens with
a rate infinite in time but finite in space. Infinitely many spots
appear (and disappear) during any time interval (because of locally
infinite variation of a Brownian path); only a finite number of them
survive till the end of the interval. In fact, at every instant the
spots are (distributed like) a Poisson point process of rate $ 1/\la $
on $ (-a_{0,t},\infty) $.

Being discrete in space, the spots may seem to appear suddenly in time
(like jumps of the usual Poisson process). However. this is a delusion
(similarly to Sect.~\ref{4d}).

\emph{The noise of stickiness is predictable.}

A spot can appear at an instant $ s $ only if $ s $ is `visible from
the right' in the sense that $ a_{s,t} > 0 $ for all $ t $ close
enough to $ s $ (that is, $ \exists \eps>0 \; \forall t \in (s,s+\eps)
\;\; a_{s,t}>0 \, $).
\[
\begin{gathered}\includegraphics[scale=0.8]{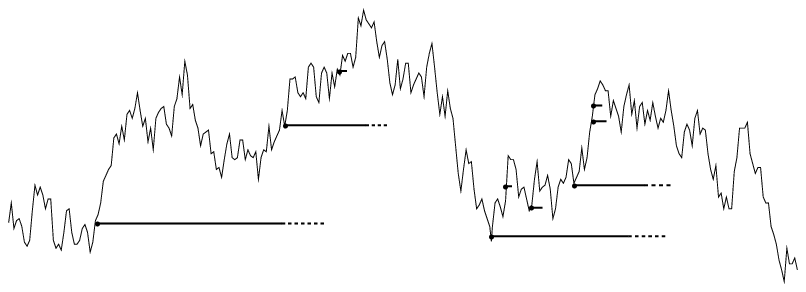}\end{gathered}
\quad\parbox{5cm}{\small A random dense countable subset of the
 continuum of points visible from the right. (Few chords are shown,
 others are too short.)}
\]
Points visible from the right are (a.s.) a dense Borel set of
cardinality continuum but Lebesgue measure zero. Knowing the past
(according to $ \F_{-\infty,s} $) but not the future ($ \F_{s,\infty}
$) we cannot guess that $ s $ is (or rather, will appear to be)
visible from the right. (Compare it with Sect.~\ref{4d}: knowing the
past we cannot guess that $ s $ is a local minimizer.)

In contrast, knowing the future ($ \F_{s,\infty} $) but not the past
($ \F_{-\infty,s} $) we know, whether $ s $ is visible from the right
or not. (This asymmetry reminds me that we often know the date of
death of a great man but not the date of birth\dots)

\emph{The time-reversed noise of stickiness is not predictable.}

In other words, the filtration $ (\F_{-\infty,t})_{t\in\R} $ admits
of continuous martingales only, but the filtration $
(\F_{-t,\infty})_{t\in\R} $ admits of some discontinuous martingales.

All the birth instants (when new spots appear) are (a.s.) a dense
countable subset of the set of points visible from the
right. Conditionally, given the Brownian path $ (B_t)_{t\in\R} =
(a_{0,t})_{t\in\R} $, birth instants are a Poisson random subset of $
\R $ whose intensity measure is a singular \sifinite\ measure $
(\D B)^+ $ concentrated on points visible from the right. Such a measure
$ (\D f)^+ $ may be defined for \emph{every} continuous function $ f $
(not just a Brownian path); note that $ f $ need not be of locally
finite variation. Namely, $ (\D f)^+ $ is the supremum (over $
t \in \R $) of images of Lebesgue measure on $ (-\infty,f(t)) $ under
the maps $ x \mapsto \max \{ s \in (-\infty,t) : f(s)=x \} $, provided
that $ \inf \{ f(s) : s \in (-\infty,t) \} = -\infty $ (which holds
a.s.\ for Brownian paths); otherwise $ (-\infty,f(t)) $ should be
replaced with $ \( \inf \{ f(s) : s \in (-\infty,t) \}, \, f(t) \) $.
The measure $ (\D f)^+ $ is always positive and \sifinite, but need
not be locally finite. That is, $ \R $ can be decomposed into a
sequence of Borel subsets of finite measure $ (\D f)^+ $. However, it
does not mean that all (or even, some) intervals are of finite measure
$ (\D f)^+ $.

The \sifinite\ positive measure $ (\D B)^+ $ is infinite on every
interval (because of locally infinite variation of the Brownian
path). Such measures are a singular space (recall Sect.~\ref{4d}); a
\emph{random} element of such a space should be treated with great
care. Interestingly, the singular space of Sect.~\ref{4d} is naturally
embedded into the singular space considered here. Indeed, every dense
countable set may be identified with its counting measure (consisting
of atoms of mass $ 1 $); the measure is \sifinite, but infinite
on every interval. In contrast, the measure $ (\D B)^+ $ is non-atomic.

Similarly to Sect.~\ref{4d} it should be possible to define a
stationary local random \sifinite\ positive measure, infinite
on every interval, over (say) the white noise. One example is $ (\D B)^+
$. Replacing $ B_t $ with $ -B_t $, or $ B_{-t} $, or $ -B_{-t} $ we
get three more examples. Similarly to Question \ref{4d4} we may ask:
do these four examples and their linear combinations exhaust
\emph{all} possible cases? New examples could lead to new noises.

See also \cite{Wa} and \cite[Sect.~4]{Ts03}.

\mysection{From convolution semigroups to continuous products of
 probability spaces}
\label{sec:2}
\mysubsection{Preliminaries: probability spaces, morphisms etc.}
\label{2aa}

\emph{Throughout, either by assumption or by construction, all
probability spaces are standard. All claims and constructions are
invariant under $ \modO $ isomorphisms.}

Recall that a
\emph{standard probability space}\index{standard!probability space}
(known also as a
Lebesgue-Rokhlin space) is a probability space isomorphic $ (\modO) $
to an interval with the Lebesgue measure, a finite or countable
collection of atoms, or a combination of both (see
\cite[17.41]{Ke}). Nonseparable $ L_2 $ spaces of random variables
are thus disallowed!

A \sif\ $ \F $ is sometimes shown in the notation $ (\Om,\F,P) $,
sometimes suppressed in the shorter notation $ (\Om,P) $.

Every function on any probability space is treated $ \modO $. That is,
I write $ f : \Om \to \R $ for convenience, but I mean that $ f $ is
an equivalence class. The same for maps $ \Om_1 \to \Om_2 $ etc. A
\emph{morphism}
$ \Om_1 \to \Om_2 $ is a measure preserving
(not just non-singular) measurable map ($ P_1,P_2 $ are suppressed in
the notation). An\index{morphism!of probability spaces}
\emph{isomorphism}\index{isomorphism!of probability spaces}
(known also as
`$ \modO $ isomorphism') is an invertible morphism whose inverse is
also a morphism. An
\emph{automorphism}\index{automorphism}
is an isomorphism from $ \Om $ to
itself. Every sub-\sif\ is assumed to contain all negligible
sets. Every morphism $ \al : \Om \to \Om' $ generates a sub-\sif\ $
\Ec \subset \F $, and every sub-\sif\ $ \Ec \subset \F $ is generated
by a morphism $ \al : \Om \to \Om' $, determined by $ \Ec $ uniquely
up to isomorphism ($ \Om' \leftrightarrow \Om'' $, making the diagram
commutative\dots); it is the
\emph{quotient space}\index{quotient space}
$ (\Om',P') = (\Om,P)/\Ec $.

A \emph{standard measurable space}\index{standard!measurable space}
(or `standard Borel space') is a set $ E $ equipped with a \sif\ $ \B
$ such that the measurable space $ (E,\B) $ is isomorphic either to $
\R $ (with the Borel \sif) or to its finite or countable subset. See
\cite[12.B and 15.B]{Ke}.

Equivalence classes of all measurable functions $ \Om \to \R $ are a
topological linear space $ L_0(\Om,\F,P) $%
\index{zzl@$ L_0 $, space}
(denoted also $ L_0(\Om) $, $ L_0(\F) $, $ L_0(P) $ etc.); its metrizable
topology corresponds to convergence in probability. Any Borel function
$ \phi : \R \to \R $ leads to a (nonlinear) map $ L_0(\Om) \to
L_0(\Om) $, $ X \mapsto \phi \circ X $, discontinuous (in general),
but Borel measurable. (Hint: if $ \phi_n(x) \to \phi(x) $ for all $ x
\in \R $ then $ \phi_n \circ X \to \phi \circ X $ for all $ X \in
L_0(\Om) $.) Given a standard measurable space $ E $, the set $
L_0(\Om \to E) $ of equivalence classes of all measurable maps $ \Om
\to E $ carries a natural Borel \sif\ and is a standard measurable
space (neither linear nor topological, in general).

A stochastic flow (and any random process) is generally treated as a
family of equivalence classes (rather than functions). The distinction
is essential when dealing with uncountable families of random
variables.\index{random variables!uncountable families}
The phrase (say)
\[
f_t = g_t \quad \text{a.s. for all $ t $}
\]
is interpreted as
\[
\inf_t \, \Prob \{ \om : f_t(\om) = g_t(\om) \} = 1
\]
rather than $ \Pr{ \cap_t \{ \om : f_t(\om) = g_t(\om) \} } = 1 $. In
spite of that, when dealing with (say) a Brownian motion $ (B_t)_t $
and writing (say) $ \max_{t\in[0,1]} B_t $, we rely on path
continuity. Here the Brownian motion is treated as a random continuous
function rather than a family of random variables.

\emph{I stop writing `standard' \textup{(}probability space\textup{),}
`$ \modO $' and `measure preserving', but I still mean it!}

\mysubsection{From convolution systems to flow systems}
\label{2a}

Recall the convolution semigroups in $ \R $, mentioned in \ref{4aa}.
The convolution relation $ \mu_s * \mu_t = \mu_{s+t} $ means that the
map $ \R^2 \to \R $, $ (x,y) \mapsto x+y $ sends the product measure $
\mu_s \times \mu_t $ into $ \mu_{s+t} $. More generally, each $ \mu_t
$ may sit on its own space $ G_t $, in which case some measure
preserving maps $ G_s \times G_t \to G_{s+t} $ should be given
(instead of the binary operation $ G \times G \to G $). Another generalization
is, abandoning stationarity (that is, homogeneity in time).

\begin{definition}\label{2a1}
A \emph{convolution system}\index{convolution system}
consists of probability spaces $ (G_{s,t}, \linebreak[0]
\mu_{s,t}) $ given for all $ s,t \in \R $, $ s<t $, and
morphisms $ G_{r,s} \times G_{s,t} \to G_{r,t} $ given for all $
r,s,t \in \R $, $ r<s<t $, satisfying the associativity condition:
\[
(xy)z = x(yz) \quad \text{for almost all $ x \in G_{r,s} $,  $ y \in
G_{s,t} $,  $ z \in G_{t,u} $}
\]
whenever $ r,s,t,u \in \R $, $ r<s<t<u $.
\end{definition}

Here and henceforth the given map $ G_{r,s} \times G_{s,t} \to G_{r,t}
$ is
denoted simply $ (x,y) \mapsto xy $. Any linearly ordered set (not
just $ \R $) may be used as the time set.

Every convolution semigroup $ (\mu_t) $ in $ \R $ leads to a
convolution system; namely, $ (G_{s,t}, \mu_{s,t}) = ( \R, \mu_{t-s} ) $, and
the map $ G_{r,s} \times G_{s,t} \to G_{r,t} $ is $ (x,y) \mapsto x+y
$. Another example: $ (G_{s,t}, \mu_{s,t}) = ( \Z_m, \mu ) $ for all $
s,t $; here $ m \in \{2,3,\dots\} $ is a parameter, and $ \mu $ is the
uniform distribution on the finite cyclic group $ \Z_m = \Z / m\Z $;
the map $ G_{r,s} \times G_{s,t} \to G_{r,t} $ is $ (x,y) \mapsto x+y
\pmod m $, and the time set is still $ \R $. The latter example is much worse
than the former; indeed, the former is separable (see Definition \ref{2a5}
below), and the latter is not.

Here is a generalization of the classical transition from convolution
semigroups to independent increments.

\begin{definition}\label{2a3}
Let $ (G_{s,t}, \mu_{s,t})_{s<t;s,t\in T} $ be a convolution system
over a linearly ordered set $ T $. A 
\emph{flow system}\index{flow system}
(corresponding to the given convolution system) consists of a
probability space $ (\Om,\F,P) $ and morphisms
$ X_{s,t} : \Om \to G_{s,t} $ (for $ s<t $; $ s,t\in T $) such that $
\F $ is generated by all $ X_{s,t} $ (`non-redundancy'), and
\begin{gather}
X_{t_1,t_2}, X_{t_2,t_3}, \dots, X_{t_{n-1},t_n} \text{ are
 independent for } t_1 < t_2 < \dots < t_n \, ; \tag{a} \\
X_{r,t} = X_{r,s} X_{s,t} \quad \text{(a.s.) for } r<s<t \,
 . \tag{b}
\end{gather}
\end{definition}

\begin{sloppypar}
The non-redundancy can be enforced by taking the quotient space $
(\Om,P) / \F_{-\infty,\infty} $, where $ \F_{-\infty,\infty} $ is the
\sif\ generated by all $ X_{s,t} $.
\end{sloppypar}

\begin{proposition}\label{2a2}
For every convolution system over a finite or countable
$ T $, the corresponding flow system exists and is unique up to
isomorphism.
\end{proposition}

By an isomorphism between flow
systems\index{isomorphism!of flow systems}
$ (X_{s,t})_{s<t} $, $ X_{s,t}
: \Om \to G_{s,t} $ and $ (X'_{s,t})_{s<t} $, $ X'_{s,t} : \Om' \to
G_{s,t} $ we mean an isomorphism $ \al : \Om \to \Om' $ such that $
X_{s,t} = X'_{s,t} \circ \al $ for $ s<t $.

\beginproof
Existence: if $ T $ is finite, $ T = \{ t_1, \dots, t_n \} $, $ t_1 <
\dots < t_n
$, we just take $ \Om = G_{t_1,t_2} \times \dots \times
G_{t_{n-1},t_n} $ with the product measure. If $ T $ is countable, we
have a consistent family of finite-dimensional distributions on the
product $ \prod_{s<t;s,t\in T} G_{s,t} $ of countably many probability
spaces.

Uniqueness follows from the fact that the joint distribution of all
$ X_{s,t} $ is uniquely determined by the measures $ \mu_{s,t} $ according to
\ref{2a3}(a,b).
\proofend

Given a convolution system $ (G_{s,t}, \mu_{s,t})_{s<t;s,t\in T} $ and
a subset $ T_0 \subset T $, the restriction $ (G_{s,t},
\mu_{s,t})_{s<t;s,t\in T_0} $ is also a convolution system. If $ T $
is countable and $ T_0 \subset T $, we get two flow systems, $
(X_{s,t})_{s<t;s,t\in T} $ on $ (\Om,P) $ and $
(X^0_{s,t})_{s<t;s,t\in T_0} $ on $ (\Om_0,P_0) $ related via a
morphism $ \al : \Om \to \Om_0 $ such that $ X_{s,t} = X^0_{s,t} \circ
\al $ (a.s.) for $ s,t\in T_0 $, $ s<t $. It may happen that $ \al $
is an isomorphism, in which case we say that $ T_0 $ is \emph{total}
in $ T $ (with respect to the given convolution system).

\begin{definition}\label{2a5}
A convolution system $ (G_{s,t}, \mu_{s,t})_{s<t;s,t\in\R} $ is
\emph{separable,}\index{separable!convolution system}
if there exists a countable set $ T_0 \subset \R $
such that for every countable $ T \subset \R $ satisfying $ T_0
\subset T $, the subset $ T_0 $ is total in $ T $ with respect to the
restriction $ (G_{s,t}, \mu_{s,t})_{s<t;s,t\in T} $ of the given
convolution system.
\end{definition}

When checking separability, one may restrict himself to the case when
the difference $ T \setminus T_0 $ is a single point.

Here is a counterpart of Proposition \ref{2a2} for the uncountable
time set $ \R $.

\begin{proposition}\label{2a6}
\begin{sloppypar}
The following two conditions on a convolution system $ (G_{s,t},
\mu_{s,t})_{s<t;s,t\in\R} $ are equivalent:
\end{sloppypar}

(a) there exists a flow system corresponding to the given convolution
system;

(b) the given convolution system is separable.
\end{proposition}

\beginproof
(a) \imp (b):
the \sif\ $ \F $, being generated by the uncountable set $ \{ X_{s,t} : s<t \}
$ of random variables, is necessarily generated by some countable subset
(since the probability space is standard).

(b) \imp (a):
we take the flow system for $ T_0 $; separability implies that each $
X_{s,t} $ (for $ s,t\in\R $) is equal (a.s.) to a function of $
(X_{s,t})_{s<t;s,t\in T_0} $.
\proofend

\mysubsection{From flow systems to continuous products, and back}
\label{2b}

\begin{definition}\label{2b1}
A \emph{continuous product of probability spaces}%
\index{continuous product!of probability spaces}
consists of a
probability space $ (\Om,\F,P) $ and sub-\sif s $ \F_{s,t}
\subset \F $ (given for all $ s,t \in \R $, $ s<t $) such that $ \F $
is generated by the union of all $ \F_{s,t} $ (`non-redundancy'), and
\begin{equation}\label{2b2}
\F_{r,s} \otimes \F_{s,t} = \F_{r,t} \quad \text{whenever } r<s<t \,
.
\end{equation}
\end{definition}

\index{zzo@$ \otimes $ (for sub-\sif s)}%
The latter means that $ \F_{r,s} $ and $ \F_{s,t} $ are independent
and generate $ \F_{r,t} $. See also Def.~\ref{2b6} below. The
non-redundancy can be enforced by taking the quotient space $ (\Om,P)
/ \F_{-\infty,\infty} $. Any linearly ordered set (not just $ \R $)
may be used as the time set. It is convenient to enlarge the time set
from $ \R $ to $ [-\infty,\infty] $ defining $ \F_{-\infty,t} $ as the
\sif\ generated by the union of all $ \F_{s,t} $ for $ s \in
(-\infty,t) $; the same for $ \F_{s,\infty} $ and $
\F_{-\infty,\infty} $.

\begin{proposition}\label{2b3}
Let $ (X_{s,t})_{s<t} $ be a flow system, and $ \F_{s,t} $ be defined
(for $ s<t $) as the sub-\sif\ generated by $ \{ X_{u,v} : s \le u < v
\le t \} $. Then sub-\sif s $ \F_{s,t} $ form a continuous product of
probability spaces.
\end{proposition}

\beginproof
$ \F_{r,s} $ and $ \F_{s,t} $ generate $ \F_{r,t} $ by \ref{2a3}(b)
and are independent by \eqref{2a3}(a) (and (b)).
\proofend

Having a continuous product of probability spaces $ (\F_{s,t})_{s<t} $
we may introduce quotient spaces 
\begin{equation}\label{2b4}
(\Om_{s,t},P_{s,t}) = (\Om,P) / \F_{s,t} \, .
\end{equation}
The relation \eqref{2b2} becomes
\begin{equation}\label{2b5}
(\Om_{r,s},P_{r,s}) \times (\Om_{s,t},P_{s,t}) = (\Om_{r,t},P_{r,t})
\, ;
\end{equation}
the equality is treated here via a canonical isomorphism. It is not
unusual; for example, the evident equality $ (A\times B)\times C =
A\times (B\times C) $ for Cartesian products of (abstract) sets is
also treated not literally but via a canonical bijection $ ((a,b),c)
\mapsto (a,(b,c)) $ between the two sets. The canonical isomorphisms
implicit in \eqref{2b5} satisfy associativity (stipulated by Definition
\ref{2a1}); indeed, for $ r<s<t<u $ we have $ (\Om_{r,s},P_{r,s})
\times (\Om_{s,t},P_{s,t}) \times (\Om_{t,u},P_{t,u}) =
(\Om_{r,u},P_{r,u}) $. Thus, $ (\Om_{s,t},P_{s,t}) $ form a
convolution system (as defined by \ref{2a1}) satisfying an
additional condition: the morphisms $ G_{r,s} \times G_{s,t} \to
G_{r,t} $ become isomorphisms. This is another approach to continuous
products of probability spaces.

\begin{definition}\label{2b6}
A \emph{continuous product of probability spaces}%
\index{continuous product!of probability spaces}
consists of
probability spaces $ (\Om_{s,t}, P_{s,t}) $ (given for all $ s,t \in
\R $, $ s<t $), and isomorphisms $ \Om_{r,s} \times \Om_{s,t} \to
\Om_{r,t} $ given for all $ r,s,t \in \R $, $ r<s<t $, satisfying the
associativity condition:
\[
(\om_1 \om_2) \om_3 = \om_1 (\om_2 \om_3) \quad \text{for almost all $
\om_1 \in \Om_{r,s} $, $ \om_2 \in \Om_{s,t} $, $ \om_3 \in \Om_{t,u}
$} 
\]
whenever $ r,s,t,u \in \R $, $ r<s<t<u $.
\end{definition}

(As before, the given map $ \Om_{r,s} \times \Om_{s,t} \to \Om_{r,t} $
is denoted simply $ (\om_1,\om_2) \mapsto \om_1 \om_2 $.)
Having $ (\F_{s,t})_{s<t} $ as in Definition \ref{2b1} we get the
corresponding $ (\Om_{s,t},P_{s,t}) $ as in Definition \ref{2b6} by
means of \eqref{2b4}. And conversely, each $ (\Om_{s,t},P_{s,t}) $ as
in Definition \ref{2b6} leads to the corresponding $ (\F_{s,t})_{s<t}
$ of Definition \ref{2b1}. Namely, we may take $ (\Om,P) =
\prod_{k\in\Z} (\Om_{k,k+1}, P_{k,k+1}) $, define $ X_{k,k+1} : \Om
\to \Om_{k,k+1} $ as coordinate projections, use the relation $
\Om_{k,k+1} = \Om_{k,k+\theta} \times \Om_{k+\theta,k+1} $ for
constructing $ X_{k,k+\theta} : \Om \to \Om_{k,k+\theta} $ and so
forth. Alternatively, we may treat $ (\Om_{s,t},P_{s,t})_{s<t} $ as a
(special, necessarily separable) convolution system and use relations
discussed below.

A separable convolution system leads to a flow system by \ref{2a6}; a
flow system leads to a continuous product of probability spaces by
\ref{2b3}; and a continuous product of probability spaces is a special
case of a separable convolution system.
\begin{equation}
\begin{gathered}
\xymatrix{
 \txt{separable convolution systems}
  \ar[r] &
 \txt{flow systems}
  \ar[dl]
 \\
\txt{continuous products\\of probability spaces}
  \ar@{^{(}->}[u]
}
\end{gathered}
\end{equation}
For example, a weakly continuous convolution semigroup $
(\mu_t)_{t\in(0,\infty)} $ in $ \R $ is a separable convolution
system. (Any dense countable subset of $ \R $ may be used as $ T_0 $
in Def.~\ref{2a5}.) The corresponding flow system consists of the
increments of
the L\'evy process corresponding to $ (\mu_t)_t $. It leads to a
continuous product of probability spaces $ (\Om_{s,t}, P_{s,t})
$. Namely, $ (\Om_{0,t}, P_{0,t}) $ may be treated as the space of
sample paths of the L\'evy process on $ [0,t] $; $ (\Om_{s,t},
P_{s,t}) $ is a copy of $ (\Om_{0,t-s}, P_{0,t-s}) $; and the
composition $ \Om_{r,s} \times \Om_{s,t} \ni (\om_1,\om_2) \mapsto
\om_3 \in \Om_{r,t} $ is
\[
\om_3 (u) = \begin{cases}
 \om_1(u) &\text{for $ u \in [r,s] $},\\
 \om_1(s)+\om_2(u-s) &\text{for $ u \in [s,t] $}.
\end{cases}
\]
Note that $ (\Om_{s,t}, P_{s,t}) $ is much larger than $ (G_{s,t},
\mu_{s,t}) $. We may treat $ (\Om_{s,t}, P_{s,t})_{s<t} $ as another
convolution system; the two convolution systems, $
(\R,\mu_{t-s})_{s<t} $ and $ (\Om_{s,t}, P_{s,t}) $, lead to the same
(up to isomorphism) continuous product of probability spaces. It holds in
general:
\[
\xymatrix{
 \txt{separable convolution system}
  \ar@{|->}[r] &
 \txt{flow system}
  \ar@{|->}[d]
\\
 \txt{\emph{another} separable\\convolution system} &
 \txt{continuous product\\of probability spaces}
  \ar@{|->}[l]
}
\]
but
\[
\xymatrix{
 \txt{continuous product\\of probability spaces}
  \ar@{|->}[r] \ar@{=}[d] &
 \txt{separable convolution system}
  \ar@{|->}[d]
\\
 \txt{\emph{the same} continuous product\\of probability spaces} &
 \txt{flow system}
  \ar@{|->}[l]
}
\]

\mysubsection{Stationary case: noise}
\label{2c}

Returning to stationarity (that is, homogeneity in time), abandoned in
Sections \ref{2a}, \ref{2b}, we add time shifts to Def.~\ref{2b1},
after a general discussion of one-parameter groups of automorphisms.

\smallskip
\textsc{digression: measurable action}
\smallskip

In the spirit of our conventions (Sect.~\ref{2aa}), an action of $ \R
$ on a probability space $ \Om = (\Om,\F,P) $ is treated as a
homomorphism of $ (\R,+) $ to the group of automorphisms $ \Aut(\Om)
$, each automorphism being an equivalence class rather than a map $
\Om \to \Om $. Thus, an action is not quite a map $ \R \times \Om
\to \Om $. The group $ \Aut(\Om) $ is topological (in fact, Polish)
\cite[17.46]{Ke}, and a homomorphism $ \R \to \Aut(\Om) $ is Borel
measurable if and only if it is continuous (which is a special case of
a well-known general theorem \cite[9.10]{Ke}). Such a homomorphism
will be called a
\emph{measurable action}\index{measurable!action}
of $ \R $ on $ \Om $. Every such action $ T : \R \to \Aut(\Om) $
corresponds to some (non-unique) measurable map $ \R \times \Om \to
\Om $. Moreover, the map can be chosen to satisfy \emph{everywhere}
the relation $ T_r(T_s(\om)) = T_{r+s}(\om) $ (Mackey, Varadarajan
and Ramsy), but this fact will not be used. More detailed discussion can be
found in \cite[Introduction]{GTW}.

\smallskip
\textsc{end of digression}
\smallskip

\begin{definition}\label{2c1}
A \emph{noise,}\index{noise}
or a \emph{homogeneous%
\index{homogeneous continuous product!of probability spaces}
continuous product of
probability spaces,} consists of a probability space $
(\Om,\F,P) $, sub-\sif s $ \F_{s,t} \subset \F $ given for all $
s,t \in \R $, $ s<t $, and a measurable action $ (T_h)_h $ of $ \R $
on $ \Om $, having the following properties:
\begin{gather}
\F_{r,s} \otimes \F_{s,t} = \F_{r,t} \quad \text{whenever } r<s<t \,
 , \tag{a} \\
T_h \text{ sends } \F_{s,t} \text{ to } \F_{s+h,t+h} \quad
 \text{whenever } s<t \text{ and } h \in \R \, , \tag{b} \\
\F \text{ is generated by the union of all } \F_{s,t} \, . \tag{c}
\end{gather}
\end{definition}

The time set $ \R $ may be enlarged to $ [-\infty,\infty] $, as noted
after Def.~\ref{2b1}. Of course, the index $ h $ of $ T_h $ runs over
$ \R $ only, and $ -\infty+h = -\infty $, $ \infty+h = \infty $.

Weakly continuous convolution semigroups in $ \R $ (unlike convolution
systems in general) lead to noises.

Measurability of the action does not follow from other conditions; here is a
counterexample.

\begin{example}\label{counterexample}
We choose a weakly continuous semigroup $ (\mu_t)_{t\in(0,\infty)} $ in $ \R^2
$ and assume that it is \emph{isotropic} in the sense that each $ \mu_t $ is
invariant under rotations $ (x,y) \mapsto (x\cos\phi-y\sin\phi, \linebreak[0]
x\sin\phi+y\cos\phi) $ of $ \R^2 $. (For instance, $ (\mu_t)_t $ may
correspond to the standard Brownian motion in $ \R^2 $. One may also use
Poissonian jumps in random directions\dots) The corresponding noise is also
isotropic in the sense that we have a measurable action of the rotation group
$ \T $ on $ \Om $, commuting with the time shifts $ T_h $ and preserving the
sub-\sif s $ \F_{s,t} $. We take some \emph{non-measurable} additive function
$ \phi : \R \to \R $ (that is, $ \phi(s+t) = \phi(s) + \phi(t) $ for all $
s,t\in\R $) and define a new (`spoiled') action $ T' : \R \to \Aut(\Om) $;
namely, $ T'_h $ is the composition of $ T_h $ and the rotation by the angle $
\phi(h) $. The structure $ \( (\Om,\F,P), (\F_{s,t})_{s<t}, (T'_h)_h \) $ is
not a noise only because the action $ (T'_h)_h $ is not measurable.
\end{example}

\begin{proposition}\label{2c2}
Every noise satisfies the
`upward continuity'\index{upward continuity!for probability spaces}
condition
\begin{equation}\label{2c2a}
\F_{s,t} \text{ is generated by } \bigcup_{\eps>0} \F_{s+\eps,t-\eps}
\quad \text{for all $ s,t \in \R $, $ s < t $} \, .
\end{equation}
\end{proposition}

Using the enlarged time set $ [-\infty,\infty] $ we interprete $
-\infty+\eps $ as $ -1/\eps $ and $ \infty-\eps $ as $ 1/\eps $.

\beginproof
In the Hilbert space $ H = L_2(\Om,\F,P) $ we consider projections $
Q_{s,t} : f \mapsto \cE{f}{\F_{s,t}} $. They commute, and $ Q_{s,t} =
Q_{-\infty,t} Q_{s,\infty} $. The \emph{monotone} operator-valued
function $ t \mapsto Q_{-\infty,t} $ must be continuous (in the strong
operator
topology) at every $ t \in \R $ except for an at most countable set,
since $ H $ is separable. By shift invariance, continuity at a single
$ t $ implies continuity at all $ t $. Thus, $ \| Q_{-\infty,t} f -
Q_{-\infty,t-\eps} f \| \to 0 $ when $ \eps \to 0 $ for every $ f \in
H $. Similarly, $ \| Q_{s,\infty} f - Q_{s+\eps,\infty} f \| \to 0
$. For $ f \in Q_{s,t} H $ we have $ \| f - Q_{s+\eps,t-\eps} f \| \le
\| f - Q_{-\infty,t-\eps} f \| + \| Q_{-\infty,t-\eps} ( f -
Q_{s+\eps,\infty} f ) \| \to 0 $.
\proofend

\begin{corollary}\label{2c3}
Every noise satisfies the
`downward continuity'\index{downward continuity!for probability spaces}
condition
\begin{equation}\label{2c3a}
\F_{s,t} = \bigcap_{\eps>0} \F_{s-\eps,t+\eps} \quad \text{for all }
s,t \in \R, \, s \le t \, ;
\end{equation}
here $ \F_{t,t} $ is the trivial \sif.
\end{corollary}

Using the enlarged time set $ [-\infty,\infty] $ we interprete $
-\infty-\eps $ as $ -\infty $ and $ \infty+\eps $ as $ \infty $.

\begin{sloppypar}
\beginproof
The \sif\ $ \F_{s-,t+} = \cap_{\eps>0} \F_{s-\eps,t+\eps} $ is
independent of $ \F_{-\infty,s-\eps} \vee \F_{t+\eps,\infty} $ for
every $ \eps $, therefore (using the proposition), also of $
\F_{-\infty,s} \vee \F_{t,\infty} $. We have $ \F_{s,t} \subset
\F_{s-,t+} $ and $ \F_{-\infty,s} \vee \F_{s,t} \vee \F_{t,\infty} =
\F_{-\infty,s} \vee \F_{s-,t+} \vee \F_{t,\infty} $; in combination
with the independence it implies $ \F_{s,t} = \F_{s-,t+} $.
\proofend
\end{sloppypar}

The two continuity conditions (`upward' and `downward') make sense
also for (non-homogeneous) continuous products of probability
spaces. The time set may be $ \R $, or $ [-\infty,\infty] $, or (say)
$ [0,1] $, etc. Of course, using $ [0,1] $ we interprete $ 0-\eps $ as
$ 0 $ and $ 1+\eps $ as $ 1 $. Still, the upward continuity implies
the downward continuity. (Indeed, the proof of \ref{2c3} does not use
the homogeneity. See also \ref{6d15}.) The converse does not hold. For
example, the \flow{\T} 
$ Y^{\text{\ref{1b}}} $ (the \flow{\T} $ Y $ of Sect.~\ref{1b}) leads to a
continuous product of probability spaces, continuous downwards but not
upwards. Namely, $ Y^{\text{\ref{1b}}}_{0,1} $ is \measurable{\F_{0,1}}
but independent of the \sif\ $ \F_{0+,1} $ generated by $
\cup_{\eps>0} \F_{\eps,1} $. On the other hand, the \sif\ $
\F_{-\infty,0+} = \F_{0-,0+} $ is trivial. See also \cite[3d6 and
3e3]{Ts03}, \cite[2.1]{Ts98}, \cite[1.5(2)]{TV}.

Similarly to Proposition \ref{2c2}, upward continuity holds for
(non-homogeneous) continuous products of probability spaces, provided
however that $ r,s $ do not belong to a finite or countable set of
discontinuity points.

The notions discussed so far are roughly as follows:
\begin{equation*}
\begin{gathered}
\xymatrix{
 \txt{convolution\\systems}
  \ar[r] &
 \txt{flow\\systems}
  \ar[r] &
 \txt{cont.\ products\\of prob.\ spaces}
  &
 \Big\} \txt{ \emph{non-stationary}}
\\
 \txt{convolution\\semigroups}
  \ar@{^{(}->}[u]
  \ar@{-->}[r] &
 \txt{ \large ? \ }
  \ar@{-->}[r] &
 \txt{noises}
  \ar@{^{(}->}[u]
  &
 \Big\} \txt{ \rlap{\emph{stationary}}\phantom{\emph{non-stationary}}}
}
\end{gathered}
\end{equation*}
It is natural to try a `stationary route' from convolution semigroups to
noises. We may define a \emph{stationary convolution system} as consisting of
probability spaces $ (G_t,\mu_t) $ for $ t \in (0,\infty) $ and morphisms $
G_s \times G_t \to G_{s+t} $ satisfying the evident associativity
condition. Assuming separability we get a flow system $ (X_{s,t})_{s<t} $
stationary in the sense that the distribution of $ X_{s,t} $ depends only on $
t-s $. It leads to time shifts $ T_h $ such that $ X_{s,t} \circ T_h =
X_{s+h,t+h} $ a.s. However, it does not mean that we get a noise, since the
measurability of the action $ T : \R \to \Aut(\Om) $ is not guaranteed (see
below).

In the other direction, the transition is simple; every noise corresponds to a
stationary convolution system. Just take $ (G_t,\mu_t) = (\Om,P) / \F_{0,t} $
and construct the morphism (in fact, isomorphism) $ G_s \times G_t \to G_{s+t}
$, combining the natural isomorphism $ ( (\Om,P) / \F_{0,s} ) \times ( (\Om,P)
/ \F_{s,s+t} ) \leftrightarrow (\Om,P) / ( \F_{0,s} \otimes \F_{s,s+t} ) =
(\Om,P) / \F_{0,s+t} $ and the isomorphism $ (\Om,P) / \F_{s,s+t}
\leftrightarrow (\Om,P) / \F_{0,t} $ given by the time shift $ T_s $.

This construction does not rely on the measurability of the action $ T
$. Therefore it may be applied to the counterexample of
\ref{counterexample} (it is in fact stationary). The `spoiled' isomorphism $
G_s \times G_t \leftrightarrow G_{s+t} $ rotates the second component
(belonging to $ G_t $) by the angle $ \phi(s) $. This (pathologic) stationary
convolution system does not lead to a noise. A more economical example: $
(G_t,\mu_t) = (\R^2,\mu_t) $, and the map $ G_s \times G_t \to G_{s+t} $ is
the vector addition `spoiled' by rotating the second vector (belonging to $
G_t $) by the angle $ \phi(s) $.

Adding stationarity to Def.~\ref{2b1} we get Def.~\ref{2c1}. It is
less simple, to add stationarity to Def.~\ref{2b6}. It should mean
isomorphisms $ \Om_s \times \Om_t \to \Om_{s+t} $. However, measurability of
the action $ T : \R \to \Aut(\Om) $ must be required; here $ \Om $ is the
`global' probability space constructed out of the given `local' probability
spaces $ \Om_t $. A simpler formulation is suggested by Liebscher's idea of
cyclic time (see Sect.~\ref{3d}). We define automorphisms $ \al_t : \Om_1 \to
\Om_1 $ for $ 0 < t < 1 $ by
\[
\al_t ( \om \om' ) = \om' \om \quad \text{for } \om \in \Om_t, \, \om' \in
\Om_{1-t} \, ;
\]
$ \om \om' $ stands for the image of $ (\om,\om') $ under the given map $
\Om_t \times \Om_{1-t} \to \Om_1 $; similarly, the map $ \Om_{1-t} \times
\Om_t \to \Om_1 $ is used in the right-hand side. These automorphisms are an
action of the circle $ \T = \R / \Z $ on $ \Om_1 $. Its measurability appears
to be equivalent to measurability of the action $ T $ of $ \R $ on $ \Om $.

It is not clear, whether the time dependence of $ (G_t,\mu_t) $ is really
essential, or not. The question may be formalized as follows.

By a \emph{Borel semigroup}\index{Borel semigroup}
we mean a semigroup $ G $ equipped with a \sif\ $
\B $ such that the measurable space $ (G,\B) $ is standard, and the
binary operation $ (x,y) \mapsto xy $ is a measurable map $ (G\times
G, \B \otimes \B) \to (G,\B) $.

Given a Borel semigroup $ G $, by a \flow{G} we mean a family $
(X_{s,t})_{s<t} $ of \valued{G} random variables $ X_{s,t} $ (given
for all $ s,t \in \R $, $ s<t $ on some probability space), satisfying
\ref{2a3}(a,b).

\begin{question}\label{2c4}\index{question}
(a) Does every noise correspond to a \flow{G} (for some Borel
semigroup $ G $)?

\begin{sloppypar}
(b) More specifically, is the following statement true?
For every noise $ \( (\Om,P), (\F_{s,t})_{s<t}, (T_h)_h \) $ there
exists a Borel semigroup $ G $ and a \flow{G} $ (X_{s,t})_{s<t} $
such that for all $ s,t,h \in \R $, $ s<t $,
\begin{gather*}
\F_{s,t} \text{ is the \sif\ generated by $ X_{s,t} $} \, , \\
X_{s,t} \circ T_h = X_{s+h,t+h} \quad \text{a.s,}
\end{gather*}
and in addition, there exists a Borel map $ L : G \to (0,\infty) $
(`graduation') such that
\begin{gather*}
L (X_{s,t}) = t-s \quad \text{a.s. whenever } s<t \, , \\
L (g_1 g_2) = L(g_1) + L(g_2) \quad \text{for all } g_1,g_2 \in G \,
 .
\end{gather*}
\end{sloppypar}
\end{question}

A straightforward idea is, to use the disjoint union $ G =
\uplus_{t\in(0,\infty)} \Om_t = \{ (t,\om) : 0<t<\infty, \, \om \in \Om_t \} $,
put $ L(t,\om) = t $ and combine the given maps $ \Om_s \times \Om_t \to
\Om_{s+t} $ into a binary operation $ G \times G \to G $. However, the matter
is more subtle than it may seem. The binary operation must be defined and
associative everywhere (rather than almost everywhere). Does it constraint the
noise somehow?

The convolution $ \mu*\nu $ of two probability measures $ \mu, \nu $
on a Borel semigroup $ G $ is defined evidently (as the image of $ \mu
\times \nu $ under $ (x,y) \mapsto xy $). A (one-parameter)
convolution semigroup $ (\mu_t)_{t\in(0,\infty)} $ in $ G $ is defined
accordingly and may be treated as a special (stationary) case of a
convolution system (as defined by \ref{2a1}). If it is separable, we
get a \valued{G} flow system $ (X_{s,t})_{s<t} $ (recall \ref{2a5},
\ref{2a6}), and in addition, automorphisms $ T_h $ of the
corresponding probability space, satisfying $ X_{s,t} \circ T_h =
X_{s+h,t+h} $.

\begin{question}\index{question}
Is $ T_h(\om) $ jointly measurable in $ \om $ and $ h $?
In other words: does every separable convolution semigroup (in every
Borel semigroup) lead to a noise, or not?
\end{question}

Finally, we have to prove two claims formulated without proof in
Sect.~\ref{4aa}, namely, Proposition \ref{2c7} and the last phrase of
Sect.~\ref{4aa}.

\begin{proof}[Proof \textup{(sketch)} of Prop.~\textup{\ref{2c7}}]
Treating the given convolution semigroup as a (stationary) convolution system,
we see that it is separable, since the corresponding flow system $
(X_{s,t})_{s<t;s,t\in T} $ on any countable $ T \subset \R $ satisfies 
$ X_{s_n,t_n} \to X_{s,t} $ in probability whenever $ s < t $, $ s_n \uparrow
s $ and $ t_n \downarrow t $ (which follows from \ref{2c6}(c) and
\eqref{4aa*}), provided that $ s,t,s_n,t_n \in T $. The latter restriction
becomes unnecessary after applying \ref{2a6}.
\proofend

\smallskip
\textsc{digression: stationary process}
\smallskip

Stationarity of a random process $ X : \R \to L_0(\Om) $ is usually defined by
the probability theory in terms of (finite-dimensional) joint distributions;
they must be invariant under time shifts. On the other hand, the ergodic
theory usually defines stationarity in terms of a measurable action $ T : \R
\to \Aut(\Om) $ such that $ X_t \circ T_h = X_{t+h} $. The two definitions are
equivalent under appropriate conditions, formulated below. However, stochastic
flows are random processes on the two-dimensional domain $ \{ (s,t) : -\infty
< s < t < \infty \} $ rather than $ \R $, which is taken into account by our
general formulation.

\begin{lemma}\label{2cdig}
Let $ E $ be a set, $ (S_h)_{h\in\R} $ a one-parameter group of
transformations $ S_h : E \to E $, and $ X : E \to L_0(\Om,\F,P) $ a random
process satisfying two conditions:

(a) for all $ n = 1,2,\dots $ and all $ e_1,\dots,e_n \in E $, the joint
distribution of $ X(S_h(e_1)), \dots, X(S_h(e_n)) $ does not depend on $ h \in
\R $;

(b) for all $ e \in E $, the map $ h \mapsto X(S_h(e)) $ from $ \R $ to $
L_0(\Om,\F,P) $ is Borel measurable.

Then there exists one and only one measurable action $ T : \R \to \Aut(\Om) $
such that
\[
X(e) \circ T_h = X(S_h(e)) \quad \text{a.s.}
\]
for all $ e \in E $, $ h \in \R $.
\end{lemma}

\beginproof
For each $ h $ separately, existence and uniqueness of $ T_h $ follow from
\cite[15.11 and 17.46]{Ke}. In order to prove Borel measurability of $ T $ it
is sufficient to prove that $ h \mapsto X \circ T_h $ is a Borel measurable
map $ \R \to L_0(\Om) $ for every $ X \in L_0(\Om) $. The set of all $ X $
possessing this property is closed, and contains $ \phi(X_1,\dots,X_n) $
whenever it contains $ X_1,\dots,X_n $, for every Borel measurable $ \phi :
\R^n \to \R $. Also, the set generates the \sif\ $ \F $, since it contains all
$ X(e) $ (because the map $ h \mapsto X(S_h(e)) = X(e) \circ T_h $ is
Borel). Therefore the set is the whole $ L_0(\R) $.
\proofend

The group $ \R $ may be replaced with an arbitrary Borel group (which will not
be used).

\smallskip
\textsc{end of digression}
\smallskip

Recall that the notion `topo-semigroup' was defined in \ref{2c6}. Every
topo-semigroup is also a Borel semigroup.

\begin{proof}[Proof \textup{(sketch) of the last phrase of Sect.~\ref{4aa}}]
  \emph{``Every weakly continuous convolution semigroup in a topo-semigroup
    leads to a stationary flow and further to a noise''.}

The flow is given by \ref{2c7}; it remains to prove measurability of the time
shift action. By \ref{2cdig} it is sufficient to prove that the map $ h
\mapsto X_{s+h,t+h} $ from $ \R $ to $ L_0(\Om \to G) $ is Borel measurable. A
stronger claim will be proven: $ X_{s,t} $ is jointly measurable in $ s,t
$. First, for every $ s $ the function $ t \mapsto X_{s,t} $ from $ (s,\infty)
$ to $ L_0(\Om \to G) $ is right-continuous, therefore Borel
measurable. Second, for every $ t $ the function $ s \mapsto X_{s,t} $ from $
(-\infty,t) $ to $ L_0(\Om \to G) $ is left-continuous, therefore Borel
measurable. It follows that $ X_{r,t} = X_{r,s} X_{s,t} $ is a Borel
measurable function of $ r,t $.
\proofend

In fact, Condition \ref{2c7}(c) may be weakened to $ x_n y \to y $ and $ y z_n
\to y $ (instead of $ x_n y z_n \to y $); it still leads to a noise (but
\ref{2c7}(d) could fail). Separability is achieved by $ X_{s,t} = \lim_m
\lim_n X_{s_m,t_n} $ (rather than $ \lim_n X_{s_n,t_n} $).

\mysection{Stability and sensitivity}
\label{sec:5}
\mysubsection{Morphism, joining, maximal correlation}
\label{5a}

The idea, presented in Sect.~\ref{1c}, is formalized below.

\begin{definition}
Let $ \( (\Om_1,P_1), (\F_{s,t}^{(1)})_{s<t} \) $ and $ \(
(\Om_2,P_2), (\F_{s,t}^{(2)})_{s<t} \) $ be two continuous products of
probability spaces.

(a)
A \emph{morphism}%
\index{morphism!of continuous products!of probability spaces}
from the first product to the second is a morphism
of probability spaces $ \al : \Om_1 \to \Om_2 $, measurable from $
(\Om_1, \F_{s,t}^{(1)}) $ to $ (\Om_2, \F_{s,t}^{(2)}) $ whenever $
s<t $.

(b)
An \emph{isomorphism}%
\index{isomorphism!of continuous products!of probability spaces}
from the first product to the second is a
morphism $ \al $ such that the inverse map $ \al^{-1} $ exists and is
also a morphism (of the products).
\end{definition}

If a morphism of products is an isomorphism of probability spaces then
it is an isomorphism of products.

\begin{example}\label{5a2}
Let $ (B_t^{(1)}, B_t^{(2)})_{t\in[0,\infty)} $ be the standard
Brownian motion in $ \R^2 $, and $ \( (\Om_1,P_1),
(\F_{s,t}^{(1)})_{s<t} \) $ be the continuous product of probability
spaces generated by the (two-dimensional) increments $
( B_t^{(1)}-B_s^{(1)}, B_t^{(2)}-B_s^{(2)} ) $. Let $ \(
(\Om_2,P_2), (\F_{s,t}^{(2)})_{s<t} \) $ correspond in the same way to
the standard Brownian motion $ (B_t)_{t\in[0,\infty)} $ in $ \R
$. Then for every $ \phi \in \R $ the formula
\[
B_t = B_t^{(1)} \cos\phi + B_t^{(2)} \sin\phi
\]
defines a morphism (not an isomorphism, of course) from the first
product to the second.
\end{example}

\begin{definition}
A \emph{morphism}\index{morphism!of noises}
from a noise to another noise is a morphism $ \al $
between the corresponding continuous products of probability spaces
that intertwines the corresponding shifts:
\[
\al \circ T_h^{(1)} = T_h^{(2)} \circ \al \quad \text{a.s.}
\]
for every $ h \in \R $.
\end{definition}

Similarly to Example \ref{5a2} we have for each $ \phi $ a morphism
from the two-dimensional white noise to the one-dimensional white
noise.

\begin{example}\label{5a4}
The homomorphism $ f_{a,b,c} \mapsto a $ from the semigroup $
G^{\text{\ref{4c}}} $ ($ =G $ of Sect.~\ref{4c}) to $ (\R,+) $ leads to a
morphism  (not an isomorphism) from the noise of stickiness to the
(one-dimensional) white noise.
The same holds for the noise of splitting.
\end{example}

\begin{definition}
A \emph{joining}%
\index{joining!of continuous products!of probability spaces}%
\index{coupling|see{joining}}
(or \emph{coupling}) of two continuous products of
probability spaces $ \( (\Om_1,P_1), (\F_{s,t}^{(1)})_{s<t} \) $ and $
\( (\Om_2,P_2), (\F_{s,t}^{(2)})_{s<t} \) $ consists of a third
continuous product of probability spaces $ \( (\Om,P),
(\F_{s,t})_{s<t} \) $ and two morphisms $ \al : \Om \to \Om_1 $, $ \be
: \Om \to \Om_2 $ of these products such that $ \F_{-\infty,\infty} $
is generated by $ \al, \be $ (that is, by inverse images of $
\F^{(1)}_{-\infty,\infty} $ and $ \F^{(2)}_{-\infty,\infty} $).
\end{definition}

Each joining leads to a measure on $ \Om_1 \times \Om_2 $ with given
projections $ P_1, P_2 $; namely, the image of $ P $ under the
(one-to-one) map $ \om \mapsto (\al(\om),\be(\om)) $. Two joinings
that lead to the same measure (on $ \Om_1 \times \Om_2 $) will be
called isomorphic.%
\index{isomorphic joinings!of continuous products!of probability spaces}

A joining of a continuous product of probability spaces with itself
will be called a
\emph{self-joining.}%
\index{self-joining!of continuous products!of probability spaces}
A \emph{symmetric}%
\index{self-joining!symmetric}\index{symmetric!self-joining}
self-joining
is a self-joining $ (\al,\be) $ isomorphic to $ (\be,\al) $. For
example, every pair of angles $ \phi, \psi $ leads to a symmetric
self-joining of the (one-dimensional) white noise,
\begin{equation}\label{5a6}
\begin{aligned}
B_t \circ \al &= B_t^{(1)} \cos\phi + B_t^{(2)} \sin\phi \, , \\
B_t \circ \be &= B_t^{(1)} \cos\psi + B_t^{(2)} \sin\psi \, .
\end{aligned}
\end{equation}
Only the difference $ |\phi-\psi| $ matters (up to isomorphism), thus, we have
a one-parameter family of self-joinings; $ \rho = \cos(\phi-\psi) \in [-1,1] $
is a natural parameter. In fact, this family exhausts \emph{all} self-joinings
of the one-dimensional white noise; therefore, they all are
symmetric. However, the two-dimensional white noise admits also asymmetric
self-joinings. For example, the self-joining satisfying
\begin{equation}\label{5a*}
\begin{aligned}
B_t^{(1)} \circ \be &= ( B_t^{(1)} \circ \al ) \cos\phi - ( B_t^{(2)} \circ
  \al ) \sin\phi \, , \\
B_t^{(2)} \circ \be &= ( B_t^{(1)} \circ \al ) \sin\phi + ( B_t^{(2)} \circ
  \al ) \cos\phi
\end{aligned}
\end{equation}
is symmetric if and only if $ \sin \phi = 0 $.

Every joining $ (\al,\be) $ of two continuous products of probability
spaces has its \emph{maximal correlation}%
\index{maximal correlation!for probability spaces}
\[
\rho^{\max} (\al,\be) = \sup | \Ex (f\circ\al) (g\circ\be) | \, ,
\]
where the supremum is taken over all $ f \in L_2(\Om_1,P_1) $, $ g \in
L_2(\Om_2,P_2) $ such that $ \Ex f = 0 $, $ \Ex g = 0 $, $ \Var f \le
1 $, $ \Var g \le 1 $. (All $ L_2 $ spaces are real, not complex.) The
product structure is irrelevant to the
`global' correlation $ \rho^{\max} (\al,\be) $, but relevant to `local'
correlations $ \rho^{\max}_{s,t} (\al,\be) $; here the supremum is taken
under an additional condition: $ f $ is \measurable{\F_{s,t}^{(1)}},
and $ g $ is \measurable{\F_{s,t}^{(2)}}. Surprisingly, the global
correlation is basically the supremum of local correlations over
infinitesimal time intervals.

\begin{proposition}\label{5a7}
$ \rho^{\max}_{r,t} (\al,\be) = \max \( \rho^{\max}_{r,s} (\al,\be), \,
\rho^{\max}_{s,t} (\al,\be) \) $ whenever $ r<s<t $.
\end{proposition}

\begin{sloppypar}
\beginproof
Generally, $ L_2(\F_{r,t}) = L_2(\F_{r,s}) \otimes L_2(\F_{s,t}) = \(
\R \oplus L_2^0(\F_{r,s}) \) \otimes \( \R \oplus L_2^0(\F_{s,t}) \) =
\R \otimes \R \oplus \R \otimes L_2^0(\F_{s,t}) \oplus L_2^0(\F_{r,s})
\otimes \R \oplus L_2^0(\F_{r,s}) \otimes L_2^0(\F_{s,t}) $, where $
\R $ is the one-dimensional space of constants, and $ L_2^0(\dots) $
--- its orthogonal complement (the zero-mean space). We apply the
argument to $ \F^{(1)} $ and $ \F^{(2)} $, decompose $ f $ and $ g $
into three orthogonal summands each ($ \R \otimes \R $ does not
appear), and get the maximum of $ \rho^{\max}_{r,s} $, $
\rho^{\max}_{s,t} $ and $ \rho^{\max}_{r,s} \rho^{\max}_{s,t} $.
\proofend
\end{sloppypar}

In fact, $ \rho^{\max} (\al,\be) = \cos(\phi-\psi) $ for the
self-joining \eqref{5a6}; however, $ \rho^{\max} (\al,\be) \linebreak[0]
 = 1 $ for the self-joining \eqref{5a*}, irrespective of $ \phi $.

\mysubsection{A generalization of the Ornstein-Uhlenbeck semigroup}
\label{5b}

Here is the `best' self-joining for a given maximal correlation. (See
also \cite[Lemma~2.1]{WW03}.)

\begin{proposition}\label{5a8}
Let $ \( (\Om,P), (\F_{s,t})_{s<t} \) $ be a continuous product of
probability spaces, and $ \rho \in [0,1] $. Then there exists a
symmetric self-joining $ (\al_\rho,\be_\rho) $%
\index{zza@$ (\al_\rho, \be_\rho) $, self-joining!for probability spaces}
of the given product such that
\[
\rho^{\max} (\al_\rho,\be_\rho) \le \rho
\]
and
\[
| \Ex ( f \circ \al ) ( f \circ \be ) | \le \Ex ( f \circ \al_\rho ) (
  f \circ \be_\rho ) 
\]
for all $ f \in L_2(\Om) $ and all self-joinings $ (\al,\be) $
satisfying $ \rho^{\max} (\al,\be) \le \rho $.

The self-joining $ (\al_\rho,\be_\rho) $ is unique up to isomorphism.
\end{proposition}

\smallskip
\textsc{digression: the compact space of joinings}
\smallskip

Here we forget about continuous products and deal with joinings of two
probability spaces. Let them be just $ [0,1] $ with Lebesgue
measure; the non\-atomic case is thus covered. (Atoms do not invalidate
the results and are left to the reader.) Up to isomorphism, joinings%
\index{joining!of probability spaces}
are probability measures $ \mu $ on the square $ [0,1] \times [0,1] $
with given (Lebesgue) projections to both coordinates. They are a
closed subset of the compact metrizable space of all probability
measures on the square, equipped with the weak topology (generated by
integrals of continuous functions). Surprisingly, the topology of $
[0,1] $ plays no role.

\begin{lemma}
Let $ \mu, \mu_1, \mu_2, \dots $ be probability measures on $ [0,1]
\times [0,1] $ with Lebesgue projections to both coordinates. Then the
following conditions are equivalent:

(a) $ \int f \, \D \mu_n \to \int f \, \D \mu $ for all continuous
functions $ f : [0,1] \times [0,1] \to \R $;

(b) $ \mu_n (A \times B) \to \mu(A \times B) $ for all measurable sets
$ A,B \subset [0,1] $;

(c) $ \int f(x) g(y) \, \mu_n(\D x \D y) \to \int f(x) g(y) \, \mu(\D
x \D y) $ for all $ f,g \in L_2[0,1] $.
\end{lemma}

\beginproof
(c) \imp (b): trivial; (b) \imp (a): we approximate $ f $ uniformly by
functions constant on each $ \(\frac k n,\frac{k+1}n\) \times \(\frac
l n,\frac{l+1}n\) $.

(a) \imp (c): we choose continuous $ f_\eps, g_\eps $ such that $ \| f
- f_\eps \| \le \eps $ and $ \| g - g_\eps \| \le \eps $ (the norms
are in $ L_2[0,1] $). We apply (a) to the continuous function $ (x,y)
\mapsto f_\eps(x) g_\eps(y) $ and note that $ \iint | f(x)g(y) -
f_\eps(x) g_\eps(y) | \, \mu_n(\D x \D y) \le \| f \| \| g - g_\eps \|
+ \| f - f_\eps \| \| g_\eps \| $.
\proofend

(It is basically Slutsky's lemma, see \cite[Lemma 0.5.7]{RY} or
\cite[Th.~1]{BEKSY}; see also \cite[Lemma B3]{TV}.)

We see that all joinings (of two given standard probability spaces)
are (naturally) a compact metrizable space.

Here is another (unrelated to the compactness) general fact about
joinings, used in the sequel. It is formulated for $ [0,1] $ but holds
for all probability spaces.

\begin{lemma}\label{5a10}
Let $ \mu_1, \mu_2, \dots $ be probability measures on $ [0,1] \times
[0,1] $ with Lebesgue projections to both coordinates. Then there
exists a sub-\sif\ $ \F $ of the Lebesgue \sif\ on $ [0,1] $ such that
\[
\iint |f(x)-f(y)|^2 \, \mu_n(\D x \D y) \to 0 \quad \text{if and only
if $ f $ is \measurable{\F}}
\]
for all $ f \in L_2[0,1] $.
\end{lemma}

\beginproof
For each $ n $ the quadratic form $ \Ec_n(f) = \int |f(x)-f(y)|^2 \,
\mu_n(\D x \D y) $ satisfies $ \Ec_n(T\circ f) \le \Ec_n(f) $ for
every $ T : \R \to \R $ such that $ |T(a)-T(b)| \le |a-b| $ for all $
a,b \in \R $. (Thus, $ \Ec_n $ is a continuous symmetric Dirichlet
form.) The set $ E $ of all $ f $ such that $ \Ec_n(f) \to 0 $ is a
(closed) linear subspace of $ L_2[0,1] $, and $ f \in E $ implies $
T\circ f \in E $. Such $ E $ is necessarily of the form $ L_2(\F) $,
see for instance \cite[Problem IV.3.1]{Ne}.
\proofend

Note also that every joining $ (\al,\be) $ leads to a bilinear form $
(f,g) \mapsto \Ex (f\circ\al) (g\circ\be) $ and the corresponding
operator $ U_{\al,\be} : L_2[0,1] \to L_2[0,1] $, $ \ip{U_{\al,\be}f}g
= \Ex (f\circ\al) (g\circ\be) $. Generally $ U_{\al,\be} $ maps one $
L_2 $ space into another, but for a self-joining we deal with a single
space. Clearly, $ U_{\be,\al} = (U_{\al,\be})^* $; $ U_{\al,\be} $ is
Hermitian if and only if the joining is symmetric.

\smallskip
\textsc{end of digression}

\begin{proof}[Proof \textup{(sketch)} of Proposition \textup{\ref{5a8}}]
Uniqueness: a self-joining $ (\al,\be) $ is uni\-quely determined by its
bilinear form $ (f,g) \mapsto \Ex(f\circ\al)(g\circ\be) $, therefore
a symmetric self-joining is uniquely determined by its quadratic form
$ f \mapsto \Ex(f\circ\al)(f\circ\be) $.

Existence. First, on the space $ \Om\times\Om $ we consider the maps $
\al_\rho (\om_1,\om_2) = \om_1 $, $ \be_\rho (\om_1,\om_2) = \om_2 $
and the measure $ \ti P_0 = \rho P_{\text{diag}} + (1-\rho) P \times P
$; here $ P_{\text{diag}} $ is the image of $ P $ under the map $ \om
\mapsto (\om,\om) $. We get a symmetric self-joining of the
probability space $ (\Om,P) $, but not of the continuous product $
(\F_{s,t})_{s<t} $. Note that $ | \Ex(f\circ\al)(f\circ\be) | \le (\Ex
f)^2 + \rho \Var f = \int (f\circ\al_\rho) (f\circ\be_\rho)
\, \D \ti P_0 $ for all $ f \in L_2(\Om) $ and all self-joinings $
(\al,\be) $ satisfying $ \rho^{\max} (\al,\be) \le \rho $.

Second, we apply the same construction to $ \Om_{-\infty,0} $ and $
\Om_{0,\infty} $ separately, and consider $ \ti P_1 = \ti
P_{-\infty,0} \times \ti P_{0,\infty} $. Similarly to the proof of
Proposition \ref{5a7} we see that $ | \Ex (f\circ\al) (f\circ\be) |
\le \int (f\circ\al_\rho) (f\circ\be_\rho) \, \D \ti P_1 \le \int
(f\circ\al_\rho) (f\circ\be_\rho) \, \D \ti P_0 $ for all $ f \in
L_2(\Om) $ and all self-joinings $ (\al,\be) $ satisfying $
\rho^{\max} (\al,\be) \le \rho $.

Third, we do it for every decomposition $ \Om = \Om_{-\infty,t_1}
\times \Om_{t_1,t_2} \times \dots \times \Om_{t_{n-1},t_n} \times
\Om_{t_n,\infty} $ and get a net of measures, symmetric self-joinings,
and their quadratic forms
\begin{equation}\label{5b3a}
\begin{gathered}
\ti U^\rho_{t_1,\dots,t_n} f = \bigotimes_{k=0}^n \( \rho f_k +
 (1-\rho) \Ex f_k \) \, , \\
\ip{ \ti U^\rho_{t_1,\dots,t_n} f }{ f } = \prod_{k=0}^n \( (\Ex
 f_k)^2 + \rho \Var f_k \)
\end{gathered}
\end{equation}
for $ f = f_0 \otimes \dots \otimes f_n $, $ f_0 \in
L_2(\Om_{-\infty,t_1}), \, f_1 \in L_2(\Om_{t_1,t_2}), \, \dots, \,
f_n \in L_2(\Om_{t_n,\infty}) $. The net converges (in the compact space of
joinings) due to monotonicity of the net of quadratic forms. The limit
is a symmetric self-joining $ (\al_\rho,\be_\rho) $ of the continuous
product of probability spaces. It majorizes $ | \Ex (f\circ\al)
(f\circ\be) | $, since every element of the net does.
\end{proof}

See also \cite[5b4]{Ts03}.
Basically, each infinitesimal element of the data set is replaced with
a fresh copy, independently of others, with probability $ 1-\rho $.
Doing it twice with parameters $ \rho_1 $ and $ \rho_2 $ is equivalent
to doing it once with parameter $ \rho = \rho_1 \rho_2 $. In terms of
operators $ U^\rho = U_{\al_\rho,\be_\rho} : L_2(\Om) \to L_2(\Om) $%
\index{zzu@$ U^\rho $, operator!for probability spaces}
it means $ U^{\rho_1} U^{\rho_2} = U^{\rho_1 \rho_2} $; a
one-parameter semigroup! It seems to lead to an \valued{\Om}
stationary Markov process $ (X_u)_{u\in\R} $, $ X_u : \ti\Om \to \Om
$, such that for every $ u>0 $ the pair $ (X_0,X_u) $ is distributed
like the pair $ (\al_\rho,\be_\rho) $ where $ \rho = \E^{-u}
$. However, $ L_2(\ti\Om) $ is separable only in classical cases.

\begin{sloppypar}
If the given continuous product of probability spaces is the white
noise then the Markov process is the well-known Ornstein-Uhlenbeck
process (infinite-dimensional, over the Gaussian measure that
describes the white noise). 
\end{sloppypar}

The proof of the relation $ U^{\rho_1} U^{\rho_2} = U^{\rho_1 \rho_2}
$ is an easy supplement to the proof of Proposition \ref{5a8}; an
elementary check for each element $ \ti U^\rho_{t_1,\dots,t_n} $ of
the net, and a passage to the
limit. In the same way we prove that the spectrum of the Hermitian
operator $ U^\rho $ is contained in $ \{ 1, \rho, \rho^2, \dots \}
\cup \{ 0 \} $. The spectral theorem gives the following.

\begin{proposition}\label{5a11}
Let $ \( (\Om,P), (\F_{s,t})_{s<t} \) $ be a continuous product of
probability spaces, $ (\al_\rho,\be_\rho) $ the self-joinings given by
Prop.~\ref{5a8}, and $ U^\rho $ the corresponding operators, that is,
$ \Ex (f\circ\al_\rho) (g\circ\be_\rho) = \ip{ U^\rho f}g $ for all $
f,g \in L_2(\Om) $. Then there exist (closed linear) subspaces $ H_0,
H_1, H_2, \dots $ and $ H_\infty $%
\index{zzh@$ H_n $, subspace!for probability spaces}
of $ L_2(\Om) $ such that
\[
L_2(\Om) = ( H_0 \oplus H_1 \oplus H_2 \oplus \dots ) \oplus H_\infty
\]
(that is, the subspaces are orthogonal and span the whole $ L_2(\Om)
$), and
\begin{gather*}
U^\rho f = \rho^n f \quad \text{for $ f \in H_n $, $ \rho \in [0,1] $}
 \, , \\
U^\rho f = 0 \quad \text{for $ f \in H_\infty $, $ \rho \in [0,1) $}
 \, .
\end{gather*}
\end{proposition}

Of course, $ U^1 f = f $ for all $ f $. The semigroup $ (U^\rho)_\rho
$ is strongly continuous if and only if $ \dim H_\infty = 0 $. Note
also that $ H_0 $ is the one-dimensional space of constants.

The spaces $ H_n $ may be called
\emph{chaos spaces,}\index{chaos spaces}
since for the
white noise $ H_n $ is the $ n $-th Wiener chaos space (and $ \dim
H_\infty = 0 $).

\mysubsection{The stable $ \sigma $-field; classical and nonclassical}
\label{5c}

\begin{definition}
Let $ \( (\Om,P), (\F_{s,t})_{s<t} \) $ be a continuous product of
probability spaces.

(a)
A random variable $ f \in L_2(\Om) $ is
\emph{stable}\index{stable!random variable}
if there exist
symmetric self-joinings $ (\al_n,\be_n) $ of the continuous product
such that
\begin{gather*}
\rho^{\max} (\al_n,\be_n) < 1 \quad \text{for every } n \, , \\
\Ex | f\circ\al_n - f\circ\be_n |^2 \to 0 \quad \text{as } n \to
\infty \, .
\end{gather*}

(b)
A random variable $ f \in L_2(\Om) $ is
\emph{sensitive}\index{sensitive random variable}
if $
\Ex(f\circ\al)(g\circ\be) = 0 $ for all $ g \in L_2(\Om) $ and all
symmetric self-joinings $ (\al,\be) $ of the continuous product
such that $ \rho^{\max} (\al,\be) < 1 $.
\end{definition}

\begin{sloppypar}
\begin{theorem}\index{theorem}
(Tsirelson \cite[2.5]{Ts99}, \cite[5b11]{Ts03})
For every continuous product of probability spaces $ \( (\Om,P),
(\F_{s,t})_{s<t} \) $ there exists a sub-\sif\ $ \F^\stable \subset
\F_{-\infty,\infty} $%
\index{stable!sub-$\sigma$-field}%
\index{zzf@$ \F^\stable $, sub-$\sigma$-field!for probability spaces}
such that
\begin{gather*}
\text{$ f $ is stable if and only if $ f $ is \measurable{\F^\stable}}
\, , \\
\text{$ f $ is sensitive if and only if $ \cE{f}{\F^\stable} = 0 $}
\end{gather*}
for all $ f \in L_2(\Om) $.
\end{theorem}
\end{sloppypar}

\beginproof
Let $ (\al,\be) $ be a symmetric self-joining, $ \rho^{\max} (\al,\be)
\le \rho $. Rewriting the inequality $ \Ex (f\circ\al) (f\circ\be) \le
\Ex (f\circ\al_\rho) (f\circ\be_\rho) $ as $ \Ex | f\circ\al -
f\circ\be |^2 \ge \Ex | f\circ\al_\rho - f\circ\be_\rho |^2 $ we see
that $ f $ is stable iff $ \Ex | f\circ\al_\rho - f\circ\be_\rho |^2
\to 0 $ as $ \rho \to 1 $. Lemma \ref{5a10} gives us a \sif\ $
\F^\stable $ such that $ f $ is stable iff $ f $ is
\measurable{\F^\stable}. Also, $ f $ is stable iff $ \ip{U^\rho f}f
\to \| f \|^2 $ as $ \rho \to 1- $, that is, $ f $ is orthogonal to $
H_\infty $.

We have $ | \ip{U_{\al,\be}f}f | \le \ip{U^\rho f}f $ for all $ f $,
therefore $ | \ip{U_{\al,\be}f}f | \le \sqrt{\ip{U^\rho f}f} \linebreak[0]
\sqrt{\ip{U^\rho g}g} $. Rewriting sensitivity of $ f $ in the form $
\forall \al,\be \, \forall g \; \ip{U_{\al,\be}f}g = 0 $ we see that $
f $ is sensitive iff $ U^\rho f = 0 $ for all $ \rho < 1 $, that is, $
f \in H_\infty $. \qed
\proofendnoqed

\begin{theorem}\label{5c4}\index{theorem}
For every continuous product of probability spaces $ \(
(\Om,P),\linebreak[0]
(\F_{s,t})_{s<t} \) $ there exists a symmetric self-joining $
(\al_{1-}, \be_{1-}) $
\index{zza@$ (\al_{1-}, \be_{1-}) $, self-joining}
of the given product such that
\[
\Ex (f\circ\al_{1-}) (g\circ\be_{1-}) = \Ex(fg)
\]
if $ f,g \in L_2(\Om) $ are stable, but
\[
\Ex (f\circ\al_{1-}) (g\circ\be_{1-}) = 0
\]
if $ f \in L_2(\Om) $ is sensitive (and $ g \in L_2(\Om) $ is
arbitrary).
The self-joining $ (\al_{1-}, \be_{1-}) $ is unique up to
isomorphism.
\end{theorem}

\beginproof
Just take the limit of $ (\al_\rho,\be_\rho) $ in the (compact!) space
of joinings, as $ \rho\to1- $.
\proofend

See `the $ 1^- $-joining' in \cite[Def.~2.2]{WW03}; see also
\cite[Sect.~1 (for $ p=0 $)]{Wa99}.

\begin{definition}
A continuous product of probability spaces is
\emph{classical,}%
\index{classical (part of)!continuous product!of probability spaces}
if it satisfies the following equilavent conditions:

(a) all random variables are stable;

(b) no random variable is sensitive;

(c) the stable sub-\sif\ $ \F^\stable $ is the whole \sif\ $ \F $.

A noise is classical\index{classical (part of)!noise}
if the underlying continuous product of
probability spaces is classical.
\end{definition}

Equivalent definitions in terms of \flow{\R}s (L\'evy processes)
exist, see \ref{6a13}, \ref{6c2}. See also \cite[Sect.~5b, especially
Def.~5b5]{Ts03}, and \cite[2.5]{Ts99}.

\begin{remark}\label{5c5}
The continuous product of probability spaces, corresponding to a flow
system $ (X_{s,t})_{s<t} $, is classical if and only if random
variables $ \phi(X_{s,t}) $ are stable for all $ s<t $ and all bounded
Borel functions $ \phi : G_{s,t} \to \R $ (a single $ \phi $ is enough
if it is one-to-one).
\end{remark}

\beginproof
If each $ \phi(X_{s,t}) $ is \measurable{\F^\stable} then $ \F^\stable
= \F $ since $ \F $ is generated by these $ \phi(X_{s,t}) $.
\proofend

\mysubsection{Examples}
\label{5d}

The time set implicit in Sections \ref{5a}--\ref{5c} is not
necessarily $ \R $; a subset of $ \R $ (or any linearly ordered set)
is also acceptable. In particular, the theory is applicable to the
`singularity concentrated in time' cases of Sect.~\ref{sec:1}.

The \flow{\Z_m} $ X^{\text{\ref{1b}}} $ ($ = X $ of Sect.~\ref{1b})
generates a continuous product of probability spaces over the time set
$ \{ 0,1,2,\dots \} \cup \{\infty\} $. Random variables $ X_{s,s+1} $
are stable; indeed, a single (indivisible) element of the data set is
replaced with probability $ 1-\rho $, therefore $ \Pr{ X_{s,s+1} \circ
\al_\rho \ne X_{s,s+1} \circ \be_\rho } \le 1-\rho $. It follows that
$ X_{s,t} = X_{s,s+1} \dots X_{t-1,t} $ is stable whenever $
s<t<\infty $. Of course, $ X_{s,t} $, being a \valued{\Z_m} random
variable, is not an element of $ L_2(\Om) $. By stability of $
X_{s,t} $ we mean stability of $ \phi(X_{s,t}) $ for every $ \phi : \Z_m
\to \R $.

In contrast, the random variable $ X_{0,\infty} $ is sensitive by
Prop.~\ref{1c1}. The same holds for $ X_{s,\infty} $. More exactly, $
\psi(X_{s,\infty}) $ is sensitive for every $ \psi : \Z_m \to \R $ such that
$ \Ex \psi(X_{s,\infty}) = 0 $. \emph{Sketch of the proof} (see also
Sect.~\ref{1c} for $ m=2 $):
\begin{multline*}
X_{s,\infty} \circ \al_\rho - X_{s,\infty} \circ \be_\rho =
( X_{s,s+1} \circ \al_\rho - X_{s,s+1} \circ \be_\rho ) + \dots + \\
+ ( X_{t-1,t} \circ \al_\rho - X_{t-1,t} \circ \be_\rho ) + 
( X_{t,\infty} \circ \al_\rho - X_{t,\infty} \circ \be_\rho ) \, ,
\end{multline*}
the summands being independent. For large $ t $ the sum from $ s $ to
$ t $ is distributed on $ \Z_m $ approximately uniformly, therefore $ 
X_{s,\infty} \circ \al_\rho - X_{s,\infty} \circ \be_\rho $ is
uniform. The same holds conditionally, given $ \al_\rho $ (that is, $
X_{r,t} \circ \al_\rho $ for all $ r,t $ including $ t=\infty $).

We see that random variables of the form $ \phi(X_{0,1},X_{1,2},\dots) $
are stable, and random variables of the form $
\phi(X_{0,1},X_{1,2},\dots) \psi(X_{0,\infty}) $ are sensitive (as before, $
\sum_{x\in\Z_m} \psi(x) = 0 $). Their sums exhaust $ L_2(\Om)
$. Therefore $ \F^\stable $ is generated by $ X_{0,1},X_{1,2},\dots $;
random variables $ X_{s,\infty} $ are independent of $ \F^\stable $
(each one separately).

The \flow{\T} $ Y^{\text{\ref{1b}}} $ (over the time set $ [0,\infty) $)
behaves similarly: $ \F^\stable $ is generated by $ Y_{s,t} $ for $
0<s<t<\infty $; random variables $ Y_{0,t} $ are independent of $
\F^\stable $ (each one separately).

We turn to the noises of Sect.~\ref{sec:4}: splitting and
stickiness. These two may be treated uniformly. Below, $ G $ is either
$ G^{\text{\ref{4b}}} $ or $ G^{\text{\ref{4c}}} $. The Brownian motion $ (B_t)_t
= (a_{0,t})_t $ generates (via increments) sub-\sif s $
\F_{s,t}^\white \subset \F_{s,t} $. It will be shown that $ \F^\stable
= \F_{-\infty,\infty}^\white $.

The Brownian motion $ (B_t)_t $ has the \emph{predictable
representation property}%
\index{predictable representation property}
w.r.t.\ the filtration $ (\F_{-\infty,t})_t
$. That is, every local martingale $ (M_t)_t $ in this filtration is
of the form $ M_t = M_{-\infty} + \int_{-\infty}^t h_s \, \D B_s $ for
some predictable process $ (h_t)_t $ (in the considered filtration);
see \cite[Def.~V.4.8]{RY}. Note that $ M_t $ and $ h_t $ need not be
\measurable{\F_{-\infty,t}^\white}.

\begin{proof}[Proof \textup{(sketch)} of the predictable representation
property]
We may restrict ourselves to a dense set of martingales, namely,
$ M_s = \cE{ \phi(X_{t_0,t_1},\dots,X_{t_{n-1},t_n})\linebreak[0]
 }{ \F_{-\infty,s} } $ where $ \phi : G^n \to \R $ is a bounded measurable (or
even smooth) function, $ -\infty < t_0 < \dots < t_n < \infty $, and $
(X_{s,t})_{s<t} $ stands for the given \flow{G}. When $ s \in [t_{k-1},t_k] $,
we deal effectively with the case $ M_s = \cE{ \psi(X_{r,t}) }{ \F_{r,s}
} $, $ -\infty < r < t < \infty $, to which we may restrict ourselves.
By independence of $ X_{r,s} $ and $ X_{s,t} $,
\begin{align*}
& M_s = u (X_{r,s}, t-s) \, , \quad \text{where $ u : G \times \R \to
 \R $ is defined by} \\
& u(x,t) = \int_G \psi(xy) \, \mu_t(\D y) \, ,
\end{align*}
and $ (\mu_t)_t $ is the given convolution semigroup in $ G $.

The semigroup $ G $ is in fact a smooth manifold with boundary, and
the function $ u $ is smooth up to the boundary (which can be checked
using explicit formulas for $ \mu_t $ and the binary operation in $ G
$). The random process $ (X_{r,s})_{s\in[r,t]} $ is a diffusion
process on the smooth manifold $ G $; it is a weak solution of a
stochastic differential equation driven by $ (B_s)_s $. It\^o's
formula gives the needed representation.
\end{proof}

By a Brownian motion adapted%
\index{adapted to a continuous product}
to a continuous product of probability
spaces $ \( (\Om,P), (\F_{s,t})_{s<t} \) $ we mean a family $
(B_t)_{t\in\R} $ of random variables $ B_t $ such that
\[
B_t - B_s \text{ is \measurable{\F_{s,t}}, and distributed normally
$ \N(0,t-s) $}
\]
whenever $ -\infty < s < t < \infty $; and in addition, $ B_0 = 0 $.
Note that the $ n $-th Wiener chaos space over $ (B_t)_t $ is included into
the chaos space $ H_n $ over the continuous product.

\begin{proposition}
Let $ \( (\Om,P), (\F_{s,t})_{s<t} \) $ be a continuous product of
probability spaces, and $ (B_t)_t $ a Brownian motion adapted to the
continuous product. If $ (B_t)_t $ has the predictable representation
property w.r.t.\ the filtration $ (\F_{-\infty,t})_t $, then the
sub-\sif\ generated by $ (B_t)_t $ is equal to $ \F^\stable $.
\end{proposition}

\beginproof
The sub-\sif\ $ \F^\white $ generated by $ (B_t)_t $ is contained in $
\F^\stable $, since Wiener chaos spaces (with finite indices) exhaust
the corresponding $ L_2 $ space. We have to prove that $ \F^\white
\supset \F^\stable $.

Every $ f \in L_2^0 (\Om) $ is of the form $ f = \int_{-\infty}^\infty
h_t \, \D B_t $. We have
\[
\ip{ U^\rho f} f = \rho \int_{-\infty}^\infty \ip{ U^\rho h_t}{ h_t }
\, \D t \, ,
\]
since\footnote{%
 Usually $ \ip{\cdot}{\cdot} $ stands for the scalar product, but now
 it will denote the predictable quadratic (co)variation.}
\begin{multline*}
\textstyle
\ip{ (\int h_t \, \D B_t) \circ \al_\rho }{ (\int h_t \, \D B_t) \circ
 \be_\rho } =
\ip{ \int (h_t\circ\al_\rho) \, \D (B_t\circ\al_\rho) }{ \int
 (h_t\circ\be_\rho) \, \D (B_t\circ\be_\rho) } = \\
\textstyle = \int (h_t\circ\al_\rho) (h_t\circ\be_\rho) \,
 \underbrace{ \D \ip{B_t\circ\al_\rho}{B_t\circ\be_\rho} }_{=\rho \,
 dt} \, .
\end{multline*}
In particular, if $ f \in H_1 $ (the first chaos) then
\[
\rho \int \underbrace{ \ip{ U^\rho h_t }{ h_t } }_{\le\|h_t\|^2} \, \D t =
\ip{ U^\rho f} f = \rho \|f\|^2 = \rho \int \| h_t \|^2 \, \D t \, ,
\]
that is, $ \| h_t \|^2 = \ip{ U^\rho h_t }{ h_t } $, which means
that $ h_t \in H_0 $ is a constant (non-random) for almost every $ t
$. Therefore $ f = \int h_t \, \D B_t $ is \measurable{\F^\white},
and we get
\[
H_1 \subset L_2 (\F^\white) \, .
\]
Further, let $ f \in H_2 $, then $ h_t $ are orthogonal to $ H_0 $
(since $ f $ is orthogonal to $ H_1 $), therefore $ \ip{ U^\rho h_t }{
h_t } \le \rho \| h_t \|^2 $. On the other hand,
\[
\rho \int \ip{ U^\rho h_t }{ h_t } \, \D t = \ip{ U^\rho f} f = \rho^2
\|f\|^2 = \rho^2 \int \| h_t \|^2 \, \D t \, ,
\]
that is, $ \rho \| h_t \|^2 = \ip{ U^\rho h_t }{ h_t } $, which
means that $ h_t \in H_1 $ for almost all $ t $. It follows that $
h_t $ is \measurable{\F^\white}; therefore $ f $ is
\measurable{\F^\white}, and we get
\[
H_2 \subset L_2 (\F^\white) \, .
\]
And so on. Finally,
\[
\qquad\qquad\qquad
L_2 (\F^\stable) = H_0 \oplus H_1 \oplus H_2 \oplus \dots \subset L_2
(\F^\white) \, . \qquad\qquad \qed
\]
\proofendnoqed

\mysection{Continuous products: from probability spaces to Hilbert
 spa\-ces}
\label{sec:3}
\mysubsection{Continuous products of spaces $ L_2 $}
\label{3a}

If $ (\Om_1,P_1) $, $ (\Om_2,P_2) $ are probability spaces and $
(\Om,P) = (\Om_1,P_1) \times (\Om_2,P_2) $ is their product, then
Hilbert spaces $ H_1 = L_2 (\Om_1,P_1) $, $ H_2 = L_2 (\Om_2,P_2) $, $
H = L_2 (\Om,P) $ are related via tensor product,%
\index{tensor product}\index{zzo@$ \otimes $ (for Hilbert spaces)}
\[
H = H_1 \otimes H_2 \, .
\]
In terms of bases it means that, having orthonormal bases $
(f_i)_{i\in I} $ in $ H_1 $ and $ (g_j)_{j\in J} $ in $ H_2 $, we get
an orthonormal basis $ (f_i \otimes g_j)_{(i,j)\in I\times J} $ in $ H
$; here
\[
(f \otimes g) (\om_1,\om_2) = f (\om_1) g (\om_2) \quad
\text{for } \om_1 \in \Om_1, \, \om_2 \in \Om_2 \, .
\]
Complex spaces $ L_2 \( (\Om,P) \to \C \) $ and real spaces $ L_2 \(
(\Om,P) \to \R \) $ may be used equally well.

In other words: having a probability space $ (\Om,\F,P) $ and two
sub-\sif s $ \F_1, \F_2 \subset \F $ such that $ \F_1 \otimes \F_2 =
\F $ (recall \eqref{2b2}), we introduce Hilbert spaces $ H_1 = L_2
(\F_1) $ (that is, $ H_1 = L_2 (\Om,\F_1,P) $), $ H_2 = L_2(\F_2) $, $
H = L_2(\F) $ and get $ H = H_1 \otimes H_2 $. This time, $ f \otimes
g $ is just the (pointwise) product of the two functions $
f $, $ g $ on $ \Om $; note that these are \emph{independent}
random variables. In addition we have $ H_1 \subset H_1 \otimes H_2 $
and $ H_2 \subset H_1 \otimes H_2 $, which does not happen in
general. Here it happens because of a special vector $ \One $ (the
constant function on $ \Om $) of $ H_1 $ (and $ H_2 $); $ H_2 $ is
identified with $ \One \otimes H_2 \subset H_1 \otimes H_2 $.

Given a continuous product of probability spaces $ (\Om,P) $, $
(\F_{s,t})_{s<t} $ (as defined by \ref{2b1}), we introduce Hilbert
spaces
\begin{gather*}
H_{s,t} = L_2 (\F_{s,t}) \quad \text{for $ s<t $} \, ; \\
H_{r,t} = H_{r,s} \otimes H_{s,t} \quad \text{for $ r<s<t $} \, .
\end{gather*}

A \emph{unitary operator}\index{unitary operator}
is a linear isometric invertible operator between Hilbert spaces (over
$ \R $ or $ \C $). The group of all unitary operators $ H \to H $ will
be denoted $ \U(H) $. Here is a counterpart of Def.~\ref{2b6}.

\begin{definition}\label{3a1}
A \emph{continuous product of Hilbert spaces}%
\index{continuous product!of Hilbert spaces}
consists of separable
Hilbert spaces $ H_{s,t} $ (given for all $ s,t \in [-\infty,\infty]
$, $ s<t $; possibly finite-dimensional, but not zero-dimensional),
and unitary operators $ H_{r,s} \otimes H_{s,t} \to H_{r,t} $ (given
for all $ r,s,t \in [-\infty,\infty] $, $ r<s<t $), satisfying the
associativity condition:
\[
(fg)h = f(gh) \quad \text{for all $ f \in H_{r,s} $, $ g \in H_{s,t}
$, $ h \in H_{t,u} $} 
\]
whenever $ r,s,t,u \in [-\infty,\infty] $, $ r<s<t<u $.
Here $ fg $ stands for the image of $ f \otimes g $ under the given
operator $ H_{r,s} \otimes H_{s,t} \to H_{r,t} $.
\end{definition}

Note the time set $ [-\infty,\infty] $ rather than $ \R $. Enlarging $
\R $ to $ [-\infty,\infty] $ is easy when dealing with probability
spaces (as noted after Def.~\ref{2b1}) but not Hilbert spaces. Any
linearly ordered set could be used as the time set in Def.~\ref{3a1};
however, existence of the least and greatest elements ($ \pm\infty $)
will be used in Sect.~\ref{3b}. The time set $ \R $ will be treated in
Sects.~\ref{3c}, \ref{3d} in the stationary setup. Homeomorphic time
sets $ [-\infty,\infty] $ and $ [0,1] $ are the same in the general
setup (\ref{3a}, \ref{3b}) but quite different in the stationary setup
(\ref{3c}, \ref{3d}).

Every continuous product of probability spaces leads to a continuous
product of Hilbert spaces.

Given a continuous product of Hilbert spaces $ (H_{s,t})_{s<t} $, we
may consider the
disjoint union\index{disjoint union}
$ \Ec $ of all $ H_{s,t} $,
\[
\Ec = \biguplus_{s<t} H_{s,t} = \{ (s,t,f) : -\infty \le s < t \le
\infty, \, f \in H_{s,t} \}
\]
and a partial binary operation
\[
\( (r,s,f), (s,t,g) \) \mapsto (r,t,fg)
\]
from a subset of $ \Ec \times \Ec $ to $ \Ec $; namely, a pair $ \(
(s_1,t_1,f_1), (s_2,t_2,f_2) \) $ belongs to the subset iff $ t_1 =
s_2 $. The operation is associative.

If the continuous product of Hilbert spaces corresponds to a
continuous product of probability spaces, then all $ H_{s,t} $ are
embedded into $ H = H_{-\infty,\infty} $, therefore $ \Ec $ is a subset
of $ \R \times \R \times H $. It is a Borel subset. \emph{Sketch of
the proof:} the function $ (s,t,f) \mapsto \dist(f,H_{s,t}) $ is
Borel measurable, since it is continuous unless $ s $ or $ t $ belong
to a finite or countable set of discontinuity points (recall \ref{2c}).

The set $ \Ec $ inherits from $ \R \times \R \times H $ the structure
of a standard measurable space. The domain of the binary operation is
evidently Borel measurable. And the binary operation is (jointly)
Borel measurable. \emph{Sketch of the proof:} the (pointwise) product
$ (f,g) \mapsto fg $ is a continuous map $ L_2(\Om,P) \times
L_2(\Om,P) \to L_1(\Om,P) $.

\smallskip
\textsc{digression: measurable family of hilbert spaces}
\smallskip

Dealing with a Hilbert space that depends on a (non-discrete)
parameter, one should bother about measurability in the parameter. To
this end we choose a single model of an infinite-dimensional separable
Hilbert space, say, the space $ l_2 $ of sequences; and for each $ n
$, a single model of an \dimensional{n} Hilbert space, say, the
space $ l_2^{(n)} $ of $ n $-element sequences. These are our
favourites. Given a standard measurable space $ (X,\X) $, we have a
favourite model $ (l_2)_{x\in X} $ of a family $ (H_x)_{x\in X} $ of
infinite-dimensional separable Hilbert spaces. The disjoint union $
\biguplus_{x\in X} l_2 $, being just $ X \times l_2 $, is a standard
measurable space. More generally, given a measurable function $ n : X \to
\{ 0,1,2,\dots \} \cup \{\infty\} $, we consider $ \( l_2^{(n(x))}
\)_{x\in X} $; here $ l_2^{(\infty)} = l_2 $. Still, $ \biguplus_{x\in
X} l_2^{(n(x))} $ is a standard measurable space; indeed, it is $
\cup_k \( \{ x : n(x)=k \} \times l_2^{(k)} \) $. The general case,
defined below, is the same up to measurable, fiberwise unitary maps.

\begin{definition}
A \emph{standard measurable family of Hilbert spaces}%
\index{standard!measurable family}%
\index{measurable!family of Hilbert spaces}
(over a standard
measurable space $ (X,\X) $) consists of separable Hilbert spaces $
H_x $, given for all $ x \in X $, and a \sif\ on the disjoint union $
\biguplus_{x\in X} H_x = \{ (x,h) : x\in X, h\in H_x \} $ satisfying
the condition:

There exist a measurable function $ n : X \to \{ 0,1,2,\dots \} \cup
\{\infty\} $ and unitary operators $ U_x : l_2^{(n(x))} \to H_x $ (for
all $ x \in X $) such that the map $ (x,h) \mapsto (x,U_x h) $ is a
Borel isomorphism of $ \biguplus_{x\in X} l_2^{(n(x))} $ onto $
\biguplus_{x\in X} H_x $.
\end{definition}

Such a \sif\ on $ \biguplus_{x\in X} H_x $ will be called a
\emph{measurable structure}\index{measurable!structure}
on the family $ (H_x)_{x\in X} $ of Hilbert spaces.

Instead of unitary operators $ U_x $ one may use vectors $ e_k(x) =
U_x e_k $ where $ e_1, e_2, \dots $ are the basis vectors of $ l_2
$. For each $ x $ vectors $ e_k(x) $ are an orthonormal basis of $ H_x
$ provided that $ \dim H_x = \infty $; otherwise the first $ n = \dim
H_x $ vectors are such a basis, and other vectors vanish. Also, $ x
\mapsto \( x, e_k(x) \) $ is a measurable map $ X \to \biguplus_{x\in
X} H_x $ (for each $ k $). These properties ensure that the map $
(x,h) \mapsto (x,U_x h) $ is a Borel measurable bijective map $
\biguplus_{x\in X} l_2^{(n(x))} \to \biguplus_{x\in X} H_x $. The map
is a Borel isomorphism if and only if $ \biguplus_{x\in X} H_x $ is a
\emph{standard} measurable space.

\begin{sloppypar}
Given two standard measurable families of Hilbert spaces $
(H'_x)_{x\in X} $, $ (H''_x)_{x\in X} $ over the same base $ (X,\X) $,
the family of tensor products $ (H'_x \otimes H''_x)_{x\in X} $ is
also a standard measurable family of Hilbert spaces (according to $
U'_x \otimes U''_x : l_2^{(n'(x) n''(x))} = l_2^{(n'(x))} \otimes
l_2^{(n''(x))} \to H'_x \otimes H''_x $).
\end{sloppypar}

\begin{lemma}\label{3a2a}
Let $ h'_x \in H'_x $ and $ h''_x \in H''_x $ be such that $ h'_x
\otimes h''_x $ is measurable in $ x $ (that is, the map $ x \mapsto (
x, h'_x \otimes h''_x ) $ from $ X $ to $ \biguplus_{x\in H} H'_x
\otimes H''_x $ is measurable). Then there exists a function $ c : X
\to \C \setminus \{0\} $ such that both $ c(x) h'_x $ and $ (1/c(x))
h''_x $ are measurable in $ x $.
\end{lemma}

\beginproof
We may assume that $ \| h'_x \| = 1 $ and $ \| h''_x \| = 1 $ (since
the norm is a measurable function of a vector). Also we may assume
that $ H'_x = l_2 $ and $ H''_x = l_2 $ (finite dimensions are left to
the reader). Consider the sphere $ S(l_2) = \{ h \in l_2 : \|h\|=1 \}
$, and the map $ (h_1,h_2) \mapsto h_1 \otimes h_2 $ from $ S(l_2)
\times S(l_2) $ to $ S(l_2 \otimes l_2) $. Inverse image of each point
of $ S(l_2 \otimes l_2) $ is either empty or a \emph{compact} subset
of $ S(l_2) \times S(l_2) $ of the form $ \{ (ch_1,(1/c)h_2) : c \in
\C, \, |c|=1 \} $. There exists a Borel function (`selector') on the
set of factorizing vectors of $ S(l_2 \otimes l_2) $ that chooses a
point from each inverse image. Applying the selector to $ h'_x \otimes
h''_x $ we get $ c(x) h'_x $ and $ (1/c(x)) h''_x $.
\proofend

Assume now that $ (H'_x)_{x\in X} $, $ (H''_x)_{x\in X} $ are just
families (not `measurable'!) of Hilbert spaces, $ H_x = H'_x \otimes
H''_x $, and a measurable structure $ \B $ is given on $ (H_x)_{x\in
X} $. We say that $ \B $ is
\emph{factorizing,}\index{factorizing!measurable structure}
if it results from some measurable
structures $ \B' $, $ \B'' $ on $ (H'_x)_{x\in X} $, $ (H''_x)_{x\in
X} $. Generally, this is not the case. Indeed, a family $ (U_x)_{x\in
X} $ of unitary operators $ U_x \in \U(l_2\otimes l_2) $ in general is
not of the form $ U_x = V_x ( U'_x \otimes U''_x ) $ where $ U'_x,
U''_x \in \U(l_2) $ are arbitrary, but $ V_x \in \U(l_2\otimes l_2) $
is a measurable function of $ x $.

Assume that $ \B $ is factorizing, that is, $ \B $ results from
some $ \B', \B'' $. Does $ \B $ determine $ \B', \B'' $ uniquely? No,
it does not. Indeed, let 
$ c : X \to \C $, $ |c(\cdot)| = 1 $, be a non-measurable
function. The transformation $ (x,h) \mapsto (x,c(x)h) $ of $
\biguplus_{x\in X} H'_x $ sends $ \B' $ to another \sif. Combining it
with the transformation $ (x,h) \mapsto (x,(1/c(x))h) $  of $
\biguplus_{x\in X} H''_x $ we get the trivial transformation of $
\biguplus_{x\in X} H'_x \otimes H''_x $, since $ (c(x) h'_x) \otimes
((1/c(x)) h''_x) = h'_x \otimes h''_x $.

\begin{lemma}\label{3a3}
Let $ \B'_1, \B'_2 $ be two measurable structures on $ (H'_x)_{x\in X}
$ and $ \B''_1, \B''_2 $ --- on
$ (H''_x)_{x\in X} $. Assume that the corresponding structures $ \B_1,
\B_2 $ on $ (H'_x \otimes H''_x)_{x\in X} $ coincide, $ \B_1 = \B_2
$. (Here $ \B_1 $ results from $ \B'_1, \B''_1 $ and $ \B_2 $ --- from
$ \B'_2, \B''_2 $.) Then there exists a function $ c : X \to \C $ such
that $ |c(\cdot)| = 1 $, the map $ (x,h) \mapsto (x,c(x)h) $ sends $
\B'_1 $ to $ \B'_2 $, and the map $ (x,h) \mapsto (x,(1/c(x))h) $
sends $ \B''_1 $ to $ \B''_2 $.
\end{lemma}

\beginproof
If vectors $ \psi_x, \xi_x \in l_2 $ are such that $ \psi_x \otimes
\xi_x $ is a measurable function of $ x $, then $ c(x) \psi_x $ and $
(1/c(x)) \xi_x $ are measurable functions of $ x $ for some choice of
$ c(\cdot) $. Thus, if $ U'_x \otimes U''_x $ is a measurable function
of $ x $ then $ c(x) U'_x $ and $ (1/c(x)) U''_x $ are measurable
functions of $ x $ for some choice of $ c(\cdot) $.
\proofend

\textsc{end of digression}
\smallskip

In the light of these general notions, we return to the continuous product of
spaces $ L_2 $ (equipped with a Borel structure before the digression) and see
that $ \Ec = \biguplus_{s<t} H_{s,t} $ is a standard measurable family of
Hilbert spaces.
\emph{Sketch of the proof:}\footnote{%
 Similarly one can prove a more general fact: the family of \emph{all}
 (closed linear) subspaces of a separable Hilbert space is a standard
 measurable family of Hilbert spaces.}
let $ e_1, e_2, \dots $ span $ H $ and $
e_{s,t}^{(k)} $ be the projection of $ e_k $ to $ H_{s,t} \subset H $,
then $ e_{s,t}^{(1)}, e_{s,t}^{(2)}, \dots $ span $ H_{s,t} $, and $
e_{s,t}^{(k)} $ is measurable in $ s,t $ (being continuous
outside a countable set). Using orthogonalization (for each $ (s,t) $
separately; zero vectors, if any, are skipped) we turn $ e_{s,t}^{(k)}
$ into an orthonormal basis of $ H_{s,t} $.

In terms of basis vectors $ e_{s,t}^{(k)} $, measurability of the
partial binary operation means that its matrix element
\[
\ip{ e_{r,s}^{(k)} e_{s,t}^{(l)} }{ e_{r,t}^{(m)} }
\]
is a Borel measurable function of $ r,s,t \in [-\infty,\infty] $, $
r<s<t $, for any $ k,l,m $.

\mysubsection{Continuous product of Hilbert spaces}
\label{3b}

The measurable structure, introduced in Sect.~\ref{3a} on $
\biguplus_{s<t} L_2(\F_{s,t}) $, exists also on $ \biguplus_{s<t}
H_{s,t} $ in general.

\begin{theorem}\label{3b1}\index{theorem}
For every continuous product of Hilbert spaces $ (H_{s,t})_{s<t} $ (as
defined by \ref{3a1}) there exists a measurable structure on the
family $ (H_{s,t})_{s<t} $ of Hilbert spaces that makes the given map
$ H_{r,s} \otimes H_{s,t} \to H_{r,t} $ Borel measurable in $ r,s,t $.
\end{theorem}

In other words, there exist orthonormal bases $ \( e_{s,t}^{(k)} \)_k
$ in the spaces $ H_{s,t} $ such that $ \ip{ e_{r,s}^{(k)}
e_{s,t}^{(l)} }{ e_{r,t}^{(m)} } $ is Borel measurable in $ r,s,t
$. Recall that $ e_{r,s}^{(k)} e_{s,t}^{(l)} $ is the image of $
e_{r,s}^{(k)} \otimes e_{s,t}^{(l)} $ inder the given map $ H_{r,s}
\otimes H_{s,t} \to H_{r,t} $.

In Sect.~\ref{3a}, spaces $ H_{s,t} $ are both subspaces and factors
of $ H = H_{-\infty,\infty} $; now they are only factors (in the sense
that $ H $ may be treated as $ H_{-\infty,s} \otimes H_{s,t} \otimes
H_{t,\infty} $), which means that a different technique is needed.

\smallskip
\textsc{digression: factors}
\smallskip

The algebra $ \B(l_2 \otimes l_2) $ of all (bounded linear) operators
on the Hilbert space $ l_2 \otimes l_2 $ contains two special
subalgebras, $ \B(l_2) \otimes \One = \{ A \otimes \One : A \in
\B(l_2) \} $ and $ \One \otimes \B(l_2) = \{ \One \otimes A : A \in
\B(l_2) \} $. Recall that $ (A \otimes B) (x \otimes y) = Ax \otimes
By $, thus, $ (A \otimes \One) (x \otimes y) = Ax \otimes y $ and $
(\One \otimes A) (x \otimes y) = x \otimes Ay $. The two subalgebras
are commutants to each other: $ \One \otimes \B(l_2) = \{ A \in \B(l_2
\otimes l_2) : \forall B \in \B(l_2) \otimes \One \;\; AB=BA \, \} $.

A unitary operator $ U \in \U(l_2 \otimes l_2) $ transforms the two
subalgebras in two other subalgebras, $ U ( \B(l_2) \otimes \One )
U^{-1} $ and $ U ( \One \otimes \B(l_2) ) U^{-1} $; still, they are
commutants to each other. Of course, $ U ( \B(l_2) \otimes \One )
U^{-1} = \{ U A U^{-1} : A \in \B(l_2) \otimes \One \} $. If $ U $ is
factorizing, that is, $ U = U_1 U_2 $ for some unitary $ U_1, U_2 \in
\B(l_2) $ then $ U ( \B(l_2) \otimes \One ) U^{-1} = \B(l_2) \otimes
\One $ and $ U ( \One \otimes \B(l_2) ) U^{-1} = \One \otimes \B(l_2)
$. And conversely, these two (mutually equivalent) relations imply
factorizability of $ U $.

The set of all subalgebras $ \A $ of the form $ U ( \B(l_2) \otimes
\One ) U^{-1} $ may be turned into a measurable space as follows. The
ball $ \{ A \in \B(l_2 \otimes l_2) : \| A \| \le 1 \} $ equipped with
the weak operator topology is a metrizable compact topological space,
and $ \{ A \in \A : \| A \| \le 1 \} $ is its closed subset. The set
of all closed subsets of a metrizable compact space is a standard
measurable space, known as Effros space, see
\cite[Sect.~12.C]{Ke}. Thus, each algebra $ \A = U ( \B(l_2) \otimes
\One ) U^{-1} $ may be treated as a point of the Effros space.

The set $ \U (l_2 \otimes l_2) $ of all unitary operators, being
a subset of the ball, is also a measurable space. It is well-known to
be a standard measurable space (and in fact, a non-closed $ G_\de
$-subset of the ball), see \cite[9.B.6]{Ke}.

\begin{lemma}\label{3b2}
(a) The set $ \AAA $ of all subalgebras $ \A $ of the form $ U (
\B(l_2) \otimes \One ) U^{-1} $ is a standard measurable space.

(b) There exists a Borel map $ \A \mapsto U_\A $ from $ \AAA $
to the space of unitary operators on $ l_2 \otimes l_2 $ such that
\[
\A = U_\A ( \B(l_2) \otimes \One ) U_\A^{-1} \quad \text{for all } \A
\in \AAA \, .
\]
\end{lemma}

\beginproof
The group $ G = \U(l_2 \otimes l_2) $ is a Polish group, and
factorizing operators are its closed subgroup $ G_0 = \U(l_2) \times
\U(l_2) $. Left-cosets $ gG_0 = \{ gh : h \in G_0 \} $ (for $ g \in G
$) are a Polish space $ G/G_0 $ \cite[1.2.3]{BK}, and by a theorem of
Dixmier (see \cite[12.17]{Ke} or \cite[1.2.4]{BK}) there exists a
Borel function (`selector') $ s : G/G_0 \to G $ such that $ s (gG_0)
\in gG_0 $ for all $ g $.

The map $ U \mapsto U ( \B(l_2) \otimes \One ) U^{-1} $ is a Borel map
$ G \to \AAA $; indeed, for each $ A \in \B(l_2) $ the map $ U
\mapsto U ( A \otimes \One ) U^{-1} $ is Borel, and $ U ( \B(l_2)
\otimes \One ) U^{-1} $ is the closure of the sequence of $ U ( A_k
\otimes \One ) U^{-1} $ where $ A_k $ are a dense sequence in $
\B(l_2) $. Being constant on each $ gG_0 $, the Borel map $ G \to
\AAA $ leads to a Borel map $ G/G_0 \to \AAA $. The latter
map is bijective, and $ \AAA $ is a part of a standard
measurable space. By a Lusin-Souslin theorem \cite[15.2]{Ke}, $ \AAA
$ is a Borel subset, which proves (a), and the inverse map $ \AAA \to
G/G_0 $ is Borel. The map $ \AAA \to G/G_0 \xrightarrow{s} G $ ensures
(b).
\proofend

\textsc{end of digression}
\smallskip

We return to a continuous product of Hilbert spaces $ (H_{s,t})_{s<t}
$ and assume for simplicity that all $ H_{s,t} $ are
infinite-dimensional. The family $ (H_{-\infty,t} \otimes
H_{t,\infty})_{t\in\R} $ of Hilbert spaces evidently carries a
measurable structure (according to the given unitary operators $
H_{-\infty,t} \otimes H_{t,\infty} \to H_{-\infty,\infty} $). We will
see that the measurable structure is factorizing,\footnote{%
 Similarly one can prove a more general fact: the natural measurable
 structure on $ \biguplus_{\A\in\AAA} H'_\A \otimes H''_\A $ is
 factorizing. Here the disjoint union is taken, roughly speaking,
 over all possible decompositions of $ l_2 $ (or $ l_2 \otimes l_2 $)
 into the tensor product of two infinite-dimensional Hilbert
 spaces. The exact formulation is left to the reader.}
which is close to Theorem \ref{3b1}. Indeed, it means existence of
measurable structures on $ (H_{-\infty,t})_{t\in\R} $ and $
(H_{t,\infty})_{t\in\R} $ that make the given map $ H_{-\infty,t}
\otimes H_{t,\infty} \to H_{-\infty,\infty} $ Borel measurable in $ t
$. Note that such measurable structures on $ (H_{-\infty,t})_{t\in\R}
$ and $ (H_{t,\infty})_{t\in\R} $ are unique up to scalar factors $
(c_t)_{t\in\R} $ according to Lemma \ref{3a3}.

For convenience we let $ H_{-\infty,\infty} = H = l_2 \otimes l_2
$. For any $ t \in \R $ the given unitary operator $ W_t :
H_{-\infty,t} \otimes H_{t,\infty} \to H $ sends $ \B(H_{-\infty,t})
\otimes \One $ to an algebra $ \A_{-\infty,t} \in \AAA $. The function
$ t \mapsto \A_{-\infty,t} $ is increasing ($ s<t $ implies $
\A_{-\infty,s} \subset \A_{-\infty,t} $), therefore Borel measurable
(and in fact, continuous outside a finite or countable set).

Lemma \ref{3b2} gives us unitary operators $ V_t $ on $ H = l_2
\otimes l_2 $ such that $ \A_{-\infty,t} = V_t ( \B(l_2) \otimes \One
) V_t^{-1} $ for all $ t \in \R $, and the map $ t \mapsto V_t $ is
Borel measurable. On the other hand, $ \A_{-\infty,t} = W_t (
\B(H_{-\infty,t})
\otimes \One ) W_t^{-1} $. Thus, $ \B(H_{-\infty,t}) \otimes \One =
(V_t^{-1} W_t)^{-1} ( \B(l_2) \otimes \One ) V_t^{-1} W_t $, which
means that $ (V_t^{-1} W_t)^{-1} $ is a factorizing operator $ l_2
\otimes l_2 \to H_{-\infty,t} \otimes H_{t,\infty} $;
\[
W_t^{-1} V_t = U_{-\infty,t} \otimes U_{t,\infty}
\]
for some unitary operators $ U_{-\infty,t} : l_2 \to H_{-\infty,t} $
and $ U_{t,\infty} : l_2 \to H_{t,\infty} $. Operators $ U_{-\infty,t}
$ define a measurable structure on $ (H_{-\infty,t})_{t\in\R} $. The
same for $ U_{t,\infty} $ and $ (H_{t,\infty})_{t\in\R} $. The 
partial binary operation $ \( (t,x), (t,y) \) \mapsto xy $ becomes
Borel measurable, since $ xy = W_t ( x \otimes y ) = V_t (
U_{-\infty,t}^{-1} x \otimes U_{t,\infty}^{-1} y ) $ and $ V_t $ is
measurable in $ t $.

The proof of Theorem \ref{3b1} is similar. Algebras $ \A_{s,t} \in
\AAA $, corresponding to $ H_{s,t} $, are used. Joint measurability of
$ \A_{s,t} $ in $ s $ and $ t $ follows from the formula $ \A_{s,t} =
\A_{-\infty,t} \cap \A_{s,\infty} $ and a general fact: on a
\emph{compact} metric space, the intersection of two closed subsets is
a jointly Borel measurable function of these two subsets
\cite[27.7]{Ke}.

Non-uniqueness of the measurable structure on $ (H_{s,t})_{s<t} $ is
described by scalar factors $ (c_{s,t})_{s<t} $, $ c_{s,t} \in \C $, $
|c_{s,t}| = 1 $ such that
\[
c_{r,s} c_{s,t} = c_{r,t} \quad \text{whenever } r < s < t
\]
(which means that $ c_{s,t} = c_t / c_s $ for some $ (c_t)_{t\in\R} $;
for example, one may take $ c_t = c_{0,t} $ for $ t>0 $, $ c_t =
1/c_{t,0} $ for $ t<0 $, and $ c_0 = 1 $).
The transformation $ (s,t,h) \mapsto (s,t,c_{s,t}h) $ of $
\biguplus_{s<t} H_{s,t} $ preserves the given maps $ H_{r,s} \otimes
H_{s,t} \to H_{r,t} $ but changes the measurable structure (unless $
c_{s,t} $ is measurable in $ s,t $).

See also \cite[Sect.~1]{Ts02}.

\mysubsection{Stationary case; Arveson systems}
\label{3c}

Let $ (\Om,\F,P) $, $ (\F_{s,t})_{s<t} $, $ (T_h)_{h\in\R} $ be a
noise (as defined by \ref{2c1}), then $ (\F_{s,t})_{s<t} $, being a
continuous product of probability spaces, leads to a continuous
product of Hilbert spaces $ (H_{s,t})_{s<t} $, while each $ T_h $,
being a measure preserving transformation of $ (\Om,\F,P) $, leads to
a unitary operator $ \theta^h : H \to H $ (where $ H =
H_{-\infty,\infty} = L_2(\Om,\F,P) $); namely, $ \theta^h f = f \circ
T_h $ for $ f \in H $. The one-parameter group $ (\theta^h)_{h\in\R}
$, being measurable (in $ h $), is of the form $ \theta^h = \exp(\I hX)
$, where $ X $ (the generator) is a self-adjoint operator. Knowing
that $ T_h $ sends $ \F_{s,t} $ to $ \F_{s+h,t+h} $ we get unitary
operators $ \theta^h_{s,t} : H_{s,t} \to H_{s+h,t+h} $ satisfying $
\theta^h_{r,s} \otimes \theta^h_{s,t} = \theta^h_{r,t} $ and $
\theta^{h_2}_{s+h_1,t+h_1} \theta^{h_1}_{s,t} = \theta^{h_1+h_2}_{s,t}
$.

The property \ref{2c1}(c) ensures that the global algebra $
\A_{-\infty,\infty} = \B(H_{-\infty,\infty}) $ is generated by (the
union of all) local algebras $ \A_{s,t} $, $ -\infty < s < t < \infty
$. See the proof of \ref{6e1}(c \imp a); the same argument works
here. As before, $ \A_{s,t} $ is the image of $ \One \otimes \B(H_{s,t})
\otimes \One $ under the given map $ H_{-\infty,s} \otimes H_{s,t}
\otimes H_{t,\infty} \to H_{-\infty,\infty} = H $; `generated by'
means here `is the closure of' (in the weak operator topology).

\begin{definition}\label{3c2}
A \emph{homogeneous continuous product of Hilbert spaces}%
\index{homogeneous continuous product!of Hilbert spaces}
consists of
a continuous product of Hilbert spaces $ (H_{s,t})_{s<t} $ and unitary
operators $ \theta^h_{s,t} : H_{s,t} \to H_{s+h,t+h} $ (given for all
$ h \in \R $ and $ s,t \in [-\infty,\infty] $, $ s<t $; of course, $
(-\infty) + h = (-\infty) $ and $ (+\infty) + h = (+\infty) $)
satisfying

(a) $ \theta^h_{r,s} \otimes \theta^h_{s,t} = \theta^h_{r,t} $ for $
-\infty \le r < s < t \le \infty $ and $ h \in \R $;

(b) $ \theta^{h_2}_{s+h_1,t+h_1} \theta^{h_1}_{s,t} =
\theta^{h_1+h_2}_{s,t} $ for $ -\infty \le s < t \le \infty $ and $
h_1, h_2 \in \R $;

(c) there exists a self-adjoint operator $ X $ such that $
\theta^h_{-\infty,\infty} = \exp (\I hX) $ for $ h \in \R $;

(d) $ \A_{-\infty,\infty} $ is the weak closure of the union of all $
\A_{s,t} $ for $ -\infty < s < t <\nolinebreak[4] \infty $.
\end{definition}

Every noise leads to a homogeneous continuous product of Hilbert
spaces.

Here are counterparts of Proposition \ref{2c2} and Corollary
\ref{2c3}.

\begin{proposition}\label{3c23}
(Liebscher \cite[Prop.~3.4]{Li}; see also \cite[4.2.1]{Ar}).
Every homogeneous continuous product of Hilbert spaces satisfies the
`upward continuity'\index{upward continuity!for Hilbert spaces}
condition
\begin{equation}\label{3c*}
\A_{s,t} \text{ is generated by } \bigcup_{\eps>0} \A_{s+\eps,t-\eps}
\quad \text{for $ -\infty \le s < t \le \infty $} \, ;
\end{equation}
here $ -\infty+\eps $ is interpreted as $ -1/\eps $ and $ \infty-\eps
$ as $ 1/\eps $.
\end{proposition}

\begin{sloppypar}
\beginproof
Assume $ s,t \in \R $ (other cases, $ s=-\infty $ and $ t=\infty $,
follow via \ref{3c2}(d)).
It is enough to prove that $ \A_{r,s} \A_{s,t} $ is generated by $
\bigcup \A_{r+\eps,t-\eps} $. For any $ A \in \A_{r,s} $ and $ B \in
\A_{s,t} $,
\[
\underbrace{ \E^{\I\eps X} A \E^{-\I\eps X} }_{\in\A_{r+\eps,s+\eps}}
\underbrace{ \E^{-\I\eps X} B \E^{\I\eps X} }_{\in\A_{s-\eps,t-\eps} }
\to AB \quad \text{as } \eps \to 0
\]
weakly and even strongly, since $ \| \E^{\I\eps X} A \E^{-2\I\eps X} B
\E^{\I\eps X} f - ABf \| \le \| A \| \| B \| \| \E^{\I\eps X} f - f \|
+ \| A \| \| \E^{-2\I\eps X} B f - Bf \| + \| \E^{\I\eps X} ABf - ABf
\| \to 0 $.
\proofend
\end{sloppypar}

\begin{corollary}\label{3c24}
(Liebscher \cite[Prop.~3.4]{Li}; see also \cite[4.2.1]{Ar}).
Every homogeneous continuous product of Hilbert spaces satisfies the
`downward continuity'\index{downward continuity!for Hilbert spaces}
condition
\begin{equation}\label{3c**}
\A_{s,t} = \bigcap_{\eps>0} \A_{s-\eps,t+\eps} \quad \text{for all }
s,t \in \R, \, s \le t \, ;
\end{equation}
here $ \A_{t,t} $ is the trivial subalgebra, and $ -\infty-\eps =
-\infty $, $ \infty+\eps = \infty $.
\end{corollary}

\beginproof
Every operator $ A \in \bigcap \A_{s-\eps,t+\eps} $ commutes with $
\bigcup \A_{-\infty,s-\eps} $ and $ \bigcup \A_{t+\eps,\infty} $,
therefore (using Prop.~\ref{3c23}) with $ \A_{-\infty,s} $ and $
\A_{t,\infty} $, which means $ A \in \A_{s,t} $.
\proofend

The two continuity conditions (`upward' and `downward') make sense
also for (non-homogeneous) continuous products of Hilbert
spaces. Still, the upward continuity implies the downward
continuity. (Indeed, the proof of \ref{3c24} does not use the
homogeneity.) Unlike Sect.~\ref{2c}, the converse is true. Systems of
Sect.~\ref{sec:1} do not lead to a counterexample! Especially, for the
system of \ref{1d}, triviality of the limiting \sif\ $
\F_{\infty-,\infty} $ means that the limiting operator algebra $
\A_{\infty-,\infty} $ contains no multiplication operators; but still,
it contains projections to the `superselection sectors' $ H_0, \dots,
H_{m-1} $. See also \ref{6e1}--\ref{6e2}, and \cite{TV}, Lemma 1.5
(and Example 4.3): (mu) $ \Longrightarrow $ (md), but (Hu)\nolinebreak
$ \Longleftrightarrow $ (Hd).

Definition \ref{3c2} may seem to be unsatisfactory, since it does not
stipulate measurability of $ \theta^h_{s,t} $ in $ h $ for finite $
s,t $. Recall however the non-uniqueness of the measurable structure
on $ (H_{s,t})_{s<t} $.

\begin{theorem}\label{3c3}\index{theorem}
For every homogeneous continuous product of Hilbert spaces $
(H_{s,t})_{s<t} $, $ (\theta^h_{s,t})_{s<t;h} $ there exists a
measurable structure on the family $ (H_{s,t})_{s<t} $ of Hilbert
spaces that makes the given map $ H_{r,s} \otimes H_{s,t} \to H_{r,t}
$ Borel measurable in $ r,s,t $, and also makes $ \theta^h_{s,t} $
Borel measurable in $ h,s,t $.
\end{theorem}

\beginproof
By Theorem \ref{3b1} (restricted to $ (s,t) = (-\infty,t) $ or $
(s,\infty) $), there exist unitary $ V_{-\infty,t} : l_2 \to
H_{-\infty,t} $ and $ V_{t,\infty} : l_2 \to H_{t,\infty} $ such that
the unitary operator $ W_t ( V_{-\infty,t} \otimes V_{t,\infty} ) :
l_2 \otimes l_2 \to H $ is a Borel function of $ t $; here, as before,
$ H = H_{-\infty,\infty} $ and $ W_t $ is the given unitary operator $
H_{-\infty,t} \otimes H_{t,\infty} \to H $. The equality $
\theta^h_{-\infty,t} \otimes \theta^h_{t,\infty} =
\theta^h_{-\infty,\infty} $ (a special case of \eqref{3c2}(a)) means in
fact $ W_{t+h} ( \theta^h_{-\infty,t} \otimes \theta^h_{t,\infty} ) = 
\theta^h_{-\infty,\infty} W_t $. We define unitary operators $
\al_{s,t} : l_2 \otimes l_2 \to l_2 \otimes l_2 $ for $ s,t \in \R $
by
\begin{multline*}
\al_{s,t} = V^{-1}_{-\infty,t} \theta^{t-s}_{-\infty,s} V_{-\infty,s}
 \otimes V^{-1}_{t,\infty} \theta^{t-s}_{s,\infty} V_{s,\infty} = \\
= ( V^{-1}_{-\infty,t} \otimes V^{-1}_{t,\infty} ) (
 \theta^{t-s}_{-\infty,s} \otimes \theta^{t-s}_{s,\infty} ) ( V_{-\infty,s}
 \otimes V_{s,\infty} ) = \\
= ( V^{-1}_{-\infty,t} \otimes V^{-1}_{t,\infty} ) W^{-1}_{t}
\theta^{t-s}_{-\infty,\infty} W_s ( V_{-\infty,s} \otimes V_{s,\infty} )
 \, .
\end{multline*}
We see that $ \al_{s,t} $ is a Borel function of $ s $ and $ t $, and
for every $ s,t $ it is a factorizing operator, $ \al_{s,t} =
\be_{s,t} \otimes \ga_{s,t} $ for some unitary $ \be_{s,t}, \ga_{s,t}
: l_2 \to l_2 $. These $ \be_{s,t}, \ga_{s,t} $ are unique up to a
coefficient: $ \al_{s,t} = (c\be_{s,t}) \otimes ( (1/c) \ga_{s,t} ) $,
$ c \in \C $, $ |c| = 1 $. Similarly to the proof of Lemma \ref{3b2},
we use a Borel selector $ G/G_0 \to G $, but for $ G=\U(l_2) \times
\U(l_2) $ and $ G_0 = \{ (c,1/c) : c\in\C, |c|=1 \} $. This way we
make $ \be_{s,t},\ga_{s,t} $ Borel measurable in $ s $ and $ t
$. Also, $ \be_{s,t} = c_{s,t} V^{-1}_{-\infty,t}
\theta^{t-s}_{-\infty,s} V_{-\infty,s} $, $ c_{s,t} \in \C $, $
|c_{s,t}| = 1 $. The product $ c_{r,s} c_{s,t} c_{t,r} $ is Borel
measurable in $ r,s,t $ since $ \be_{t,r} \be_{s,t} \be_{r,s} =
c_{t,r} c_{s,t} c_{r,s} \cdot \One $. Multiplying each $ V_{-\infty,t}
$ by $ 1/c_{0,t} $ we get Borel measurability in $ s,t $ of $
V^{-1}_{-\infty,t} \theta^{t-s}_{-\infty,s} V_{-\infty,s} = (c_{0,s}
c_{s,t} c_{t,0})^{-1} \be_{s,t} $. That is, we get measurable
structures on $ (H_{-\infty,t})_t $ and $ (H_{t,\infty})_t $ that
conform to the shifts. It remains to use the relation $ H_{-\infty,s}
\otimes H_{s,t} = H_{-\infty,t} $; two terms ($ H_{-\infty,s} $ and $
H_{-\infty,t} $) are understood, the third ($ H_{s,t} $) comes out.
\proofend

Waiving the infinite points $ \pm\infty $ on the time axis we get a
\emph{local} homogeneous continuous product of Hilbert spaces.%
\index{local!homogeneous continuous product}
In this
case we may treat $ H_{s,t} $ as a copy of $ H_{0,t-s} $, forget about
shift operators $ \theta^h_{s,t} $, and stipulate unitary operators $
H_{0,s-r} \otimes H_{0,t-s} \to H_{0,t-r} $ instead of $ H_{r,s}
\otimes H_{s,t} \to H_{r,t} $. See also \cite[Prop.~4.1.8]{Za}.

\begin{definition}\label{3c4}
An \emph{algebraic product system of Hilbert spaces}%
\index{algebraic product system}
consists of
separable Hilbert spaces $ H_t $ (given for all $ t \in (0,\infty) $;
possibly finite-dimensio\-nal, but not zero-dimensional), and unitary
operators $ H_s \otimes H_t \to H_{s+t} $ (given for all $ s,t \in
(0,\infty) $), satisfying the associativity condition:
\[
(fg)h = f(gh) \quad \text{for all $ f \in H_r $, $ g \in H_s $, $ h
\in H_t $}
\]
whenever $ r,s,t \in (0,\infty) $. Here $ fg $ stands for the image of
$ f \otimes g $ under the given operator $ H_r \otimes H_s \to H_{r+s}
$.
\end{definition}

All spaces $ H_t $ are infinite-dimensional, unless they all are
one-dimensional; indeed, $ \dim H_{s+t} = \dim H_s \cdot \dim H_t $.

Algebraic product systems are in a natural one-to-one correspondence
with local homogeneous continuous products of Hilbert spaces.

Every noise leads to a homogeneous continuous product of Hilbert
spaces, therefore to a local homogeneous continuous product of
Hilbert spaces, therefore to an algebraic product system of Hilbert
spaces. In particular, every L\'evy process in $ \R $ (or $ \R^n $)
does.

Absence of measurability conditions opens the door to pathologies. An
example is suggested by the pathologic stationary convolution system of
Sect.~\ref{2c}. We start with an isotropic L\'evy process in $ \R^2 $, as in
\ref{counterexample}. Rotating sample paths we get (measure preserving)
automorphisms of the `global' probability space $ (\Om,P) $, as well as
`local' probability spaces $ (\Om_{s,t}, P_{s,t}) $. These automorphisms lead
to unitary operators $ U_{s,t}^\phi $ on $ H_{s,t} = L_2(\Om_{s,t}, P_{s,t})
$; note that
\[
U^\phi_{r,s} \otimes U^\phi_{s,t} = U^\phi_{r,t} \quad \text{and}
\quad U^\phi_{s,t} U^\psi_{s,t} = U^{\phi+\psi}_{s,t} \, .
\]
Being a group of automorphisms of the homogeneous continuous product
of Hilbert spaces, they lead to a group of automorphisms of the
corresponding algebraic product system of Hilbert spaces: $ H_t =
H_{0,t} $; $ U^\phi_t = U^\phi_{0,t} $;
\begin{gather*}
U^\phi_{s+t} (fg) = (U^\phi_s f)(U^\phi_t g) \quad \text{for } f \in
 H_s, \, g \in H_t \, ; \\
U^\phi_t U^\psi_t = U^{\phi+\psi}_t \, .
\end{gather*}
No doubt, $ U^\phi_t $ is a Borel function of $ \phi $ and $ t $. We
spoil the algebraic product system of Hilbert spaces, replacing the
given operators $ W_{s,t} : H_s \otimes H_t \to H_{s+t} $ with
operators $ \ti W_{s,t} $ defined by
\[
\ti W_{s,t} (f \otimes g) = W_{s,t} ( f \otimes U_t^{\phi(s)} g )
\quad \text{for } f \in H_s, \, g \in H_t \, ;
\]
here $ \phi : \R \to \R $ is some non-measurable additive function
(that is, $ \phi(s+t) = \phi(s) + \phi(t) $ for all $ s,t\in\R $). The
associativity condition is still satisfied:
\begin{multline*}
\ti W_{r+s,t} \( \ti W_{r,s} (f \otimes g) \otimes h )
 = f \cdot ( U^{\phi(r)}_s g ) \cdot ( U_t^{\phi(r+s)} h ) = \\
= f \cdot (
 U^{\phi(r)}_s g ) \cdot ( U_t^{\phi(r)} U_t^{\phi(s)} h ) = f \cdot
 U^{\phi(r)}_{s+t} ( g \cdot U_t^{\phi(s)} h )
 = \ti W_{r,s+t} ( f \otimes \ti W_{s,t} (g \otimes h) )
\end{multline*}
for $ f \in H_r $, $ g \in H_s $, $ h \in H_t $; here $ f \cdot g $
means $ W_{r,s} (f \otimes g) $ rather than $ \ti W_{r,s} (f \otimes
g) $.

We will see in Sect.~\ref{3d} that the `spoiled' binary operation is
not Borel measurable, no matter which measurable structure is chosen
on the family $ (H_t)_{t>0} $ of Hilbert spaces.

\begin{definition}\label{3c5}
A \emph{product system of Hilbert spaces,}%
\index{product system}
or \emph{Arveson system,}%
\index{Arveson system}
is a family $ (H_t)_{t>0} $ of Hilbert spaces, equipped with two
structures: first, an algebraic product system of Hilbert spaces, and
second, a standard measurable family of Hilbert spaces, such that the
binary operation $ (f,g) \mapsto fg $ on $ \biguplus_{t>0} H_t $ is
Borel measurable.
\end{definition}

\begin{corollary}\label{3c9}
(From Theorem \ref{3c3}.)
Every homogeneous continuous product of Hilbert spaces leads to an
Arveson system.
\end{corollary}

Existence of a good measurable structure was derived in Theorem
\ref{3c3} from measurability of a unitary group of shifts on the
`global' Hilbert space $ H_{-\infty,\infty} $. Arveson systems in
general seem to need a different idea, since no `global' Hilbert space
is stipulated. Nevertheless the same idea (group of shifts) works,
being combined with another idea: cyclic time.

See also \cite[Chap.~3]{Ar}, \cite[Sect.~3.1 and 7]{Li},
\cite[Sect.~1]{Ts02}.

\mysubsection{Cyclic time; Liebscher's criterion}
\label{3d}

Till now, our time set was $ \R $, or $ [-\infty,\infty] $, or a
subset of $ \R $; in every case it was a linearly ordered set. Now we
want to use the circle $ \T = \R / \Z $ as the time set. It makes no
sense for processes with independent increments (every periodic
process with independent increments on $ \R $ is deterministic), but
it makes sense for convolution systems, flow systems, continuous
products of probability spaces or Hilbert spaces, noises and product
systems.\index{cyclic time}
Definitions \ref{2a1}, \ref{2a3}, \ref{2b1},
\ref{2b6}, \ref{3a1} may be transferred to $ \T $. To this end we just
replace `$ r,s,t \in \R $' (or `$ r,s,t \in [-\infty,\infty] $') with
`$ r,s,t \in \T $' and interprete `$ r<s<t $' according to the cyclic
order on $ \T $. More formally, $ t_1 < \dots < t_n $ means (for $
t_1, \dots, t_n \in \T $) that there exist $ \ti t_1, \dots, \ti t_n
\in \R $ such that $ t_k = \ti t_k \bmod 1 $ for $ k = 1,\dots,n $ and
$ \ti t_1 < \dots < \ti t_n \le \ti t_1 + 1 $. Special cases $ n=2,3,4
$ give us relations $ s<t $, $ r<s<t $, $ r<s<t<u $.

The general (non-homogeneous) case is described by probability spaces
$ (G_{s,t},\linebreak[0]
 \mu_{s,t}) $, \valued{G_{s,t}} random variables $ X_{s,t}
$, sub-\sif s $ \F_{s,t} $, probability spaces $ (\Om_{s,t},P_{s,t}) $
and finally, Hilbert spaces $ H_{s,t} $. The degenerate
case $ \ti t_n = \ti t_1 + 1 $ is allowed, and leads to $ G_{t,t} $,
\dots, $ H_{t,t} $ ($ t \in \T $). Note that the interval from $ t $
to $ t $ is of length $ 1 $ (zero length intervals are excluded by the
strict inequalities $ \ti t_1 < \dots < \ti t_n $); one could prefer
the notation $ G_{t,t+1} $, \dots, $ H_{t,t+1} $ (taking into account
that $ t+1 = t $ in $ \T $). For a flow system $ (X_{s,t})_{s<t} $,
random variables $ X_{0,0} $ and $ X_{t,t} $ are generally different;
$ X_{0,0} = X_{0,t} X_{t,0} $ but $ X_{t,t} = X_{t,0} X_{0,t} $. (Also
$ G_{0,0} $ and $ G_{t,t} $ are generally different.) For \flow{G}s
in a group $ G $ these random variables are conjugate: $ X_{t,t} =
X^{-1}_{0,t} X_{0,0} X_{0,t} $. If $ G $ is commutative then $ X_{t,t}
= X_{0,0} $, but generally $ X_{t,t} \ne X_{0,0} $. Nevertheless $
\F_{0,0} = \F_{t,t} $ (it is the \sif\ generated by the whole flow),
which leads to $ (\Om_{0,0},P_{0,0}) = (\Om_{t,t},P_{t,t}) $ and $
H_{0,0} = H_{t,t} $ where $ H_{s,t} = L_2(\F_{s,t}) = L_2
(\Om_{s,t},P_{s,t}) $. Transferring Definition \ref{2b1} to the time
set $ \T $ we get $ \F_{0,0} = \F_{0,t} \otimes \F_{t,0} = \F_{t,0}
\otimes \F_{0,t} = \F_{t,t} $. Using the approach of Definition
\ref{2b6} we identify $ \Om_{0,0} $ and $ \Om_{t,t} $ according to $ 
\Om_{0,0} = \Om_{0,t} \times \Om_{t,0} = \Om_{t,0}
\times \Om_{0,t} = \Om_{t,t} $. Similarly, when transferring
Definition \ref{3a1} to $ \T $ we identify $ H_{0,0} $ and $ H_{t,t} $
according to $ H_{0,0} = H_{0,t} \otimes H_{t,0} = H_{t,0} \otimes
H_{0,t} = H_{t,t} $. We may denote $ H_{0,0} $ by $ H_\T $ and write $
H_{t,t} = H_\T $ for all $ t \in \T $; similarly, $ \Om_{t,t} = \Om_\T
$ etc. (However, $ X_\T $ makes sense only in commutative semigroups.)

Cyclic-time systems (of various kinds) correspond naturally to
\emph{periodic} linear-time systems. Here `periodic' means, invariant
under the discrete group of time shifts $ t \mapsto t+n $, $ n \in \Z
$.

Homogeneous linear-time systems correspond to homogeneous cyclic-time
systems. Here homogeneity is defined as before (in Definitions
\ref{2c1}, \ref{3c2}) via shifts of the cyclic time set $ \T $.

Given a (linear-time) algebraic product system of Hilbert spaces (or
equivalently, a local homogeneous continuous product of Hilbert
spaces), we may consider the corresponding cyclic-time system. The
latter (in contrast to the former) stipulates the `global' Hilbert
space $ H_\T $, and a group $ (\theta^h_\T)_{h\in\T} $ of unitary
operators on $ H_\T $. In terms of the local homogeneous continuous
product of Hilbert spaces, $ H_\T = H_{0,1} $ and $ \theta^t_\T (fg) =
(\theta^{-t}_{t,1} g) (\theta^{1-t}_{0,t} f) $ for $ f \in H_{0,t} $,
$ g \in H_{t,1} $, $ t \in (0,1) $. In terms of the algebraic product
system of Hilbert spaces, $ H_\T = H_{0,1} $ and $ \theta^t_\T (fg) =
gf $ for $ f \in H_t $, $ g \in H_{1-t} $, $ t \in (0,1) $. 

\begin{theorem}\index{theorem}
(Liebscher \cite[Th.~7]{Li})
A (linear-time) algebraic product system of Hilbert spaces can be
upgraded to an Arveson system if and only if the corresponding
cyclic-time shift operators $ \theta^h_\T $ are a Borel measurable
(therefore continuous) function of $ h \in \T $.
\end{theorem}

Let us apply Liebscher's criterion to the pathologic example of
Sect.~\ref{3c}. We have $ \ti\theta^t_\T \( \ti W_{t,1-t}(f \otimes g) \) =
\ti W_{1-t,t} (g \otimes f) $ for $ f \in H_t $, $ g \in H_{1-t} $, $ t \in
(0,1) $. That is, $ \ti\theta^t_\T ( f \cdot U_{1-t}^{\phi(t)} g ) = g
\cdot U_t^{\phi(1-t)} f $; as before, $ f \cdot g $ means $ W_{t,1-t}
(f \otimes g) $ rather than $ \ti W_{t,1-t} (f \otimes g) $. We have $
\ti\theta^t_\T
( f \cdot g ) = ( U_{1-t}^{\phi(-t)} g ) \cdot ( U_t^{\phi(1-t)} f ) $,
which means that $ \ti \theta^t_\T $ is not a measurable function of $
t $. Indeed, we may take $ f = \exp ( \I X_{0,t}^{(1)} ) $ and $ g =
\exp ( \I X_{0,1-t}^{(1)} ) $; here $ (X_{s,t}^{(1)}, X_{s,t}^{(2)}) $
are the increments of the underlying isotropic two-dimensional L\'evy
process. Then $ f \cdot g = \exp ( \I X_{0,t}^{(1)} ) \exp (
\I X_{t,1}^{(1)} ) = \exp ( \I X_{0,1}^{(1)} ) $ does not depend on $ t
$, but $ \ti\theta^t_\T (f\cdot g) = \exp \( \I ( X_{0,1-t}^{(1)}
\cos\phi(-t) - X_{0,1-t}^{(2)} \sin\phi(-t) ) \) \exp \( \I (
X_{1-t,1}^{(1)} \cos\phi(1-t) - X_{1-t,1}^{(2)} \sin\phi(1-t) ) \)
$. Even in the special case $ \phi(1) = 0 $ we get $ \exp \( \I
( X_{0,1}^{(1)} \cos\phi(t) + X_{0,1}^{(2)} \sin\phi(t) ) \) $,
which is not measurable in $ t $.

\mysection{Classical part of a continuous product}
\label{sec:6}
\mysubsection{Probability spaces: additive flows}
\label{6a}

By the \emph{classical part}%
\index{classical (part of)!continuous product!of probability spaces}
of a continuous product of probability spaces $ \( (\Om,P), \linebreak[0]
(\F_{s,t})_{s<t} \) $ we mean the quotient space
$ (\Om,P) / \F^\stable $ equipped with $ (\F^\stable_{s,t})_{s<t} $
where $ \F^\stable_{s,t} $%
\index{zzf@$ \F^\stable_{s,t} $, sub-\sif}
is the stable part $ \F_{s,t} \cap
\F^\stable $ of $ \F_{s,t} $ transferred to the quotient space. The
classical part is a continuous product of probability spaces; indeed,
\[
\F^\stable_{r,t} = \F^\stable_{r,s} \otimes \F^\stable_{s,t} \quad
\text{for } r<s<t \, ,
\]
since local versions $ U^\rho_{s,t} $%
\index{zzu@$ U^\rho_{s,t} $, operator!for probability spaces}
of the operators $ U^\rho $ satisfy
\[
U^\rho_{r,t} = U^\rho_{r,s} \otimes U^\rho_{s,t} \quad \text{for }
r<s<t \, .
\]

Recall that $ L_2 (\F^\stable) = H_0 \oplus H_1 \oplus H_2 \oplus
\dots $, the chaos spaces $ H_n $ being defined by $ U^\rho f = \rho^n
f $ for $ f \in H_n $; also, $ H_0 $ is the one-dimensional space of
constant functions. Similarly, $ L_2 (\F^\stable_{s,t}) =
\bigoplus_{n<\infty} H_n(s,t) $, $ U^\rho_{s,t} f = \rho^n f $ for $ f
\in H_n(s,t) = H_n \cap L_2(\F_{s,t}) $.%
\index{zzh@$ H_n (s,t) $, subspace!for probability spaces}

\begin{proposition}\label{6a1}
(\cite[2.9]{Ts99})
The following conditions are equivalent for every $ f \in L_2(\Om) $:

(a) $ f \in H_1 $;

(b) $ f = \cE f {\F_{-\infty,t}} + \cE f { \F_{t,\infty}} $ for all $
t \in \R $;

\begin{sloppypar}
(c) $ \cE f {\F_{r,t}} = \cE f {\F_{r,s}} + \cE f {\F_{s,t}} $
whenever $ -\infty \le r < s < t \le \nolinebreak[4] \infty $.
\end{sloppypar}
\end{proposition}

\beginproof
(b) \imp (c): in terms of the projections $ Q_{s,t} : f \mapsto \cE f
{\F_{s,t}} $ we have $ Q_{r,t} f = Q_{r,t} Q_{-\infty,s} f + Q_{r,t}
Q_{s,\infty} f = Q_{r,s} f + Q_{s,t} f $.

(c) \imp (a): $ U^\rho f = \rho f $, since it holds for each element
$ \ti U^\rho_{t_1,\dots,t_n} $ of the net converging to $ U^\rho $
(recall \eqref{5b3a}).

(a) \imp (b): eigenvalues of the operator $ U^\rho =
U^\rho_{-\infty,t} \otimes U^\rho_{t,\infty} $ are products of
eigenvalues, $ \rho^k
\rho^l = \rho^{k+l} $, $ k,l \in \{0,1,2,\dots\} \cup \{\infty\} $. We
have $ H_n = \bigoplus_{k=0}^n H_k (-\infty,t) \otimes H_{n-k}
(t,\infty) $. Especially, $ H_1 = H_0 (-\infty,t) \otimes H_1
(t,\infty) \oplus \linebreak[4]
 H_1 (-\infty,t) \otimes H_0 (t,\infty) = H_1
(-\infty,t) \oplus H_1 (t,\infty) $. Thus, $ f = g+h $ for some $ g
\in H_1 (-\infty,t) $ and $ h \in H_1 (t,\infty) $. However, $
\cE{h}{\F_{-\infty,t}} = \Ex h = 0 $, therefore $ \cE{ f }{
\F_{-\infty,t}} = g $; similarly $ \cE{ f }{ \F_{t,\infty}} = h $ and
we get (b).
\proofend

\begin{corollary}\label{6a1b}
Every square integrable \flow{\R} adapted to a continuous product of
probability spaces is adapted to its classical part.
\end{corollary}

\begin{corollary}\label{6a1a}
A continuous product of probability spaces generated by square
integrable \flow{\R}s is classical.
\end{corollary}

The integrability condition can be removed, see \ref{6a12},
\ref{6a13}.

\begin{theorem}\label{6a2}\index{theorem}
(Tsirelson \cite[Th.~2.12]{Ts99}, \cite[Th.~6a3]{Ts03})
The sub-\sif\ generated by $ H_1 $ is equal to $ \F^\stable $.
\end{theorem}

\beginproof
Clearly, $ \F_1 \subset \F^\stable $ ($ \F_1 $ being generated by $
H_1 $); the other inclusion, $ \F^\stable \subset \F_1 $, follows from
the next lemma.
\proofend

\begin{lemma}\label{6a4a}
The space $ L_2(\F^\stable) $ is the closure of the union of all
subspaces of the form
\[
\bigotimes_{k=0}^n \( H_0(t_k,t_{k+1}) \oplus H_1(t_k,t_{k+1}) \)
\]
where $ -\infty = t_0 < t_1 < \dots < t_n < t_{n+1} = \infty $, $
n=0,1,2,\dots $
\end{lemma}

\beginproof
$ U^\rho $ is the limit of the \emph{decreasing} net of
\emph{commuting} operators $ \ti U^\rho_{t_1,\dots,t_n} $ (recall
\eqref{5b3a}). Therefore for each $ n $ the spectral subspace $ H_0
\oplus \dots \oplus H_n $ of $ U^\rho $ corresponding to the upper
part $ \{ \rho^n,\rho^{n-1},\dots,1 \} $ of its spectrum, is the limit
(that is, the intersection) of the decreasing net of the corresponding
subspaces for $ \ti U^\rho_{t_1,\dots,t_n} $. Similarly, the subspace
$ ( H_n \oplus H_{n+1} \oplus \dots ) \oplus H_\infty $ is the limit
(that is, the closure of the union) of the increasing net of the
corresponding subspaces for $ \ti U^\rho_{t_1,\dots,t_n} $. The latter
subspace, being intersected with $ H_n $, gives a subspace of $
\bigotimes_{k=0}^n \( H_0(t_k,t_{k+1}) \oplus H_1(t_k,t_{k+1}) \) $.
\proofend

By the way, it follows from the lemma above that
\begin{equation}\label{6a4aa}
H_n \text{ is the closed linear span of } \bigcup_{t\in\R} H_1
(-\infty,t) \otimes H_{n-1} (t,\infty)
\end{equation}
for each $ n $.

Each $ f \in H_1 $ leads to an \flow{\R} $ (f_{s,t})_{s<t} $, $
f_{s,t} = \cE f {\F_{s,t}} $, adapted to $ (\F_{s,t})_{s<t} $ in the
sense that $ f_{s,t} $ is \measurable{\F_{s,t}} whenever $ s<t
$. Choosing a sequence $ (f_k)_k $ that spans $ H_1 $ we get the
following.

\begin{corollary}
For every continuous product of probability spaces, its classical part
is generated by (a finite or countable collection of) square
integrable adapted \flow{\R}s.
\end{corollary}

These \flow{\R}s may be combined into a single vector-valued flow,
say, \flow{l_2}. Assuming the downward continuity (recall
\ref{2c3a}) we may use the infinite-dimensio\-nal L\'evy-It\^o theorem
\cite[4.1]{Fe} for representing the \flow{l_2} via a Gaussian
process and (compensated, nonstationary) Poisson processes, those
processes being independent \cite[5.1]{Fe}. In fact, the whole
Poissonian component can be generated by a single \flow{\R}
\cite[6.1]{Fe}, in contrast to the Gaussian component. The framework
of Feldman \cite{Fe} is different from ours, but the difference is
inessential for the classical part, as explained below.

Spaces $ H_1(s,t) $ satisfy the additive relation
\[
H_1(r,t) = H_1(r,s) \oplus H_1(s,t) \quad \text{for } r < s < t \, ,
\]
much simpler than the multiplicative relation $ L_2(\F_{r,t}) =
L_2(\F_{r,s}) \otimes L_2(\F_{s,t}) $. Orthogonal projections $
Q_{s,t} : H_1 \to H_1 $, $ Q_{s,t} H_1 = H_1(s,t) $, lead to a
projection\nobreakdash-\hspace{0pt}valued measure $ (Q_A)_A $; $ Q_A :
H_1 \to H_1 $ for Borel
sets $ A \subset \R $, $ Q_{(s,t]} = Q_{s,t} $ for $ s<t $. To this
end, however, we must assume right-continuity of $ Q_{-\infty,t} $ in
$ t $. Otherwise we should split each point $ t $ of a finite
of countable set in two, $ t_{\text{left}} $ and $ t_{\text{right}}
$. (Alternatively, we could replace the time set $ [-\infty,\infty] $
with an arbitrary, not just connected, compact subset of $ \R $, thus
making $ Q_{-\infty,t} $ continuous in $ t $.) Assume for simplicity
the right-continuity (for a while; the assumption expires before
Prop.~\ref{6a6}). We get (closed linear) subspaces $ H_1(A) = Q_A H_1
\subset H_1 $ satisfying
\begin{equation}\label{6a4b}
\begin{gathered}
H_1 (A \cup B) = H_1(A) \oplus H_1(B) \quad \text{when } A \cap B =
 \emptyset \, , \\
H_1 ( A_1 \cap A_2 \cap \dots ) = H_1(A_1) \cap H_1(A_2) \cap \dots \,
 , \\
H_1 \( (s,t] \) = H_1 (s,t) \quad \text{for } s<t \, .
\end{gathered}
\end{equation}
Defining $ \F^\stable_A $ as the sub-\sif\ generated by $ H_1(A) $ we
get (\cite[6c4]{Ts03})
\begin{gather}
\F_{A\cup B}^\stable = \F_A^\stable \otimes \F_B^\stable \quad
 \text{whenever } A \cap B = \emptyset \, , \label{6a5} \\
A_n \uparrow A \quad \text{implies} \quad \F_{A_n}^\stable \uparrow
 \F_A^\stable \, , \label{6a3} \\
A_n \downarrow A \quad \text{implies} \quad \F_{A_n}^\stable
 \downarrow \F_A^\stable \, , \label{6a4} \\
\F^\stable_{(s,t]} = \F^\stable_{s,t} \quad \text{for } s<t \,
 . \label{6a5a}
\end{gather}
It means that the classical part of any continuous product of
probability spaces is a 
\emph{factored probability space}\index{factored probability space (Feldman)}
as defined
by Feldman \cite[1.1]{Fe}, which cannot be extended beyond the
classical part, see Theorem \ref{11a2}.

\beginproof
\eqref{6a3}: if $ A_n \uparrow A $ then $ H_1(A_n) \uparrow H_1(A) $.

\eqref{6a5}: for each $ f \in H_1 $ the equality $ \Ex \E^{\I f} = \(
\Ex \E^{\I Q_A f} \) \( \Ex \E^{\I Q_{\R\setminus A} f} \) $
holds (by independence) if the set $ A \subset \R $ is an interval or
the union of a finite number of intervals. The monotone class theorem
extends the equality to all Borel sets $ A $. Thus, $ \F_A^\stable $
and $ \F_{\R\setminus A}^\stable $ are independent; \eqref{6a5} follows.

\eqref{6a4} follows from \eqref{6a3} and \eqref{6a5} similarly to
\ref{2c3}.

\eqref{6a5a}: see \ref{6a2}.
\proofend

\begin{proposition}\label{6a6}
The following conditions are equivalent for a \emph{classical}
continuous product of probability spaces:

(a) upward continuity \eqref{2c2a};

(b) downward continuity \eqref{2c3a};

(c) the subspace $ \bigcap_{\eps>0} H_1 (t-\eps, t+\eps) $ is trivial
for every $ t \in \R $.
\end{proposition}

\beginproof
By \eqref{6a3}, (a) is equivalent to triviality of $
H_1(\{s_{\text{right}}\}) $ and $ H_1(\{t_{\text{left}}\}) $ for $ s<t
$. By \eqref{6a4}, (b) is equivalent to triviality of $
H_1(\{s_{\text{left}}\}) $ and $ H_1(\{t_{\text{right}}\}) $ for $ s
\le t $.
Also, (c) is equivalent to triviality of $ H_1(\{t_{\text{left}}\}) $
and $ H_1(\{t_{\text{right}}\}) $ for all $ t $. (At $ \pm\infty $ use
the non-redundancy stipulated by Def.~\ref{2b1}.)
\proofend

\mysubsection{Probability spaces: multiplicative flows}
\label{6b}

We turn to \flow{\T}s; the circle $ \T $ will now be treated as the
complex circle $ \T = \{ z \in \C : |z| = 1 \} $ rather than $ \R / \Z
$. Accordingly, $ L_2 $ spaces over $ \C $ are used. Corollaries
\ref{6a1b}, \ref{6a1a} fail for \flow{\T}s. For counterexamples see
Sect.~\ref{1b}; the singular time point must be finite (not $
\pm\infty $), since the upward continuity at $ \pm\infty $ is ensured
by Def.~\ref{2b1}. Results presented below (\ref{6a10}, \ref{6a12},
\ref{6a13}) are close to \cite[1.7]{TV}.
The time set is $ \R $, but may be enlarged to $ [-\infty,\infty] $.

\begin{proposition}\label{6a10}
If a  continuous product of probability spaces satisfies the upward
continuity condition \eqref{2c2a}, then every \flow{\T}
adapted to the continuous product is adapted to the classical part.
\end{proposition}

\beginproof
Every \flow{\T} $ (X_{s,t})_{s<t} $ satisfies the inequality
\begin{equation}\label{6a11}
\ip{ U^\rho X_{s,t} }{ X_{s,t} } \ge | \Ex X_{s,t} |^{2(1-\rho)} \quad
\text{for $ s<t $ and $ \rho \in [0,1] $} \, ,
\end{equation}
since it holds for each element $ \ti U^\rho_{t_1,\dots,t_n} $ of the
net converging to $ U^\rho $ (recall \eqref{5b3a}): $ \ip{ \ti
U^\rho_{t_1,\dots,t_n} X_{s,t} }{ X_{s,t} } = \prod_{k=0}^n \( \rho
+ (1-\rho) | \Ex X_{t_k,t_{k+1}} |^2 \) $, the logarithm of each
factor being concave in $ \rho $. Stability of $ X_{s,t} $ is thus
ensured, if $ \Ex X_{s,t} \ne 0 $. The latter follows from the upward
continuity: $ | \Ex X_{s-\eps,s+\eps} |^2 = \Ex | \cE{ X_{r,t} }{
\F_{r,s-\eps} \vee \F_{s+\eps,t} } |^2 \to 1 $ as $ \eps \to 0 $; we
cover the compact interval $ [r,t] $ by a finite number of open
intervals $ (s-\eps,s+\eps) $ such that $ \Ex X_{s-\eps,s+\eps} \ne 0
$ and get $ \Ex X_{r,t} \ne 0 $.
\proofend

A stable (that is, adapted to the classical part) \flow{\T} $
(X_{s,t})_{s<t} $ can satisfy $ \Ex X_{s,t} = 0
$ for some $ s<t $; indeed, $ H_1(t-,t+) $ can contain a random
variable $ X_{t-,t+} = \pm 1 $ such that $ \Ex X_{t-,t+} = 0 $. On the
other hand, it must be $ \Ex X_{s,t} \ne 0 $ for some $ s,t $
(irrespective of stability), and moreover, the equivalence relation $
s \sim t \Longleftrightarrow \Ex X_{s,t} \ne 0 $ divides $ \R $ into
at most countable number of intervals (maybe, sometimes degenerate),
since $ L_2(\Om) $ is separable.

\begin{corollary}\label{6a12}
Every \flow{\R} adapted to a continuous product of probability
spaces is adapted to its classical part.
\end{corollary}

\beginproof
\flow{\T}s $ (\E^{\I\la X_{s,t}})_{s<t} $ corresponding to the given
\flow{\R} $ (X_{s,t})_{s<t} $ satisfy $ \Ex \E^{\I\la X_{s,t}} \to 1
$ as $ \la \to 0 $. By \eqref{6a11}, $ \E^{\I\la X_{s,t}} $ is
\measurable{\F^\stable} for all $ \la $ small
enough. Therefore $ X_{s,t} $ is \measurable{\F^\stable}.
\proofend

See also \cite[proof of Th.~1.7]{TV}.

\begin{corollary}\label{6a13}
A continuous product of probability spaces is classical if and only if
it is generated by (a finite or countable collection of) \flow{\R}s.
\end{corollary}

\begin{sloppypar}
The same holds for \flow{(\C,\cdot)}s (valued in the multiplicative
semigroup of complex numbers) if the upward continuity is assumed, see
\cite[1.7]{TV}. Relations between \flow{(\C,\cdot)}s and \flow{(\C,+)}s
described below appear in different forms in \cite[Appendix A]{TV} and earlier
works cited there.
\end{sloppypar}

\begin{proposition}\label{6a17}
A stable square integrable \flow{(\C,\cdot)} $ (X_{s,t})_{s<t} $ is
uniquely determined by the projections $ Q_0 X_{s,t} = \Ex X_{s,t} $
and $ Q_1 X_{s,t} $ of each $ X_{s,t} $ to $ H_0(s,t) $ and $ H_1(s,t)
$.
\end{proposition}

\begin{sloppypar}
\beginproof
The projection $ \prod_{k=0}^n (Q_0+Q_1) X_{t_k,t_{k+1}} $ of $
X_{t_0,t_{n+1}} $ to $ \bigotimes_{k=0}^n \( H_0(t_k,t_{k+1}) \oplus
H_1(t_k,t_{k+1}) \) $ is uniquely determined. It remains to use Lemma
\ref{6a4a}.
\proofend
\end{sloppypar}

\begin{sloppypar}
Clearly, $ \Ex X_{r,t} = (\Ex X_{r,s}) (\Ex X_{s,t}) $ and $ Q_1
X_{r,t} = (Q_1 X_{r,s}) (\Ex X_{s,t}) + (\Ex X_{r,s}) (Q_1 X_{s,t})
$. In particular, if $ \Ex X_{s,t} = 1 $ for all $ s<t $, then $ Q_1
X_{r,t} = Q_1 X_{r,s} + Q_1 X_{s,t} $, that is, $ (Q_1 X_{s,t})_{s<t}
$ is a \flow{(\C,+)}.
\end{sloppypar}

\begin{sloppypar}
\begin{proposition}\label{6a18}
For every square integrable, zero-mean \flow{(\C,+)} $
(Y_{s,t})_{s<t} $ there exists a square integrable \flow{(\C,\cdot)}
$ (X_{s,t})_{s<t} $ such that
\[
\Ex X_{s,t} = 1 \quad \text{and} \quad Q_1 X_{s,t} = Y_{s,t} \quad
\text{for } s<t \, .
\]
\end{proposition}
\end{sloppypar}

\beginproof
Similarly to the proof of \ref{6a17} we calculate the projection of the
desired $ X_{s,t} $ to subspaces of the form $ \bigotimes_{k=0}^n \(
H_0(t_k,t_{k+1}) \oplus H_1(t_k,t_{k+1}) \) $. The subspaces are an
increasing net. The projections are consistent, and bounded in $ L_2
$:
\begin{multline*}
\Big\| \bigotimes_k (1+Y_{t_k,t_{k+1}}) \Big\|^2 = \prod_k
 (1+\|Y_{t_k,t_{k+1}}\|^2) \le \\
\le \exp \Big( \sum_k \|Y_{t_k,t_{k+1}}\|^2
 \Big) = \exp \( \| Y_{s,t} \|^2 \) \, .
\end{multline*}
Thus, they are a net converging in $ L_2 $; its limit is the desired $
X_{s,t} $.
\proofend

The relation between the flows $ X $ and $ Y $ as in \ref{6a18} will
be denoted by
\begin{equation}\label{6a19}
X = \Exp Y \, ; \qquad Y = \Log X \, .
\end{equation}%
\index{zze@$ \Exp $ (over cont.\ products)}%
\index{zzl@$ \Log $ (over cont.\ products)}
It is a one-to-one correspondence between (the set of all) square
integrable \flow{(\C,\cdot)}s $ (X_{s,t})_{s<t} $ satisfying $ \Ex
X_{s,t} = 1 $ for $ s<t $ (which implies stability by \eqref{6a11}),
and (the set of all) square integrable \flow{(\C,+)}s $
(Y_{s,t})_{s<t} $ satisfying $ \Ex Y_{s,t} = 0 $ for $ s<t $ (these
are stable by \ref{6a1b}).

Relations \eqref{6a19} do not mean that $ X_{s,t} = \exp(Y_{s,t})
$. In fact, if $ Y $ is sample continuous, therefore Gaussian, then
$ X_{s,t} = \exp(Y_{s,t}-\frac12 \|Y_{s,t}\|^2) $ for $ s<t $.
Of course, `$ \exp $' is the usual exponential function $ \C \to \C
$, while `$ \Exp $' is introduced by \eqref{6a19}. If the
limit $ Y_{-\infty,\infty} = \lim_{s\to-\infty,t\to\infty} Y_{s,t} $
exists (in $ L_2 $) then the limit $ X_{-\infty,\infty} $ exists, and
we may write
\[
X_{-\infty,\infty} = \Exp (Y_{-\infty,\infty}) \, ; \qquad \Exp : H_1
\to L_2(\F^\stable) \, .
\]

\begin{proposition}\label{6a20}
The space $ L_2(\F^\stable) $ is the closed linear span of $ \{ \Exp f
: f \in H_1 \} $.
\end{proposition}

\beginproof
The projection of $ \Exp f $ to $ H_0 \oplus H_1 $ is $ 1+f $, thus
(using the upper bound from the proof of \ref{6a18})
\[
\| (\Exp f) - (1+f) \|^2 = \| \Exp f \|^2 - \| 1+f \|^2 \le ( \exp
\|f\|^2 ) - (1+\|f\|^2) \, .
\]
It follows that
\[
f = \lim_{\eps\to0} \frac1{\eps} \( \Exp(\eps f) - 1 \) \quad \text{in
} L_2 \, ,
\]
and we see that $ H_1 \subset E $, where $ E $ is the closed linear
span of $ \{ \Exp f : f \in H_1 \} $. Similarly, $ H_1(s,t) \subset
E(s,t) $, where $ E(s,t) $ is the closed linear span of $ \{ \Exp f
: f \in H_1(s,t) \} $. However, $ E(r,t) \supset E(r,s) \otimes E(s,t)
$ for $ r<s<t $. We see that $ E $ contains each subspace of the form
$ \bigotimes_{k=0}^n \( H_0(t_k,t_{k+1}) \oplus H_1(t_k,t_{k+1}) \) $;
it remains to use \ref{6a4a}.
\proofend

\begin{proposition}\label{6a21}
If $ Y_{t-,t+} = 0 $ for all $ t $, then $ \| X_{s,t} \|^2 = \exp ( \|
Y_{s,t} \|^2 ) $ for $ s<t $.
\end{proposition}

\beginproof
Recall the proof of \ref{6a18} and note that $ \exp \( \|
Y_{t_k,t_{k+1}} \|^2 \) = 1 + \| Y_{t_k,t_{k+1}} \|^2 + o \( \|
Y_{t_k,t_{k+1}} \|^2 \) $.
\proofend

\begin{proposition}\label{6b10}
If a classical continuous product of probability spaces satisfies the
equivalent continuity conditions \ref{6a6}(a--c), then the map $ \Exp
: H_1 \to L_2(\Om) $ has the property
\[
\ip{ \Exp f }{ \Exp g } = \exp \ip f g \quad \text{for } f,g \in H_1
\, .
\]
\end{proposition}

\begin{sloppypar}
\beginproof
Similarly to \ref{6a21}, $ \ip{ 1+f_{t_k,t_{k+1}} }{ 1+g_{t_k,t_{k+1}} }
= 1 + \ip{ f_{t_k,t_{k+1}} }{ g_{t_k,t_{k+1}} } \approx \exp \ip{
f_{t_k,t_{k+1}} }{ g_{t_k,t_{k+1}} } $.
\proofend
\end{sloppypar}

It means that $ L_2(\Om) $ is nothing but the Fock space%
\index{Fock space}
$ \E^{H_1} $,
see \cite[Sect.~2.1.1, especially (2.7)]{Ar}. More generally,
\[
L_2 (\F^\stable) = \E^{H_1}
\]
for all continuous products of probability spaces satisfying the
downward continuity condition.

\mysubsection{Noises}
\label{6c}

Given a noise $ (\F_{s,t})_{s<t} $, $ (T_h)_h $, we may consider the
classical part of the continuous product of probability spaces $
(\F_{s,t})_{s<t} $. It consists of sub-\sif s $ \F_{s,t}^\stable =
\F_{s,t} \cap \F^\stable $ (transferred to the quotient space $
(\Om,P) / \F^\stable $, which does not matter now). Time shifts $ T_h
$ leave $ \F^\stable $ invariant (since operators $ U^\rho $ evidently
commute with time shifts), thus $ T_h $ sends $ \F_{s,t}^\stable $ to
$ \F_{s+h,t+h}^\stable $. It means that the classical part of a noise
is a (classical) noise.%
\index{classical (part of)!noise}

A classical noise is generated by \flow{\R}s, like any other
classical continuous product of probability spaces (see
\ref{6a13}). However, we want these \flow{\R}s $ (X_{s,t})_{s<t} $
to be \emph{stationary} in the sense that $ X_{s+h,t+h} = X_{s,t}
\circ T_h $ (where $ T_h : \Om \to \Om $ are time shifts). The
following result is proven in \cite[2.9]{Ts98} under assumptions
excluding the Poisson component, but the argument works in general.

\begin{theorem}\label{6c1}\index{theorem}
Every classical noise is generated by (a finite or countable
collection of) square integrable stationary \flow{\R}s.  
\end{theorem}

\beginproof
Time shifts $ T_h $ on $ \Om $ induce unitary operators $ U_h $ on $
L_2(\Om) $, commuting with $ U^\rho $ and therefore leaving invariant
the first chaos space $ H_1 $; we will treat $ U_h $ as operators $
H_1 \to H_1 $. They are connected with the projections $ Q_{s,t} : H_1
\to H_1 $ by the relation $ U_h^{-1} Q_{s,t} U_h = Q_{s+h,t+h}
$. Integrating the function $ t \mapsto \E^{\I\la t} $ by the
projection-valued measure $ (Q_A)_A $ (recall \eqref{6a4b}) we get
unitary operators $ V_\la
: H_1 \to H_1 $ satisfying Weyl relations $ U_h V_\la = \E^{\I\la h}
V_\la U_h $. By the well-known theorem of von Neumann (see
\cite[Th.~VIII.14]{RS}), $ H_1 $ decomposes into the direct sum of a
finite or countable number of irreducible components, --- subspaces,
each carrying an irreducible representation of Weyl relations. Each
irreducible representation is unitarily equivalent to the standard
representation in $ L_2(\R) $, where $ U_h $ acts as the shift by $ h
$, and $ V_\la $ acts as the multiplication by $ t \mapsto \E^{\I\la
t} $. Defining $ X^{(k)}_{s,t} $ as the vector that corresponds to the
indicator function of the interval $ (s,t) $ in the $ k $-th
irreducible component of $ H_1 $ we get the needed stationary
\flow{\R}s $ (X^{(k)}_{s,t})_{s<t} $.
\proofend

\begin{corollary}\label{6c2}
A noise is classical if and only if it is generated by (a finite or
countable collection of) stationary \flow{\R}s.
\end{corollary}

These \flow{\R}s may be combined into a single vector-valued flow
(\flow{\R^n} or \flow{l_2}). The infinite-dimensional L\'evy-It\^o 
theorem \cite[4.1]{Fe} may be used for representing the stationary
\flow{l_2} via Brownian motions and (stationary, compensated) Poisson
processes.

We may treat $ H_1 $ as the tensor product, $ H_1 = L_2(\R) \otimes
\cH = L_2(\R,\cH) $, where $ L_2(\R) $ carries the standard
representation of Weyl relations, and $ \cH $ is the
Hilbert space of all square integrable, zero mean, stationary
\flow{\R}s. Further, the space $ \cH $ decomposes in two orthogonal
subspaces, the Brownian part and the Poissonian part. The (finite or
infinite) dimension of the Brownian part is the maximal number of
independent Brownian motions adapted to the given classical noise. The
Poissonian part may be identified with the $ L_2 $ space over the
corresponding L\'evy-Khinchin measure.

\mysubsection{Pointed Hilbert spaces}
\label{6d}

\begin{definition}
Let $ (H_{s,t})_{s<t} $ be a continuous product of Hilbert spaces. A
vector $ f \in H_{r,t} $ is \emph{decomposable,}%
\index{decomposable!vector}
if $ f \ne 0 $ and
for every $ s \in (r,t) $ there exist $ g \in H_{r,s} $ and $ h \in
H_{s,t} $ such that $ f = gh $.
\end{definition}

(As before, $ gh $ is the image of $ g \otimes h $ under the given
unitary operator $ H_{r,s} \otimes H_{s,t} \to H_{r,t} $.)

\begin{lemma}
If $ g \in H_{r,s} $ and $ h \in H_{s,t} $ are such that the vector $
gh \in H_{r,t} $ is decomposable then $ g $ and $ h $ are
decomposable.
\end{lemma}

See also \cite[6.0.2 and 6.2.1]{Ar}.

\beginproof
We may assume $ \|g\| = 1 $, $ \|h\| = 1 $. Consider the
one-dimensional orthogonal projection $ Q_g : H_{r,s} \to H_{r,s} $, $
Q_g \psi = \ip{\psi}g g $. Note that $ \ip{ (Q_g \otimes \One_{s,t}) f
}{ f } = \ip{ Q_g g }{ g } \ip h h = 1 $. Let $ s' \in (r,s) $. Then $
f = g' h' $ for some vectors $ g' \in H_{r,s'} $ and $ h \in H_{s',t}
$ (of norm $ 1 $). As before, $ \ip{ ( Q_{g'} \otimes \One_{s',t} ) f
}{ f } = 1 $. Therefore
\[
\underbrace{ \ip{ ( Q_{g'} \otimes \One_{s',s} \otimes \One_{s,t} ) f
}{ f } }_{=1} = \ip{ ( Q_{g'} \otimes \One_{s',s} ) g }{ g }
\underbrace{ \ip{ \One_{s,t} h }{ h } }_{=1} \, .
\]
The equality $ \ip{ ( Q_{g'} \otimes \One_{s',s} ) g }{ g } = 1 $
means that $ g = g' \psi $ for some $ \psi \in H_{s',s} $. Thus,
$ g $ is decomposable.
\proofend

Theorem \ref{3b1} gives us a measurable structure on the family $
(H_{s,t})_{s<t} $ of Hilbert spaces. The structure is non-unique, but
we can adapt factor-vectors to any given structure, as stated below.
See also \cite[Corollary 5.2]{Li}.

\begin{proposition}\label{6d3}
For every measurable structure as in Theorem \ref{3b1} and every
decomposable vector $ f \in H_{-\infty,\infty} $ there exists a family
$ (f_{s,t})_{s<t} $ of vectors $ f_{s,t} \in H_{s,t} $ (given for all
$ s,t $ such that $ s < t $) satisfying the conditions
\begin{gather*}
f_{r,s} f_{s,t} = f_{r,t} \quad \text{whenever } r < s < t \, , \\
f_{s,t} \quad \text{is measurable in $ s,t $} \, .
\end{gather*}
\end{proposition}

\beginproof
For every $ t \in \R $ we choose $ g_t \in H_{-\infty,t} $ and $ h_t
\in H_{t,\infty} $ such that $ g_t h_t = f $. Lemma \ref{3a2a} gives
us complex numbers $ c_t $ such that the vectors $ f_{-\infty,t} = c_t
g_t $ and $ f_{t,\infty} = (1/c_t) h_t $ are measurable in $ t $. Now
$ f_{s,t} $ are uniquely determined by requiring $ f_{-\infty,s}
f_{s,t} f_{t,\infty} = \nolinebreak[4] f $.
\proofend

Applying \ref{6d3} to a continuous product of spaces $ L_2 $ (recall
Sect.~\ref{3a}) we may get the following.

\begin{corollary}\label{6c4}
Every square integrable \flow{(\C,\cdot)} $ (X_{s,t})_{s<t} $
adapted to a continuous product of probability spaces can be written
as $ X_{s,t} = c_{s,t} Y_{s,t} $ where $ c_{s,t} $ are complex numbers
satisfying $ c_{r,s} c_{s,t} = c_{r,t} $ for $ r<s<t $, and $
(Y_{s,t})_{s<t} $ is a \flow{(\C,\cdot)} such that the map $ (s,t)
\mapsto Y_{s,t} $ from $ \{ (s,t) \in \R^2 : s<t \} $ to $ L_2(\Om) $
is Borel measurable.
\end{corollary}

\begin{corollary}\label{6c5}
If a square integrable \flow{(\C,\cdot)} $ (X_{s,t})_{s<t} $
satisfies $ \Ex X_{s,t} = 1 $ whenever $ s<t $, then the map $ (s,t)
\mapsto X_{s,t} $ from $ \{ (s,t) \in \R^2 : s<t \} $ to $ L_2(\Om) $
is Borel measurable.
\end{corollary}

One may prove \ref{6c5} via \ref{6c4} or, alternatively, via
\eqref{6a19} and the proof of \ref{6a17}.

Decomposable vectors need not exist in a continuous product of Hilbert
spaces in general, but they surely exist in every continuous product
of spaces $ L_2 $, since constant functions are decomposable vectors.

\begin{definition}\label{6d6}
A \emph{continuous product of pointed Hilbert spaces}%
\index{continuous product!of pointed Hilbert spaces}%
\index{pointed Hilbert spaces}
consists of a
continuous product of Hilbert spaces $ (H_{s,t})_{s<t} $ and vectors $
u_{s,t} \in H_{s,t} $ such that
\begin{gather*}
u_{r,s} u_{s,t} = u_{r,t} \quad \text{whenever } -\infty\le
 r<s<t\le\infty \, , \\
\| u_{s,t} \| = 1 \quad \text{whenever } -\infty\le s<t\le\infty \, .
\end{gather*}
\end{definition}

Such family $ (u_{s,t})_{s<t} $ will be called a \emph{unit}%
\index{unit!of a cont.\ prod.\ of Hilbert spaces}
(of the
given continuous product of Hilbert spaces). It is basically the same
as a decomposable vector of $ H = H_{-\infty,\infty} $.
Each $ H_{s,t} $ may be identified with a subspace of $ H $, namely
(the image of) $ u_{-\infty,s} \otimes H_{s,t} \otimes u_{t,\infty} $.

Waiving the infinite points $ \pm\infty $ on the time axis we get a
\emph{local}%
\index{local!cont.\ prod.\ of (pointed) Hilbert spaces}
continuous product of pointed Hilbert spaces. The
embeddings $ H_{-1,1} \subset H_{-2,2} \subset \dots $ may be used for
enlarging the time set $ \R $ to $ [-\infty,\infty] $; to this end $
H_{-\infty,\infty} $ is constructed as the completion of the union of
$ H_{-n,n} $ (see also \eqref{6d13}). In this respect (and many
others), continuous products of pointed Hilbert spaces are closer to
continuous products of probability spaces than Hilbert spaces.

Every continuous product of spaces $ L_2 $ is naturally a continuous
product of pointed Hilbert spaces, $ u_{s,t} $ being the function
that equals to $ 1 $ everywhere.

\begin{question}\label{6d6a}\index{question}
Does every continuous product of pointed Hilbert spaces emerge from
some continuous product of probability spaces (that is, is isomorphic
to some continuous product of $ L_2 $ spaces with `probabilistic'
units)?
\end{question}

Many results and arguments of Sections \ref{5a}, \ref{5b}, \ref{6a},
\ref{6b} may be generalized to continuous products of pointed Hilbert
spaces.

\begin{definition}\label{6d7}
(a)
Let $ (H_{s,t}^{(1)})_{s<t} $ and $ (H_{s,t}^{(2)})_{s<t} $ be two
continuous products of Hilbert spaces. An \emph{embedding}%
\index{embedding!of continuous products!of Hilbert spaces}\footnote{%
 Not `morphism' for not contradicting \cite[3.7.1]{Ar}.}
of the first product to the second is a family $ (\al_{s,t})_{s<t} $
of isometric linear embeddings $ \al_{s,t} : H_{s,t}^{(1)} \to
H_{s,t}^{(2)} $ such that
\[
(\al_{r,s} f) (\al_{s,t} g) = \al_{r,t} (fg) \quad \text{for } f \in
H_{r,s}^{(1)}, \, g \in H_{s,t}^{(1)}
\]
(as before, $ fg $ is the image of $ f \otimes g $ in $ H_{r,t} $).

(b)
Let $ (H_{s,t}^{(1)},u_{s,t}^{(1)})_{s<t} $ and $
(H_{s,t}^{(2)},u_{s,t}^{(2)})_{s<t} $ be two continuous products of
pointed Hilbert spaces. An \emph{embedding}%
\index{embedding!of continuous products!of pointed Hilbert spaces}
of the first product to the second is an embedding $ (\al_{s,t})_{s<t}
$ of $ (H_{s,t}^{(1)})_{s<t} $ to $ (H_{s,t}^{(2)})_{s<t} $ as in
(a) satisfying the additional condition
\[
\al_{s,t} u_{s,t}^{(1)} = u_{s,t}^{(2)} \quad \text{for } s<t \, .
\]
\end{definition}

If $ \al_{s,t} (H_{s,t}^{(1)}) $ is the whole $ H_{s,t}^{(2)} $ for $
s<t $, then $ (\al_{s,t})_{s<t} $ is an \emph{isomorphism}.%
\index{isomorphism!of continuous products!of Hilbert spaces}

Every morphism between continuous products of probability spaces leads
to an embedding of the corresponding continuous product of pointed
Hilbert spaces (in the opposite direction). See Examples \ref{5a2},
\ref{5a4}.

\begin{definition}
A \emph{joining}%
\index{joining!of continuous products!of pointed Hilbert spaces}
(or \emph{coupling}) of two continuous products of pointed Hilbert
spaces $ (H_{s,t}^{(1)},u_{s,t}^{(1)})_{s<t} $ and $
(H_{s,t}^{(2)},u_{s,t}^{(2)})_{s<t} $ consists of a third
continuous product of pointed Hilbert spaces $ (H_{s,t},u_{s,t})_{s<t}
$ and two embeddings $ (\al_{s,t})_{s<t} $, $ (\be_{s,t})_{s<t} $, $
\al_{s,t} : H_{s,t}^{(1)} \to H_{s,t} $, $ \be_{s,t} : H_{s,t}^{(2)}
\to H_{s,t} $ of these products such that $ H_{-\infty,\infty} $ is
the closed linear span of $ \al_{-\infty,\infty}
(H_{-\infty,\infty}^{(1)}) \cup \be_{-\infty,\infty}
(H_{-\infty,\infty}^{(2)}) $.
\end{definition}

Each joining leads to bilinear forms $ (f,g) \mapsto \ip{ \al_{s,t}f
}{ \be_{s,t}g } $ for $ f \in H_{s,t}^{(1)} $, $ g \in H_{s,t}^{(2)}
$. Two joinings that lead to the same bilinear form will be
called isomorphic.%
\index{isomorphic joinings!of continuous products!of pointed Hilbert spaces}
A joining with itself will be called a \emph{self-joining.}%
\index{self-joining!of continuous products!of Hilbert spaces}
A \emph{symmetric}%
\index{self-joining!symmetric}\index{symmetric!self-joining}
self-joining is a self-joining $ (\al,\be) $ isomorphic to $ (\be,\al)
$.

Every joining of two continuous products of probability spaces leads
to a joining of the corresponding continuous products of pointed
Hilbert spaces. The same holds for self-joinings and symmetric
self-joinings.

Every joining $ (\al,\be) $ of two continuous products of pointed
Hilbert spaces has its \emph{maximal correlation}%
\index{maximal correlation!for Hilbert spaces}
\begin{gather*}
\rho^{\max} (\al,\be) = \rho_{-\infty,\infty}^{\max} (\al,\be) \, , \\
\rho_{s,t}^{\max} (\al,\be) = \sup | \ip{ \al_{s,t}f }{ \be_{s,t}g } |
\, ,
\end{gather*}
where the supremum is taken over all $ f \in H_{s,t}^{(1)} $, $ g \in
H_{s,t}^{(2)} $ such that $ \| f \| \le 1 $, $ \| g \| \le 1 $, $ \ip
f {u_{s,t}^{(1)}} = 0 $, $ \ip g {u_{s,t}^{(2)}} = 0 $.

The maximal correlation defined in Sect.~\ref{5a} for a joining of
continuous products of probability spaces is equal to the maximal
correlation of the corresponding joining of continuous products of
pointed Hilbert spaces.

\begin{proposition}\label{6d9}
$ \rho^{\max}_{r,t} (\al,\be) = \max \( \rho^{\max}_{r,s} (\al,\be), \,
\rho^{\max}_{s,t} (\al,\be) \) $ whenever $ r<s<t $.
\end{proposition}

\begin{sloppypar}
\beginproof
Similar to \ref{5a7}; each $ H_{s,t} $ decomposes into the
one-dimensional subspace spanned by $ u_{s,t} $ and its orthogonal
complement $ H_{s,t}^0 $.
\proofend
\end{sloppypar}

\begin{proposition}
For every continuous product of pointed Hilbert spaces and every $
\rho \in [0,1] $ there exists a symmetric self-joining $
(\al_\rho,\be_\rho) $%
\index{zza@$ (\al_\rho, \be_\rho) $, self-joining!for pointed Hilbert spaces}
of the given product such that
\[
\rho^{\max} (\al_\rho,\be_\rho) \le \rho
\]
and
\[
\ip{ \al_{s,t} f }{ \be_{s,t} f } \le \ip{ (\al_\rho)_{s,t} f }{
(\be_\rho)_{s,t} f }
\]
for all $ s<t $, $ f \in H_{s,t} $ and all self-joinings $ (\al,\be) $
satisfying $ \rho^{\max} (\al,\be) \le \rho $.

The self-joining $ (\al_\rho,\be_\rho) $ is unique up to isomorphism.
\end{proposition}

\beginproof
Similar to \ref{5a8} but simpler; we just take the limit of the
decreasing net of (commuting) Hermitian operators (or their quadratic
forms)
\begin{multline*}
\ti U^\rho_{t_1,\dots,t_n} f = \bigotimes_{k=0}^n \( \rho f_k +
(1-\rho) \ip{ f_k }{ u_{t_k,t_{k+1}} } u_{t_k,t_{k+1}} \) \quad
 \text{for } f = f_0 \otimes \dots \otimes f_n, \\
f_0 \in H_{-\infty,t_1}, \, f_1 \in H_{t_1,t_2}, \, \dots,
 \, f_n \in H_{t_n,\infty} \, . \qquad \qed
\end{multline*}
\proofendnoqed

Operators $ U^\rho_{s,t} $%
\index{zzu@$ U^\rho_{s,t} $, operator!for pointed Hilbert spaces}
satisfy
\[
U^\rho_{r,s} \otimes U^\rho_{s,t} = U^\rho_{r,t} \, , \quad
U^{\rho_1}_{s,t} U^{\rho_2}_{s,t} = U^{\rho_1 \rho_2}_{s,t} \, ,
\]
and the spectrum of $ U^\rho_{s,t} $ is contained in $ \{ 1, \rho,
\rho^2, \dots \} \cup \{ 0 \} $. Similarly to \ref{5a11} we have the
following.

\begin{proposition}\label{6d11}
For every continuous product of pointed Hilbert spaces there exist
(closed linear) subspaces $ H_0, H_1, H_2, \dots $%
\index{zzh@$ H_n $, subspace!for pointed Hilbert spaces}
and $ H_\infty $ of
$ H = H_{-\infty,\infty} $ such that 
\begin{gather*}
H = ( H_0 \oplus H_1 \oplus H_2 \oplus \dots ) \oplus H_\infty \, , \\
U^\rho f = \rho^n f \quad \text{for $ f \in H_n $, $ \rho \in [0,1] $}
 \, , \\
U^\rho f = 0 \quad \text{for $ f \in H_\infty $, $ \rho \in [0,1) $}
 \, .
\end{gather*}
\end{proposition}

The space $ H_0 $ is one-dimensional, spanned by $ u_{-\infty,\infty}
$. Similarly we introduce subspaces $ H_n (s,t) $.%
\index{zzh@$ H_n (s,t) $, subspace!for pointed Hilbert spaces}
The relation $
U^\rho_{r,s} \otimes U^\rho_{s,t} = U^\rho_{r,t} $ implies
\begin{equation}\label{6d11a}
H_n (r,t) = \bigoplus_{k=0}^n H_k (r,s) \otimes H_{n-k} (s,t) \, .
\end{equation}

We recall the embedding of each $ H_{s,t} $ into $ H =
H_{-\infty,\infty} $ by $ f \mapsto u_{-\infty,s} f u_{t,\infty} $ for
$ f \in H_{s,t} $, and introduce for $ s<t $ the orthogonal projection
$ Q_{s,t} $ of $ H $ onto $ H_{s,t} \subset H $; clearly,
\[
Q_{s,t} ( fgh ) = \ip f {u_{-\infty,s}} u_{-\infty,s} g \ip h
{u_{t,\infty}} u_{t,\infty}
\]
for $ f \in H_{-\infty,s} $, $ g \in H_{s,t} $, $ h \in H_{t,\infty}
$.

\begin{proposition}
The following conditions are equivalent for every $ f \in H $:

(a) $ f \in H_1 $;

(b) $ f = Q_{-\infty,t} f + Q_{t,\infty} f $ for all $ t \in \R $;

(c) $ Q_{r,t} f = Q_{r,s} f + Q_{s,t} f $ whenever $ -\infty \le r < s
< t \le \infty $.
\end{proposition}

\beginproof
Similar to \ref{6a1}, with $ Q_{s,t} $ instead of $
\cE{\cdot}{\F_{s,t}} $ and projection to the unit instead of
expectation.
\proofend

\begin{lemma}\label{6d12}
The space $ H_0 \oplus H_1 \oplus H_2 \oplus \dots $ is the closure of
the union of all subspaces of the form
\[
\bigotimes_{k=0}^n \( H_0(t_k,t_{k+1}) \oplus H_1(t_k,t_{k+1}) \)
\]
where $ -\infty = t_0 < t_1 < \dots < t_n < t_{n+1} = \infty $, $
n=0,1,2,\dots $
\end{lemma}

\beginproof
Similar to \ref{6a4a}.
\proofend

The formula \eqref{6a4aa} holds as well. Similarly to \eqref{6a4b},
spaces $ H_1(A) $ may be defined for all Borel sets $ A \subset \R $.

Given a continuous product of pointed Hilbert spaces $ (H_{s,t},
u_{s,t})_{s<t} $, we introduce the `upward
continuity'\index{upward continuity!for pointed Hilbert spaces}
condition, similar to \eqref{2c2a},
\begin{equation}\label{6d13}
H_{s,t} \text{ is the closure of } \bigcup_{\eps>0} H_{s+\eps,t-\eps}
\quad \text{for $ -\infty \le s < t \le \infty $}
\end{equation}
(here $ -\infty+\eps $ means $ -1/\eps $, $ \infty-\eps $ means $
1/\eps $), and the `downward
continuity'\index{downward continuity!for pointed Hilbert spaces} 
condition, similar to \eqref{2c3a},
\begin{equation}\label{6d14}
H_{s,t} = \bigcap_{\eps>0} H_{s-\eps,t+\eps} \quad \text{for }
-\infty \le s \le t \le \infty
\end{equation}
(here $ H_{t,t} $ is the one-dimensional subspace spanned by the unit,
$ -\infty-\eps $ means $ -\infty $, and $ \infty+\eps $ means $ \infty
$).

Upward continuity \eqref{2c2a} for continuous products of probability
spaces is evidently equivalent to upward continuity \eqref{6d13}  of
the corresponding continuous products of pointed Hilbert spaces. The
same holds for downward continuity. As noted in Sect.~\ref{2c} (after
\ref{2c3}), downward continuity does not imply upward continuity. The
argument of \ref{2c3} may be generalized as follows.

\begin{proposition}\label{6d15}
Upward continuity implies downward continuity.
\end{proposition}

\beginproof
It is sufficient to prove that $ \bigcap_{\eps>0} H_{s,s+\eps} $ is
one-dimensional (spanned by the unit), since $ \bigcap_{\eps>0}
H_{-\infty,s+\eps} = H_{-\infty,s} \otimes \bigcap_{\eps>0}
H_{s,s+\eps} $. Assuming the contrary, we take $ f \in
\bigcap_{\eps>0} H_{s,s+\eps} $, $ \| f \| = 1 $, orthogonal to the
unit. Using upward continuity we approximate $ f $ by $ g \in
H_{s+\eps,\infty} $, $ \| g \| = 1 $. We have $ f = f_{s,s+\eps}
u_{s+\eps,\infty} $, $ g = u_{s,s+\eps} g_{s+\eps,\infty} $; thus, $ |
\ip{ f_{s,s+\eps} }{ u_{s,s+\eps} } | \cdot | \ip{ u_{s+\eps,\infty}
}{ g_{s+\eps,\infty} } | = | \ip f g | $ is close to $ 1 $, while $ |
\ip{ f_{s,s+\eps} }{ u_{s,s+\eps} } | $ is small; a contradiction.
\proofend

\begin{proposition}\label{6d16}
(Zacharias \cite[Lemma 2.2.1]{Za}, Arveson \cite[Th.~6.2.3]{Ar}.)
Let $ (H_{s,t},u_{s,t})_{s<t} $ be a continuous product of pointed
Hilbert spaces, satisfying the upward continuity condition
\eqref{6d13}. Let $ -\infty \le r < t \le \infty $ be given, and $ f
\in H_{r,t} $ be a decomposable vector. Then $ \ip{f}{u_{r,t}} \ne 0
$.
\end{proposition}

\beginproof
Similar to the last argument of the proof of \ref{6a10}. Namely, $ f $
is the limit of the projection $ g_\eps = f_{r,s-\eps} \ip{
f_{s-\eps,s+\eps} }{ u_{s-\eps,s+\eps} } u_{s-\eps,s+\eps}
f_{s+\eps,t} $ of $ f $ to $ H_{r,s-\eps} \otimes u_{s-\eps,s+\eps}
\otimes H_{s+\eps,t} $. Therefore (see also \cite[Sect.~6.1]{Ar})
\[
\frac{ | \ip{ f_{s-\eps,s+\eps} }{ u_{s-\eps,s+\eps} } | }{ \|
f_{s-\eps,s+\eps} \| } = \frac{ \| f_{r,s-\eps} \ip{ f_{s-\eps,s+\eps}
}{ u_{s-\eps,s+\eps} } u_{s-\eps,s+\eps} f_{s+\eps,t} \| }{ \|
f_{r,s-\eps} \| \| f_{s-\eps,s+\eps} \| \| f_{s+\eps,t} \| } = \frac{
\| g_\eps \| }{ \| f \| } \to 1 \, .
\]
We cover the compact interval $ [r,t] $ by a finite number of open
intervals $ (s-\eps,s+\eps) $ such that $ \ip{ f_{s-\eps,s+\eps} }{
u_{s-\eps,s+\eps} } \ne 0 $ and get $ \ip{f}{u_{r,t}} \ne 0 $.
(The argument may be adapted to the compact interval $
[-\infty,\infty] $.)
\proofend

Now we generalize \ref{6a10}, \ref{6a17} and \ref{6a18}.

\begin{proposition}\label{6d17}
Let $ (H_{s,t},u_{s,t})_{s<t} $ be a continuous product of pointed
Hilbert spaces, satisfying the upward continuity condition
\eqref{6d13}. Then all decomposable vectors of $ H_{s,t} $ belong to $
H_0(s,t) \oplus H_1(s,t) \oplus H_2(s,t) \oplus \dots $ (that is, are
orthogonal to $ H_\infty(s,t) $), for $ -\infty \le s < t \le \infty
$.
\end{proposition}

\beginproof
According to \ref{6d11}, it is sufficient to prove that $ \ip{ U^\rho
f }{ f } \to \| f \|^2 $ as $ \rho \to 1 $. We use \ref{6d16} and the
inequality
\begin{equation}\label{6d18}
\ip{ U^\rho f }{ f } \ge \| f \|^{2\rho} | \ip{ f }{ u_{s,t} }
|^{2(1-\rho)}
\end{equation}
proven similarly to \eqref{6a11} (irrespective of the upward
continuity).
\proofend

Let $ f \in H_{-\infty,\infty} $ be a decomposable vector such that $
\ip{ f }{ u_{-\infty,\infty} } = 1 $. Then $ f_{s,t} $ as in \ref{6d3}
may be chosen such that $ \ip{ f_{s,t} }{ u_{s,t} } = 1 $ for $ s<t
$. Using \eqref{6d18} we see that $ f $ is orthogonal to $ H_\infty $.
Similarly to \ref{6a17} (using \ref{6d12}), $ f $ is uniquely
determined by its projection $ g $ to $ H_1 $. Similarly to
\ref{6a18}, every $ g \in H_1 $ is the projection of some decomposable
$ f $. Similarly to \eqref{6a19} we denote the relation between $ f $
and $ g $ by
\begin{equation}
f = \Exp g \, ; \qquad g = \Log f \, .
\end{equation}%
\index{zze@$ \Exp $ (over cont.\ products)}%
\index{zzl@$ \Log $ (over cont.\ products)}
It is a one-to-one correspondence between decomposable vectors $ f \in
H_{-\infty,\infty} $ satisfying $ \ip{ f }{ u_{-\infty,\infty} } = 1
$, and vectors $ g \in H_1 $. Still, $ \| \Exp g \|^2 \le \exp ( \| g
\|^2 ) $ and $ g = \lim_{\eps\to0} \frac1\eps \( \Exp(\eps g) -
u_{-\infty,\infty} \) $. Similarly to \ref{6a20}, the space $ H_0
\oplus H_1 \oplus H_2 \oplus \dots $ is the closed linear span of $ \{
\Exp g : g \in H_1 \} $. Similarly to \ref{6a21}, the equality $ \|
\Exp g \|^2 = \exp ( \| g \|^2 ) $ is ensured if $ g_{t-,t+} = 0 $ for
all $ t $, which in turn is ensured by the downward continuity
condition \eqref{6d14}, since $ H_1(t-,t+) \subset H_{t-,t+} $ (and $
H_1 $ is orthogonal to the unit). Similarly to \ref{6b10} we get the
following.

\begin{proposition}\label{6d21}
If $ (H_{s,t},u_{s,t})_{s<t} $ is a continuous product of pointed
Hilbert spaces satisfying the downward continuity condition
\eqref{6d14}, then the map $ \Exp : H_1 \to H_0 \oplus H_1 \oplus H_2
\oplus \dots $ has the property 
\[
\ip{ \Exp f }{ \Exp g } = \exp \ip f g \quad \text{for } f,g \in H_1
\, .
\]
\end{proposition}

See also \cite[Sect.~6.4]{Ar}, \cite[Th.~2.2.4]{Za}, \cite[Appendix
A]{TV}. We conclude that
\[
H_0 \oplus H_1 \oplus H_2 \oplus \dots = \E^{H_1}
\]
is the Fock space.\index{Fock space}

\mysubsection{Hilbert spaces}
\label{6e}

We want to know, to which extent results of Sect.~\ref{6d} depend on
the choice of a unit $ (u_{s,t})_{s<t} $ in a given continuous product
of Hilbert spaces $ (H_{s,t})_{s<t} $ (if a unit exists). As before,
the time set is $ [-\infty,\infty] $. First, note that the enlargement
of $ \R $ to $ [-\infty,\infty] $ mentioned after \ref{6d6} depends
heavily on the choice of a unit. Second, we compare two kinds of
continuity, one being unit-independent \eqref{3c*}, \eqref{3c**}, the
other unit-dependent \eqref{6d13}, \eqref{6d14}.

\begin{proposition}\label{6e1}
Let $ (H_{s,t})_{s<t} $ be a continuous product of Hilbert spaces, and
$ (u_{s,t})_{s<t} $ a unit. Then the following three conditions are
equivalent:

(a) the upward continuity \eqref{3c*} of $ (H_{s,t})_{s<t} $;

(b) the downward continuity \eqref{3c**} of $ (H_{s,t})_{s<t} $;

(c) the upward continuity \eqref{6d13} of $ (H_{s,t},u_{s,t})_{s<t}
$.
\end{proposition}

\beginproof
(a) \imp (b): see the proof of \ref{3c24}.

(b) \imp (c): \cite[Prop.~3.4]{Li} Projections $ Q_\eps : H \to H $
defined by $ Q_\eps ( fgh ) = f \ip g {u_{-\eps,\eps}} u_{-\eps,\eps}
h $ for $ f \in H_{-\infty,-\eps} $, $ g \in H_{-\eps,\eps} $, $ h \in
H_{\eps,\infty} $ belong to algebras $ \A_{-\eps,\eps} $. Their limit
$ Q_{0+} = \lim_{\eps\to0} Q_\eps $ belongs to the trivial algebra $
\A_{0,0} $ by \eqref{3c**}, and $ Q_{0+} u_{-\infty,\infty} =
u_{-\infty,\infty} $, therefore $ Q_{0+} = \One $. It follows that $
H_{-\infty,-\eps} \uparrow H_{-\infty,0} $ and $ H_{\eps,\infty}
\uparrow H_{0,\infty} $.

(c) \imp (a): we take $ \eps_n \downarrow 0 $ and introduce
projections $ Q_n $ of $ H_{s,t} $ onto $ H_{s+\eps_n,t-\eps_n}
\subset H_{s,t} $, then $ Q_n \uparrow \One $ by \eqref{6d13}. For any
operator $ A \in \A_{s,t} $ we define $ A_n \in \A_{s+\eps_n,t-\eps_n}
$ by $ A_n f = Q_n A f $ for $ f \in H_{s+\eps_n,t-\eps_n} \subset
H_{s,t} $ and observe that $ A_n \to A $ strongly, since $ \| A_n \|
\le \| A \| $ and $ A_n f \to A f $ for all $ f \in \cup_{\eps>0}
H_{s+\eps,t-\eps} $.
\proofend

See also \cite[proof of 6.1.1]{Ar}.

\begin{question}\label{6e2}\index{question}
What about a counterpart of \ref{6e1} for the downward continuity
\eqref{6d14} of $ (H_{s,t},u_{s,t})_{s<t} \, $?
\end{question}

Operators $ U^\rho $ and subspaces $ H_n $ (recall \ref{6d11}) depend
on the choice of a unit.

\begin{theorem}\label{6e3}\index{theorem}
For every continuous product of Hilbert spaces, containing at least
one unit and satisfying the (equivalent) continuity conditions
\ref{6e1}(a,b), the subspaces $ H_0 \oplus H_1 \oplus \dots $ and $
H_\infty $ do not depend on the choice of the unit.
\end{theorem}

The theorem follows immediately from the next result (or
alternatively, from \ref{6e3b}).

\begin{proposition}\label{6e3a}
\begin{sloppypar}
Let $ (H_{s,t})_{s<t} $ be a continuous product of Hilbert spaces
satisfying \ref{6e1}(a,b) and
containing at least one unit. Then $ H_0(s,t) \oplus H_1(s,t) \oplus
\nolinebreak
\dots $ is the closed linear span of (the set of all) decomposable
vectors of $ H_{s,t} $, for $ -\infty \le s < t \le \infty $.
\end{sloppypar}
\end{proposition}

\beginproof
Combine Proposition \ref{6d17} and the generalization of \ref{6a20}
(mentioned in Sect.~\ref{6d}).
\proofend

Assuming \ref{6e1}(a,b) we define $ H^\cls_{s,t} $%
\index{zzh@$ H^\cls_{s,t} $, subspace}
as the closed
linear span of decomposable vectors of $ H_{s,t} $ and get for $
H^\cls = H^\cls_{-\infty,\infty} $%
\index{zzh@$ H^\cls $, subspace}
\begin{equation}\label{6e4a}
H^\cls = H_0 \oplus H_1 \oplus \dots = H \ominus H_\infty
\end{equation}
if at least one unit exists; otherwise $ \dim H^\cls = 0 $. Clearly,
\[
H^\cls_{r,t} = H^\cls_{r,s} \otimes H^\cls_{s,t} \quad \text{for }
-\infty \le r < s < t \le \infty \, ,
\]
and we get the classical part%
\index{classical (part of)!continuous product!of Hilbert spaces}
$ (H^\cls_{s,t})_{s<t} $ of a continuous
product of Hilbert spaces $ (H_{s,t})_{s<t} $ provided that $ \dim
H^\cls \ne 0 $. Proposition \ref{6d21} shows that the classical part
is the Fock space,\index{Fock space}
\[
H^\cls = \E^{H_1} \, .
\]
Of course, the space $ H_1 $ and the map $ \Exp : H_1 \to H^\cls $
depend on the choice of a unit $ u $. We make the dependence explicit by
writing $ uH_1 $\index{zzu@$ uH_n $, subspace}
instead of $ H_1 $ and $ u\Exp $\index{zzu@$ u\Exp $}
instead of $ \Exp
$. Recall that an \emph{affine}\index{affine}
operator between two Hilbert spaces $ H', H'' $ is an operator of the
form $ x \mapsto y_0 + Lx $ (for $ x \in H' $) where $ L : H' \to H''
$ is a linear operator and $ y_0 \in H'' $. Here and in the following
proposition, Hilbert spaces over $ \R $ or $ \C $ are acceptable.

\begin{proposition}\label{6e6a}
Let $ (H_{s,t})_{s<t} $ be a continuous product of Hilbert spaces
satisfying \ref{6e1}(a,b), and $ (u_{s,t})_{s<t} $, $ (v_{s,t})_{s<t}
$ two units. Then there exists an isometric affine invertible map $ A
: uH_1 \to vH_1 $ such that the following conditions are equivalent
for all $ f \in uH_1 $, $ g \in vH_1 $:

(a) $ A(f) = g $;

(b) tho vectors $ u\Exp f $, $ v\Exp g $ span the same one-dimensional
subspace.
\end{proposition}

\beginproof
It follows from \ref{6d21} that
\begin{multline}\label{6e7}
\frac{ \ip{ u\Exp f_1 }{ u\Exp f_2 } }{ \| u\Exp f_1 \| \| u\Exp f_2
 \| } = \\
\frac{ \ip{ u\Exp f_1 }{ u } }{ | \ip{ u\Exp f_1 }{ u } | }
 \frac{ \ip{ u }{ u\Exp f_2 } }{ | \ip{ u }{ u\Exp f_2 } | }
 \exp \( \ip{ f_1 }{ f_2 } - \tfrac12 \| f_1 \|^2 - \tfrac12 \| f_2
\|^2 \)
\end{multline}
for all $ f_1, f_2 \in uH_1 $. Therefore
\begin{equation}\label{6e8}
\frac{ | \ip{ u\Exp f_1 }{ u\Exp f_2 } | }{ \| u\Exp f_1 \| \| u\Exp
f_2 \| } =
\exp \( -\tfrac12 \| f_1 - f_2 \|^2 \) \, ,
\end{equation}
which expresses the distance between $ f_1 $ and $ f_2 $ in terms of
the one-dimensional subspaces spanned by $ u\Exp f_1 $ and $ u\Exp f_2
$. We get an isometric invertible map $ A : uH_1 \to vH_1 $; it
remains to prove that $ A $ is affine. In the real case (over $ \R $)
it is well-known (and easy to see) that every isometry is affine. In
the complex case (over $ \C $) one more implication of \eqref{6e7} is
used:
\begin{multline}
\exp \( \I \Im \ip{ f_2-f_1 }{ f_3-f_1 } \) = \\
= \frac{ \ip{ u\Exp f_1 }{ u\Exp f_2 } \ip{ u\Exp f_2 }{ u\Exp f_3 }
 \ip{ u\Exp f_3 }{ u\Exp f_1 } }
 { | \ip{ u\Exp f_1 }{ u\Exp f_2 } \ip{ u\Exp f_2 }{ u\Exp f_3 }
 \ip{ u\Exp f_3 }{ u\Exp f_1 } | } \, . \qquad \qed
\end{multline}
\proofendnoqed

The dependence of the operators $ U^\rho $ on the choice of a unit is
estimated below (which gives us another proof of Theorem \ref{6e3}).
For convenience we write $ U^\rho $ and $ V^\rho $ rather than $ uU^\rho $ and
$ vU^\rho $.

\begin{proposition}\label{6e3b}
Let $ (H_{s,t})_{s<t} $ be a continuous product of Hilbert spaces
satisfying \ref{6e1}(a,b), and $ (u_{s,t})_{s<t} $, $ (v_{s,t})_{s<t}
$ two units. Then the operators $ U^\rho_{s,t} $, $ V^\rho_{s,t} $
corresponding to these units satisfy the inequality
\[
U^\rho_{s,t} \ge ( V^\rho_{s,t} )^2 \exp ( -4 (\ln\rho) \ln | \ip{
u_{s,t} }{ v_{s,t} } | )
\]
for $ -\infty \le s < t \le \infty $ and $ \rho \in (0,1] $.
\end{proposition}

(Inequalities between operators are treated as inequalities between
their quadratic forms, of course.)

\beginproof
First, a general inequality
\begin{equation}\label{6e4}
( 2 | \ip x z |^2 - 1 ) ( 2 | \ip y z |^2 - 1 ) \le | \ip x y |^2
\end{equation}
is claimed for any three unit vectors $ x,y,z $ of a Hilbert space satisfying
$ 2 | \ip x z |^2 \ge 1 $. We prove it introducing $ \al, \be, \ga
\in [0,\pi/2] $ by $ \cos\al = | \ip x y | $, $ \cos\be = | \ip x z |
$,  $ \cos\ga = | \ip y z | $. We have $ \al \le \be + \ga $. However,
$ \be \le \pi/4 $; also $ \ga \le \pi/4 $, otherwise there is nothing
to prove. Therefore $ \cos 2\al \ge \cos 2(\be+\ga) $, $
2(\cos2\be) (\cos2\ga) = \cos 2(\be-\ga) + \cos 2(\be+\ga) \le 1 +
\cos 2\al $ and $ (2\cos^2 \be - 1) (2\cos^2 \ga - 1) \le \cos^2 \al
$, which is \eqref{6e4}.

Second, for a given $ \rho \in [\frac12,1] $ and $ s<t $ we define
operators $ U,V $ by
\[
Uf = \rho f + (1-\rho) \ip{ f }{ u_{s,t} } u_{s,t} \, , \quad
Vf = (2\rho-1) f + 2(1-\rho) \ip{ f }{ v_{s,t} } v_{s,t} \, .
\]
Assuming $ \| f \| = 1 $ and applying \eqref{6e4} to $ x = u_{s,t} $,
$ y = f $, $ z = v_{s,t} $ we get
\[
\( 2 | \ip{ u_{s,t} }{ v_{s,t} } |^2 - 1 \) \( 2 | \ip{ f }{ v_{s,t} }
|^2 - 1 \) \le | \ip{ f }{ u_{s,t} } |^2
\]
provided that $ 2 | \ip{ u_{s,t} }{ v_{s,t} } |^2 \ge 1 $. It may be
written as an operator inequality
\begin{gather*}
\( 2 | \ip{ u_{s,t} }{ v_{s,t} } |^2 - 1 \) \bigg( \frac{ V -
 (2\rho-1) \One }{ 1-\rho } - \One \bigg) \le \frac{ U - \rho \One }{
 1 - \rho } \, ; \\
\( 2 | \ip{ u_{s,t} }{ v_{s,t} } |^2 - 1 \) (V-\rho\One) \le U - \rho
 \One \, ; \\
\( 2 | \ip{ u_{s,t} }{ v_{s,t} } |^2 - 1 \) V + 2\rho \( 1 - | \ip{
 u_{s,t} }{ v_{s,t} } |^2 \) \One = \qquad\qquad \\
\qquad\qquad = \( 2 | \ip{ u_{s,t} }{ v_{s,t} } |^2 - 1 \)
(V-\rho\One) + \rho\One \le U \, .
\end{gather*}
Taking into account that $ V \le \One $ we have
\begin{equation}\label{6e5}
\( 2\rho - 1 + 2(1-\rho) | \ip{ u_{s,t} }{ v_{s,t} } |^2 \) V \le U \,
.
\end{equation}

Third, using the upward continuity condition (similarly to
\ref{6d16}), for any $ \eps $ we can choose $ s = t_0 < t_1 < \dots <
t_n = t $ such that $ | \ip{ u_{t_{k-1},t_k} }{ v_{t_{k-1},t_k} } | >
1-\eps $ for $ k = 1,\dots,n $. We apply \eqref{6e5} on each interval
$ (t_{k-1},t_k) $ rather than $ (s,t) $ and multiply the inequalities:
\[
\bigg( \prod_{k=1}^n \( 2\rho - 1 + 2(1-\rho) | \ip{ u_{t_{k-1},t_k}
}{ v_{t_{k-1},t_k} } |^2 \bigg) \ti V_{t_1,\dots,t_n}^{2\rho-1} \le
\ti U_{t_1,\dots,t_n}^\rho \, ;
\]
taking the limit of the net we get for $ \rho \in [\frac12,1] $
\begin{equation}\label{6e6}
| \ip{ u_{s,t} }{ v_{s,t} } |^{4(1-\rho)} V_{s,t}^{2\rho-1} \le
 U_{s,t}^\rho \, ,
\end{equation}
since $ \ln ( 2\rho - 1 + 2(1-\rho) a^2 ) \sim -2(1-\rho) (1-a^2) \sim
4(1-\rho) \ln a $ as $ a \to 1 $.

Fourth, we apply \eqref{6e6} to $ \rho = 1 - \frac c n $ (for an
arbitrary $ c>0 $) and raise to the $ n $-th power, using the
semigroup property of $ U^\rho $, $ V^\rho $:
\[
| \ip{ u_{s,t} }{ v_{s,t} } |^{\frac{4c}n} V_{s,t}^{1-\frac{2c}n} \le
 U_{s,t}^{1-\frac c n} \, ; \qquad
| \ip{ u_{s,t} }{ v_{s,t} } |^{4c} V_{s,t}^{(1-\frac{2c}n)^n} \le
 U_{s,t}^{(1-\frac c n)^n} \, ;
\]
the limit (as $ n \to \infty $) gives
\[
| \ip{ u_{s,t} }{ v_{s,t} } |^{4c} V_{s,t}^{\exp(-2c)} \le
 U_{s,t}^{\exp(-c)} \, .
\]
Substituting $ \rho = \exp(-c) $ we get $ | \ip{ u_{s,t} }{ v_{s,t} }
|^{-4\ln\rho} V_{s,t}^{\rho^2} \le U_{s,t}^\rho $.
\proofend

Definition \ref{3a1} of a continuous product of Hilbert spaces $
(H_{s,t})_{s<t} $ stipulates the global space $ H_{-\infty,\infty}
$. However, all said about local spaces (for $ s,t \ne \pm\infty $)
holds for \emph{local} continuous products of Hilbert spaces (defined
similarly to \ref{3a1} but waiving the infinite points $ \pm\infty $
on the time axis).

\mysubsection{Homogeneous case; Arveson systems of type $ I $}
\label{6f}

Given a homogeneous continuous product of Hilbert spaces $
(H_{s,t})_{s<t} $, $ (\theta_{s,t}^h)_{s<t;h} $, we may consider the
classical part of the continuous product of Hilbert spaces $
(H_{s,t})_{s<t} $. It consists of spaces $ H^\cls_{s,t} $ spanned by
decomposable vectors. The time shift $ \theta_{s,t}^h $ sends $
H^\cls_{s,t} $ to $ H^\cls_{s+h,t+h} $, since decomposable vectors go to
decomposable vectors. We see that the classical part of a homogeneous
continuous product of Hilbert spaces%
\index{classical (part of)!homogeneous continuous product!of Hilbert spaces}
is a \emph{homogeneous}
continuous product of Hilbert spaces, provided that $ \dim H^\cls > 0
$.

If the homogeneous continuous product of Hilbert spaces corresponds to
a noise then surely $ \dim H^\cls > 0 $, since constant functions are
decomposable vectors. They are also shift-invariant, that is,
invariant under the group of shifts $ (\theta^h)_h $ (that is, $
(\theta^h_{-\infty,\infty})_h $).

As before, the time set is $ [-\infty,\infty] $, which is crucial
below.

\begin{proposition}\label{6f1}
For every homogeneous continuous product of Hilbert spaces $
(H_{s,t})_{s<t} $, $ (\theta_{s,t}^h)_{s<t;h} $, the subspace spanned
by all shift-invariant decomposable vectors is either \dimensional{0}
or \dimensional{1}.
\end{proposition}

\beginproof
Let $ u, v $ be two such vectors, $ \| u \| = 1 $, $ \| v \| = 1 $. By
\ref{6e1}, the upward continuity \eqref{3c*} of $ (H_{s,t})_{s<t} $,
ensured by \ref{3c23}, implies the upward continuity \eqref{6d13} of $
(H_{s,t},u_{s,t})_{s<t} $. By \ref{6d16}, $ \ip u v \ne 0 $.

Clearly, $ | \ip u v | \le \prod_{k=1}^n | \ip{ u_{k-1,k} }{ v_{k-1,k}
} | $. However, $ | \ip{ u_{k-1,k} }{ v_{k-1,k} } | $ does not depend
on $ k $ by the shift invariance. Thus, $ 0 < | \ip u v | \le | \ip{
u_{0,1} }{ v_{0,1} } |^n $ for all $ n $, which means that $ | \ip{
u_{0,1} }{ v_{0,1} } | = 1 $.

Similarly to the proof of \ref{6d16}, $ | \ip{ u_{-\infty,-n} }{
v_{-\infty,-n} } | \to 1 $ and $ | \ip{ u_{n,\infty} }{ v_{n,\infty} }
| { \to 1 } $ as $ n \to \infty $. Therefore $ | \ip u v | = \lim_n |
\ip{ u_{-n,n} }{ v_{-n,n} } | = \lim_n | \ip{ u_{0,1} }{ v_{0,1} }
|^{2n} { = 1 } $.
\proofend

\begin{sloppypar}
\begin{example}
It can happen that decomposable vectors exist, but no one of them is
shift-invariant. An example will be constructed from a
(non-stationary) \flow{\R} $ (X_{s,t})_{s<t} $ such that each $
X_{s,t} $ is distributed normally, $ \Var(X_{s,t}) = t-s $ (just like
Brownian increments) and $ \Ex X_{s,t} = 0 $ for $ -\infty < s <
t \le 0 $, but $ \Ex X_{s,t} = t-s $ for $ 0 \le s < t < \infty $ (a
drift after $ 0 $). Time shifts $ T_h : \Om \to \Om $ defined by $
X_{s,t} \circ T_h = X_{s+h,t+h} $ do not preserve the measure $ P $
(thus, our object is not a noise), however, they transform $ P $ into
an equivalent (that is, mutually absolutely continuous) measure. These
time shifts lead to unitary operators $ \theta_h $ on $ H = L_2(\Om,P)
$; unitarity is achieved by multiplying the given function at $ T_h
\om $ by the square root of the corresponding density (that is,
Radon-Nikodym derivative), see Sect.~\ref{10a} for details. The
absolute continuity of measures is not uniform in $ h $. Moreover, $
\theta_h \to 0 $ in the weak operator topology (as $ h \to \pm\infty
$), which excludes any non-zero shift-invariant vector (decomposable
or not). On the other hand, waiving time shifts we get a
\emph{classical} system; decomposable vectors span $ H $.
\end{example}
\end{sloppypar}

Shift-invariant decomposable vectors are scarce, as far as global
vectors are meant. Think about $ u = \exp
\(\I\la(B_\infty-B_{-\infty})\) $ for a Brownian motion $ B $; $ u $
is ill-defined (unless $ \la=0 $), however $ u_{s,t} = \exp \( \I\la
(B_t-B_s) \) $ is well-defined (for each $ \la $).

Waiving the global space $ H_{-\infty,\infty} $ we get \emph{local}
homogeneous continuous products of Hilbert spaces, or equivalently,
algebraic product systems of Hilbert spaces (recall \ref{3c4}). In
order to exclude pathologies we impose a natural condition of
measurability, thus turning to Arveson systems (recall \ref{3c5}).

\begin{definition}\label{6f3}
A \emph{unit}\index{unit!of Arveson system}
of an Arveson system $ (H_t)_{t>0} $ is a family $
(u_t)_{t>0} $ of vectors $ u_t \in H_t $ such that

(a) $ u_s u_t = u_{s+t} $ for $ s,t \in (0,\infty) $,

(b) the map $ t \mapsto u_t $ from $ (0,\infty) $ to $ \biguplus_{t>0}
H_t $ is Borel measurable,

(c) $ \| u_t \| = 1 $ for $ t \in (0,\infty) $.
\end{definition}

Arveson \cite[3.6.1]{Ar} admits $ \| u_t \| = \E^{ct} $ for any $ c
\in \R $, but I prefer $ \| u_t \| = 1 $ for compatibility with
\ref{6d6}. 

Translating \ref{6f3} into the language of local homogeneous continuous
products of Hilbert spaces we get units $
(u_{s,t})_{-\infty<s<t<\infty} $ satisfying $ \theta^h_{s,t} u_{s,t} =
u_{s+h,t+h} $ and measurability. The latter can be enforced by
appropriate scalar coefficients, see \ref{6d3}. Accordingly,
\ref{6f3}(b) can be enforced by replacing $ u_t $ with $ c_t u_t $ for
some $ (c_t)_{t>0} $ such that $ c_s c_t = c_{s+t} $. Continuity
conditions \ref{6e1}(a--c) are satisfied locally (on every finite time
interval).

\begin{lemma}
Let $ (H_t)_{t>0} $ be an Arveson system such that $ H_t $ contains at
least one decomposable vector, for at least one $ t $. Then the system
contains a unit.
\end{lemma}

\beginproof
Assuming that $ H_1 $ contains a decomposable vector, we get a
decomposable vector $ v \in H_\T $ of the cyclic-time system
corresponding to the given system (recall Sect.~\ref{3d}). The space $
vH_1 $ need not be shift-invariant, however, the set $ \{ c \cdot
v\Exp f : c \in \C \setminus \{0\}, \, f \in vH_1 \} $ of all
decomposable vectors is shift-invariant. Applying \ref{6e6a} we see
that $ \T $ acts on $ vH_1 $ by affine transformations $ A^h : vH_1
\to vH_1 $; namely, $ v\Exp (A^h f) = c \theta^h_\T (v\Exp f) $ for
some $ c \in \C \setminus \{0\} $ (that depends on $ f $ and $ h
$). The action is continuous (recall \eqref{6e8}), and has at least
one fixed point, since the point $ \int_\T A^h f \, dh $ evidently is
fixed, irrespective of $ f \in vH_1 $. (Roughly speaking, the
geometric mean of all shifts of a decomposable vector is a
shift-invariant decomposable vector.) Let $ f $ be a fixed point, then
the decomposable vector $ u = v\Exp f $ satisfies $ \theta^h_\T u =
c(h) u $, therefore $ \theta^h_\T u = \E^{2\pi\I nh} u $ for some $ n
\in \Z $.

Now it will be shown that $ n=0 $, by checking the equality $
\theta^{1/n}_\T u = u $ for all $ n=1,2,\dots $ Instead of $ u $ we
use here another vector of the same one-dimensional subspace, namely,
$ \bigotimes_{k=0}^{n-1} \theta_{0,\frac1n}^{\frac k n} u_{0,\frac1n}
$. Applying to it the operator $ \theta^{\frac1n}_\T =
\bigotimes_{k=0}^{n-1} \theta_{\frac k n,\frac{k+1}n}^{\frac 1 n} $ we
get $ \bigotimes_{k=0}^{n-1} \theta_{\frac k n,\frac{k+1}n}^{\frac 1 n}
\theta_{0,\frac1n}^{\frac k n} u_{0,\frac1n} = \bigotimes_{k=0}^{n-1}
\theta_{0,\frac1n}^{\frac{k+1}n} u_{0,\frac1n} =
\bigotimes_{k=0}^{n-1} \theta_{0,\frac1n}^{\frac k n} u_{0,\frac1n} $.

Having $ \theta^h_\T u = u $ for $ h \in \T $, we return to the linear
time and construct a unit $ (u_t)_{t>0} $, namely, $ u_{n+t} =
\underbrace{ u_\T \otimes \dots \otimes u_\T }_n \otimes u_{0,t} $ for
$ n = 0,1,2,\dots $ and $ { 0 \le t < 1 } $.
\proofend

The classical part $ H_t^\cls $ of $ H_t $ may be defined as the
closed linear span of all decomposable vectors of $ H_t $. Clearly, $
H_s^\cls \otimes H_t^\cls = H_{s+t}^\cls $, and we get another Arveson
system $ (H_t^\cls)_{t>0} $, the classical part of the given Arveson
system, provided that $ \dim H_t^\cls \ne 0 $.

\begin{definition}
\cite[6.0.3]{Ar}
An Arveson system $ (H_t)_{t>0} $ is \emph{decomposable,}%
\index{decomposable!Arveson system}
if for every
$ t > 0 $, the space $ H_t $ is the closed linear span of its
decomposable vectors.
\end{definition}

The classical part of an Arveson system is decomposable. Here is a
counterpart of Theorem \ref{6c1}.

\begin{theorem}\label{6f6}\index{theorem}
Every decomposable Arveson system is generated by its units.
\end{theorem}

\begin{sloppypar}
That is, $ H_t $ is the closed linear span of vectors of the form $
(u_1)_{\frac t n} (u_2)_{\frac t n} \dots \linebreak[3]
(u_n)_{\frac t n} $, where $ u_1,u_2, \dots, u_n $ are units. See also
\cite[6.0.5]{Ar}, \cite[Cor.~6.6]{Li}. In order to prove the theorem we
first translate everything into the language of a local homogeneous
continuous product of Hilbert spaces $ (H_{s,t})_{s<t} $, $
(\theta^h_{s,t})_{s<t;h} $ and its classical part $
(H^\cls_{s,t})_{s<t} $, $ (\theta^h_{s,t}|_{H_{s,t}^\cls})_{s<t;h}
$. We know from Sect.~\ref{6e} that $ H_{s,t}^\cls = uH_0(s,t) \oplus
uH_1(s,t) \oplus \dots = H_{s,t} \ominus uH_\infty(s,t) $ for any unit
$ u $ (not necessarily shift-invariant; and if $ H_{s,t} $ contains no
units then $ \dim H^\cls_{s,t} = 0 $). Also, we have the map $
u_{s,t}\Exp : u_{s,t}H_1(s,t) \to H_{s,t}^\cls $ satisfying $ \ip{
u_{s,t}\Exp f }{ u_{s,t}\Exp g } = \exp \ip f g $. All decomposable
vectors of $ H_{s,t} $ are of the form $ u_{s,t}\Exp f $ up to a
coefficient. Thus, $ H_{s,t}^\cls $ is the Fock space,
\[
H_{s,t}^\cls = \E^{u_{s,t}H_1(s,t)} \, .
\]
We have no global space $ H^\cls_{-\infty,\infty} $, but still, the
global first chaos space $ uH_1 = uH_1(-\infty,\infty) $ is
well-defined, according to \emph{additive} relations $ u_{r,s}H_1(r,s)
\oplus u_{s,t}H_1(s,t) = u_{r,t}H_1(r,t) $.
\end{sloppypar}

Assume now that the unit $ u $ is shift-invariant. Then the subspace $
uH_1 $ is shift-invariant. We have no global exponential map $ u\Exp :
uH_1 \to H^\cls $, but we have a family of local exponential maps $
u_{s,t} \Exp : u_{s,t}H_1(s,t) \to H^\cls_{s,t} $ shift-invariant in
the sense that $ u_{s+h,t+h} \Exp (\theta^h_{s,t} f) = \theta^h_{s,t}
u_{s,t} \Exp f $ for $ f \in u_{s,t}H_1(s,t) $.

Similarly to the proof of Theorem \ref{6c1}, time shifts induce
unitary operators $ U_h : uH_1 \to uH_1 $ such that $ U_h^{-1} Q_{s,t}
U_h = Q_{s+h,t+h} $, where $ Q_{s,t} $ is the projection of $ uH_1 $
onto $ u_{s,t}H_1(s,t) \subset uH_1 $. They lead to Weyl relations $
U_h V_\la = \E^{\I\la h} V_\la U_h $, and so, $ H_1 $ decomposes into
the direct sum of a finite or countable number of irreducible
components, unitarily equivalent to the standard representation of
Weyl relations in $ L_2(\R) $. Similarly to Sect.~\ref{6c} we may
treat $ uH_1 $ as the tensor product, $ uH_1 = L_2(\R) \otimes \cH =
L_2(\R,\cH) $, where $ L_2(\R) $ carries the standard representation
of Weyl relations, and $ \cH $ is the Hilbert space of all families $
(g_{s,t})_{s<t} $ of vectors $ g_{s,t} \in u_{s,t} H_1 $ satisfying $
g_{r,t} = g_{r,s} + g_{s,t} $ and shift-invariant in the sense that $
\theta^h_{s,t} g_{s,t} = g_{s+h,t+h} $.

\begin{sloppypar}
Given a decomposable vector $ u_{0,t} \Exp g \in H_{0,t} $, we approximate $ g
\in u_{0,t} H_1(0,t) = L_2((0,t),\cH) $ by step functions $ g_n : (0,t) \to
\cH $ constant on $ (0,\frac t n), (\frac t n, \frac{2t}n), \dots,
(\frac{(n-1)t}n, t) $. Applying $ u_{0,t} \Exp $ to this step function we
complete the proof of Theorem \ref{6f6}.
\end{sloppypar}

At the same time we classify all \emph{classical} Arveson systems up
to isomorphism. They consist of Fock spaces,
\[
H_t = \E^{L_2((0,t),\cH)}
\]
with the natural multiplication (and Borel structure). The classifying
parameter is $ \dim \cH \in \{0,1,2,\dots\} \cup \{\infty\} $. See
Arveson \cite[Th.~6.7.1, Def.~3.1.6 and Prop.~3.5.1]{Ar} and Zacharias
\cite[Th.~4.1.10]{Za}.

\mysubsection{Examples}
\label{6g}

The two noises of Sect.~\ref{sec:4} (splitting and stickiness) are
nonclassical noises; both satisfy $ \F^\stable =
\F^\white_{-\infty,\infty} \subsetneqq \F_{-\infty,\infty} $ (recall
Sect.~\ref{5d}). The corresponding continuous products of spaces $ L_2
$ are nonclassical continuous products of \emph{pointed} Hilbert
spaces (recall Sect.~\ref{6d}), their classical parts being $ L_2
(\F^\stable) $. According to Sect.~\ref{6e}, they are also
nonclassical continuous products of Hilbert spaces; still, $ L_2
(\F^\stable) $ are their classical parts. According to Sect.~\ref{6f},
these nonclassical products lead to Arveson systems. See also Question
\ref{9d10}. In both cases
(splitting and stickiness), the classical part of the Arveson system
is the Fock space $ H^\cls = L_2 (\F^\white_{-\infty,\infty}) $. We
see that the classical part is neither trivial nor the whole system;
such Arveson systems (or rather, the corresponding $ E_0 $-semigroups)
are known as type $ II $ systems, see \cite[2.7.6]{Ar}. Type $ I $
means a classical system, while type $ III $ means a system with a
trivial classical part.%
\index{type $I,II,III$}\index{Arveson system!of type $I,II,III$}

\mysection{Distributed singularity, black noise (according to Le Jan and
 Raimond)}
\label{sec:7}
\mysubsection{Black noise in general}
\label{7a}

\begin{definition}
A noise is \emph{black}\index{black noise}
if its classical part is trivial, but the
whole noise is not.
\end{definition}

In other words: all stable random variables are constant, but some
random variables are not constant. There exist nontrivial centered
(that is, zero-mean) random variables, and they all are sensitive.
The self-joinings $ (\al_\rho,\be_\rho) $ and the operators $ U^\rho $
introduced in Sect.~\ref{5b} are trivial, irrespective of $ \rho \in
[0,1) $, if the noise is black. (See also \cite[Remark 2.1]{WW03}.)

Existence of some black noise is proven by Vershik and the author
\cite[Sect.~5]{TV}, but the term `black noise' appeared in
\cite{Ts98}. Why `black'? Well, the white noise is called `white'
since its spectral density is constant. It excites harmonic
oscillators of all frequencies to the same extent. For a black noise,
however, the response of any linear sensor is zero!

What could be a physically reasonable nonlinear sensor able to sense a
black noise? Maybe a fluid could do it, which is hinted at by the
following words of Shnirelman \cite[p.~1263]{Shn} about a paradoxical
motion of an ideal incompressible fluid: `\dots\ very strong external
forces are present, but they are infinitely fast oscillating in space
and therefore are indistinguishable from zero in the sense of
distributions. The smooth test functions are not ``sensitive'' enough
to ``feel'' these forces.'

The very idea of black noises, nonclassical continuous products etc.\
was suggested to me by Anatoly Vershik in 1994.

Two black noises are presented in Sections \ref{7f}, \ref{7j}. Two
more ways of constructing black noises are available, see \cite{Wa01}
and \cite[8b]{Ts03}.

By Theorem \ref{6a2}, a noise is black if and only if $ \dim H_1 = 0
$, that is, the first chaos space is trivial. More generally, a
continuous product of probability spaces has no classical part (that
is, its classical part is trivial) if and only if $ \dim H_1 = 0 $.

\begin{lemma}\label{7a2}
(a)
For every continuous product of probability spaces, the first chaos
space $ H_1 $ is equal to the intersection of spaces of the form
\[
L_2^0 (\F_{-\infty,t_0}) + L_2^0 (\F_{t_0,t_1}) + \dots + L_2^0
(\F_{t_{n-1},t_n}) + L_2^0 (\F_{t_n,\infty})
\]
over all finite sets $ \{ t_0,\dots,t_n \} \subset \R $, $ -\infty <
t_0 < \dots < t_n < \infty $. Here $ L_2^0(s,t) \subset
L_2(\F_{-\infty,\infty}) $ consists of all $ f \in L_2(\F_{s,t}) $
such that $ \Ex f = 0 $.

(b)
If the continuous product of probability spaces satisfies the downward
continuity condition \eqref{2c3a}, then the intersection over rational
$ t_0, \dots, t_n $ is also equal to $ H_1 $. The same holds for every
dense subset of $ \R $.
\end{lemma}

\beginproof
(a): Follows easily from \ref{6a1}.

\begin{sloppypar}
(b):
Given an irrational $ t \in \R $, we choose rational $ r_k \uparrow t
$. The downward continuity gives $ \cE{ f }{ \F_{r_k,\infty} } \to
\cE{ f }{ \F_{t,\infty} } $. If $ f $ belongs to the intersection over
rationales then $ f = \cE{ f }{ \F_{-\infty,r_k} } + \cE{ f }{
\F_{r_k,\infty} } $, thus $ \cE{ f }{ \F_{-\infty,r_k} } \to f - \cE{
f }{ \F_{t,\infty} } $. Taking $ \cE{ \dots }{ \F_{-\infty,t} } $ we
get $ \cE{ f }{ \F_{-\infty,r_k} } \to \cE{ f }{ \F_{-\infty,t} } $,
therefore $ f = \cE{ f }{ \F_{-\infty,t} } + \cE{ f }{ \F_{t,\infty} }
$. \qed
\end{sloppypar}
\proofendnoqed

\begin{corollary}
(a)
For every continuous product of probability spaces, the orthogonal
projection of $ f \in L_2 (\F_{-\infty,\infty}) $ such that $ \Ex f =
0 $ to $ H_1 $ is the limit (in $ L_2 $) of the net of random
variables
\[
\cE{ f }{ \F_{-\infty,t_0} } + \cE{ f }{ \F_{t_0,t_1} } + \dots + \cE{
f }{ \F_{t_{n-1},t_n} } + \cE{ f }{ \F_{t_n,\infty} }
\]
indexed by all finite sets $ \{ t_0,\dots,t_n \} \subset \R $, $
-\infty < t_0 < \dots < t_n < \infty $. 

(b)
For every noise, the orthogonal projection of $ f \in L_2
(\F_{0,1}) $ such that $ \Ex f = 0 $ to $ H_1(0,1) $ is equal to
\[
\lim_{n\to\infty} \sum_{k=1}^{2^n} \cE{ f }{ \F_{(k-1)2^{-n},k2^{-n}}
} \, .
\]

(c) A noise is black if and only if
\[
\lim_{n\to\infty} \sum_{k=1}^{2^n} \cE{ f }{ \F_{(k-1)2^{-n},k2^{-n}}
} = 0 \quad \text{for all } f \in L_2 (\F_{0,1}), \, \Ex f = 0
\]
(or equivalently, for all $ f $ of a dense subset of $ L_2^0(\F_{0,1})
$).
\end{corollary}

\beginproof
(a), (b):
Consider the projections to the spaces treated in \ref{7a2}; use
\ref{2c3}.

(c):
If $ \dim H_1(0,1) = 0 $ then $ \dim H_1(n,n+1) = 0 $ by homogeneity,
therefore $ \dim H_1 = 0 $.
\proofend

See also \cite[6a4]{Ts03}.

\begin{corollary}\label{7a4}
For every continuous product of probability spaces and every function
$ f \in L_2(\Om) $, if
\[
\Var \( \cE{ f }{ \F_{s,s+\eps} } \) = o(\eps) \quad \text{as } \eps
\to 0
\]
uniformly in $ s \in [r,t] $, then $ f $ is orthogonal to $ H_1(r,t)
$.
\end{corollary}

\beginproof
Assuming $ \Ex f = 0 $ we see that
\[
\Big\| \sum_{k=1}^n \cE{ f }{ \F_{s_{k-1},s_k} } \Big\|^2 =
\sum_{k=1}^n \Var \( \cE{ f }{ \F_{s_{k-1},s_k} } \)
\]
is much smaller than $ \sum_{k=1}^n (s_k-s_{k-1}) = t-r $ whenever $
r = s_0 < s_1 < \dots < s_n = t $ are such that $ \max (s_k-s_{k-1}) $
is small enough. Thus, the limit of the net vanishes.
\proofend

\mysubsection{Black noise and flow system}
\label{7b}

\begin{proposition}\label{7b1}
Let a flow system $ (X_{s,t})_{s<t} $, $ X_{s,t} : \Om \to G_{s,t} $
be such that for every interval $ (r,t) \subset \R $ and bounded
measurable function $ f : G_{r,t} \to \R $,
\begin{equation}\label{7b2}
\Var \( \cE{ f(X_{r,t}) }{ \F_{s,s+\eps} } \) = o(\eps) \quad
\text{as } \eps \to 0
\end{equation}
uniformly in $ s \in [r,t] $. Assume also that the classical part of
the corresponding continuous product of probability spaces satisfies
the equivalent continuity conditions \ref{6a6}(a--c). Then the
classical part is trivial.
\end{proposition}

\beginproof
By \ref{7a4} the random variable $ f(X_{r,t}) $ is orthogonal to $
H_1(r,t) $. However, $ \( H_{r,s} \ominus H_1(r,s) \) \otimes \( H_{s,t}
\ominus H_1(s,t) \) \subset H_{r,t} \ominus H_1(r,t) $, since $
H_{r,s} \ominus H_1(r,s) = \( H_0 \oplus H_2(r,s) \oplus H_3(r,s)
\oplus \dots \) \oplus H_\infty(r,s) $, and $ H_k(r,s) \otimes
H_l(s,t) \subset H_{k+l}(r,t) $ for $ k,l \in \{ 0,1,2,\dots \} \cup
\{\infty\} $. Therefore random variables of the form $ f_1(X_{r,s})
f_2(X_{s,t}) $ are orthogonal to $ H_1(r,t) $. Similarly, random
variables of the form $ f_1(X_{t_0,t_1}) f_2(X_{t_1,t_2}) \dots
f_n(X_{t_{n-1},t_n}) $ for $ t_0 < t_1 < \dots < t_n $ are
orthogonal to $ H_1 $. These random variables being dense in $
L_2(\Om) $, we get $ \dim H_1 = 0 $.
\proofend

We have
\begin{multline*}
\cE{ f(X_{r,t}) }{ \F_{s,s+\eps} } =
 \cE{ f(X_{r,s} X_{s,s+\eps} X_{s+\eps,t}) }{ \F_{s,s+\eps} } = \\
= \cE{ f(X_{r,s} X_{s,s+\eps} X_{s+\eps,t}) }{ X_{s,s+\eps} } =
 f_{s,s+\eps} (X_{s,s+\eps}) \, ,
\end{multline*}
where $ f_{s,s+\eps} : G_{s,s+\eps} \to \R $ is defined by
\begin{equation}\label{7b3}
f_{s,s+\eps} (y) = \int f(xyz) \, \mu_{r,s}(\D x)
\mu_{s+\eps,t}(\D z) \, ;
\end{equation}
as before, the measure $ \mu_{s,t} $ on $ G_{s,t} $ is the
distribution of $ X_{s,t} $. Therefore
\[
\Var \( \cE{ f(X_{r,t}) }{ \F_{s,s+\eps} } \) = \int |
f_{s,s+\eps} |^2 \, \D\mu_{s,s+\eps} - \bigg| \int f_{s,s+\eps}
\, \D\mu_{s,s+\eps} \bigg|^2 \, ,
\]
and \eqref{7b2} becomes
\begin{equation}\label{7b4}
\int | f_{s,s+\eps} |^2 \, \D\mu_{s,s+\eps} - \bigg| \int
f_{s,s+\eps} \, \D\mu_{s,s+\eps} \bigg|^2 = o(\eps)
\end{equation}
uniformly in $ s $, for all $ f \in L_\infty (G_{r,t},\mu_{r,t}) $. It
is sufficient to check the condition for all $ f $ of a subset of $
L_\infty (G_{r,t},\mu_{r,t}) $ dense in $ L_2 (G_{r,t},\mu_{r,t}) $
(recall the last argument of the proof of \ref{7b1}).

\mysubsection{About random maps in general}
\label{7c}

It is often inconvenient to treat a random process as a random
function, that is, a map from $ \Om $ to the space of functions. Say,
a Poisson process has a right-continuous modification, a
left-continuous modification and a lot of other modifications, but the
choice of a modification is often irrelevant. It is already
stipulated in Sect.~\ref{2aa} that ``a stochastic flow (and any random
process) is generally treated as a family of equivalence classes
(rather than functions)'', but now we have to be more explicit.

\begin{definition}\label{7c1}
(a)
An \emph{\Smap}\index{Smap@\Smap}
from a set $ A $ to a standard measurable space
$ B $ consists of a probability space $ (\Om,\F,P) $ and a family $
\Xi = (\Xi_a)_{a\in A} $ of equivalence classes $ \Xi_a $ of
measurable maps $ \Om \to B $ such that $ \F $ is generated by all $
\Xi_a $ (`non-redundancy').

(b)
The \emph{distribution}\index{distribution!of \Smap}
$ \La = \La_\Xi $ of an \Smap\ $ \Xi $
is the family of its finite-dimensional distributions, that is, the
joint distributions $ \la_{a_1,\dots,a_n} $ of \valued{B} random
variables $ \Xi_{a_1}, \dots, \Xi_{a_n} $ for all $ n = 1,2,\dots $
and $ a_1,\dots,a_n \in A $.
\end{definition}

Of course, two maps $ \Om \to B $ are called equivalent iff they are
equal almost everywhere on $ \Om $. As always, $ \Om $ is a standard
probability space. The non-redundancy can be enforced by replacing $
(\Om,\F,P) $ with its quotient space.

One may say that an \Smap\ $ A \to B $ is a \valued{B} random
process on $ A $, provided that all modifications are treated as the
same process.

Two \Smap s $ \Xi, \Xi' $ are identically distributed ($
\La_{\Xi} = \La_{\Xi'} $) if and only if they are isomorphic in the
following sense: there exists an isomorphism $ \al $ between the
corresponding probability spaces such that $ \Xi_a = \Xi'_a \circ \al
$ for $ a \in A $.

\begin{example}
A stationary Gaussian random process on $ \R $, continuous in
probability, may be treated as a special curve (`helix') in a Hilbert
space of Gaussian random variables. Depending on the covariance
function, sometimes it has continuous sample paths, sometimes not. In
the latter case we have no idea of a `favorite' modification, but
anyway, the corresponding \Smap\ from $ \R $ to $ \R $ is well-defined
(and continuous in probability), and its distribution is uniquely
determined by the covariance function.
\end{example}

\begin{example}
Skorokhod \cite[Sect.~1.1.1]{Sk84} defines a strong random operator on
a Hilbert space $ H $ as a continuous linear map from $ H $ into the
space of \valued{H} random variables. It may be treated as a linear
continuous \Smap\ $ H \to H $, but generally not a random linear
continuous operator $ H \to H $.

My term `\Smap' is derived from Skorokhod's `strong random
operator'; `S' may allude to `stochastic', `Skorokhod', or `strong'.
\end{example}

Every \Smap\ $ \Xi $ from $ A $ to $ B $ leads to a linear
operator $ \Tl_1^\Xi $\index{zzt@$ \Tl $, operator}
from the space of all bounded measurable
functions on $ B $ to the space of all bounded functions on $ A $;
namely,
\[
(\Tl_1^\Xi \phi) (a) = \Ex \phi(\Xi_a) \, .
\]
However, $ \Tl_1^\Xi $ involves only one-dimensional distributions of
$ \Xi $. Joint distributions are involved by operators $ \Tl_n^\Xi $
from bounded measurable functions on $ B^n $ to bounded functions on $
A^n $; here $ n \in \{1,2,\dots\} \cup \{\infty\} $, and $ A^\infty,
B^\infty $ consist of infinite sequences:
\begin{equation}\label{7c3a}
\begin{gathered}
(\Tl_n^\Xi \phi) (a_1,\dots,a_n) = \Ex \phi (\Xi_{a_1}, \dots,
 \Xi_{a_n}) \quad \text{for } n < \infty \, , \\
(\Tl_\infty^\Xi \phi) (a_1,a_2,\dots) = \Ex \phi (\Xi_{a_1},
 \Xi_{a_2}, \dots) \, .
\end{gathered}
\end{equation}
Clearly, the operator $ \Tl_n^\Xi $ determines uniquely (and is
determined by) the \dimensional{n} distributions $ \la_{a_1,\dots,a_n}
$ of $ \Xi $; thus, the distribution $ \La_\Xi $ of $ \Xi $
(determines and) is uniquely determined by the operators $ \Tl_1^\Xi,
\Tl_2^\Xi, \Tl_3^\Xi, \dots $ together (or equivalently, the operator
$ \Tl_\infty^\Xi $ alone). See also Sections \ref{7g}, \ref{7h} for a
description of the class of all operators of the form $ \Tl_\infty^\Xi
$.

Assume now that $ A $ is also a standard measurable space (like $ B
$).

\begin{proposition}\label{7c4}
The following two conditions are equivalent for every \Smap\ $
\Xi $ from $ A $ to $ B $:

(a)
the map $ a \mapsto \Xi_a $ is measurable from $ A $ to the space $
L_0(\Om\to B) $;

(b)
there exists a measurable function $ \xi : A \times \Om \to B $ such
that for every $ a \in A $ the function $ \om \mapsto \xi(a,\om) $
belongs to the equivalence class $ \Xi_a $.
\end{proposition}

\beginproof
We may assume that $ \Om = A = B = (0,1) $.

(b) \imp (a):
the set of all $ \xi $ satisfying (a) is closed under pointwise (on $
A \times \Om $) convergence.

(a) \imp (b):
we may take
\[
\xi(\om,a) = \limsup_{\eps\to0} \frac1{2\eps}
\int_{\om-\eps}^{\om+\eps} \Xi_a(\om) \, \D\om
\]
(which is just one out of many appropriate $ \xi $).
\proofend

See also \cite[Introduction (the proof of (II)\imp(I))]{GTW}.

\Smap s satisfying the equivalent conditions \ref{7c4}(a,b) will
be called \emph{measurable}.\index{measurable!\Smap}

Every measurable \Smap\ $ \Xi $ from $ A $ to $ B $ leads to an
operator $ \Tr_1^\Xi $\index{zzt@$ \Tr $, operator}
from probability measures on $ A $ to
probability measures on $ B $ (or rather, a linear operator on finite
signed measures), namely,
\[
\int_B \phi \, \D (\Tr_1^\Xi \nu) = \int_A \Ex \phi(\Xi_a) \, \nu(\D
a) = \int_A (\Tl_1^\Xi \phi) \, \D\nu
\]
for all bounded measurable $ \phi : B \to \R $. In other words, $
\Tr_1^\Xi \nu $ is the image of $ \nu \times P $  under the map $ \xi
: A \times \Om \to B $ corresponding to $ \Xi $ as in
\ref{7c4}(b). In fact, $ \Tr_1^\Xi \nu = \Ex \Xi_\nu $ where $ \Xi_\nu
$ is a \emph{random} measure, the image of $ \nu $ under the map $
\xi_\om : A \to B $, $ \xi_\om(a) = \xi(a,\om) $. Another choice of $
\xi $ (for the given $ \Xi $) may change $ \Xi_\nu $ only on a set of
zero probability. Similarly, for any measure $ \nu $ on $ A^n $, 
\[
\int_{B^n} \phi \, \D (\Tr_n^\Xi \nu) = \int_{A^n} \Ex
\phi(\Xi_{a_1},\dots,\Xi_{a_n}) \, \nu(\D a_1 \dots \D a_n) =
\int_{A_n} (\Tl_n^\Xi \phi) \, \D\nu
\]
for $ \phi : B^n \to \R $. The same holds for $ n=\infty $.

Let $ A,B,C $ be three standard measurable spaces and $ \Xi', \Xi'' $
be measurable \Smap s, $ \Xi' $ from $ A $ to $ B $, $ \Xi'' $
from $ B $ to $ C $, on probability spaces $ \Om', \Om'' $
respectively. Their \emph{composition}%
\index{composition!of \Smap s}
$ \Xi = \Xi' \Xi'' $ is a
measurable \Smap\ from $ A $ to $ C $ on $ \Om = \Om' \times \Om'' $
(or rather its quotient space, for non-redundancy), defined as
follows:
\[
\xi \( a, (\om',\om'') \) = \xi'' \( \xi'(a,\om'), \om'' \) \quad
\text{for } \om' \in \Om', \, \om'' \in \Om'' \, ,
\]
where $ \xi : A \times (\Om'\times\Om'') \to C $, $ \xi' : A \times
\Om' \to B $, $ \xi'' : B \times \Om'' \to C $ correspond to $ \Xi,
\Xi', \Xi'' $ as in \ref{7c4}(b). The composition is well-defined by
the next lemma.

\begin{lemma}\label{7c5}
\cite[Prop.~1.2.2/1.1]{LJR}
The composition $ \Xi = \Xi' \Xi'' $ does not depend on the choice of
$ \xi', \xi'' $.
\end{lemma}

\beginproof
For a given $ a $, a change of $ \xi' $ influences $
\Xi_a(\cdot,\cdot) $ on a subset of $ \Om_1 \times \Om_2 $,
negligible, since its first projection is negligible.
A change of $ \xi'' $ influences $
\Xi_a(\cdot,\cdot) $ on a subset of $ \Om_1 \times \Om_2 $,
negligible, since all its sections ($ \om' = \const $) are
negligible.
\proofend

The distribution $ \La $ of $ \Xi = \Xi' \Xi'' $ is uniquely
determined by the distributions $ \La', \La'' $ of $ \Xi' $ and $
\Xi'' $ (since distributions correspond to isomorphic classes), and
will be called the \emph{convolution}%
\index{convolution!for \Smap s}
of these two distributions, $
\La = \La' * \La'' $. It is easy to see that $ \Xi = \Xi' \Xi'' $
implies
\[
\Tl^\Xi_n \phi = \Tl^{\Xi'}_n ( \Tl^{\Xi''}_n \phi) \quad \text{and}
\quad \Tr^\Xi_n \nu = \Tr^{\Xi''}_n ( \Tr^{\Xi'}_n \nu) \, .
\]

We may treat $ \Xi' $ and $ \Xi'' $ as independent \Smap s on the same
probability space $ \Om $; namely, $ \Xi'_a (\om',\om'') = \Xi'_a
(\om') $ and $ \Xi''_b (\om',\om'') = \Xi''_b (\om'') $. Conditional
expectations are given by
\begin{gather}
\cE{ \phi(\Xi_a) }{ \Xi' } = \( \Tl_1^{\Xi''} \phi \) ( \Xi'_a ) \, ,
 \label{7c7} \\
\cE{ \phi(\Xi_a) }{ \Xi'' } = \int_B \phi(\Xi''_b) \, \( \Tr_1^{\Xi'}
\de_a \) (\D b) \label{7c8}
\end{gather}
for all bounded measurable $ \phi : C \to \R $, where $ \de_a $ is the
probability measure concentrated at $ a $, and $ \Xi''_b $ means $
\xi''(b,\cdot) $; the choice of $ \xi'' $ does not matter (similarly
to \ref{7c5}). More generally,
\begin{align}
\cE{ \phi(\Xi_{a_1},\dots,\Xi_{a_n}) }{ \Xi' } &= \( \Tl_n^{\Xi''} \phi
 \) ( \Xi'_{a_1},\dots,\Xi'_{a_n} ) \, , \label{7c9} \\
\cE{ \phi(\Xi_{a_1},\dots,\Xi_{a_n}) }{ \Xi'' } &= \label{7c10} \\
 = \int_{B^n} \phi(\Xi''_{b_1},\dots,\Xi''_{b_n}) & \, \(
 \Tr_n^{\Xi'} \de_{a_1,\dots,a_n} \) (\D b_1 \dots \D b_n) \notag
\end{align}
for all bounded measurable $ \phi : C^n \to \R $ and $ n < \infty $;
similar formulas hold for $ n=\infty $.

See also \cite[Sect.~1]{LJR} and \cite[Sect.~8d]{Ts03}.

\subsection[Flow systems of S-maps]{Flow systems of \Smap s}
\label{7d}

Let $ \X $ be a compact metrizable space (mostly, the circle will be
used). Then $ L_0(\Om\to\X) $ is equipped with the (metrizable)
topology of convergence in probability.

\begin{definition}\label{7d1}
An \Smap\ $ \Xi $ from $ \X $ to $ \X $ is \emph{continuous in
probability,}%
\index{continuous in probability!\Smap}
if the map $ x \mapsto \Xi_x $ is continuous from $ \X
$ to $ L_0(\Om\to\X) $.
\end{definition}

Clearly, \ref{7d1} is stronger than \ref{7c4}(a). Continuity in
probability is preserved by the composition, which is made clear by
the next lemma.

\begin{lemma}\label{7d2}
The following three conditions are equivalent for every \Smap\ $ \Xi $
from $ \X $ to $ \X $:

(a) $ \Xi $ is continuous in probability;

(b) $ \Tl_2^\Xi \phi \in C(\X^2) $ for all $ \phi \in C(\X^2) $;

(c) $ \Tl_n^\Xi \phi \in C(\X^n) $ for all $ \phi \in C(\X^n) $ and all $
n = 1,2,3,\dots $
\end{lemma}

\beginproof
(a) \imp (c): the map $ (x_1,\dots,x_n) \mapsto
(\Xi_{x_1},\dots,\Xi_{x_n}) $ from $ \X^n $ to $ L_0(\Om\to\X^n) $
is continuous, therefore $ \Ex \phi (\Xi_{x_1},\dots,\Xi_{x_n}) $ is
continuous in $ x_1,\dots,x_n $.

(c) \imp (b): trivial.

(b) \imp (a): let $ \phi $ be the metric, $ \phi(x_1,x_2) =
\dist(x_1,x_2) $, then $ (\Tl_2^\Xi \phi) (x_1,x_2) = \Ex \dist
(\Xi_{x_1},\Xi_{x_2}) = 0 $ on the diagonal $ x_1 = x_2 $; by (b), $
\Ex \dist (\Xi_{x_1},\Xi_{x_2}) \to 0 $ as $ \dist(x_1,x_2) \to 0 $,
which is (a).
\proofend

A Lipschitzian version, given below, will be used in Sect.~\ref{7f}.
By $ \Lip (\phi) $ we denote the least $ C \in [0,\infty] $ such that $ \dist
( \phi(x_1), \phi(x_2) ) \le C \dist (x_1,x_2) $ for all $ x_1, x_2 $; $ \phi
$ is a Lipschitz function iff $ \Lip(\phi) < \infty $.

\begin{lemma}\label{7d3}
The following three conditions are equivalent for every \Smap\ $ \Xi $
from $ \X $ to $ \X $:

(a)
there exists $ C $ such that $ \Ex \dist (\Xi_{x_1},\Xi_{x_2}) \le C
\dist(x_1,x_2) $ for all $ x_1,x_2 \in \X $;

(b)
there exists $ C_2 $ such that $ \Lip (\Tl_2^\Xi \phi) \le C_2 \Lip(\phi) $
for all Lipschitz functions $ \phi : \X^2 \to \R $;

(c)
for each $ n=1,2,3,\dots $ there exists $ C_n $ such that $ \Lip
(\Tl_n^\Xi \phi) \le C_n \Lip(\phi) $ for all Lipschitz functions $ \phi :
\X^n \to \R $.
\end{lemma}

\begin{sloppypar}
(Any `reasonable' metric on $ \X^n $ may be used, say, $ \dist \(
(x'_1,\dots,x'_n), \linebreak[0]
(x''_1,\dots,x''_n) \) = \sum_k \dist (x'_k,x''_k) $, or $ \max_k
\dist (x'_k,x''_k) $.)
\end{sloppypar}

\beginproof
(c) \imp (b) \imp (a): trivial;

(a) \imp (c): the same as `(b) \imp (a) \imp (c)' in the proof of
\ref{7d2}, but quantitative.
\proofend

A \emph{convolution system of \Smap s}%
\index{convolution system!of \Smap s}
(over $ \X $) may be
defined as a family $ (\La_{s,t})_{s<t} $, where each $ \La_{s,t} $ is
the distribution of an \Smap\ from $ \X $ to $ \X $,
continuous in probability, and
\[
\La_{r,t} = \La_{r,s} * \La_{s,t}
\]
whenever $ r<s<t $.

Every convolution system of \Smap s leads to a convolution
system as defined by \ref{2a1}. Namely, each $ \La_{s,t} $ leads to a
probability space $ (G_{s,t}, \mu_{s,t}) $\footnote{%
 Not $ (\Om_{s,t},P_{s,t}) $ for conformity to Sect.~\ref{2a}.}
carrying an \Smap\ $ \Xi_{s,t} = ( \Xi^{s,t}_x)_{x\in\X} $, $
\Xi^{s,t}_x \in L_0(G_{s,t}\to\X) $ and unique up to isomorphism
(between \Smap s). Given $ r<s<t $ we have $ \Xi_{r,t} =
\Xi_{r,s} \Xi_{s,t} $ up to isomorphism, which gives us a
representation of $ G_{r,t} $ as a quotient space of $ G_{r,s} \times
G_{s,t} $, that is, a morphism $ G_{r,s} \times G_{s,t} \to G_{r,t}
$. The convolution system $ (G_{s,t},\mu_{s,t})_{s<t} $ is determined
by $ (\La_{s,t})_{s<t} $ uniquely up to isomorphism. If it is
separable (as defined by \ref{2a5}) then it leads to a flow system,
that is, all \Smap s $ \Xi_{s,t} = (\Xi_x^{s,t})_x $ may be
defined on a single probability space $ (\Om,P) $, satisfying (recall
\ref{2a3}(a,b))
\begin{equation}\label{7d4}
\begin{gathered}
\Xi_{t_1,t_2}, \Xi_{t_2,t_3}, \dots, \Xi_{t_{n-1},t_n} \text{ are
 independent for } t_1 < t_2 < \dots < t_n \, ; \\
\Xi_{r,t} = \Xi_{r,s} \Xi_{s,t} \quad \text{for } r<s<t \, .
\end{gathered}
\end{equation}
According to Sect.~\ref{2b}, the flow system leads to a continuous
product of probability spaces.

A sufficient condition for the separability is, \emph{temporal
continuity in probability}\index{temporal continuity}
(in addition to the spatial continuity in
probability assumed before for each $ \Xi_{s,t} $):
\begin{equation}\label{7d5}
\text{both } \Xi_x^{s-\eps,s} \text{ and } \Xi_x^{s,s+\eps} \text{
converge to $ x $ in probability as $ \eps\to0+ $} \, ,
\end{equation}
for all $ s\in\R $ and $ x \in \X $. It involves only one-dimensional
distributions, $ \la_x^{s,t} $, and may be reformulated in terms of
the operators $ T_1^{s,t} = \Tl_1^{\Xi_{s,t}} : C(\X) \to C(\X) $,
namely,
\[
\text{both } T_1^{s-\eps,s} (\phi) \text{ and } T_1^{s,s+\eps} (\phi)
\text{ converge to $ \phi $ pointwise as $ \eps\to0+ $} \, ,
\]
for all $ s \in \R $ and $ \phi \in C(\X) $.

\begin{lemma}\label{7d6}
Condition \eqref{7d5} implies separability.
\end{lemma}

\beginproof
It is sufficient (recall the end of Sect.~\ref{sec:2}) to prove that $
\Xi_x^{s-\eps,t} \to \Xi_x^{s,t} $ and $ \Xi_x^{s,t+\eps} \to \Xi_x^{s,t} $ in
probability as $ \eps \to 0 $, for $ x \in \X $ and $ s<t $. We have
\[
\Ex \dist ( \Xi_x^{s-\eps,t}, \Xi_x^{s,t} ) = \int_\X \( \Ex \dist (
\Xi_{x'}^{s,t}, \Xi_x^{s,t} ) \) \, \( \Tr_1^{\Xi_{s-\eps,s}} \de_x \)
(\D x') \, .
\]
By \eqref{7d5}, $ \Tr_1^{\Xi_{s-\eps,s}} \de_x \to \de_x $. By
continuity in probability, $ \Ex \dist ( \Xi_{x'}^{s,t}, \Xi_x^{s,t} )
\to 0 $ as $ x' \to x $. Therefore $ \Ex \dist ( \Xi_x^{s-\eps,t},
\Xi_x^{s,t} ) \to 0 $.

Also,
\[
\Ex \dist ( \Xi_x^{s,t+\eps}, \Xi_x^{s,t} ) = \int_\X \( \Ex \dist (
\Xi_{x'}^{t,t+\eps}, x' ) \) \, \( \Tr_1^{\Xi_{s,t}} \de_x \) (\D x')
\to 0 \, ,
\]
since the integrand converges to $ 0 $ for every $ x' $.
\proofend

\begin{remark}\label{7d7}
The continuous product of probability spaces, corresponding to a flow
system $ (\Xi_{s,t})_{s<t} $ of \Smap s from $ \X $ to $ \X $, is
classical if and only if random variables $ \phi(\Xi_x^{s,t}) $ are
stable for all $ s<t $ and all bounded Borel functions $ \phi : \X \to
\R $ (a single $ \phi $ is enough if it is one-to-one). Proof:
similar to \ref{5c5}.
\end{remark}

See also \cite[Sect.~1]{LJR}.

\subsection[From S-maps to black noise]{From \Smap s to black noise}
\label{7e}

Let $ \X $, $ (G_{s,t},\mu_{s,t})_{s<t} $ and $ (\Xi_{s,t})_{s<t} $
be as in Sect.~\ref{7d}, satisfying \eqref{7d5}.

Sect.~\ref{7b} gives us a sufficient condition
\eqref{7b2}-\eqref{7b3}-\eqref{7b4} (in combination with
\ref{6a6}(a--c)) for triviality of the classical part of the
corresponding continuous product of probability spaces.

Instead of all $ f \in L_\infty(G_{s,t}) $ we may use a subset of $
L_\infty(G_{s,t}) $ dense in $ L_2(G_{s,t}) $.

The \sif\ on $ G_{r,t} $ is generated by \valued{\X} random
variables $ \Xi_x^{r,t} $ for $ x \in \X $. Therefore functions of
the form
\begin{equation}\label{7e1}
\phi ( \Xi_{x_1}^{r,t}, \dots, \Xi_{x_n}^{r,t} ) \quad \text{for $ \phi \in
C(\X^n) $, $ (x_1,\dots,x_n) \in \X_n $}
\end{equation}
and $ n=1,2,\dots $ are dense in $ L_2 (G_{r,t}) $.

\begin{lemma}\label{7e2}
Assume that $ n \in \{1,2,\dots\} $ is given, a linear subset $ F_n $
of $ C(\X^n) $, dense in the norm topology, and a linear subset $ N_n
$ of the space of (finite, signed) measures on $ \X^n $, dense in the
weak topology; and an interval $ (r,t) \subset \R $. Then functions of
the form
\[
\int_{\X^n} \phi ( \Xi_{x_1}^{r,t}, \dots, \Xi_{x_n}^{r,t} ) \, \nu
(\D x_1 \dots \D x_n) \quad \text{for } \phi \in F_n \text{ and } \nu
\in N_n
\]
are $ L_2 $-dense among functions of the form \eqref{7e1} for the
given $ n $.
\end{lemma}

\beginproof
For every $ \phi \in F $, the map $ \X^n \to L_2(G_{r,t}) $ given by $
(x_1,\dots,x_n) \linebreak[0]
 \mapsto \phi ( \Xi_{x_1}^{r,t}, \dots, \Xi_{x_n}^{r,t}
) $ is continuous (since $ \Xi_{r,t} $ is continuous in
probability). Therefore, $ \nu_k \to \de_{x_1,\dots,x_n} $ (as $ k
\to \infty $) implies
\[
\quad
\int_{\X^n} \phi ( \Xi_{x'_1}^{r,t}, \dots, \Xi_{x'_n}^{r,t} ) \,
\nu_k (\D x'_1 \dots \D x'_n) \to \phi ( \Xi_{x_1}^{r,t}, \dots,
\Xi_{x_n}^{r,t} ) \quad \text{in } L_2(G_{r,t}) \, . \quad \qed
\]
\proofendnoqed

Substituting $ f = \int \phi ( \Xi_{x_1}^{r,t}, \dots, \Xi_{x_n}^{r,t}
) \, \nu (\D x_1 \dots \D x_n) $ to \eqref{7b3} we get
\begin{equation}\label{7e3}
f_{s,s+\eps} = \int_{\X^n} \( \Tl_n^{\Xi_{s+\eps,t}} \phi \) (
\Xi_{x_1}^{s,s+\eps}, \dots, \Xi_{x_n}^{s,s+\eps} ) \, \(
\Tr_n^{\Xi_{r,s}} \nu \) (\D x_1 \dots \D x_n) \, .
\end{equation}

\begin{proposition}\label{7e4}
(Le Jan and Raimond; implicit in \cite{LJR1}.)
Let a flow system $ (\Xi_{s,t})_{s<t} $ of \Smap s from $ \X $ to
$ \X $ and a probability measure $ \nu_0 $ on $ \X $ satisfy the
conditions

(a) (stationarity)
the distribution of $ \Xi_{s,s+h} $ does not depend on $ s $;

(b) (invariant measure) $ \Tr_1^{\Xi_{s,t}} \nu_0 = \nu_0 $ for $ s<t
$;

(c) (Lipschitz boundedness)
if $ \phi \in C(\X^n) $ is a Lipschitz function, then $
\Tl_n^{\Xi_{s,t}} \phi $ is also a Lipschitz function, with a Lipschitz
constant
\[
\Lip (\Tl_n^{\Xi_{s,t}} \phi) \le C_n \Lip(\phi) \quad \text{for } 0 \le s <
t \le 1 \, ,
\]
where $ C_n < \infty $ depends only on $ n $;

(d)
$\displaystyle
 \sup_{\phi,\nu} \Var \bigg( \int_{\X^n} \phi \(
 \Xi^{0,\eps}_{x_1},\dots,\Xi^{0,\eps}_{x_n} \) \, \nu (\D x_1 \dots \D
 x_n) \bigg) = o(\eps) \quad \text{as } \eps \to 0
$,\linebreak
where the supremum is taken over all $ \phi \in C(\X^n) $ such that $
\Lip(\phi) \le 1 $ and all positive measures $ \nu $ on $ \X^n $ such
that $ \nu_1 \le \nu_0, \dots, \nu_n \le \nu_0 $; here $
\nu_1,\dots,\nu_n $ are coordinate projections of $ \nu $, that is, $
\int \phi(x_k) \, \nu(\D x_1\dots\D x_n) = \int \phi(x) \, \nu_k(\D x) $ for
$ \phi \in C(\X) $.

Then the corresponding noise is black.
\end{proposition}

\beginproof
The continuity condition \ref{6a6}(a--c) is ensured by (a). By Lemma
\ref{7e2} it is sufficient to verify \eqref{7b2} in the form
\eqref{7b4} when $ f_{s,s+\eps} $ is given by \eqref{7e3}, assuming
that $ \phi $ and $ \nu $ have the properties formulated in (d). In order
to apply (d) to the function $ \phi_{s+\eps} = \Tl_n^{\Xi_{s+\eps,t}} \phi $
(instead of $ \phi $) and the measure $ \nu_s = \Tr_n^{\Xi_{r,s}} \nu $
(instead of $ \nu $) we only need to check that these properties of $
\phi $ and $ \nu $ are inherited by $ \phi_{s+\eps} $ and $ \nu_s $. The
Lipschitz property of $ \phi_{s+\eps} $ follows from (c), up to the
(harmless) constant $ C_n $. The property of $ \nu_s $ (majorization
of its coordinate projections) follows from (b).
\proofend

\mysubsection{Example: Arratia's coalescing flow, or the Brownian web}
\label{7f}

\[
\setlength{\unitlength}{0.8cm}\nopagebreak[4]
\begin{picture}(12,3)
\put(0,0.5){\includegraphics[scale=0.8]{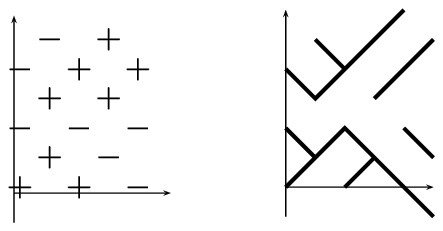}}
\put(1,0){\makebox(0,0){(a)}}
\put(3.8,0){\makebox(0,0){(b)}}
\put(6,0.5){\includegraphics[scale=0.8]{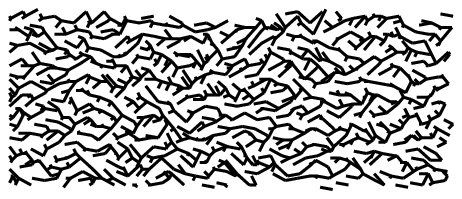}}
\put(8.5,0){\makebox(0,0){(c)}}
\end{picture}
\]
A two-dimensional array of random signs (a) produces a system of
coalescing random walks (b) that converges in the scaling limit (c) to
a flow system of \Smap s, introduced by Arratia in 1979 \cite{Arr} and
investigated further by T\'oth, Werner, Soucaliuc \cite{TW},
\cite{STW}, Fontes, Isopi, Newman, Ravishankar \cite{FINR}, Le Jan and
Raimond \cite[Sect.~3]{LJR1}. It consists of infinitely many
coalescing Brownian motions, independent before coalescence.
Our approach, based on \Smap s, deals with equivalence classes $
\Xi_x^{s,t} $ rather than sample functions $ (x,s,t) \mapsto
\xi_{s,t}(\om,x) $ (recall \ref{7c4}(b)); fine properties of sample
functions, examined in some of the works cited above, are irrelevant
here.

The space $ \X $ is the circle $ \T = \R/\Z $. The distribution of $
\Xi_{s,s+t} $ does not depend on $ s $ (stationarity). The
(one-dimensional) distribution of $ \Xi_x^{0,t} $ is the normal
distribution $ N(x,t) $ (or rather, the distribution of $ ( x + \zeta )
\bmod 1 $ where $ \zeta \sim N(0,t) $). It may be thought of as the
distribution of $ x + B_t $ where $ (B_s)_s $ is the (standard)
Brownian motion in $ \T $. The (two-dimensional) joint distribution of
$ \Xi^{0,t}_{x_1} $ and $ \Xi^{0,t}_{x_2} $ is the joint distribution
of $ x_1 + B_t^{(1)} $ and $ x_k + B_t^{(k)} $ where $ (B_s^{(1)})_s,
(B_s^{(2)})_s $ are two independent Brownian motions in $ \T $, and $
k : \Om \to \{1,2\} $ is a random variable, defined as follows:
\begin{equation}\label{7f1}
k = \begin{cases}
 1 &\text{if $ \min \{ s : x_1 + B_s^{(1)} = x_2 + B_s^{(2)} \} \le t
  $}, \\
 2 &\text{otherwise}.
\end{cases}
\end{equation}
That is, the second Brownian motion $ \( x_{k(s)} + B_s^{(k(s))} \)_s
$ is independent of the first one, as long as they do not meet.
Afterwards they are equal. In spite of the asymmetry (the second 
motion joins the first), the resulting distribution does not depend on
the order of initial points. Joint distributions of higher dimensions
are defined similarly. We have
\begin{equation}\label{7f2}
\Ex \dist \( \Xi^{0,t}_{x_1}, \Xi^{0,t}_{x_2} \) \le \dist(x_1,x_2) \,
\end{equation}
since the process $ t \mapsto \dist \( \Xi^{0,t}_{x_1},
\Xi^{0,t}_{x_2} \) $ is a supermartingale. (On $ \R $ it would be a
martingale, but we are on the circle.) Thus, for every $ t $ the
\Smap\ $ \Xi_{0,t} $ is continuous in probability. The temporal
continuity in probability \eqref{7d5} is evident. The uniform
distribution $ \nu_0 $ on the circle evidently is an invariant measure
in the sense of \ref{7e4}(b).

By \eqref{7f2}, $ \Xi_{0,t} $ satisfies \ref{7d3}(a), which implies
\ref{7d3}(c), the constants $ C_n $ not depending on $ t $. Thus,
\ref{7e4}(c) holds (as well as \ref{7e4}(a,b)). In order to get a
black noise, it remains to verify \ref{7e4}(d).

\begin{sloppypar}
Let us start with the case $ n=1 $. We consider $ \Var \( \int
\phi(\Xi^{0,\eps}_x) \, \nu(\D x) \) $ assuming $ \Lip(\phi) \le 1 $ and $ 0
\le \nu \le \nu_0 $. Note that
\[
\Var \bigg( \int \dots \bigg) = \iint \Cov \( \phi(\Xi^{0,\eps}_{x_1}),
\phi(\Xi^{0,\eps}_{x_2}) \) \, \nu(\D x_1) \nu(\D x_2) \, .
\]
On one hand,
\[
| \Cov \( \phi(\Xi^{0,\eps}_{x_1}), \phi(\Xi^{0,\eps}_{x_2}) \) | \le \sqrt{
\Var \phi(\Xi^{0,\eps}_{x_1}) } \sqrt{ \Var \phi(\Xi^{0,\eps}_{x_2}) } \le
\eps
\]
since $ | \phi(\Xi^{0,\eps}_x) - \phi(x) | \le | \Xi^{0,\eps}_x - x | $ and
$ \Ex | \Xi^{0,\eps}_x - x |^2 \le \eps $. (The latter would be `$
=\eps $' on $ \R $, but we are on the circle.) On the other hand, we
may assume that $ \| \phi \|_{C(\X)} \le \const \cdot \Lip(\phi) \le \const
$; using \eqref{7f1},
\begin{multline*}
| \Cov \( \phi(\Xi^{0,\eps}_{x_1}), \phi(\Xi^{0,\eps}_{x_2}) \) | = 
 | \Cov \( \phi(x_1 + B_\eps^{(1)}), \phi(x_k + B_\eps^{(k)}) \) | = \\
= | \Cov \( \phi(x_1 + B_\eps^{(1)}), \phi(x_k + B_\eps^{(k)}) - \phi(x_2 +
 B_\eps^{(2)}) \) | \le \| \phi \|_{C(\X)}^2 \Pr{k=1} \, .
\end{multline*}
However, the probability of meeting, $ \Pr{k=1} $, is (exponentially)
small for $ \dist(x_1,x_2) \gg \sqrt\eps $, therefore $
\iint_{\dist(x_1,x_2)\ge\de} |\Cov(\dots)| \, \nu(\D x_1) \nu(\D x_2)
$ is $ o(\eps) $ (in fact, exponentially small) as $ \eps \to 0 $, for
every $ \de > 0 $. Also, $ \iint_{\dist(x_1,x_2)\le\de} |\Cov(\dots)|
\, \nu(\D x_1) \nu(\D x_2) \le \iint_{\dist(x_1,x_2)\le\de} \eps \,
\nu_0(\D x_1) \nu_0(\D x_2) = { \eps \cdot 2\de } $. It follows that $
\Var \( \int \dots \) = o(\eps) $ as $ \eps \to 0 $, uniformly in $
\phi $ and $ \nu $.
\end{sloppypar}

Generalization for $ n=2,3,\dots $ is straightforward. One estimates
the\linebreak
$ { (\nu \times \nu) } $\nobreakdash-\hspace{0pt}measure of the set $ \{
\( (x'_1,\dots,x'_n), (x''_1,\dots,x''_n) \) \in \X^n \times \X^n :
\dist(x'_k,x''_l)\linebreak[0]
\le\de \} $ for each pair $ (k,l) $ separately
(taking into account that $ \nu_k \le \nu_0 $, $ \nu_l \le \nu_0 $),
and consider the union of these $ n^2 $ sets.

By Proposition \ref{7e4}, the noise corresponding to Arratia's
coalescing flow is black.

The proof presented above, due to Le Jan and Raimond
\cite[Sect.~3]{LJR1}, is simpler than \cite[Sect.~7]{Ts03}.

\mysubsection{Random kernels}
\label{7g}

By a \emph{kernel}\index{kernel}
from a set $ A $ to a standard measurable space $ B
$ we mean a map $ A \to \Pc(B) $. Here $ \Pc(B) $%
\index{zzp@$ \Pc(B) $, space}
is the standard
measurable space of all probability measures on $ B $, equipped with
the \sif\ generated by the functions $ \Pc(B) \to \R $ of the form $
\mu \to \int \phi \, \D\mu $, where $ \phi $ runs over bounded
measurable functions $ B \to \R $, see \cite[Sect.~17.E]{Ke}.

\begin{definition}\label{7g1}
An \emph{\Skernel}\index{Skernel@\Skernel}
from a set $ A $ to a standard measurable space $ B
$ is an \Smap\ from $ A $ to $ \Pc(B) $.
\end{definition}

This idea was introduced by Le Jan and Raimond in order to describe
``turbulent evolutions where [\dots] two points thrown initially at
the same place separate'' \cite[Introduction]{LJR}.

Note that $ B $ is naturally embedded into $ \Pc(B) $ (by $ b \mapsto
\de_b $, the measure concentrated at $ b $). Accordingly, a map $ A
\to B $ may be treated as a special case of a kernel, $ A \to B
\subset \Pc(B) $. Similarly, an \Smap\ $ \Xi $ (from $ A $ to $ B $)
may be treated as a special case of an \Skernel\ $ K $ (from $ A $ to
$ B $); namely, $ K_a = \de_{\Xi_a} $.

By the distribution $ \La_K $ of an \Skernel\ $ K $ from $ A $ to $ B
$ we mean the distribution of $ K $ as an \Smap\ from $ A $ to $
\Pc(B) $; it consists (recall \ref{7c1}(b)) of the joint distributions
$ \la_{a_1,\dots,a_n} $ of \valued{\Pc(B)} random variables $ K_{a_1},
\dots, K_{a_n} $.

Usually it is difficult to construct an \Skernel\ (or a flow of
\Skernel s) directly, by specifying joint distributions of
\emph{measures} (or corresponding
infinitesimal\nobreakdash-\hspace{0pt}time data). It is easier to do
it indirectly, by specifying joint distributions of \emph{points} and
using a moment method described below.

Let $ K $ be an \Skernel\ from $ A $ to $ B $. Combining formally
\ref{7g1} and \eqref{7c3a} one could treat $ \Tl_1^K $ as defined on
the (huge) space of functions on $ \Pc(B) $, but we prefer it to be
defined on the same (modest) space as $ \Tl_1^\Xi $ in
Sect.~\ref{7c}. Namely, we define a linear operator $ \Tl_1^K $ from
the space of all bounded measurable functions on $ B $ to the space of
all bounded functions on $ A $ by
\[
(\Tl_1^K \phi) (a) = \Ex \int_B \phi \, \D K_a \, .
\]
In other words, we restrict ourselves to \emph{linear} functions on $
\Pc(B) $, $ \mu \mapsto \int \phi \, \D\mu $. For an \Smap\ $ \Xi $
from $ A $ to $ B $, treated as (a special case of) an \Skernel\ $ K
$, we have $ \Tl_1^\Xi = \Tl_1^K $, since $ \int \phi \, \D\de_{\Xi_a}
= \phi(\Xi_a) $. Generally, for $ n \in \{ 1,2,\dots \} \cup
\{\infty\} $ we define a linear operator $ \Tl_n^K $ from bounded
measurable functions on $ B^n $ to bounded functions on $ A^n $ by 
\begin{equation}\label{7g2}
\begin{gathered}
(\Tl_n^K \phi) (a_1,\dots,a_n) = \Ex \int_{B^n} \phi \, \D ( K_{a_1}
 \times\dots\times K_{a_n} ) \quad \text{for } n < \infty \, , \\
(\Tl_\infty^K \phi) (a_1,a_2,\dots) = \Ex \int_{B^\infty} \phi \, \D
 (K_{a_1} \times K_{a_2} \times \dots ) \, .
\end{gathered}
\end{equation}
(Of course, $ \int \phi \, \D ( K_{a_1} \times\dots\times K_{a_n} ) $
means $ \int \phi(b_1,\dots,b_n) \, K_{a_1}(\D b_1) \dots K_{a_n}(\D
b_n) $.) These expectations of multilinear functions of $ K_{a_1},
\dots, K_{a_n} $ are sometimes called the moments%
\index{moments of \Skernel s}
of $ K $. A solution
of the corresponding moment problem is given below, see \ref{7g3}
(uniqueness) and \ref{7g6}, \ref{7h3} (existence).

\begin{lemma}\label{7g3}
The distribution $ \La_K $ is uniquely
determined by the operators $ \Tl_n^K $, $ n=1,2,3,\dots $ (or
equivalently, by a single operator $ \Tl^K_\infty $).
\end{lemma}

\beginproof
The $ n $-th moment of the (bounded) random variable $ \int \phi \, \D
K_a $ is equal to
\[
\Ex \int_{B^n} \phi(b_1) \dots \phi(b_n) \, K_a (\D b_1) \dots
K_a (\D b_n) = \( \Tl_n^K (\phi\otimes\dots\otimes\phi) \)
(a,\dots,a) \, ;
\]
thus, the distribution of $ \int \phi \, \D K_a $ is determined
uniquely via its moments. Similarly, the joint distribution of $ \int
\phi_1 \, \D K_{a_1} $ and $ \int \phi_2 \, \D K_{a_2} $ is determined
via its (mixed) moments; and so on.
\proofend

The moments are basically the same as (non-random) kernels $ T_n^K $
from $ A^n $ to $ B^n $ defined by
\[
T_n^K (a_1,\dots,a_n) = \Ex (K_{a_1}\times\dots\times K_{a_n}) \, ,
\]
that is,
\[
\int \phi \, \D \( T_n^K (a_1,\dots,a_n) \) = \Ex \int \phi \,
\D(K_{a_1}\times\dots\times K_{a_n}) = ( \Tl_n^K \phi) (a_1,\dots,a_n)
\]
for bounded measurable $ \phi : B^n \to \R $. Also $ n=\infty $ is
admitted,
\[
T_\infty^K (a_1,a_2,\dots) = \Ex (K_{a_1}\times K_{a_2}\times\dots) \,
.
\]
Each $ T_n^K $ is a marginal of $ T_\infty^K $ in the sense that
\begin{multline*}
\int \phi(b_1,\dots,b_n) \, \( T_\infty^K (a_1,a_2,\dots) \) (\D b_1
 \D b_2 \dots) = \\
= \int \phi(b_1,\dots,b_n) \, \( T_n^K (a_1,\dots,a_n) \) (\D b_1
 \dots \D b_n) \, ;
\end{multline*}
similarly, $ T_n^K $ is a marginal of $ T_{n+1}^K $
(\emph{consistency}). Everyone knows that a probability distribution on $
B^\infty $ is basically the same as a consistent family of probability
distributions on $ B^n $, $ n<\infty $. Accordingly, a consistent
family of kernels from $ A^n $ to $ B^n $ is basically the same as a
kernel from $ A^\infty $ to $ B^\infty $ satisfying the condition
\begin{equation}\label{7g4}
\int \phi(b_1,\dots,b_n) \, T(a_1,a_2,\dots) (\D b_1 \D b_2 \dots)
\quad \text{depend on $ a_1,\dots,a_n $ only.}
\end{equation}

Measures $ T_2^K (a_1,a_2) $ and $ T_2^K (a_2,a_1) $ are mutually
symmetric: $ T_2^K (a_2,a_1) \linebreak[0]
(\D b_2 \D b_1) = T_2^K (a_1,a_2) (\D b_1
\D b_2) $; more formally, $ \int \phi(b_2,b_1) \, T_2^K (a_2,a_1) (\D
b_1 \D b_2) = \int \phi(b_1,b_2) \, T_2^K (a_1,a_2) (\D b_1 \D b_2)
$. Similarly,
\[
T_n^K ( a_{\si(1)}, \dots, a_{\si(n)} ) ( \D b_{\si(1)} \dots \D
b_{\si(n)} ) \quad \text{does not depend on $ \si $} \, ,
\]
$ \si $ being a permutation, $ \si : \{1,\dots,n\} \to \{1,\dots,n\} $
bijectively. Also $ n=\infty $ is admitted,
\begin{equation}\label{7g5}
T_\infty^K ( a_{\si(1)}, a_{\si(2)}, \dots ) ( \D b_{\si(1)} \D
b_{\si(2)} \dots ) \quad \text{does not depend on $ \si $} \, ,
\end{equation}
$ \si : \{1,2,\dots\} \to \{1,2,\dots\} $ bijectively.

\begin{lemma}\label{7g6}
Let $ A $ be a finite or countable set, and $ B $ a standard
measurable space. A kernel $ T $ from $ A^\infty $ to $ B^\infty $ is
of the form $ T = T_\infty^K $ for some \Skernel\ $ K $ from $ A $ to
$ B $ if and only if $ T $ satisfies \eqref{7g4} and \eqref{7g5}.
\end{lemma}

\beginproof
We know that $ T_\infty^K $ satisfies \eqref{7g4}, \eqref{7g5}.

Assume that $ T $ satisfies \eqref{7g4}, \eqref{7g5}. For every $ a
\in A $ the measure $ T(a,a,\dots) $ on $ B^\infty $ is invariant
under permutations. The general form of such a measure is well-known
(de Finetti type theorem on exchangeability, see \cite[Th.~4.2]{Link}),
it is $ \int (\nu\times\nu\times\dots) \, \mu_a(\D \nu) $, the mixture
of product measures $ \nu\times\nu\times\dots $ over $ \nu \in \Pc(B)
$ distributed according to some (uniquely determined) measure $ \mu_a
\in \Pc(\Pc(B)) $. The distribution $ \mu_a $ of a single
\valued{\Pc(B)} random variable $ K_a $ is thus constructed.

Given $ a_1, a_2 \in A $, the measure $ T(a_1,a_2,a_1,a_2,\dots) $ on
$ (A\times A)^\infty = A^\infty \times A^\infty $ is invariant under
(the product of) two permutation groups, each acting on only one of
the two $ A^\infty $ factors. The general form of such a measure is
also well-known (de Finetti type theorem on partial exchangeability
\cite[Th.~4.1]{Link}), it is a mixture of products, namely, $ \int 
(\nu_1\times\nu_2\times\nu_1\times\nu_2\dots) \, \mu_{a_1,a_2}(\D
\nu_1 \D \nu_2) $ for some (uniquely determined) measure $
\mu_{a_1,a_2} \in \Pc(\Pc(B)\times\Pc(B)) $; this is the joint
distribution of $ K_{a_1} $ and $ K_{a_2} $. And so on.
\proofend

See also \cite[Sect.~2.5.1]{LJR} and \cite[8d3]{Ts03}.

The statement of Lemma \ref{7g6} does not hold for uncountable sets $
A $ (unless separability is stipulated in the spirit of \ref{2a5});
here is a counterexample. Let $ B $ contain only two points. We define
the measure $ T(a_1,a_2,\dots) $ as the uniform distribution on the
set of all sequences $ (b_1,b_2,\dots) $ such that $ \forall k,l \(
a_k = a_l \imply b_k = b_l \) $. Thus, if the sequence $
(a_1,a_2,\dots) $ contains a finite number $ n $ of different points,
then $ T(a_1,a_2,\dots) $ consists of $ 2^n $ equiprobable atoms (and
if $ n=\infty $ then the measure is continuous). Such $ T $ is of the
form $ T_\infty^K $ if and only if $ A $ is finite or countable.

Separability is naturally treated via continuity in probability, see
Sect.~\ref{7h}.

Assume now that $ A $ is also a measurable space. An \Skernel\ from $
A $ to $ B $ will be called \emph{measurable,}%
\index{measurable!\Skernel}
if it is a measurable \Smap\ (from $ A $ to $ \Pc(B) $).

Every measurable \Skernel\ $ K $ from $ A $ to $ B $ leads to an
operator $ \Tr^K_1 : \Pc(A) \to \Pc(B) $ (or rather, a linear operator
on finite signed measures), namely,
\[
\int_B \phi \, \D (\Tr^K_1 \nu) = \int_A \bigg( \Ex \int_B \phi \, \D
K_a \bigg) \, \nu(\D a) = \int_A \Tl^K_1 \phi \, \D\nu
\]
for all bounded measurable $ \phi : B \to \R $. In other words, $
\Tr^K_1 \nu $ is the barycenter of the measure $ \Tr^\Xi_1 \nu $ on $
\Pc(B) $ (that is, in $ \Pc(\Pc(B)) $), where $ \Xi $ is the same as $
K $ but treated as an \Smap\ from $ A $ to $ \Pc(B) $ (thus, $
\Tr^\Xi_1 $ is defined according to Sect.~\ref{7c}). The well-known
`barycenter' map $ \Pc(\Pc(B)) \to \Pc(B) $ is used, $ \mu \mapsto
\int \nu(\cdot) \, \mu(\D\nu) $. Similarly, for any measure $ \nu $ on
$ A^n $,
\[
\int\limits_{B^n} \phi \, \D(\Tr^K_n \nu) = \int\limits_{A^n} \bigg(
\Ex \int\limits_{B^n} \phi \, \D (K_{a_1}\times\dots\times K_{a_n})
\bigg) \nu(\D a_1 \dots \D a_n) = \int\limits_{A^n} (\Tl^K_n \phi) \,
\D\nu
\]
for $ \phi : B^n \to \R $. The same for $ n=\infty $. For an \Smap\ $
\Xi $ from $ A $ to $ B $, treated as an \Skernel\ $ K $, we have $
\Tr_n^\Xi = \Tr_n^K $.

Integrating out $ a $ while keeping $ \om $ one may get a
\emph{random} measure $ K_\nu = \int K_a \, \nu(\D a) $ on $ B $. To
this end we consider a measurable function $ \xi : A \times \Om \to
\Pc(B) $ related to $ \Xi $ as in \ref{7c4}(b), $ \Xi $ being related
to $ K $ as before. For almost every $ \om \in \Om $ we have a
measurable function $ \xi_\om : A \to \Pc(B) $, $ \xi_\om(a) =
\xi(a,\om) $. The function $ \xi_\om $ sends $ \nu $ into a measure on
$ \Pc(B) $; its barycenter is $ K_\nu $. The choice of $ \xi $ does
not matter, since $ K_\nu \in L_0(\Om\to\Pc(B)) $ is treated $ \modO
$. For an \Smap\ $ \Xi $ from $ A $ to $ B $, treated as an \Skernel\
$ K $, we have $ \Xi_\nu = K_\nu $. In general,
\[
\int_B \phi \, \D K_\nu = \int_A \bigg( \int_B \phi \, \D K_a \bigg)
\, \nu(\D a) \quad \text{a.s.}
\]
for every bounded measurable $ \phi : B \to \R $, and
\[
\Ex K_\nu = \Tr^K_1 \nu \, .
\]
The map $ \nu \mapsto K_\nu $ is a linear map $ \Pc(A) \to
L_0(\Om\to\Pc(B)) $. The family $ (K_\nu)_{\nu\in\Pc(A)} $ of
\valued{\Pc(B)} random variables $ K_\nu $ is an \Smap\ from $ \Pc(A)
$ to $ \Pc(B) $. Unlike an arbitrary \Smap\ from $ \Pc(A) $ to $
\Pc(B) $, the \Smap\ $ (K_\nu)_\nu $ is linear (in $ \nu $).

Let $ A,B,C $ be three standard measurable spaces and $ K', K'' $
be measurable \Skernel s, $ K' $ from $ A $ to $ B $, $ K'' $
from $ B $ to $ C $, on probability spaces $ \Om', \Om'' $
respectively. In order to define the composition%
\index{composition!of \Skernel s}
$ K = K'K'' $ of
\Skernel s we may turn to the corresponding \Smap s $ (K'_\nu)_\nu $,
$ (K''_\nu)_\nu $. Their composition is a \emph{linear} \Smap\ from $
\Pc(A) $ to $ \Pc(C) $, it is of the form $ (K_\nu)_\nu $, which
defines $ K = K' K'' $, an \Skernel\ from $ A $ to $ C $ on $ \Om'
\times \Om'' $; roughly speaking,
\[
(K'K'')_a (\om',\om'') = \int_B K''_b(\om'') \, K'_a(\om')(\D b) \, ,
\]
but rigorously, $ \xi', \xi'' $ should be used (as in \ref{7c5}).
Similarly to Sect.~\ref{7c}, the composition of \Skernel s, $ K = K'
K'' $, is related to convolution of their distributions, $ \La_{\Xi} =
\La_{\Xi'} * \La_{\Xi''} $, and composition of operators,
\[
\Tl^K_n \phi = \Tl^{K'}_n ( \Tl^{K''}_n \phi) \quad \text{and}
\quad \Tr^K_n \nu = \Tr^{K''}_n ( \Tr^{K'}_n \nu) \, .
\]
For \Smap s from $ \X $ to $ \X $, treated as \Skernel s, the
composition defined here conforms to that of Sect.~\ref{7c}. In
general, treating $ K' $ and $ K'' $ as two independent \Skernel s on
the same probability space, we generalize \eqref{7c7}--\eqref{7c10},
\begin{gather}
\CE{ \int_C \phi \, \D K_a }{ K' } = \int_B (
 \Tl_1^{K''} \phi ) \, \D K'_a \, ; \\
\CE{ \int_C \phi \, \D K_a }{ K'' } = \qquad \\
\qquad = \int_B \bigg( \int_C \phi \, \D K''_b \bigg) (\Tr_1^{K'}
 \de_a) (\D b) = \int_C \phi \, \D K''_\nu \, , \quad \text{where }
 \nu = \Tr_1^{K'} \de_a \, ; \notag \\
\CE{ \int_{C^n} \phi \, \D (K_{a_1}\times\dots\times K_{a_n}) }{
  K' } = \qquad \\
\qquad = \int_{B^n} ( \Tl_n^{K''} \phi ) \, \D
 (K'_{a_1}\times\dots\times K'_{a_n}) \, ; \notag \\
\CE{ \int_{C^n} \phi \, \D (K_{a_1}\times\dots\times K_{a_n}) }{
  K'' } = \qquad \\
\qquad = \int_{B^n} \bigg( \int_{C^n} \phi \, \D
 (K''_{b_1}\times\dots\times K''_{b_n}) \bigg) \( \Tr_n^{K'}
 (\de_{a_1}\times\dots\times\de_{a_n}) \) (\D b_1 \dots \D b_n) \,
 . \notag
\end{gather}

\subsection[Flow systems of S-kernels]{Flow systems of \Skernel s}
\label{7h}

Let $ \X $ be a compact metrizable space (mostly, the circle will be
used). Then $ \Pc(\X) $, equipped with the weak topology, is also a
compact metrizable space, and $ L_0(\Om\to\Pc(\X)) $ is equipped with
the (metrizable) topology of convergence in probability.

\begin{definition}\label{7h1}
An \Skernel\ $ K $ from $ \X $ to $ \X $ is \emph{continuous in
probability,}\index{continuous in probability!\Skernel}
if the map $ x \mapsto K_x $ is continuous from $ \X
$ to $ L_0(\Om\to\Pc(\X)) $.
\end{definition}

Clearly, \ref{7h1} implies measurability of $ K $. For an \Smap\ $ \Xi
$ from $ \X $ to $ \X $, treated as an \Skernel\ $ K $, \ref{7h1}
conforms to \ref{7d1} (since the natural embedding $ x \mapsto \de_x $
of $ \X $ into $ \Pc(\X) $ is homeomorphic). Here is a
generalization of Lemma \ref{7d2}. It characterizes continuity in
probability of an \Skernel\ in terms of $ C(\X^n) $ (while a
straightforward use of \ref{7d2} would involve huge spaces $
C(\Pc^n(\X)) $). It also shows that continuity in probability is
preserved by the composition.

\begin{lemma}\label{7h2}
The following three conditions are equivalent for every \Skernel\ $ K
$ from $ \X $ to $ \X $:

(a) $ K $ is continuous in probability;

(b) $ \Tl_2^K \phi \in C(\X^2) $ for all $ \phi \in C(\X^2) $;

(c) $ \Tl_n^K \phi \in C(\X^n) $ for all $ \phi \in C(\X^n) $ and all $
n = 1,2,3,\dots $
\end{lemma}

\beginproof
(a) \imp (c): the composition $ \Tl_n^K \phi $ of a chain of
continuous maps $ (x_1,\dots,x_n) \mapsto (K_{x_1},\dots,K_{x_n})
\mapsto K_{x_1} \times \dots \times K_{x_n} \mapsto \int \phi \, \D (
K_{x_1} \times \dots \times K_{x_n} ) \mapsto \Ex \int \phi \, \D (
K_{x_1} \times \dots \times K_{x_n} ) $ between the spaces $ \X^n \to
L_0(\Om\to\Pc^n(\X)) \to L_0(\Om\to\Pc(\X^n)) \to L_0 ( \Om,
[-\|\phi\|, \|\phi\|] ) \to \R $ is continuous.

(c) \imp (b): trivial.

\begin{sloppypar}
(b) \imp (a): let $ \phi \in C(\X) $, then
\begin{multline*}
\Ex \bigg| \int \phi \, \D K_{x_1} - \int \phi \, \D K_{x_2} \bigg|^2
 = \\
= \Ex \iint \phi(x) \phi(y) \( K_{x_1}(\D x) K_{x_1}(\D y) -
 K_{x_1}(\D x) K_{x_2}(\D y) - \\
- K_{x_2}(\D x) K_{x_1}(\D y) + K_{x_2}(\D x) K_{x_2}(\D y) \) = \\
= \psi(x_1,x_1) - \psi(x_1,x_2) - \psi(x_2,x_1) + \psi(x_2,x_2) \, ,
\end{multline*}
where $ \psi = \Tl_2^K (\phi\otimes\phi) $. By (b), $ \psi $ is a
continuous function on $ \X^2 $. We see that $ (x_1,x_2) \mapsto \Ex
| \int \phi \, \D K_{x_1} - \int \phi \, \D K_{x_2} |^2 $ is a
continuous function on $ \X^2 $ vanishing on the diagonal $ x_1 = x_2
$. Therefore $ x \mapsto \int \phi \, \D K_x $ is a continuous map $
\X \to L_2(\Om) $, which is (a). \qed
\end{sloppypar}
\proofendnoqed

Condition \ref{7h2}(b) may be reformulated in terms of the kernel $
T_2^K $ (non-random, from $ \X^2 $ to $ \X^2 $),
\[
T_2^K \text{ is a continuous map $ \X^2 \to \Pc(\X^2) $} \, .
\]
Similarly, \ref{7h2}(c) means continuity of all $ T_n^K : \X^n \to
\Pc(\X^n) $, or equivalently,
\[
T_\infty^K \text{ is a continuous map $ \X^\infty \to \Pc(\X^\infty)
$} \, ,
\]
since finite-dimensional functions $ (x_1,x_2,\dots) \mapsto
\phi(x_1,\dots,x_n) $ are dense in $ C(\X^\infty) $. Of course, $
\X^\infty $ is equipped with the product topology, and is a compact
metrizable space.

\begin{proposition}\label{7h3}
The following two conditions are equivalent for every kernel $ T $
from $ \X^\infty $ to $ \X^\infty $:

(a) $ T = T^K_\infty $ for some \Skernel\ $ K $ from $ \X $ to $ \X
$, continuous in probability;

(b) $ T $ satisfies \eqref{7g4}, \eqref{7g5}, and is a continuous map
$ \X^\infty \to \Pc(\X^\infty) $.
\end{proposition}

\beginproof
(a) \imp (b): evident.

(b) \imp (a): we choose a countable dense subset $ A \subset \X $
and apply Lemma \ref{7g6} to the restriction $ T_0 $ of $ T $ to $
A^\infty $, thus getting an \Skernel\ $ K_0 $ from $ A $ to $ \X $
such that $ T_0 = T_\infty^{K_0} $. For every $ \phi \in C(\X^2) $
the function $ \Tl_2^{K_0} \phi $ is \emph{uniformly} continuous on $
A^2 $ (since the map $ (x_1,x_2,\dots) \mapsto \int \phi(x'_1,x'_2)
T_{x_1,x_2,\dots} (\D x'_1 \D x'_2 \dots) $ is continuous). Similarly
to the proof of \ref{7h2} ((b) \imp (a)) we deduce that $ K_0 $ is
\emph{uniformly} continuous in probability, that is, $ K_0 $ is a
uniformly continuous map $ A \to L_0 ( \Om \to \Pc(\X) ) $. It
remains to extend it to $ \X $ by continuity.
\proofend

We observe a natural one-to-one correspondence between
\begin{itemize}
\item[]
distributions $ \La_K $ of \Skernel s $ K $ from $ \X $ to $ \X $,
continuous in probability;
\item[]
kernels $ T $ from $ \X^\infty $ to $ \X^\infty $, satisfying
\ref{7h3}(b);
\item[]
consistent systems $ (T_n)_{n=1}^\infty $ of kernels $ T_n $ from $
\X^n $ to $ \X^n $, satisfying the finite-dimensional counterpart of
\ref{7h3}(b).
\end{itemize}
The convolution of distributions corresponds to the composition of
kernels.

See also \cite[Sect.~2.5.1]{LJR} and \cite[8d3]{Ts03}.

For an \Smap\ $ \Xi $ from $ \X $ to $ \X $, treated as an \Skernel\
$ K $, the kernel $ T_2^K $ satisfies an additional condition: for all
$ x \in \X $,
\begin{equation}\label{7h4}
\text{the measure $ T_2^K(x,x) $ is concentrated on the diagonal of $
\X^2 $} \, .
\end{equation}
It is easily reformulated in terms of $ T_\infty^K $, and we get a
natural one-to-one correspondence between
\begin{itemize}
\item[]
distributions $ \La_\Xi $ of \Smap s $ \Xi $ from $ \X $ to $ \X $,
continuous in probability;
\item[]
kernels $ T $ from $ \X^\infty $ to $ \X^\infty $, satisfying
\ref{7h3}(b) and the reformulated additional condition \eqref{7h4};
\item[]
consistent systems $ (T_n)_{n=1}^\infty $ of kernels $ T_n $ from $
\X^n $ to $ \X^n $, satisfying the finite-dimensional counterpart of
\ref{7h3}(b), and \eqref{7h4}.
\end{itemize}
The convolution of distributions corresponds to the composition of
kernels.

A \emph{convolution system of \Skernel s}%
\index{convolution system!of \Skernel s}
(over $ \X $) may be
defined as a family $ (\La_{s,t})_{s<t} $, where each $ \La_{s,t} $ is
the distribution of an \Skernel\ from $ \X $ to $ \X $,
continuous in probability, and
\[
\La_{r,t} = \La_{r,s} * \La_{s,t}
\]
whenever $ r<s<t $. An equivalent description is a family $
(T_{s,t})_{s<t} $ of kernels $ T_{s,t} $ from $ \X^\infty $ to $
\X^\infty $, satisfying \ref{7h3}(b) and
\[
T_{r,t} = T_{r,s} T_{s,t}
\]
whenever $ r<s<t $.

Similarly to Sect.~\ref{7d}, every convolution system of \Skernel s
leads to a convolution system as defined by \ref{2a1}, and if
separability (as defined by \ref{2a5}) holds then we get a flow
system, and further, a continuous product of probability spaces. For
\Smap s from $ \X $ to $ \X $, treated as \Skernel s, the new
construction conforms to that of Sect.~\ref{7d}. In general, \Skernel
s from $ \X $ to $ \X $ (or their distributions) may be treated as
\Smap s from $ \Pc(\X) $ to $ \Pc(\X) $ (or their distributions,
respectively); thus, the `old' construction (of Sect.~\ref{7d}) may be
applied, as well as the `new' construction introduced above. The `old'
and `new' flow systems are isomorphic, and the `old' and `new'
separability conditions are equivalent (since the relevant `old' and
`new' \sif s coincide). The sufficient condition \eqref{7d5} for the
separability (the \emph{temporal continuity in probability}) may be
reformulated in terms of the (non-random) kernels $ T_1^{K_{s,t}} $
from $ \X $ to $ \X $, namely,
\begin{equation}\label{7h5}
\text{both $ T_1^{K_{s-\eps,s}}(x) $ and $ T_1^{K_{s,s+\eps}}(x) $
converge to $ \de_x $ as $ \eps \to 0+ $}
\end{equation}
for $ s \in \R $ and $ x \in \X $. It follows that $ K_\nu^{s-\eps,s}
\to \nu $ and $ K_\nu^{s,s+\eps} \to \nu $ in probability, since the
expectation of the transportation distance between $ K_\nu^{s,t} $ and
$ \nu $ does not exceed
\[
\Ex \iint \dist(x,x') \, K_x^{s,t} (\D x') \nu(\D x) = \iint
\dist(x,x') \, T_1^{K_{s,t}} (x) (\D x') \nu(\D x) \, .
\]
The conclusion follows.

\begin{lemma}\label{7h6}
Condition \eqref{7h5} implies separability.
\end{lemma}

\begin{remark}\label{7h7}
The continuous product of probability spaces, corresponding to a flow
system $ (K_{s,t})_{s<t} $ of \Skernel s from $ \X $ to $ \X $, is
classical if and only if random variables $ \int \phi \, \D K_x^{s,t}
$ are stable for all $ s<t $ and all bounded Borel functions $ \phi :
\X \to \R $. Proof: similar to \ref{7d7}.
\end{remark}

\subsection[From S-kernels to black noise]{From \Skernel s to black noise}
\label{7i}

Sect.~\ref{7h} gives us $ \X $, $ (G_{s,t},\mu_{s,t})_{s<t} $ and $
(K_{s,t})_{s<t} $ in the same way as Sect.~\ref{7d} provided $ \X $,
$ (G_{s,t},\mu_{s,t})_{s<t} $ and $ (\Xi_{s,t})_{s<t} $ to
Sect.~\ref{7e}. The temporal continuity in probability \eqref{7h5} is
assumed to hold. Taking into account higher moments (like in the proof of
\ref{7g3}) we note, similarly to \eqref{7e1}, that functions of the form
\begin{equation}\label{7i1}
\int_{\X^n} \phi \, \D ( K_{x_1}^{r,t} \times \dots \times
K_{x_n}^{r,t} ) \quad \text{for $ \phi \in C(\X^n) $, $
(x_1,\dots,x_n) \in \X_n $}
\end{equation}
and $ n=1,2,\dots $ are dense in $ L_2 (G_{r,t}) $. Here are
counterparts of \ref{7e2}--\ref{7e4}.

\begin{lemma}\label{7i2}
Assume that $ n \in \{1,2,\dots\} $ is given, a linear subset $ F_n $
of $ C(\X^n) $, dense in the norm topology, and a linear subset $ N_n
$ of the space of (finite, signed) measures on $ \X^n $, dense in the
weak topology; and an interval $ (r,t) \subset \R $. Then functions of
the form
\[
\int_{\X^n} \bigg( \int_{\X^n} \phi \, \D ( K_{x_1}^{r,t} \times
\dots \times K_{x_n}^{r,t} ) \bigg) \, \nu (\D x_1 \dots \D x_n) \quad
\text{for } \phi \in F_n \text{ and } \nu \in N_n
\]
are $ L_2 $-dense among functions of the form \eqref{7i1} for the
given $ n $.
\end{lemma}

\beginproof
Similar to \ref{7e2}.
\proofend

Substituting $ f = \int_{\X^n} \( \int_{\X^n} \phi \, \D (
K_{x_1}^{r,t} \times \dots \times K_{x_n}^{r,t} ) \) \, \nu (\D x_1
\dots \D x_n) $ to \eqref{7b3} we get
\begin{multline}\label{7i3}
f_{s,s+\eps} = \\
= \int_{\X^n} \bigg( \int_{\X^n} \( \Tl_n^{K_{s+\eps,t}} \phi \) \,
\D ( K_{x_1}^{s,s+\eps} \times \dots \times K_{x_n}^{s,s+\eps} )
\bigg) \, \( \Tr_n^{K_{r,s}} \nu \) (\D x_1 \dots \D x_n) \, .
\end{multline}

\begin{proposition}\label{7i4}
(Le Jan and Raimond; implicit in \cite{LJR1}.)
Let a flow system $ (K_{s,t})_{s<t} $ of \Skernel s from $ \X $ to
$ \X $ and a probability measure $ \nu_0 $ on $ \X $ satisfy the
conditions

(a) (stationarity)
the distribution of $ K_{s,s+h} $ does not depend on $ s $;

(b) (invariant measure) $ \Tr_1^{K_{s,t}} \nu_0 = \nu_0 $ for $ s<t
$;

(c) (Lipschitz boundedness)
if $ \phi \in C(\X^n) $ is a Lipschitz function, then $
\Tl_n^{K_{s,t}} \phi $ is also a Lipschitz function, with a Lipschitz
constant
\[
\Lip (\Tl_n^{K_{s,t}} \phi) \le C_n \Lip(\phi) \quad \text{for } 0 \le s <
t \le 1 \, ,
\]
where $ C_n < \infty $ depends only on $ n $;

(d)
$\displaystyle
\sup_{\phi,\nu} \Var \bigg( \int_{\X^n} \bigg( \int_{\X^n} \phi \,
\D \( K^{0,\eps}_{x_1} \times \dots \times K^{0,\eps}_{x_n} \) \bigg)
\, \nu (\D x_1 \dots \D x_n) \bigg) = o(\eps)
$ \linebreak
as $ \eps \to 0 $,
where the supremum is taken over all $ \phi \in C(\X^n) $ such that $
\Lip(\phi) \le 1 $ and all positive measures $ \nu $ on $ \X^n $ such
that $ \nu_1 \le \nu_0, \dots, \nu_n \le \nu_0 $; here $
\nu_1,\dots,\nu_n $ are coordinate projections of $ \nu $, that is, $
\int \phi(x_k) \, \nu(\D x_1\dots\D x_n) = \int \phi(x) \, \nu_k(\D x) $ for
$ \phi \in C(\X) $.

Then the corresponding noise is black.
\end{proposition}

\beginproof
Similar to \ref{7e4}.
\proofend

\mysubsection{Example: the sticky flow (Le Jan, Lemaire, Raimond)}
\label{7j}

\vspace*{12pt}

\[
\setlength{\unitlength}{0.8cm}\nopagebreak[4]
\begin{picture}(12,3)
\put(0,0.65){\includegraphics[scale=0.8]{pic7a.eps}}
\put(1,0.1){\makebox(0,0){(a)}}
\put(3.8,0.1){\makebox(0,0){(b)}}
\put(6,0.35){\includegraphics[scale=0.8]{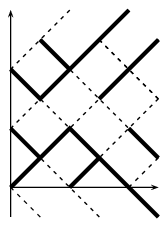}}
\put(7,0.1){\makebox(0,0){(c)}}
\put(9.5,0.65){\includegraphics[scale=0.8]{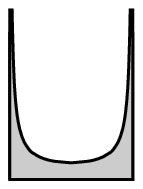}}
\put(10.25,0.1){\makebox(0,0){(d)}}
\end{picture}
\]
Recall the discrete model of Sect.~\ref{7f} (see  (a), (b)). Every
space-time lattice point $ (t,x) $ is connected with one of the two
points $ (t+1,x\pm1) $ according to a two-valued random variable. Now we
perturb the model (see (c)): $ (t,x) $ is connected strongly with one
of the two points $ (t+1,x\pm1) $ and weakly with the other, according
to a random variable $ \theta_{t,x} $ taking on values on the interval
$ [0,1] $. The case of \ref{7f} appears if $ \theta_{t,x} $ takes on two
values $ 0,1 $ only (equiprobably). Generally, $ \theta_{t,x} $ is the
strength of the connection to $ (t+1,x+1) $, and $ 1 - \theta_{t,x} $
to $ (t+1,x-1) $. Random variables $ \theta_{t,x} $ are independent,
identically distributed; their (common) distribution has a density
shown on (d), to be specified later.

Particles move on the lattice. They are conditionally independent,
given all $ \theta_{t,x} $. Conditional probabilities of transitions
from $ (t,x) $ to $ (t+1,x+1) $ and $ (t+1,x-1) $ are equal to $
\theta_{t,x} $ and $ 1-\theta_{t,x} $ respectively. Unconditionally
(when $ \theta_{t,x} $ are not given), two particles move
independently until they meet, after which they prefer moving
together, but have a chance to separate, in contrast to the model of
\ref{7f}; they are sticky rather than coalescing.

Given $ n $ particles at $ (t,x) $, the probability of $ k $ particles
to choose $ (t+1,x+1) $ (and the other $ n-k $ particles to choose $
(t+1,x-1) $) is $ \binom n k \Ex \( \theta^k (1-\theta)^{n-k} \) $. A
clever choice (specified later) of the distribution of $ \theta $
makes the expectation a product of the form
\begin{equation}\label{7j1}
\binom n k \Ex \( \theta^k (1-\theta)^{n-k} \) = \frac{ \al(k)
\al(n-k) }{ \be(n) }
\end{equation}
whenever $ k,n \in \Z $, $ 0 \le k \le n $, for some functions $ \al,
\be : \{ 0,1,\dots \} \to (0,\infty) $. This fact leads to a
stationary distribution of a simple form described below (Le Jan and
Lemaire \cite{LJL}).

Let $ x $ run over the cyclic group $ \Z_m = \Z / m\Z $; this is our
space ($ m $ is a parameter). Each configuration of $ n $ unnumbered
particles is described by a family of $ m $ \emph{occupation numbers},
$ s = (s_x)_{x\in\Z_m} $, $ \sum_x s_x = n $; namely, $ s_x $ is the
number of particles situated at $ x $. We ascribe to each
configuration a probability
\begin{equation}\label{7j2}
\mu_n(s) = \const \cdot \prod_{x\in\Z_m} \be(s_x)
\end{equation}
(`$\const$' being a normalization constant) and claim that such
probability measure $ \mu_n $ is invariant under our
dynamics. Moreover, it satisfies the \emph{detailed balance}
condition:
\[
\mu_n(s') p(s'\to s'') = \mu_n(s'') p(s''\to s')
\]
for all \particle{n} configurations $ s', s'' $; here $ p(s'\to s'') $
stands for the transition probability (the probability of $ s'' $ at $
t+1 $ given $ s' $ at $ t $).

A transition from $ s' $ to $ s'' $ may be described by \emph{edge
occupation numbers} $ s(x\to x-1) $, $ s(x\to x+1) $; these must
satisfy $ s'_x = s(x\to x-1) + s(x\to x+1) $ and $ s''_x = s(x-1\to x)
+ s(x+1\to x) $ for all $ x $. Let us call each such family of edge
occupation numbers a \emph{transition channel}. Several transition
channels from $ s' $ to $ s'' $ may exist (try increasing each $
s(x\to x-1) $ while decreasing each $ s(x\to x+1) $); the detailed
balance holds for each transition channel separately:
\begin{multline*}
\mu_n(s') p_{\text{channel}} (s'\to s'') = \prod_x \be(s_x) \cdot
 \prod_x \frac{ \al(s(x\to x+1)) \al(s(x\to x-1)) }{ \be(s_x) } = \\
= \prod_x \( \al(s(x\to x+1)) \al(s(x\to x-1)) \) \, ,
\end{multline*}
which evidently is symmetric (that is, time-reversible).

In order to get \eqref{7j1}, following \cite{LJL}, we choose for $
\theta $ the beta distribution,
\[
\theta \sim \operatorname{Beta} (\eps,\eps) \, ,
\]
its density being
\[
a \mapsto \frac{ \Gamma(2\eps) }{ \Gamma^2(\eps) } a^{\eps-1}
(1-a)^{\eps-1} \quad \text{for } 0<a<1 \, ;
\]
$ \eps \in (0,1) $ is a parameter. We have
\[
\Ex \( \theta^k (1-\theta)^{n-k} \) = \frac{ \Gamma(2\eps) }{
\Gamma^2(\eps) } \cdot \frac{ \Gamma(k+\eps) \Gamma(n-k+\eps) }{
\Gamma(n+2\eps) } \, ,
\]
thus, \eqref{7j1} holds for
\begin{gather*}
\al(k) = \frac{ \Gamma(k+\eps) }{ k! \Gamma(\eps) } = \frac1{k!} \eps
 (1+\eps) \dots (k-1+\eps) \, , \\
\be(n) = \frac{ \Gamma(n+2\eps) }{ n! \Gamma(2\eps) } = \frac1{n!}
2\eps (1+2\eps) \dots (n-1+2\eps) \, .
\end{gather*}

In order to take the scaling limit we embed the discrete space $ \Z_m
$ into a continuous space, the circle $ \T = \R / \Z $, by $ x \mapsto
\frac1m x $, and the discrete time set $ \Z $ into the continuous time
set $ \R $ by $ t \mapsto \frac1{m^2} t $. We also let $ \eps $ depend
on $ m $, namely,
\[
\eps = \frac a m \, ,
\]
$ a \in (0,\infty) $ being a parameter of the continuous
model. Convergence (as $ m \to \infty $) is proven by Le Jan and
Lemaire \cite{LJL}; the continuous model is (a special case of) the
sticky flow introduced by Le Jan and Raimond \cite{LJR}. The motion of
a single particle is the standard Brownian motion in $ \T $. Two
particles spend together a non-zero time, but never a time interval
(rather, a nowhere dense closed set of non-zero Lebesgue measure on
the time axis).

The convergence (as $ m \to \infty $) is proven at the level of the
kernels $ T_n^{K_{s,t}} $ (`moments'); the (continuous) sticky flow is
a flow system $ (K_{s,t})_{s<t} $ of \Skernel s $ K_{s,t} $ from $ \X
$ to $ \X $ (here $ \X 
= \T $) that corresponds to the constructed consistent system of
kernels $ T_n^{K_{s,t}} $ from $ \X^n $ to $ \X^n $ (recall
Proposition \ref{7h3} and the paragraph after it). A consistent system
$ (\mu_n)_n $ of invariant measures $ \mu_n $ on $ \X^n $ (or
equivalently, an invariant measure on $ \X^\infty $) is written out
explicitly by a continuous counterpart of \eqref{7j2}. The kernels $
T_n^{K_{s,t}} $ are described for infinitesimal $ t-s $ via their
Dirichlet forms,
\[
\phi \mapsto \lim_{\eps\to0+} \int_{\X^n} \bigg( \int_{\X^n}
|\phi(x)-\phi(\cdot)|^2 \, \D T_n^{K_{0,\eps}} (x) \bigg) \mu_n(\D x)
\, ,
\]
written out explicitly for smooth functions $ \phi : \X^n \to \R $.

The flow system of \Skernel s satisfies the conditions
\ref{7i4}(a--d), therefore the corresponding noise is black.

For details see \cite{LJL}, \cite{LJR}.

\mysection{Unitary flows are classical}
\label{sec:8}
\mysubsection{Can a Brownian motion be nonclassical?}
\label{8a}

We know (see \ref{6a12}) that every \flow{\R} is classical (that is,
generates a classical continuous product of probability spaces). The
same holds for every stationary \flow{\T}, and moreover, for every
\flow{\T} satisfying the upward continuity condition (see \ref{6a10}
and \ref{2c2}). However, nonclassical stationary \flow{G}s exist in
some Borel semigroups (in fact, finite-dimensional topo-semigroups) $
G $, see Sect.~\ref{5d}.

\begin{question}\index{question}\label{8a*}
Whether every stationary \flow{G} in every Borel group $ G $ (not just
semigroup) is classical, or not?
\end{question}

If $ G $ is a group (rather than semigroup) then \flow{G}s may be identified
with \valued{G} processes with independent increments. If $ G $ is a
topological (rather than Borel) group then some of these processes are
sample-continuous and stationary; these are called Brownian motions in $ G $,
classical or not, according to the corresponding noises.

\begin{question}\index{question}
(a)
Whether every Brownian motion in the group $ G = \Homeo(\T) $ of
all homeomorphisms of the circle $ \T $ is classical, or not?

(b) \cite[1.11]{Ts98}
Whether every Brownian motion in every Polish group $ G $ is
classical, or not?
\end{question}

It is enough to examine a single Polish group $ G =
\Homeo([0,1]^\infty) $ of all homeomorphisms of the Hilbert cube $
[0,1]^\infty $ rather than all Polish groups, since every Polish group
is isomorphic to a subgroup of $ \Homeo([0,1]^\infty) $ by a theorem
of Uspenskii, see \cite[9.18]{Ke} or \cite[1.4.1]{BK}.

We turn to the unitary group $ \U(H) $ of a separable
infinite-dimensional Hilbert space $ H $ (over $ \R $ or $ \C $),
equipped with the strong (or equivalently, weak) operator topology.
Every Brownian motion in $ \U(H) $ is classical \cite{Ts98}, which is
generalized by the theorem below (its proof is basically the same as in
\cite{Ts98}). We know (recall \ref{2c7}) that every \emph{weakly continuous}
convolution semigroup in $ \U(H) $ leads to an \flow{\U(H)} $ (X_{s,t})_{s<t}
$ and a measurable action $ (T_h)_h $ of $ \R $ on $ \Om $ such that $ X_{s,t}
\circ T_h = X_{s+h,t+h} $. Such $ (X_{s,t})_{s<t} $ will be called a
\emph{stationary \flow{\U(H)}, continuous in probability.}

\begin{theorem}\label{8a2}\index{theorem}
Every stationary \flow{\U(H)}, continuous in probability, is
classical.
\end{theorem}

\begin{question}\index{question}
(a) Does Theorem \ref{8a2} hold without assuming continuity in
probability? (See also \ref{8a*}.)

(b) Is the following claim true? Every stationary \flow{G} $
(X_{s,t})_{s<t} $ is of the form $ X_{s,t} = c_s^{-1} Y_{s,t} c_t $
where $ (Y_{s,t})_{s<t} $ is a stationary \flow{G} continuous in
probability, and $ c_t \in G $ for $ t \in \R $. (The map $ t \mapsto
c_t $ need not be measurable.) Here $ G = \U(H) $; but more general $ G $ may
be considered as well.
\end{question}

\begin{corollary}\label{8a3}
Let $ G $ be a topological semigroup admitting a one-to-one continuous
homomorphism to $ \U(H) $. Then every stationary \flow{G}, continuous
in probability, is classical.
\end{corollary}

\begin{example}
(a)
The group $ \Diff(M) $ of all diffeomorphisms of a smooth manifold $ M
$ acts unitarily on $ L_2(M,m) $, where $ m $ is any smooth strictly
positive measure on $ M $. Unitarity of the action is achieved
multiplying by the root of the Radon-Nikodym derivative, see
\cite[14.4.5]{Ar}. By \ref{8a3}, every stationary \flow{\Diff(M)},
continuous in probability, is classical. Brownian motions in $
\Diff(M) $ are described by Baxendale \cite{Ba}.

(b)
Every Lie group $ G $ acts by diffeomorphisms on itself. Therefore
every stationary \flow{G}, continuous in probability, is
classical. Brownian motions in Lie groups are described by Yosida
\cite{Yos}.
\end{example}

We may also consider the semigroup $ G $ of all contractions (that is,
linear operators of norm $ \le 1 $) in a separable
infinite-dimensional Hilbert space $ H $ (over $ \R $ or $ \C $),
equipped with the strong operator topology. It is a topological
semigroup, $ G \supset \U(H) $. Here is the corresponding
generalization of Theorem \ref{8a2}.

\begin{theorem}\label{8a6}\index{theorem}
Every stationary \flow{G}, continuous in probability, is
classical. (Here $ G $ is the semigroup of contractions.)
\end{theorem}

\begin{example}
The semigroup $ G $ of all conformal endomorphisms of the disc acts by
contractions on the space $ L_2 $ on the disc. Therefore every
stationary \flow{G}, continuous in probability, is classical. Some
important Brownian motions in $ G $ are now well-known as 
SLE.\index{SLE as a Brownian motion}
\end{example}

\mysubsection{From unitary flows to quantum instruments}
\label{8b}

In order to prove Theorem \ref{8a2} it is sufficient to check
stability, namely,
\begin{equation}\label{8b1}
U^\rho f_{\phi,\psi} \to f_{\phi,\psi} \quad \text{as } \rho \to 1 \,
.
\end{equation}
Here $ U^\rho $ are the operators on $ L_2(\Om) $ introduced in
Sect.~\ref{5b}, and random variables $ f_{\phi,\psi} \in L_2(\Om) $
are matrix elements of the given random unitary operators,
\[
f_{\phi,\psi} (\om) = \ip{ X_{0,1}(\om) \phi }{ \psi }_H \quad \text{for
} \om \in \Om, \> \phi,\psi \in H \, .
\]
The \valued{\U(H)} random variable $ X_{0,1} $ belongs to a given
\flow{U(H)} $ (X_{s,t})_{s<t} $; $ X_{s,t} \in L_0 \( \Om \to \U(H) \)
$. The time interval $ [0,1] $ is used just for convenience.

The operator $ U^\rho $ is the limit of (a net of) operators $ \ti
U^\rho_{t_1,\dots,t_n} $, see \eqref{5b3a}. Calculating $
\ti U^\rho_{t_1,\dots,t_n} f_{\phi,\psi} $ we come to an important
construction (quantum operations and instruments). In order to
simplify notation we do it first for a finite time set $ T = \{
0,1,\dots,n \} $ rather than a finite subset of a
continuum. Accordingly, for now our \flow{\U(H)} consists of $ n $
independent random operators $ X_{0,1}, \dots, X_{n-1,n} $ on $ H $
(and their products).

For $ n=1 $ we have a single random operator $ X_{0,1} $, and
\begin{multline*}
\ip{ U^\rho f_{\phi,\psi} }{ f_{\phi,\psi} }_{L_2(\Om)} = | \Ex
 f_{\phi,\psi} |^2  + \rho \( \Ex | f_{\phi,\psi} |^2 - | \Ex
 f_{\phi,\psi} |^2 \) = \\
= | \ip{ \Ex X_{0,1}\phi }{ \psi }_H |^2 + \rho \( \Ex | \ip{
 X_{0,1}\phi }{ \psi }_H |^2 - | \ip{ \Ex X_{0,1}\phi }{ \psi }_H |^2
 \) \, .
\end{multline*}
For a fixed $ \phi $ we have two positive\footnote{%
 Not necessarily strictly positive.}
quadratic forms of $ \psi $; they correspond to some positive
 Hermitian operators $ S^{(0)}, S^{(1)} : H \to H $,
\begin{equation}\label{8b2}
\ip{ U^\rho f_{\phi,\psi} }{ f_{\phi,\psi} }_{L_2(\Om)} = \ip{ S^{(0)}
\psi }{ \psi }_H + \rho \ip{ S^{(1)} \psi }{ \psi }_H \, ; \quad
 S^{(0)}, S^{(1)} \ge 0 \, .
\end{equation}
In terms of one-dimensional operators $ S_\phi $,
\[
S_\phi \psi = \ip{ \psi }{ \phi }_H \phi \, ,
\]
we get $ \ip{ S_\phi \psi }{ \psi }_H = | \ip\phi\psi |^2 =
\trace(S_\phi S_\psi) $ (the trace is always taken in $ H $, never in
$ L_2(\Om) $), $ S^{(0)} = S_{\Ex X_{0,1} \phi} $ and $ S^{(0)} +
S^{(1)} = \Ex S_{X_{0,1}\phi} $. Quadratic dependence of $ S^{(0)},
S^{(1)} $ on $ \phi $ means their linear dependence on $ S_\phi $,
\[
S^{(0)} = \Ec^{(0)} S_\phi \, , \quad S^{(1)} = \Ec^{(1)} S_\phi \, ;
\]
here $ \Ec^{(0)}, \Ec^{(1)} : \V(H) \to \V(H) $ are linear operators
on the space $ \V(H) $ of all Hermitian trace class operators $ H \to
H $; $ \V(H) $ is a Banach space (over $ \R $), partially ordered by
the cone $ \V^+(H) $ of all positive operators. The operators $
\Ec^{(0)}, \Ec^{(1)} $ are positive in the sense that $ S \in \V^+(H)
$ implies $ \Ec^{(0)} S, \Ec^{(1)} S \in \V^+(H) $. Also, $ \trace \(
\Ec^{(0)} S \) \le \trace(S) $ for $ S \in V^+(H) $, and the same for
$ \Ec^{(1)} $. Such operators on $ \V(H) $ are called \emph{quantum
operations} in \cite[Chap.~2]{Da}. A stronger requirement, called
complete positivity \cite[Chap.~9.2]{Da} is satisfied, but will not be
used. See \cite[Sect.~3]{Ts98} and \cite{Da} for details.

The sum $ \Ec^{(0)} + \Ec^{(1)} $ is also a quantum operation, and in
fact a \emph{nonselective}\footnote{%
 A nonselective quantum operation is also called a quantum channel.}
one, in the sense that
\begin{equation}
\trace ( \Ec^{(0)} S + \Ec^{(1)} S ) = \trace (S) \quad \text{for } S
\in \V(H) \, ,
\end{equation}
since $ \trace \( \Ec^{(0)} S_\phi + \Ec^{(1)} S_\phi \) = \trace (
\Ex S_{X_{0,1}\phi} ) = \Ex \trace (S_{X_{0,1}\phi}) = \Ex \|
X_{0,1}\phi \|_H^2 \linebreak[0]
= \| \phi \|_H^2 = \trace (S_\phi) $ for $ \phi \in
H $. A finite family of quantum operations whose sum is nonselective
is a (quantum) \emph{instrument} \cite[Chap.~4,
Def.~1.1]{Da}.\footnote{%
 In general, an instrument is a vector measure, valued in quantum
 operations. We need only the elementary case when the underlying
 measurable space is finite.}
Thus, a random unitary operator $ X_{0,1} $ leads to an instrument $
\Ec_{0,1} $ consisting of $ \Ec^{(0)} $ and $ \Ec^{(1)} $,
\[
\Ec^{(0)} S_\phi = S_{\Ex X_{0,1} \phi} \, , \quad ( \Ec^{(0)} +
\Ec^{(1)} ) S_\phi = \Ex S_{X_{0,1} \phi} \, ,
\]
which completes the case $ n=1 $.

For $ n=2 $ we have two independent random operators $ X_{0,1},
X_{1,2} $ and their product $ X_{0,2} = X_{0,1} X_{1,2} $; the latter
means $ X_{0,2} \phi = X_{1,2} ( X_{0,1} \phi ) $.\footnote{%
 It could be more convenient to write $ \phi X $ rather than $ X \phi
 $.}
Two instruments $ \Ec_{0,1} $ and $ \Ec_{1,2} $ arise as above. Their
composition \cite[Sect.~4.2]{Da} is an instrument $ \Ec_{0,2} $
consisting of four quantum operations $ \Ec_{0,2}^{(0,0)},
\Ec_{0,2}^{(0,1)}, \Ec_{0,2}^{(1,0)}, \Ec_{0,2}^{(1,1)} $ defined by
\[
\Ec_{0,2}^{(k,l)} = \Ec_{0,1}^{(k)} \Ec_{1,2}^{(l)} \, , \quad
\text{that is,} \quad \Ec_{0,2}^{(k,l)} S = \Ec_{1,2}^{(l)} (
\Ec_{0,1}^{(k)} S ) \quad \text{for } S \in \V(H) \, .
\]
Similarly to \eqref{8b2},
\begin{equation}\label{8b3}
\ip{ U^\rho f_{\phi,\psi} }{ f_{\phi,\psi} }_{L_2(\Om)} =
\sum_{k,l\in\{0,1\}} \rho^{k+l} \trace \( ( \Ec_{0,2}^{(k,l)} S_\phi)
S_\psi \) \, .
\end{equation}

\begin{proof}[Proof \textup{(}sketch\textup{)} of \textup{\eqref{8b3}}]
First we introduce one-dimensional operators $ S_{\phi_1,\phi_2} $ for
$ \phi_1,\phi_2 \in H $ by $ S_{\phi_1,\phi_2} \psi =
\ip{\psi}{\phi_1} \phi_2 $ and observe that $
\trace(AS_{\phi_1,\phi_2}) = \ip{ A\phi_2 }{ \phi_1 } $ for $ A \in
\V(H) $. Two formulas for $ n=1 $ follow from \eqref{8b2} by
bilinearity:
\begin{gather*}
\ip{ U^\rho f_{\phi,\psi_1} }{ f_{\phi,\psi_2} }_{L_2(\Om)} =
 \sum_{k=0,1} \rho^k \trace \( (\Ec^{(k)} S_\phi) S_{\psi_1,\psi_2} \)
 \, , \\
\ip{ U^\rho f_{\phi_1,\psi} }{ f_{\phi_2,\psi} }_{L_2(\Om)} =
 \sum_{k=0,1} \rho^k \trace \( (\Ec^{(k)} S_{\phi_2,\phi_1}) S_\psi \)
 \, .
\end{gather*}

Second, we choose an orthonormal basis $ (e_n)_n $ of $ H $ and note
that
\[
A = \sum_{m,n} \trace(AS_{n,m}) S_{m,n}
\]
for $ A \in \V(H) $.

Third,
\begin{multline*}
f_{\phi,\psi} = \ip{ X_{0,2}\phi }{ \psi }_H = \ip{
 X_{1,2}(X_{0,1}\phi) }{ \psi }_H = \\
= \sum_n \ip{ X_{0,1}\phi }{ e_n }_H \ip{ X_{1,2} e_n }{ \psi }_H =
 \sum_n f'_{\phi,e_n} \otimes f''_{e_n,\psi} \, ;
\end{multline*}
here $ f'_{\dots} \in L_2(\Om,\F_{0,1}) $ and $ f''_{\dots} \in
L_2(\Om,\F_{1,2}) $. Finally,
\begin{multline*}
\ip{ U^\rho f_{\phi,\psi} }{ f_{\phi,\psi} }_{L_2(\Om)} =
 \Big\langle \sum_m (U_{0,1}^\rho f'_{\phi,e_m}) \otimes (U_{1,2}^\rho
 f''_{e_m,\psi}) , \, \sum_n f'_{\phi,e_n} \otimes f''_{e_n,\psi}
 \Big\rangle_{L_2(\Om)} = \\
= \sum_{m,n} \ip{ U_{0,1}^\rho f'_{\phi,e_m} }{ f'_{\phi,e_n}
 }_{L_2(\Om)} \ip{ U_{1,2}^\rho f''_{e_m,\psi} }{ f''_{e_n,\psi}
 }_{L_2(\Om)} = \\
= \sum_{m,n} \bigg( \sum_{k=0,1} \rho^k \trace \( (\Ec_{0,1}^{(k)}
 S_\phi) S_{m,n} \) \bigg) \bigg( \sum_{l=0,1} \rho^l \trace \(
 (\Ec_{1,2}^{(l)} S_{n,m}) S_\psi \) \bigg) = \\
= \sum_{k,l\in\{0,1\}} \rho^{k+l} \trace \( S_{1,2}^{(l)}
 (\Ec_{0,1}^{(k)} S_\phi) S_\psi \) = \sum_{k,l\in\{0,1\}} \rho^{k+l}
  \trace \( ( \Ec_{0,2}^{(k,l)} S_\phi) S_\psi \) \, .
\end{multline*}
\end{proof}

Similarly, for $ T = \{0,1,\dots,n\} $, $ n=1,2,3,\dots $
\begin{equation}\label{8b4}
\ip{ U^\rho f_{\phi,\psi} }{ f_{\phi,\psi} }_{L_2(\Om)} =
\sum_{k_1,\dots,k_n\in\{0,1\}} \rho^{k_1+\dots+k_n} \trace \( (
\Ec_{0,n}^{(k_1,\dots,k_n)} S_\phi) S_\psi \)
\end{equation}
where $ \Ec_{0,n} $ is the composition of $ n $ instruments,
\[
\Ec_{0,n}^{(k_1,\dots,k_n)} = \Ec_{0,1}^{(k_1)} \dots
\Ec_{n-1,n}^{(k_n)} \, ,
\]
$ \Ec_{m-1,m} $ being defined by
\[
\Ec_{m-1,m}^{(0)} S_\phi = S_{\Ex X_{m-1,m} \phi} \, , \quad (
\Ec_{m-1,m}^{(0)} + \Ec_{m-1,m}^{(1)} ) S_\phi = \Ex S_{X_{m-1,m}}
\phi
\]
for $ \phi \in H $, $ m=1,\dots,n $.

\mysubsection{From quantum instruments to Markov chains and stopping times}
\label{8c}

We still deal with the finite time set $ T = \{ 0,1,\dots,n \} $ and
an \flow{\U(H)} $ (X_{s,t})_{s<t;s,t\in T} $. Let $ S_0 \in \V_1^+(H)
$ be given; by $ \V_1^+(H) $ we denote $ \{ S \in \V^+(H) :
\trace(S)=1 \} $. Together with the instrument $ \Ec_{0,n} $
constructed in Sect.~\ref{8b}, $ S_0 $ leads to a probability
distribution $ \mu_{0,n} $ on the set $ \{0,1\}^n $ of $ 2^n $ points
$ (k_1,\dots,k_n) $,
\[
\mu_{0,n} (k_1,\dots,k_n) = \trace ( \Ec_{0,n}^{(k_1,\dots,k_n)} S_0 )
\, .
\]
For any $ \phi, \psi \in H $ such that $ \| \phi \| = 1 $, $ \| \psi
\| = 1 $ we have $ \trace \( ( \Ec_{0,n}^{(k_1,\dots,k_n)} S_\phi )
S_\psi \) \linebreak[0]
\le \trace ( \Ec_{0,n}^{(k_1,\dots,k_n)} S_\phi ) $, thus, \eqref{8b4}
gives
\begin{equation}\label{8c1}
\ip{ (\One-U^\rho) f_{\phi,\psi} }{ f_{\phi,\psi} }_{L_2(\Om)} \le 
\sum_{k_1,\dots,k_n\in\{0,1\}} (1-\rho^{k_1+\dots+k_n}) \mu_{0,n}
(k_1,\dots,k_n) \, ;
\end{equation}
striving to prove \eqref{8b1} we will estimate this sum from above,
showing that $ k_1+\dots+k_n $ is not too large for the most of $
(k_1,\dots,k_n) $ according to $ \mu_{0,n} $.

The probability measure $ \mu_{0,n} $ turns the set $ \{0,1\}^n $ into
another probability space (in addition to the original probability
space $ \Om $ that carries $ (X_{s,t})_{s<t} $); the `coordinate'
random process $ k_1, \dots, k_n $ generates the natural filtration on
$ ( \{0,1\}^n, \mu_{0,n} ) $. We introduce a Markov chain $ S_0,
\dots, S_n $ on the (filtered) probability space $ ( \{0,1\}^n,
\mu_{0,n} ) $; each $ S_t $ is a random element of $ \V_1^+(H) $. The
initial state $ S_0 $ is chosen from the beginning (and
non-random). For $ t \in \{ 1,\dots,n \} $ the state $ S_t $ is a
function of $ k_1,\dots,k_t $, namely,
\[
S_t = \frac1{ \trace( \Ec_{0,t}^{(k_1,\dots,k_t)} S_0 ) }{
\Ec_{0,t}^{(k_1,\dots,k_t)} S_0 } \, ;
\]
of course, $ \Ec_{0,t}^{(k_1,\dots,k_t)} = \Ec_{0,1}^{(k_1)} \dots
\Ec_{t-1,t}^{(k_t)} $. It may happen that the denominator vanishes for
some $ (k_1,\dots,k_t) $, but such cases are of probability $ 0 $ and
may be ignored. We have
\begin{gather*}
\cP{ k_{t+1} = 1 }{ k_1,\dots,k_t } = \trace ( \Ec_{t,t+1}^{(1)} S_t )
 \, , \\
S_{t+1} = \frac1{ \cP{ k_{t+1} }{ k_1,\dots,k_t } }
 \Ec_{t,t+1}^{(k_{t+1})} S_t \, ,
\end{gather*}
since the joint distribution of $ k_1, \dots, k_t $ is given by $
\trace( \Ec_{0,t}^{(k_1,\dots,k_t)} S_0 ) $.

Striving to estimate $ k_1+\dots+k_n $ from above we introduce
stopping times $ \tau_1, \tau_2, \dots $ on the (filtered) probability
space $ ( \{0,1\}^n, \mu_{0,n} ) $ by
\begin{gather*}
\tau_1 = \min \big\{ t \in \{1,\dots,n\} : k_t = 1 \big\} \, , \\
\tau_{m+1} = \min \big\{ t \in \{\tau_m+1,\dots,n\} : k_t = 1 \big\}
 \quad \text{for } m=1,2,\dots
\end{gather*}
The minimum of the empty set is, by definition, infinite; thus, $
\tau_1 = \infty $ if $ k_1=\dots=k_n=0 $, and more generally, $ \tau_m
= \infty $ if (and only if) $ k_1+\dots+k_n < m $.

\begin{lemma}\label{8c1a}
\begin{gather}
1 - \Ex \exp(-\tau_1/n) \ge \frac1{3n} \sum_{t=1}^n \trace (
 \Ec_{0,t}^{(0,\dots,0)} S_0 ) \, ; \tag{a} \\
1 - \CE{ \exp\Big(-\frac{\tau_{m+1}-\tau_m}n\Big) }{
 k_1,\dots,k_{\tau_m} } \ge \frac1{3n} \sum_{t=\tau_m+1}^n \trace (
 \Ec_{\tau_m,t}^{(0,\dots,0)} S_{\tau_m} ) \tag{b}
\end{gather}
for $ m=1,2,\dots $
\end{lemma}

\beginproof
(a):
\begin{multline*}
1 - \Ex \exp(-\tau_1/n) = \sum_{t=0}^{n-1} \( \exp(-t/n) -
 \exp(-(t+1)/n) \) \Pr{ k_1=\dots=k_t=0 } \\
+ \exp(-1) \Pr{
 k_1=\dots=k_n=0 } \ge \\
\ge \frac1{3n} \sum_{t=1}^n \Pr{ k_1=\dots=k_t=0 } = \frac1{3n}
 \sum_{t=1}^n \trace ( \Ec_{0,t}^{(0,\dots,0)} S_0 ) \, ;
\end{multline*}
(b): similar.
\proofend

The next lemma holds for any increasing sequence of stopping times $
(\tau_m)_m $, irrespective of instruments etc.

\begin{lemma}\label{8c2}
For any $ \theta \in (0,1) $ and $ m = 1,2,\dots $ the probability of
the event
\[
\tau_m \le n \quad \text{and} \quad \cE{ \exp(-(\tau_{l+1}-\tau_l)/n)
}{ k_1,\dots,k_{\tau_l} } \le \theta \text{ for } l=0,1,\dots,m-1
\]
does not exceed $ \E \theta^m $. (Here $ \tau_0 = 0 $.)
\end{lemma}

\begin{sloppypar}
\beginproof
Denote by $ M $ the (random) least $ l \in \{0,\dots,m-1\} $ such that
$ \cE{ \exp(-(\tau_{l+1}-\tau_l)/n) }{ k_1,\dots,k_{\tau_l} } > \theta
$ and let $ M=m $ if there is no such $ l $. The expectation $ \Ex (
\theta^{-(M\wedge l)} \exp(-\tau_{M\wedge l}/n) ) $ decreases in $ l $,
therefore $ \Ex ( \theta^{-M} \exp(-\tau_M/n) ) \le 1 $, which implies
$ \Pr{ \tau_m \le n , \, M = m } \le \E \theta^m $.
\proofend
\end{sloppypar}

\mysubsection{A compactness argument}

We return from the finite time set to the continuum $ \R $. All said
in \ref{8b} and \ref{8c} is applicable to any finite subset of $ \R $,
but for now we only need the operation $ \Ec^{(0)}_{s,t} : \V(H) \to
\V(H) $ for $ s<t $,
\[
\Ec^{(0)}_{s,t} S_\phi = S_{\Ex X_{s,t}\phi} \quad \text{for } \phi
\in H \, .
\]
By stationarity, $ \Ec^{(0)}_{s,t} = \Ec^{(0)}_{0,t-s} $.

\begin{lemma}\label{8d1}
(a) $ \trace(\Ec^{(0)}_{0,t} S) \to \trace(S) $ as $ t \to 0+ $, for
every $ S \in \V(H) $;

(b) convergence in (a) is uniform in $ S \in K $ whenever $ K $ is a
compact subset of $ \V(H) $.
\end{lemma}

\beginproof
(a):
\[
\trace(\Ec^{(0)}_{0,t} S_\phi) = \trace(S_{\Ex X_{s,t}\phi}) = \|
\Ex X_{s,t}\phi \|^2 \to \| \phi \|^2 = \trace(S_\phi) \, ,
\]
since $ X_{0,t}\phi \to \phi $ in probability, and $ \| X_{0,t} \phi
\|^2 \le \| \phi \|^2 $.

(b): Use monotonicity of $ \trace(\Ec^{(0)}_{0,t} S) $ in $ t $. Or
alternatively, use uniform continuity: the linear functionals $ S
\mapsto \trace(\Ec^{(0)}_{0,t} S) $ are of norm $ \le 1 $.
\proofend

Given a finite set $ T \subset \R $, $ T = \{ t_0, \dots, t_n \} $, $
t_0 < t_1 < \dots < t_n $, we may consider the restricted \flow{\U(H)}
$ (X_{s,t})_{s<t;s,t\in T} $. Applying to it the construction of
Sect.~\ref{8c} we get a Markov chain $ (S_t)_{t\in T} $ provided that
an initial state $ S_{t_0} \in \V_1^+(H) $ is chosen; each $ S_t $ is
a random element of $ \V_1^+(H) $.

\begin{lemma}\label{8d2}
For every $ S_0 \in \V_1^+(H) $ and $ \eps > 0 $ there exists a
compact set $ K \subset \V_1^+(H) $ such that for every finite set $ T
= \{ t_0, \dots, t_n \} $, $ 0 = t_0 < t_1 < \dots < t_n = 1 $, 
\[
\Pr{ S_{t_0} \in K, \, \dots, \, S_{t_n} \in K } \ge 1-\eps \, ,
\]
where $ (S_t)_{t\in T} $ is the Markov chain starting with $ S_{t_0} =
S_0 $.
\end{lemma}

\beginproof
We take finite-dimensional projection operators $ Q_m : H \to H $ such
that
\[
\Ex \trace ( Q_m X_{0,1} S_0 X^*_{0,1} ) \ge 1 - \frac1{m^3}
\]
and (given $ T $ and $ m $) define a martingale $ M_m $ by
\begin{multline*}
M_m (t_k) = \cE{ \trace (Q_m S_{t_n}) }{ S_{t_0},\dots,S_{t_k} } = \\
= \trace \( Q_m \cE{ S_{t_n} }{ S_{t_k} } \) = \trace \( Q_m \Ex
 (X_{t_k,1} S_{t_k} X^*_{t_k,1}) \) = \\
= \trace \( S_{t_k} \Ex (X^*_{t_k,1} Q_m X_{t_k,1}) \) = \trace
\( S_{t_k} Q_m(t_k) \) \, ,
\end{multline*}
where
\[
Q_m(t) = \Ex (X^*_{t,1} Q_m X_{t,1}) \, .
\]
That is,
\[
\Pr{ S_{t_0} \in A_m, \, \dots, S_{t_n} \in A_m } \ge 1 - \frac1{m^2}
\, ,
\]
where
\[
A_m = \{ S \in \V_1^+(H) : \sup_{t\in[0,1]} \trace \( S Q_m(t) \) \ge
1 - \tfrac1m \} \, .
\]
(Note that $ A_m $ does not depend on $ T $.) It remains to prove that
for every $ m $ the set $ A_m \cap A_{m+1} \cap \dots $ is compact.

We have $ Q_m(t) = \cT^*_{1-t} (Q_m) $ where linear operators $
\cT^*_t : \V(H) \to \V(H) $ defined by $ \cT^*_t (S) = \Ex ( X^*_{0,t}
S X_{0,t} ) $ are a one-parameter semigroup. The semigroup is strongly
continuous, since $ X^*_{0,t} S_\phi X_{0,t} = S_{X^*_{0,t} \phi} \to
S_\phi $ a.s.\ as $ t \to 0 $. Therefore sets $ \{ Q_m(t) : t \in
[0,1] \} $ are compact in $ \V(H) $, and we may choose
finite-dimensional projections $ \ti Q_m $ such that
\[
Q_m(t) \le \ti Q_m + \frac1m \One \quad \text{for } t \in [0,1] \, ,
\]
which implies $ \trace \( SQ_m(t) \) \le \trace ( S \ti Q_m ) +
\frac1m $, thus
\[
A_m \subset \{ S \in \V_1^+(H) : \trace ( S \ti Q_m ) \ge 1-\tfrac2m
\} \, .
\]
Compactness of $ A_m \cap A_{m+1} \cap \dots $ follows easily.
\proofend

\begin{proof}[Proof \textup{(sketch)} of Theorem \textup{\ref{8a2}}]
According to \eqref{8b1} and Sect.~\ref{5b} it is sufficient to prove
that $ \ti U^\rho_{t_1,\dots,t_n} f_{\phi,\psi} \to f_{\phi,\psi} $ (as $ \rho
\to 1 $) uniformly in $ T $, where $ T $ runs over all finite
sets $ T = \{ t_1,\dots,t_n \} \subset (0,1) $. According to
\eqref{8c1} it is sufficient to prove that $ \Ex ( 1 -
\rho^{k_1+\dots+k_n} ) \to 0 $ (as $ \rho \to 1 $) uniformly in $ T $. Here $
T $ runs over all finite sets $ T = \{ t_0,\dots,t_n \} $, $ 0 = t_0 < t_1 <
\dots < t_n = 1 $, and $ (k_1,\dots,k_n) $ is the random process introduced in
Sect.~\ref{8c} (distributed $ \mu_{0,n} $). The initial state $ S_0 \in
\V_1^+(H) $ is arbitrary but fixed (that is, the convergence need not be
uniform in $ S_0 $).

Let $ \eps > 0 $ be given; we have to find $ \rho < 1 $ such that $
\Ex ( 1 - \rho^{k_1+\dots+k_n} ) \le \eps $ (or rather, $ 3\eps $) for
all $ T $ (and $ n $).

Lemma \ref{8d2} gives us a compact set $ K \subset \V_1^+(H) $ such
that $ \Pr{ S_{t_0} \in K, \, \dots, \linebreak[0]
 S_{t_n} \in K } \ge 1-\eps
$. By Lemma \ref{8d1}, $ \trace (\Ec_{0,t}^{(0)} S ) \to 1 $ (as $ t \to 0 $)
uniformly\linebreak
 in $ S \in K $. It follows that $ \inf_T \inf_{S\in K}
\int_0^1 \trace (\Ec_{0,t}^{(0)} S ) \, \D t > 0 $. Similarly, \linebreak
 $ \inf_T \inf_{S\in K} \frac1n \sum_{k=1}^n \trace (\Ec_{0,t_k}^{(0)} S
) > 0 $ provided that $ T $ is distributed on $ [0,1] $ uniformly
enough, say, $ t_{k+1} - t_k \le \frac2n $ for all $ k $, which can be
ensured by enlarging $ T $ appropriately. Using Lemma \ref{8c1a}(b) and
taking into account that $ \tau_m = \infty $ for $ m > k_1+\dots+k_n $
we get $ \theta < 1 $ such that the inequality
\[
\min_{m=1,\dots,k_1+\dots+k_n} \CE{ \exp\Big(-\frac{\tau_{m+1}-\tau_m}n\Big)
}{ k_1,\dots,k_{\tau_m} } \le \theta
\]
holds with probability $ \ge \Pr{ S_{t_0} \in K, \, \dots, \, S_{t_n}
\in K } \ge 1-\eps $; note that $ \theta $ depends on $ K $ but not $
T,n $. Combining it with Lemma \ref{8c2} we have for $ m=1,2,\dots $
\[
\Pr{ \tau_m \le n } \le \E \theta^m + \eps \, ,
\]
that is, $ \Pr{ k_1+\dots+k_n \ge m } \le \E \theta^m + \eps
$. Choosing $ m $ such that $ \E \theta^m \le \eps $ we get $ \Pr{
k_1+\dots+k_n \ge m } \le 2\eps $; note that $ m $ depends on $ K $
but not $ T,n $.

Finally, $ \Ex ( 1 - \rho^{k_1+\dots+k_n} ) \le \Pr{ k_1+\dots+k_n
\ge m } + (1-\rho^m) \le 2\eps + (1-\rho^m) $. Choosing $ \rho < 1 $
such that $ 1 - \rho^m \le \eps $ we get $ \Ex ( 1 -
\rho^{k_1+\dots+k_n} ) \le 3\eps $.
\end{proof}

The proof of Theorem \ref{8a6} is similar. However, the quantum
operation $ \Ec^{(0)} + \Ec^{(1)} $ becomes selective (recall
\eqref{8b2}), which leads to killing for the random process $ k_1,
\dots, k_n $ (introduced in \ref{8c}). Accordingly, we enlarge the
state space of the Markov chain $ S_0, \dots, S_n $ (also introduced
in \ref{8c}) by an absorbing state $ 0 $; now, each $ S_t $ is a
random element of $ \V_1^+(H) \cup \{0\} \subset \V^+(H) $. If $ S_t $
jumps to $ 0 $ then the next stopping time $ \tau_m $ is, by
definition, infinite. In the proof of Lemma \ref{8d2} each inequality
of the form $ \Ex \trace(Q_m S) \ge 1-\de $ should be first rewritten
as $ \Ex \trace \( (1-Q_m) S \) \le \de $; the latter form is
applicable in the more general situation.

\mysection{Random sets as degrees of nonclassicality}
\label{sec:9}
\mysubsection{Discrete time (toy models)}
\label{9a}

For a finite time set $ T = \{ 1,\dots,n \} $ a continuous product $
(\Om_{s,t}, P_{s,t})_{s<t;s,t\in T} $ of probability spaces is just
the usual product,
\[
(\Om_{s,t},P_{s,t}) = (\Om_{s,s+1},P_{s,s+1}) \times \dots \times
(\Om_{t-1,t},P_{t-1,t}) \, .
\]
Accordingly,
\[
H_{s,t} = H_{s,s+1} \otimes \dots \otimes H_{t-1,t}
\]
for Hilbert spaces $ H_{s,t} = L_2(\Om_{s,t},P_{s,t}) $. Each $
H_{t,t+1} $ is a direct sum,
\[
H_{t,t+1} = H^{(1)}_{t,t+1} \oplus H^{(0)}_{t,t+1} \, ,
\]
where $ H^{(1)}_{t,t+1} $ is the one-dimensional subspace of
constants, and $ H^{(0)}_{t,t+1} $ is its orthogonal complement, the
subspace of centered (zero mean) random variables. We open the
brackets:
\begin{multline*}
H = H_{1,n} = H_{1,2} \otimes \dots \otimes H_{n-1,n} = \\
= \( H^{(1)}_{1,2} \oplus H^{(0)}_{1,2} \) \otimes \dots \otimes \(
H^{(1)}_{n-1,n} \oplus H^{(0)}_{n-1,n} \) = \bigoplus_{C\in\Comp(T)}
H_C \, ,
\end{multline*}
where $ \Comp(T) = 2^{T'} $ is the set of all subsets $ C $ of the set
$ T' = \{ 1.5, 2.5, \dots, \linebreak[0]
{n-0.5} \} $, and
\[
H_C = \bigotimes_{t=1}^{n-1} H^{(k(t))}_{t,t+1} \, , \qquad
k(t) = \begin{cases}
 0, &\text{if $ t+0.5 \in C $},\\
 1, &\text{otherwise}.
\end{cases}
\]
(Later, $ \Comp(T) $ will consist of \emph{compact} sets.) In other
words, $ H_C $ is spanned by products $ \prod_{t:t+0.5\in C}
f_{t,t+1} $ for $ f_{t,t+1} \in L_2(\Om_{t,t+1},P_{t,t+1}) $, $ \Ex
f_{t,t+1} = 0 $.

The orthogonal decomposition of $ H $ leads to a
\emph{projection-valued measure} $ Q $ on $ \Comp(T) $. That is, a
Hermitian projection operator $ Q(A) : H \to H $ corresponds to every
$ A \subset \Comp(T) $, satisfying
\begin{equation}\label{9a1}
\begin{gathered}
Q ( A \cap B ) = Q(A) Q(B) \, , \\
Q ( \Comp(T) \setminus A ) = \One_H - Q(A) \, , \\
Q ( A \cup B ) = Q(A) + Q(B) \quad \text{if } A \cap B = \emptyset
\end{gathered}
\end{equation}
for $ A,B \subset \Comp(T) $. Such $ Q $ is uniquely determined by
\[
Q ( \{C\} ) H = H_C \quad \text{for } C \in \Comp(T)
\]
or alternatively, by
\[
Q ( \{ C : t+0.5 \notin C \} ) = \cE{ \cdot }{ \F_{1,t} \vee
\F_{t+1,n} } \quad \text{for } t=1,\dots,n-1
\]
(the operator of conditional expectation, given $ \om_{1,2}, \dots,
\om_{t-1,t} $ and $ \om_{t+1,t+2}, \dots, \linebreak[0]
\om_{n-1,n} $).

The operators $ U^\rho $ introduced in Sect.~\ref{5b} are easily
expressed in terms of $ H_C $,
\[
U^\rho = \int_{\Comp(T)} \rho^{|C|} \, Q(\D C) = \sum_{C\in\Comp(T)}
\rho^{|C|} Q(\{C\}) = \bigoplus_{C\in\Comp(T)} \rho^{|C|} \, ;
\]
that is, each $ H_C $ is an eigenspace, its eigenvalue being $
\rho^{|C|} $ (here $ |C| $ stands for the number of elements in $ C
$). Accordingly, the eigenspaces $ H_k $ introduced in Sect.~\ref{5b}
are
\[
H_k = Q ( \{ C : |C|=k \} ) H = \bigoplus_{C:|C|=k} H_C \, .
\]

Every $ f \in H $ leads to a measure $ \mu_f $ on the set $ \Comp(T)
$, called the \emph{spectral measure} of $ f $; namely,
\[
\mu_f (A) = \ip{ Q(A)f }{ f } = \| Q_A f \|^2 \quad \text{for } A
\subset \Comp(T) \, .
\]
Clearly, $ \mu_f(\Comp(T)) = \|f\|^2 $. Assuming $ \|f\|=1 $ we get a
probability measure $ \mu_f $, thus, $ C $ may be thought of as a
random set. However, this random set is defined on the probability
space $ (\Comp(T),\mu_f) $ rather than $ (\Om,P) $.

Let $ g \in H_1 = \bigoplus_{|C|=1} H_C $ and $ f = \Exp g $ in the
sense of Sect.~\ref{6b}; that is, $ g = g_{1,2} + \dots + g_{n-1,n} $
and $ f = (1+g_{1,2}) \dots (1+g_{n-1,n}) $. Then $ \mu_f $ is a
product measure, $ \mu_f(\{C\}) = \prod_{t:t+0.5\in C} \| g_{t,t+1}
\|^2 $. That is, the probability measure $ \mu_{f/\|f\|} $ describes a
random set $ C $ that contains each point $ t+0.5 $ with probability $
\|g_{t,t+1}\|^2 / (1+\|g_{t,t+1}\|^2) $, independently of others. In
contrast, the probability measure $ \mu_{g/\|g\|} $ describes a
single-point random set $ C $, equal to $ \{ t+0.5 \} $ with
probability $ \|g_{t,t+1}\|^2 / \|g\|^2 $.

Spaces $ L_2 $ may be replaced with arbitrary \emph{pointed Hilbert
spaces} (recall Sect.~\ref{6d}).

Now we turn to the discrete example of Sect.~\ref{1b}, the \flow{\Z_m}
$ (X_{s,t}^{\text{\ref{1b}}})_{s<t;s,t\in T} $ over the \emph{infinite} time
set $ T = \{ 0,1,2,\dots \} \cup \{\infty\} $, and the corresponding
continuous product of probability spaces. The set $ \Comp(T) $ of all
compact subsets of the space $ T' = \{ 0.5, 1.5, \dots \} \cup
\{\infty\} $ is a standard measurable space, it may be identified with
the set $ 2^{T'\setminus\{\infty\}} $ of \emph{all} subsets of $ \{
0.5, 1.5, \dots \} $. The \sif\ on $ 2^{T'\setminus\{\infty\}} $ is
generated by the algebra of cylinder subsets. On this algebra we
define a projection-valued measure (or rather, additive set function)
$ Q $ by opening (a finite number of) brackets in $ H = \(
H^{(1)}_{0,1} \oplus H^{(0)}_{0,1} \) \otimes \dots \otimes \(
H^{(1)}_{n-1,n} \oplus H^{(0)}_{n-1,n} \) \otimes H_{n,\infty} $. Its
$ \si $\nobreakdash-\hspace{0pt}additive extension to the \sif,
evidently unique, exists by Kolmogorov's theorem combined with a
simple argument \cite[3d11]{Ts03}.

The spectral measure $ \mu_f $ of the random variable $ f = \exp \(
\frac{2\pi i}m X_{0,\infty}^{\text{\ref{1b}}} \) $ is easy to calculate. For
any $ n $ we have a product, $ f = \exp \( \frac{2\pi i}m
X_{0,1}^{\text{\ref{1b}}} \) \dots \exp \( \frac{2\pi i}m
X_{n-1,n}^{\text{\ref{1b}}} \) \linebreak[0]
 \exp \( \frac{2\pi i}m
X_{n,\infty}^{\text{\ref{1b}}} \) $, therefore the random set contains each
of the points $ \frac12, \frac32, \dots, \linebreak[1]
 n-\frac12 $ with probability
$ 1 - | \Ex \exp \( \frac{2\pi i}m X_{0,1}^{\text{\ref{1b}}} \) |^2 = \sin^2
\frac\pi m $, \emph{independently of others.} It is just a Bernoulli
process, an infinite sequence of independent equiprobable events. The
random set is infinite a.s., which means that $ f \in H_\infty $, the
sensitive space orthogonal to the stable space $ H_0 \oplus H_1 \oplus
\dots $ See also \ref{1c1} and Sect.~\ref{5d}.

Continuous-time counterparts of these ideas for the commutative setup
(see \ref{9b}) are introduced by Tsirelson \cite[Sect.~2]{Ts98}; for
the noncommutative setup --- by Liebscher \cite{Li} and Tsirelson
\cite[Sect.~2]{Ts02} independently.

\mysubsection{Probability spaces}
\label{9b}

Dealing with the time set $ \R $ we introduce the set $ \Comp(\R) $%
\index{zzc@$ \Comp(\R) $, space}
of
all compact subsets $ C \subset \R $. Endowed with the \sif\ generated
by the sets of the form $ \{ C \in \Comp(\R) : C \cap U = \emptyset \}
$, where $ U $ varies over all open subsets of $ \R $, $ \Comp(\R) $
becomes a standard measurable space \cite[12.6]{Ke}. Nowhere dense
compact sets $ C \subset \R $ (that is, with no interior points) are a
measurable subset of $ \Comp(\R) $.

By a projection-valued measure\index{projection-valued measure}
on $ \Comp(\R) $ (over a Hilbert space
$ H $) we mean a family of Hermitian projection operators $ Q(A) : H
\to H $ given for all measurable $ A \subset \Comp(\R) $, satisfying
\eqref{9a1} and countable additivity:
\[
\text{if } A_n \uparrow A \text{ then } Q(A_n) \to Q(A)
\text{ strongly}
\]
(that is, $ \| Q(A_n)f - Q(A)f \| \to 0 $ for every $ f \in H $).

\begin{theorem}\label{9b1}\index{theorem}
\cite[Th.~3d12 and (3d3)]{Ts03}
Let $ (\Om,\F,P), (\F_{s,t})_{s<t} $ be a continuous product of
probability spaces, satisfying the upward continuity condition
\eqref{2c2a}. Then

(a)
there exists one and only one projection-valued measure $ Q $ on $
\Comp(\R) $ (over $ H = L_2(\Om,\F,P) $) such that
\begin{equation}\label{9b1a}
Q ( \{ C : C \cap (s,t) = \emptyset \} ) = \cE{\cdot}{ \F_{-\infty,s}
\vee \F_{t,\infty} }
\end{equation}
whenever $ -\infty < s < t < \infty $;

(b)
$ Q $ is concentrated on (the set of all) nowhere dense compact sets $
C \subset \R $;

(c)
$ Q ( \{ C : t \in C \} ) = 0 $ for every $ t \in \R $.
\end{theorem}

Throughout Sect.~\ref{9b}, the upward continuity condition is assumed
for all continuous products of probability spaces. The time set $ \R $
may be enlarged to $ [-\infty,\infty] $, but it is the same, since a
compact subset of $ [-\infty,\infty] $ not containing $ \pm\infty $
(like \ref{9b1}(c)) is in fact a compact subset of $ \R $.

Using the relation $ Q(A\cap B) = Q(A) Q(B) $ we get
\begin{equation}\label{9b1c}
Q ( \{ C : C \cap ((r,s)\cup(t,u)) = \emptyset \} ) = \cE{\cdot}{
\F_{-\infty,r} \vee \F_{s,t} \vee \F_{u,\infty} }
\end{equation}
for $ -\infty < r < s < t < u < \infty $, and the same for any finite
number of intervals.

It is convenient to express a relation of the form $ Q(A) = \One_H $
by saying that `almost all spectral sets\index{spectral sets}
belong to $ A $'. Thus, (b)
says that almost all spectral sets are nowhere dense, while (c) says
that for every $ t $, almost all spectral sets do not contain $ t
$. By the way, the latter shows that $ Q ( \{ C : C \cap (s,t) =
\emptyset \} ) = Q ( \{ C : C \cap [s,t] = \emptyset \} )$. Also,
applying Fubini theorem we see that almost every spectral set is of
zero Lebesgue (or other) measure.

As before, every $ f \in H $, $ \| f \| = 1 $, has its \emph{spectral
measure}\index{spectral measure!over a noise or cont.\ prod.}
\[
\mu_f (A) = \ip{ Q(A)f }{ f } = \| Q_A f \|^2 \, ,
\]
a probability measure on $ \Comp(\R) $. The relation $ Q(A) = \One_H $
holds if and only if $ \mu_f(A) = 1 $ for all $ f $.

As before,
\begin{equation}\label{9b1e}
U^\rho = \int_{\Comp(\R)} \rho^{|C|} \, Q(\D C) \, ,
\end{equation}
but now $ C $ may be infinite, in which case $ \rho^{|C|} = 0 $ for
all $ \rho \in [0,1) $. We have
\[
L_2 (\F_\stable) = Q ( \{ C : |C| < \infty \} ) \, ;
\]
accordingly, the sensitive subspace is $ Q ( \{ C : |C| = \infty \} )
$. A continuous product of probability spaces is classical if and
only if almost all spectral sets are finite. The classical part of a
continuous product of probability spaces is trivial if and only if
almost all nonempty spectral sets are infinite. Also,
\begin{equation}\label{9b1g}
\ip{ U^\rho f}{f} = \int_{\Comp(\R)} \rho^{|C|} \, \mu_f (\D C) \, ;
\end{equation}
$ f $ is stable if and only if $ \mu_f
$\nobreakdash-\hspace{0pt}almost all spectral sets are finite; 
$ f $ is sensitive if and only if $ \mu_f
$\nobreakdash-\hspace{0pt}almost all nonempty spectral sets are
infinite.

Let $ g \in H_1 = Q ( \{ C : |C|=1 \} ) $ and $ f = \Exp g $ (in the
sense of Sect.~\ref{6b}), then $ \mu_{f/\|f\|} $ is the distribution
of a Poisson point process; the mean number of points on $ (s,t) $ is
equal to $ \| g_{s,t} \|^2 $ where $ g_{s,t} = \cE{ g }{ \F_{s,t} }
$.

If $ Q(A) = 0 $ then $ \mu_f(A) = 0 $ for all $ f $. For some $ f $
(for instance, $ f=1 $) the relation $ \mu_f(A) = 0 $ does not imply $
Q(A) = 0 $; however, such $ f $ are exceptional (in fact, they are a
meager subset of $ H $). For other, typical $ f \in H $ the class of $
A $ such that $ \mu_f(A) = 0 $ does not depend on $ f $. It means that
all typical spectral measures are equivalent (that is, mutually
absolutely continuous). Thus, each continuous product of probability
spaces leads to a \emph{measure type} (or `class') $ \cM $, --- an
equivalence class of probability measures on $ \Comp(\R) $.

A measure belonging to the class $ \cM $ has one and only one
atom. Namely, a spectral set has a chance to be empty.

The same construction may be applied to the restriction of a given
continuous product of probability spaces to a given time interval $
(s,t) \subset \R $. We get a projection-valued measure $ Q_{s,t} $ on
the space $ \Comp(s,t) $ of all compact subsets of $ (s,t) $, over $
H_{s,t} = L_2(\F_{s,t}) $, and a measure type $ \cM_{s,t} $ on $
\Comp(s,t) $. It appears that
\[
Q_{r,t} = Q_{r,s} \otimes Q_{s,t} \quad \text{for }
-\infty \le r < s < t \le \infty
\]
in the sense that $ Q_{r,t} (A \times B) = Q_{r,s}(A) \otimes
Q_{s,t}(B) $ for all measurable $ A \subset \Comp(r,s) $ and $ B
\subset \Comp(s,t) $. Here $ A \times B = \{ C_1 \cup C_2 : C_1 \in A,
C_2 \in B \} $. Note that $ \Comp(r,s) \times \Comp(s,t) $ is a subset
of $ \Comp(r,t) $ consisting of all $ C \in \Comp(r,t) $ such that $ s
\notin C $; it does not harm, since $ Q_{r,t} (\{C:s\in C\}) = 0 $. We
get
\begin{equation}\label{9b2}
\cM_{r,t} = \cM_{r,s} \times \cM_{s,t} \quad \text{for }
-\infty \le r < s < t \le \infty
\end{equation}
in the sense that the product measure $ \mu = \mu_1 \times \mu_2 $
belongs to $ \cM_{r,t} $ for some (therefore, all) $ \mu_1 \in
\cM_{r,s} $, $ \mu_2 \in \cM_{s,t} $. Also, $ \cM_{s,t} $ is the
marginal distribution of $ \cM $ in the sense that the marginal of $
\mu $ belongs to $ \cM_{s,t} $ for some (therefore, all) $ \mu \in \cM
$.

\begin{definition}\label{9b3}
A measure type $ \cM $ on $ \Comp(\R) $ is \emph{factorizing,}%
\index{factorizing!measure type}
if $ \{
C : s \in C \} $ is $\cM$\nobreakdash-\hspace{0pt}negligible for each
$ s $, and the marginals of $ \cM $ satisfy \eqref{9b2}.
\end{definition}

See also \cite[Def.~4.1]{Li}. Every continuous product of probability
spaces (satisfying the upward continuity condition) leads to a
factorizing measure type on $ \Comp(\R) $.

\begin{question}\label{9b7aa}\index{question}
Does every factorizing measure type on $ \Comp(\R) $ correspond to
some continuous product of probability spaces? (See also \ref{9c2} and
\ref{10a1}.)
\end{question}

A fragment $ C \cap [s,t] $ of a spectral set $ C $ is also a spectral
set. More formally, for every interval $ [s,t] \subset \R $ the map $
C \mapsto C \cap [s,t] $ of $ \Comp(\R) $ to itself is
$\cM$\nobreakdash-\hspace{0pt}nonsingular. That is, the inverse image
of each negligible set is negligible.

The following three conditions are thus equivalent:
\begin{itemize}
\item[]
almost all nonempty spectral sets are infinite;
\item[]
almost all spectral sets are perfect (no isolated points; the empty
set is perfect);
\item[]
the classical part of the continuous product of probability spaces is
trivial.
\end{itemize}

\begin{example}\label{9b7a}
All classical noises (except for the trivial case, that is, assuming $
\dim L_2(\Om) > 1 $) lead to the same measure type $ \cM =
\cM_\Poisson $ on $ \Comp(\R) $. Namely, $ \cM_\Poisson $ contains the
distribution of the Poisson random subset of $ \R $ whose intensity
measure is finite and equivalent to the Lebesgue measure on $ \R
$. (Equivalent finite intensity measures lead to equivalent Poisson
measures.) Clearly, $ \cM_\Poisson $ is factorizing and
shift-invariant.
\end{example}

\begin{example}\label{9b8}
For the noise of splitting, considered in Sect.~\ref{4d}, spectral
sets are described by Warren \cite{Wa1}, see also Watanabe
\cite{Wat0}. They are at most countable. Moreover, almost all spectral
sets $ C $ satisfy $ |C'| < \infty $ (or equivalently, $ C'' =
\emptyset $); here $ C' $ stands for the set of all limit (that is,
accumulation) points of $ C $. The relation $ |C'| < \infty $ holds
also for the noise of stickiness, considered in Sect.~\ref{4e}
\cite[6b4]{Ts03}.
\end{example}

\begin{example}\label{9b9}
For the black noise of coalescence, considered in Sect.~\ref{7f},
almost all nonempty spectral sets are perfect (therefore,
uncountable). For some especially simple random variable $ f $ the
spectral measure $ \mu_f $ is described by Tsirelson (see
\cite[Sect.~7d]{Ts03} and references there); spectral sets are of
Hausdorff dimension $ 1/2 $.
\end{example}

For every noise, finite spectral sets select a subnoise (namely, the
classical part of the noise). Similarly, spectral sets $ C $ satisfying
$ C'' = \emptyset $ select a subnoise \cite[Th.~6b2(b) for $ \al=2 $,
and 6b10]{Ts03}. The same holds for $ C''' = \emptyset $ ($ \al=3 $)
and higher levels of the Cantor-Bendixson hierarchy, but for now we
have no examples.

For every $ \al \in (0,1) $ spectral sets of Hausdorff dimension at
most $ \al $ select a subnoise \cite[Th.~6b9, 6b10]{Ts03}.
(See also Sect.~\ref{9c}.)

Consider (necessarily countable or finite) spectral sets $ C $ such
that $ { \forall t \in C } \linebreak[0]
\; { \exists \eps>0 } \; \( C \cap (t,t+\eps) = \emptyset \) $; that
is, accumulation is allowed from the left but not from the right. Do
these select a subnoise? I do not know. (See also \cite[Question
6b12]{Ts03}.)

All said above holds for (not just homogeneous) continuous products of
probability spaces.

Here is a generalization of Theorem \ref{5c4}. In Item (b), $ \dim C $
stands for the Hausdorff dimension of $ C $. In Item (a), $
C^{(\kappa)} $ is defined recursively: $ C^{(0)} = C $, $
C^{(\kappa+1)} = (C^{(\kappa)})' $ and $ C^{(\kappa)} =
\cap_{\kappa_1<\kappa} C^{(\kappa_1)} $ for limit ordinals $ \kappa
$. The case $ \kappa = 1 $ of (a) returns us to \ref{5c4}.

\begin{theorem}\label{9b6}\index{theorem}
Let $ \( (\Om,P), (\F_{s,t})_{s<t} \) $ be a continuous product of
probability spaces, satisfying the upward continuity condition
\eqref{2c2a}. Then:

(a) For every finite or countable ordinal $ \kappa $ there exists a
symmetric self-joining $ (\al_\kappa,\be_\kappa) $ of the given
product such that

\begin{itemize}
\item[]
$ \Ex (f\circ\al_\kappa) (g\circ\be_\kappa) = \Ex (fg) $ for all $ f,g
\in L_2(\Om) $ such that the relation $ C^{(\kappa)} = \emptyset $
holds both for \almost{\mu_f} all $ C $ and for \almost{\mu_g} all $ C
$;
\item[]
$ \Ex (f\circ\al_\kappa) (g\circ\be_\kappa) = 0 $ for all $ f \in
L_2(\Om) $ such that $ C^{(\kappa)} \ne \emptyset $ for \almost{\mu_f}
all $ C $, and all $ g \in L_2(\Om) $.
\end{itemize}

The self-joining $ (\al_\kappa,\be_\kappa) $ is unique up to
isomorphism.

(b) For every $ \theta \in (0,1) $ there exists a symmetric
self-joining $ (\al_\theta,\be_\theta) $ of the given product such
that

\begin{itemize}
\item[]
$ \Ex (f\circ\al_\theta) (g\circ\be_\theta) = \Ex (fg) $ for all $ f,g
\in L_2(\Om) $ such that the relation $ \dim C \le \theta $ holds both
for \almost{\mu_f} all $ C $ and for \almost{\mu_g} all $ C $;
\item[]
$ \Ex (f\circ\al_\theta) (g\circ\be_\theta) = 0 $ for all $ f \in
L_2(\Om) $ such that $ \dim C > \theta $ for \almost{\mu_f} all $ C $,
and all $ g \in L_2(\Om) $.
\end{itemize}

The self-joining $ (\al_\theta,\be_\theta) $ is unique up to
isomorphism.

\end{theorem}

\beginproof
We combine the arguments of \cite[Sect.~6b]{Ts03} with the compactness
of the space of joinings. A symmetric self-joining corresponds to
every element of the set $ S $ (of some Borel functions $ \Comp(\R)
\to [0,1] $) introduced in \cite[Sect.~5b]{Ts03} and used in
\cite[Sect.~6b, especially 6b1, 6b8]{Ts03}.
\proofend

\begin{corollary}\label{9b7}
The following two conditions are equivalent for any $ \theta \in (0,1)
$ and any continuous product of probability spaces, corresponding to a
flow system $ (X_{s,t})_{s<t} $:

(a) $ \dim C \le \theta $ for almost all spectral sets $ C $;

(b) $ \dim C \le \theta $ for \almost{\mu_f} all $ C $, where $ f =
\phi(X_{s,t}) $, for all $ s<t $ and all bounded Borel functions $
\phi : G_{s,t} \to \R $ (a single $ \phi $ is enough if it is
one-to-one).
\end{corollary}

Proof: similar to \ref{5c5}, using \ref{9b6}(b). See also \cite[Proof
of Th.~1.3, especially (4.2)]{WW03}.

\begin{remark}\label{9b12}
Five more corollaries similar to \ref{9b7} are left to the
reader. Namely, each one of \ref{5c5}, \ref{7d7}, \ref{7h7} may be
generalized using each one of \ref{9b6}(a), \ref{9b6}(b), giving $ 3
\times 2 = 6 $ corollaries, \ref{9b7} being one of them.
\end{remark}

Similarly to \eqref{9b1a} we have for $ s<t $
\begin{multline}\label{9b13}
Q ( \{ C : C' \cap (s,t) = \emptyset \} ) = Q ( \{ C : | C \cap (s,t)
 | < \infty \} ) = \\
= \cE{\cdot}{ \F_{-\infty,s} \vee \F_{s,t}^\stable \vee \F_{t,\infty}
 } \, .
\end{multline}
The counterpart of \eqref{9b1c} for $ C' $ is left to the reader. The
operator \eqref{9b1a} corresponds to a self-joining (of the continuous
product of probability spaces), a combination of $ (\al_\rho,\be_\rho)
$ with different $ \rho $; namely, $ \rho = 0 $ on $ (s,t) $ and $
\rho = 1 $ on $ (-\infty,s) \cup (t,\infty) $. The same holds for the
operator \eqref{9b13}; still, $ \rho = 1 $ on $ (-\infty,s) \cup
(t,\infty) $, but on $ (s,t) $ we use $ \rho=1- $ (recall \ref{5c4})!

Similarly to \eqref{9b1e}, \eqref{9b1g} we may construct operators $
V^\rho : L_2(\Om) \to L_2(\Om) $ by
\[
\begin{gathered}
V^\rho = \int_{\Comp(\R)} \rho^{|C\mskip1mu'|} \, Q(\D C) \, , \\
\ip{ V^\rho f}{f} = \int_{\Comp(\R)} \rho^{|C\mskip1mu'|} \, \mu_f (\D
 C)
\end{gathered}
\]
(see also \cite[Th.~6b3 and before]{Ts03}). This is another
generalization of the Ornstein-Uhlenbeck semigroup; it perturbs
sensitive random variables, while stable random variables are intact.

\mysubsection{Example: nonclassical Harris flows (Warren and Watanabe)}
\label{9c}

Similarly to the Arratia flow (considered in Sect.~\ref{7f}), a Harris
flow is a flow system of \Smap s (from $ \X $ to $ \X $, where $ \X $
may be the circle $ \T = \R/\Z $ or the line $ \R $), homogeneous both
in time and in space, such that the one-point motion is the (standard)
Brownian motion in $ \X $. However, Arratia's particles are
independent before coalescence, while Harris' particles are correlated
all the time. Namely,
\[
\D X_t \D Y_t = b(X_t-Y_t) \, \D t \, ,
\]
where $ X_t, Y_t $ are coordinates of two particles at time $ t $, and
the \emph{correlation function} $ b(\cdot) $ is a given positive
definite function $ \X \to \R $; $ b(0)=1 $, $ b(-x)=b(x) $. Thus, $
(X_t-Y_t)_t $ is a diffusion process, it becomes the (standard)
Brownian motion in $ \X $ under a random time change, the new time
being $ 2 \int_0^t \( 1 - b(X_s-Y_s) \) \, \D s $, as long as $ X_s
\ne Y_s $. Three cases emerge (see also \cite[Sect.~7.4]{LJR}):

\begin{itemize}
\item[] non-coalescing case:
\[
\text{the function } x \mapsto \frac{x}{ 1-b(x) } \text{ is
non-integrable near $ 0 $} \, ;
\]
\item[] classical coalescing case:
\[
x \mapsto \frac{x}{ 1-b(x) } \text{ is integrable near $ 0 $, but } x
\mapsto \frac{1}{ 1-b(x) } \text{ is not} \, ;
\]
\item[] nonclassical case:
\[
x \mapsto \frac{1}{ 1-b(x) } \text{ is integrable near $ 0 $} \, .
\]
\end{itemize}

The non-coalescing case: $ X_t - Y_t $ never vanishes; rather, $ X_t -
Y_t \to 0 $ as $ t \to \infty $, and $ \int_0^\infty \( 1 - b(X_s-Y_s)
\) \, \D s < \infty $. Two particles cannot meet, since their
rapprochement is infinitely long. The origin is a natural boundary for
the diffusion process $ (X_t - Y_t)_t $. In particular it happens if $
b(\cdot) $ is twice continuously differentiable (which leads to a flow
of homeomorphisms, see \cite[Sect.~8]{Ha}).

The classical coalescing case: $ X_t - Y_t $ vanishes at some $ t $
and remains $ 0 $ forever. Two particles cannot diverge after meeting;
even a small divergence would take infinite time (and no wonder: it
must involve infinitely many more meetings). The origin is an exit
boundary. In particular it happens if $ 1-b(x) \sim |x|^\al $ (as $ x
\to 0 $), $ 1 \le \al < 2 $.

In the nonclassical case (in particular, $ 1-b(x) \sim |x|^\al $, $ 0
< \al < 1 $) we have an additional freedom. The origin is
a regular boundary; we may postulate it to be absorbing, sticky or
reflecting. Only absorption leads to a flow system of \Smap s
(\emph{`the coalescing nonclassical case'}); stickiness and reflection
lead rather to flow systems of \Skernel s.

The corresponding noise is classical in the non-coalescing case, as
well as in the classical coalescing case \cite[Sect.~7.4]{LJR},
\cite[Th.~1.1]{WW03}.

We turn to the coalescing nonclassical case. Here, the noise is
nonclassical \cite[Sect.~7.4]{LJR}, \cite[Th.~1.1]{WW03}. If the
correlation function $ b(\cdot) $ is continuous (and smooth outside
the origin, and strictly positive definite) then the classical part of
the noise is generated by infinitely many independent Brownian motions
\cite[Sects.~6.4, 7.4]{LJR}, \cite[Th.~1.1]{WW03}. By the way, the
Arratia flow does not fit into this framework (it needs a
discontinuous correlation function, $ b(0)=1 $, $ b(x)=0 $ for $ x \ne
0 $), and leads to a black noise (recall Sect.~\ref{7f}).

Assuming $ 1-b(x) \asymp |x|^\al $ as $ x \to 0 $ for some $ \al \in
(0,1) $ (and some additional technical conditions on $ b(\cdot) $),
Warren and Watanabe \cite[Th.~1.3]{WW03} manage to find the Hausdorff
dimension; the inequality
\begin{equation}\label{9c1}
\dim C \le \frac{1-\al}{2-\al}
\end{equation}
holds for almost all spectral sets $ C $, but the strict inequality $
\dim C < \frac{1-\al}{2-\al} $ does not.

By \ref{9b12} it is sufficient to prove \eqref{9c1} for \almost{\mu_f}
all $ C $, where $ f $ is a random variable of the form $ \Xi_x^{s,t}
$ (the coordinate at $ t $ of a particle starting at $ s,x $). The set
$ C' $ of limit points of $ C $ is related to the self-joining $
(\al_{1-},\be_{1-}) $ (recall the paragraph after \eqref{9b13}). An
explicit description of the joining in terms of diffusions, found by
Warren and Watanabe, leads them to an explicit description of the
random set $ C' $. It is (distributed like) the set of zeros of a
diffusion process. If $ C' $ happens to be nonempty then $ \dim C' =
\frac{1-\al}{2-\al} $.

\begin{question}\label{9c2}\index{question}
Can spectral sets of a noise be of Hausdorff dimension greater than $
\frac12 $? (See also \ref{9b7aa} and \ref{10a1}.)
\end{question}

\mysubsection{Hilbert spaces}
\label{9d}

Throughout Sect.~\ref{9b} the spaces $ L_2 $ may be replaced with
arbitrary \emph{pointed Hilbert spaces} (recall Sect.~\ref{6d}). I do
it explicitly for Theorem \ref{9b1}; the rest of the work is left to
the reader.

\begin{theorem}\label{9d1}\index{theorem}
Let $ (H_{s,t},u_{s,t})_{s<t} $ be a continuous product of
pointed Hilbert spaces, satisfying the upward continuity condition
\eqref{6d13}. Then

(a)
there exists one and only one projection-valued measure $ Q $ on $
\Comp(\R) $ (over $ H = H_{-\infty,\infty} $) such that
\begin{equation}
Q ( \{ C : C \cap (s,t) = \emptyset \} ) H = H_{-\infty,s} u_{s,t}
H_{t,\infty}
\end{equation}
whenever $ -\infty < s < t < \infty $;

(b,c): the same as in \ref{9b1}.
\end{theorem}

(See also Theorem \ref{9d5}.)
As before, $ H_{-\infty,s} u_{s,t} H_{t,\infty} $ stands for the image
of $ H_{-\infty,s} \otimes u_{s,t} \otimes H_{t,\infty} $ under the
given unitary operator $ H_{-\infty,s} \otimes H_{s,t} \otimes
H_{t,\infty} \to H_{-\infty,\infty} $.

The upward continuity condition will be assumed for all continuous
products of pointed Hilbert spaces.

Every continuous product of pointed Hilbert spaces leads to a
factorizing measure type on $ \Comp(\R) $. The continuous product is
classical if and only if almost all spectral sets are finite. In this
case the factorizing measure type is Poissonian (but the underlying
measure type on $ \R $ need not be shift-invariant). By Theorem
\ref{6e3}, classicality does not depend on the choice of a unit.

\begin{question}\label{9d3}\index{question}
(\cite[Notes 3.11, 10.2]{Li})
Does the factorizing measure type depend on the choice of a unit?
\end{question}

The answer `does not depend' would result from an affirmative answer to the
following question (see \cite[Def.~8.2]{Bh01} and \cite[Notes 3.6, 5.8 and
Sect.~11 (question 1)]{Li}).
Let $ (H_{s,t})_{s<t} $ be a continuous product of Hilbert spaces,
satisfying \ref{6e1}(a-c), and $ (u_{s,t})_{s<t} $, $ (v_{s,t})_{s<t}
$ two units. Does there exist an automorphism of $ (H_{s,t})_{s<t} $
sending $ (u_{s,t})_{s<t} $ to $ (v_{s,t})_{s<t} $? In other words,
are the two continuous products of pointed Hilbert spaces $
(H_{s,t},u_{s,t})_{s<t} $, $ (H_{s,t},v_{s,t})_{s<t} $ isomorphic?
(An automorphism may be defined as an invertible embedding to itself,
recall \ref{6d7}.) However, the latter question is recently answered in the
negative \cite{Ts04}.

Every continuous product of probability spaces $
(\Om_{s,t},P_{s,t})_{s<t} $ leads to a continuous product of pointed
Hilbert spaces $ (H_{s,t},u_{s,t})_{s<t} $ (namely, $ H_{s,t} = L_2
(\Om_{s,t},P_{s,t}) $ and $ u_{s,t} = 1 $ on $ \Om_{s,t} $) and
further, to a continuous product of Hilbert spaces $ (H_{s,t})_{s<t}
$. The classical part of $ (\Om_{s,t},P_{s,t})_{s<t} $ corresponds to
the classical part of $ (H_{s,t},u_{s,t})_{s<t} $, in fact, the
classical part of $ (H_{s,t})_{s<t} $.

We say that the classical part of $ (\Om_{s,t},P_{s,t})_{s<t} $ is
trivial, if all stable random variables are constant, that is, $ \dim
L_2(\F^\stable) = 1 $. On the other hand, we say that the classical
part of $ (H_{s,t})_{s<t} $ is trivial, if $ \dim H^\cls = 0 $ (recall
\eqref{6e4a}), which means, no decomposable vectors at all. The latter
never happens for $ (H_{s,t})_{s<t} $ obtained from some $
(\Om_{s,t},P_{s,t})_{s<t} $; indeed, $ \dim H^\cls = \dim
L_2(\F^\stable) \ge 1 $. However, it may happen that $ \dim H^\cls = 1
$; every black noise is an example.

If $ \dim H^\cls = 1 $ then all units are basically the same, and we
may treat the corresponding factorizing measure type on $ \Comp(\R) $
as an (isomorphic) invariant of $ (H_{s,t})_{s<t} $. In this case the
projection-valued measure $ Q $ may be attributed to the embedding $
H^\cls \subset H $ (rather than the unit). More generally, some $ Q $
may be attributed to the embedding $ H^\cls \subset H $ assuming only
$ \dim H^\cls \ge 1 $ (as explained below), which leads to invariants
of (not pointed) continuous products of Hilbert spaces.

Let $ (H_{s,t})_{s<t} $ be a continuous product of Hilbert spaces, and
$ H'_{s,t} \subset H_{s,t} $ be (closed linear) subspaces satisfying $
\dim H'_{s,t} \ge 1 $ and $ H'_{r,s} H'_{s,t} = H'_{r,t} $ for $ r<s<t
$. Then $ (H'_{s,t})_{s<t} $ is also a continuous product of Hilbert
spaces, and identical maps $ H'_{s,t} \to H_{s,t} $ are an embedding
of the latter product to the former (recall \ref{6d7}). Such a pair of
continuous products may be called an \emph{embedded pair.}

Instead of the subspaces $ H'_{s,t} $ we may consider (following
\cite[Sect.~3.2]{Li}) Hermitian projections $ P_{s,t} : H \to H $
satisfying $ P_{s,t} \in \A_{s,t} $, $ P_{s,t} \ne 0 $ and $ P_{r,s}
P_{s,t} = P_{r,t} $ for $ r<s<t $ (the algebras $ \A_{s,t} $ appear in
Sect.~\ref{3b}). Namely, $ P_{s,t} H = H_{-\infty,s} H'_{s,t}
H_{t,\infty} $. Monotonicity of $ P_{s,t} $ in $ s $ and $ t $ ensures
existence of the limit $ P_{s-,t+} = \lim_{\eps\to0+}
P_{s-\eps,t+\eps} $ in the strong operator topology. (As before, $
-\infty-\eps = -\infty $, $ \infty+\eps = \infty $.) In fact, $
P_{s-,t+} = P_{s,t} $ unless $ s $ or $ t $ belongs to a finite or
countable set of discontinuity points. Theorem \ref{9d1}
is a special case of the following fact.

\begin{theorem}\label{9d5}\index{theorem}
(Liebscher \cite[Th.~1]{Li}, Tsirelson \cite[2.9]{Ts02})\\
Let $ (H'_{s,t}\subset H_{s,t})_{-\infty\le s<t\le\infty} $ be an
embedded pair of continuous products of Hilbert spaces, and $
(P_{s,t})_{s<t} $ the corresponding family of projections. Then there
exists one and only one projection-valued measure $ Q $ on $
\Comp([-\infty,\infty]) $ (over $ H = H_{-\infty,\infty} $) such that
\begin{equation}
Q ( \{ C : C \cap [s,t] = \emptyset \} ) = P_{s-,t+}
\end{equation}
whenever $ -\infty \le s \le t \le \infty $.
\end{theorem}

\beginproof
We define $ Q $ on the algebra generated by sets of the form $ \{ C :
C \cap [s,t] = \emptyset \} $ by $ Q ( \{ C : C \cap [s,t] = \emptyset
\} ) = P_{s-,t+} $ and $ Q ( A \cap B ) = Q(A) Q(B) $ (and
additivity). The algebra generates the Borel \sif\ on $
\Comp([-\infty,\infty]) $, and we extend $ Q $ to the \sif\ using
\cite[2.4, 2.5, 2.6]{Ts02}. \proofend

The time set $ [-\infty,\infty] $ is essential. In the local case (the
time set $ \R $) we have no global space $ H_{-\infty,\infty} $ (and
no embeddings $ H_{-1,1} \subset H_{-2,2} \subset \dots $), thus, no
global $ Q $. Still, we have $ Q_{s,t} $ for $ s<t $ (since the time
set $ [s,t] $ is similar to $ [-\infty,\infty] $).

A measure type on $ \Comp([-\infty,\infty]) $ corresponds to each
embedded pair. Still, for a given $ t $ we have $ t \notin C $ for
almost all spectral sets $ C $, unless $ t $ belongs to a finite or
countable set of discontinuity points. Therefore almost all spectral
sets are nowhere dense, of Lebesgue measure zero. They are compact
subsets of $ \R $, unless $ -\infty $ or $ +\infty $ is a
discontinuity point. If there is no discontinuity points at all then
we get a factorizing measure type (as defined by \ref{9b3}).
Especially, the factorizing measure type corresponding to the
embedding of the classical part is an invariant of a continuous
product of Hilbert spaces \cite[Th.~2]{Li}. In the local case (the
time set $ \R $) we get instead a consistent family of factorizing
measure types on $ \Comp([s,t]) $, $ -\infty < s < t < \infty $.

Let $ (H_{s,t},u_{s,t})_{s<t} $ be a continuous product of pointed
Hilbert spaces (satisfying the upward continuity condition), $ Q^u $
the corresponding projection-valued measure, and $ Q^\U $ the
projection-valued measure corresponding to the embedding of the
classical part. Similarly to \eqref{9b13},
\[
Q^u ( \{ C : C' \cap (s,t) = \emptyset \} ) = Q^u ( \{ C : | C \cap
(s,t) | < \infty \} ) = P_{s+,t-} = P_{s,t} \, ,
\]
thus $ Q^u ( \{ C : C' \cap (s,t) = \emptyset \} ) = Q^\U ( \{ C : C
\cap (s,t) = \emptyset \} ) $, which implies
\begin{equation}\label{9d6}
Q^u ( \{ C : C' \in A \} ) = Q^\U (A)
\end{equation}
for all Borel sets $ A \subset \Comp(\R) $. This is a useful relation,
due to Liebscher \cite[Prop.~3.9]{Li}, between unit-dependent
spectral sets $ C $ and unit-independent spectral sets; the latter
spectral set consists of the limit points of the former!

The factorizing measure types corresponding to $ Q^u $ and $ Q^\U $
will be denoted by $ \cM^u $ and $ \cM^\U $ respectively, and the
corresponding spectral sets by $ C_u $ and $ C_\U $ (these make sense
when saying that $ C $ satisfies something a.s.). In some sense,
\eqref{9d6} allows us to say that $ C'_u = C_\U $ in distribution (for
every $ u $).

If $ (H_{s,t})_{s<t} $ is classical then $ Q^\U (\{\emptyset\}) = \One
$, $ C_u $ is finite a.s., and $ C_\U $ is empty a.s. In general, $ Q^\U
(\{\emptyset\}) H = H^\cls $; $ C_\U $ has a chance to be empty.

Dealing with the probabilistic case ($ H_{s,t} = L_2
(\Om_{s,t},P_{s,t}) $) we assume that $ u $ is `the probabilistic unit'
($ u_{s,t}(\cdot) = 1 $ on $ \Om_{s,t} $), unless otherwise stated.

\begin{example}\label{9d8}
For the noise of splitting (recall \ref{9b8}) $ C'_u $ is finite a.s.,
thus $ C_\U $ is finite a.s. In fact, $ \cM^\U = \cM_\Poisson $ (for $
\cM_\Poisson $ see \ref{9b7a}). The same holds for the noise of
stickiness.
\end{example}

For a black noise (for instance, \ref{9b9}), spectral sets are perfect
a.s.; $ C_\U = C'_u = C_u $.

\begin{question}\label{9d9}\index{question}
The two noises (of splitting and stickiness) mentioned in \ref{9d8}
lead to two continuous products of Hilbert spaces (not pointed!) $
(H_{s,t}^\splt)_{s<t} $ and $ (H_{s,t}^\stick)_{s<t} $. Are these
products isomorphic?
\end{question}

We know that $ \cM^\U_\splt = \cM^\U_\stick = \cM_\Poisson $. However,
$ \cM^u_\splt $ and $ \cM^u_\stick $ are different ($ u $ being the
probabilistic unit in both cases). Namely, accumulation of $
C^u_\splt $ is two-sided, while accumulation of $
C^u_\stick $ is one-sided (from the right only).

\begin{sloppypar}
Note that the continuous products of Hilbert spaces $
(H_{s,t}^\splt)_{s<t} $ and $ (H_{s,t}^\stick)_{s<t} $ in Question
\ref{9d9} are not treated as homogeneous, despite the fact that they
originate from noises. The desired isomorphism need not intertwine
time shifts, it need not be an isomorphism of the two homogeneous
continuous products of Hilbert spaces $ \( (H_{s,t}^\splt)_{s<t},
(\theta_{s,t}^{\splt,h})_{s<t;h} \) $ and $ \(
(H_{s,t}^\stick)_{s<t}, (\theta_{s,t}^{\stick,h})_{s<t;h} \) $ that
arise naturally from the two noises (recall Sect.~\ref{3c}).
\end{sloppypar}

In fact, the two \emph{homogeneous} (and not local!) continuous
products of Hilbert spaces are non-isomorphic. Indeed, Prop.~\ref{6f1}
shows that the `probabilistic' units $ u^\splt, u^\stick $ are
basically the only shift-invariant units. Therefore the relation $
\cM^u_\splt \ne \cM^u_\stick $ denies isomorphism.

However, Prop.~\ref{6f1} does not apply to Arveson systems (recall
Def.~\ref{6f3} and Th.~\ref{6f6}). The Arveson systems $
(H_t^\splt)_{t>0} $, $ (H_t^\stick)_{t>0} $ have units $ u_t = \exp
\( (\al+\be\I) B_t - \frac12 \al^2 + \I \ga t \) $ parametrized by $
\al,\be,\ga \in \R $; here $ (B_t)_t = (a_{0,t})_t $ is the
corresponding Brownian motion. An isomorphism must send a unit into a
unit, but may change the parameters $ \al,\be,\ga $.

\begin{question}\label{9d10}\index{question}
The two noises (of splitting and stickiness) mentioned in \ref{9d8}
lead to two Arveson systems. Are these systems isomorphic?
\end{question}

In \ref{9d9} (unlike \ref{9d10}) on one hand, the isomorphism need not
be shift-invariant; on the other hand, it must act on the global space
$ H_{-\infty,\infty} $.

Let $ (H_t)_{t>0} $ be an Arveson system and $ (u_t)_{t>0} $ a
unit. They lead to a homogeneous local continuous product of pointed
Hilbert spaces, and we may enlarge the time set $ \R $ to $
[-\infty,\infty] $, thus getting a shift-invariant factorizing measure
type $ \cM^u $ on $ \Comp(\R) $; but the enlargement depends on the
unit $ u $. Waiving the unit we have a homogeneous \emph{local}
continuous product of Hilbert spaces, thus, a consistent
shift-invariant factorizing family of measure types $
(\cM_{s,t})_{-\infty<s<t<\infty} $; they correspond to embeddings
$ H^\cls_{s,t} \subset H_{s,t} $. However, the absence of the `global'
embedding $ H^\cls_{-\infty,\infty} \subset H_{-\infty,\infty} $ does
not prevent us from introducing a `global' measure type $ \cM $ on $
\Comp(\R) $. To this ens we note that spaces $ \Comp([-n,n]) $ are
related not only by projections $ \Comp([-n-1,n+1]) \to
\Comp([-n,n]) $, $ C \mapsto C \cap [-n,n] $, but also by embeddings $
\Comp([-n,n]) \subset \Comp([-n-1,n+1]) $. The measure classes (see
\ref{10a2}) $ (\Comp([-n,n]),\cM_{-n,n}) $ form not only a projective
(inverse) system, but also an inductive (direct) system. In contrast to
probability spaces (suitable for projective but not inductive limits),
for measure classes we may take inductive limits (but not projective
limits). Thus, we may define $ \cM $ as the measure type on $
\Comp(\R) $ compatible with all $ \cM_{-n,n} $ in the sense that the
conditional distribution of $ C $ given that $ C \subset [-n,n] $
belongs to $ \cM_{-n,n} $. (Existence and uniqueness of such $ \cM $
is easy to check.)

\begin{corollary}
Every Arveson system of type $ II $ leads to a measure type $ \cM^\U $
on $ \Comp(\R) $; every unit $ u $ of the Arveson system leads to a
measure type $ \cM^u $ on $ \Comp(\R) $; and $ \cM^\U $ is the image of
$ \cM^u $ under the map $ C \mapsto C' $.
\end{corollary}

Do not think, however, that type $ II $ systems are simpler than type
$ III $ systems. In fact, invariants like $ \cM^\U $ catch only a
small part of the structure of type $ II $ systems. The rest of the
structure cannot be simpler than the whole structure of type $ III $
systems! See \cite[Sect.~6.4]{Li}.

\mysection{Continuous products of measure classes}
\label{sec:10}
\mysubsection{From random sets to Hilbert spaces}
\label{10a}

We start with a result that involves an idea of Anatoly Vershik
(private communication, 1994) of a continuous product of measure
classes as a source of a continuous product of Hilbert spaces, and an
idea of Jonathan Warren (private communication, 1999) of constructing
a continuous product of measure classes out of a given random set. See
Liebscher \cite[Prop.~4.1 and Sect.~8.2]{Li} and Tsirelson \cite[Lemma
5.3]{Ts02}. For now this is the richest source of (non-isomorphic)
continuous products of Hilbert spaces with nontrivial classical
part. (See also \ref{9b7aa}.) First, recall that every continuous
product of pointed Hilbert spaces, satisfying the upward continuity
condition, leads to a factorizing measure type on $ \Comp(\R) $, as
explained after Theorem \ref{9d1}.

\begin{theorem}\label{10a1}\index{theorem}
For every factorizing measure type $ \cM $ on $ \Comp(\R) $ (as
defined by \ref{9b3}) there exists a continuous product of pointed
Hilbert spaces, satisfying the upward continuity condition
\eqref{6d13}, such that the corresponding factorizing measure type is
equal to $ \cM $.
\end{theorem}

A proof will be sketched later. `Square roots of measures', introduced
by Accardi \cite{Ac}, are instrumental. For more definitions and basic
facts see \cite[Sect.~14.4]{Ar} and \cite[Sect.~3]{Ts02}. My
definition (below) is somewhat more restrictive than Arveson's, since
I restrict myself to measure classes generated by a single measure
(and standard probability spaces).

\begin{definition}\label{10a2}
A \emph{measure class}\index{measure class}
is a triple $ (\Om,\F,\cM) $ consisting of a set $ \Om $, a \sif\ $ \F $
on $ \Om $ and a set $ \cM $ of probability measures on $ (\Om,\F) $ such
that for some $ \mu \in \cM $ the probability space $ (\Om,\F,\mu) $ is
standard and for every probability measure $ \nu $ on $ (\Om,\F) $, 
\[
\nu \sim \mu \quad \text{if and only if} \quad \nu \in \cM \, ,
\]
$ \nu \sim \mu $ denoting mutual absolute continuity.
\end{definition}

Hilbert spaces $ L_2(\mu) $, $ L_2(\nu) $ for $ \mu,\nu \in \cM $ may
be glued together via the unitary operator $ L_2(\mu) \to L_2(\nu) $,
\[
f \mapsto \sqrt{ \frac\mu\nu } f \, ;
\]
here $ \frac\mu\nu $ is the Radon-Nikodym derivative (denoted also by $
\frac{\D\mu}{\D\nu} $). These spaces may be treated as `incarnations' of
a single Hilbert space $ L_2 (\Om,\F,\cM) $.\index{zzl@$ L_2(\cM) $, space}
The general form of an
element of $ L_2 (\Om,\F,\cM) $ is $ f \sqrt\mu $, where $ \mu \in \cM $
and $ f \in L_2(\Om,\F,\mu) $, taking into account the relation
\[
f \sqrt\mu = \bigg( \sqrt{ \frac\mu\nu } f \bigg) \sqrt\nu \, .
\]

Any isomorphism of measure classes induces naturally a unitary
operator between the corresponding Hilbert spaces.

The product of two measure classes is defined naturally, and
\[
L_2 \( (\Om,\F,\cM) \times (\Om',\F',\cM') \) = L_2 (\Om,\F,\cM) \otimes L_2
(\Om',\F',\cM') \, ;
\]
that is, we have a canonical unitary operator between these spaces,
namely, $ f \sqrt\mu \otimes f' \sqrt{\mu'} \mapsto (f\otimes f') \sqrt{ \mu
\otimes \mu' } $, where $ (f\otimes f') (\om,\om') = f(\om) f'(\om') $. I'll
write in short $ L_2 (\cM) $ instead of $ L_2 (\Om,\F,\cM) $; thus, $
L_2 ( \cM \times \cM' ) = L_2 ( \cM ) \otimes L_2 ( \cM' ) $. Everyone knows
the similar fact for measure spaces, $ L_2 ( \mu \times \mu' ) = L_2 (
\mu ) \otimes L_2 ( \mu' ) $. Here is a counterpart of Def.~\ref{2b6}.

\begin{definition} \label{10a3}
A \emph{continuous product of measure classes}%
\index{continuous product!of measure classes}
consists of measure classes $ (\Om_{s,t}, \cM_{s,t} ) $ (given for all $
s,t \in [-\infty,\infty] $, $ s<t $), and isomorphisms $
(\Om_{r,s},\cM_{r,s}) \times (\Om_{s,t},\cM_{s,t}) \to
(\Om_{r,t},\cM_{r,t}) $ (given for all $ r,s,t \in [-\infty,\infty] $,
$ r<s<t $) satisfying the associativity condition:
\[
(\om_1 \om_2) \om_3 = \om_1 (\om_2 \om_3) \quad \text{for almost all $
\om_1 \in \Om_{r,s} $, $ \om_2 \in \Om_{s,t} $, $ \om_3 \in \Om_{t,u}
$} 
\]
whenever $ -\infty\le r<s<t\le\infty $.
\end{definition}

Note the time set $ [-\infty,\infty] $ rather than $ \R $. Enlarging $
\R $ to $ [-\infty,\infty] $ is easy when dealing with probability
spaces (as noted after Def.~\ref{2b1}) but not measure classes (nor
Hilbert spaces, as noted after \ref{3a1}). Having a \emph{local}
continuous product of measure classes (over the time set $ \R $) we
may choose $ \mu_{n,n+1} \in \cM_{n,n+1} $ for each $ n \in \Z $ and
define $ \cM_{-\infty,\infty} $ as the equivalence class that contains
the product of these $ \mu_{n,n+1} $. However, another choice of $
\mu_{n,n+1} $ may lead to another $ \cM_{-\infty,\infty} $.

Given a continuous product of measure classes $ (\cM_{s,t})_{s<t} $,
we may construct the corresponding continuous product of Hilbert
spaces $ (H_{s,t})_{s<t} $; just $ H_{s,t} = L_2 (\cM_{s,t}) $.

\begin{question}\index{question}
Does every continuous product of Hilbert spaces (up to isomorphism)
emerge from some continuous product of measure classes?
\end{question}

See also \cite[Note 8.4 and Sect.~11 (question 9)]{Li}. A counterpart
of Def.~\ref{2b1} (see \ref{10a7}) needs some preparation.

\begin{definition}
Sub-\sif s $ \F_1, \dots, \F_n $ on a measure class $ (\Om,\F,\cM) $ are
\emph{independent,}%
\index{independent!sub-\sif s on a measure class}
if there exists a probability measure $ \mu \in \cM $ such that
\begin{equation}\label{16}
\mu ( A_1 \cap \dots \cap A_n ) = \mu(A_1) \dots \mu(A_n) \quad
\text{for all } A_1 \in \F_1, \dots, A_n \in \F_n \, .
\end{equation}
\end{definition}

For independent $ \F_1, \dots, \F_n $ the sub-\sif\ $ \F_1 \vee \dots
\vee \F_n $ generated by them will be denoted also by $ \F_1 \otimes \dots
\otimes \F_n $.

Given a product of two measure classes,
\[
(\Om,\F,\cM) = (\Om_1,\F_1,\cM_1) \times (\Om_2,\F_2,\cM_2) \, ,
\]
we have two independent sub-\sif s $ \ti\F_1, \ti\F_2 $ such that $ \F
= \ti\F_1 \otimes \ti\F_2 $; roughly,
\[
\ti\F_1 = \{ A \times \Om_2 : A \in \F_1 \} \, , \quad
\ti\F_2 = \{ \Om_1 \times B : B \in \F_2 \}
\]
(however, all negligible sets must be added).

And conversely, every two independent sub-\sif s $ \F_1, \F_2 $ such
that $ \F = \F_1 \otimes \F_2 $ emerge from a representation of $
(\Om,\F,\cM) $ (up to isomorphism) as a product; in fact, $
(\Om_k,\F_k,\cM_k) = (\Om,\F,\cM)/\F_k $ is the quotient space.

The following definition (in the style of \ref{2b1}) is equivalent to
\ref{10a3}.

\begin{definition}\label{10a7}
A \emph{continuous product of measure classes}%
\index{continuous product!of measure classes}
consists of a measure class $ (\Om,\F,\cM) $ and sub-\sif s $ \F_{s,t}
\subset \F $ (given for all $ s,t \in [-\infty,\infty] $, $ s<t $)
such that $ \F_{-\infty,\infty} = \F $ (`non-redundancy'), and
\begin{equation}
\F_{r,s} \otimes \F_{s,t} = \F_{r,t} \quad \text{whenever } -\infty\le
r<s<t\le\infty \, .
\end{equation}
\end{definition}

Every factorizing measure type $ \cM $ on $ \Comp(\R) $ leads to a
continuous product of measure classes $ (\Comp(s,t),\cM_{s,t})_{s<t}
$; recall \eqref{9b2}.

\begin{proof}[Proof \textup{(sketch)} of Theorem \textup{\ref{10a1}}]
\hskip 0pt plus 20pt
The factorizing measure type $ \cM $ on $ \Comp(\R) $ leads to a
continuous product of measure classes $ (\Comp(s,t),\linebreak[0]
\cM_{s,t})_{s<t}
$ and further, to a continuous product of Hilbert spaces $
(H_{s,t})_{s<t} $, $ H_{s,t} = L_2(\cM_{s,t}) $.
Measures of $ \cM $ have an atom at $ \emptyset $, that is, $
\mu(\{\emptyset\}) > 0 $. (Indeed, on small intervals $
\mu_{s,t}(\{\emptyset\}) > 0 $ since it is close to $ 1 $; cover $
[-\infty,\infty] $ by a finite number of such intervals, and
multiply.) We define $ u_{s,t} $ as the root of the probability
measure concentrated at the atom, thus getting a continuous product of
pointed Hilbert spaces $ (H_{s,t},u_{s,t})_{s<t} $. The upward
continuity follows from the fact that $ \Comp(s,t) = \bigcup_{\eps>0}
\Comp(s+\eps,t-\eps) \bmod 0 $. The projection-valued measure $ Q $
given by Theorem \ref{9d1} satisfies $ Q ( \{ C : C \cap (s,t) =
\emptyset \} ) H = H_{-\infty,s} u_{s,t} H_{t,\infty} $, the latter
being the space of all vectors $ \psi = f \sqrt\mu $ such that the
measure $ |\psi|^2 = |f|^2 \cdot \mu $ on $ \Comp(\R) $ is
concentrated on the set $ \{ C : C \cap (s,t) = \emptyset \} $. Thus,
$ \ip{ Q(A)\psi }{ \psi } = |\psi|^2 (A) $ for sets $ A \subset
\Comp(\R) $ of the form $ A = \{ C : C \cap (s,t) = \emptyset \}
$. The same holds for finite intersections of such sets, therefore,
for all measurable $ A \subset \Comp(\R) $. It means that the spectral
measure of $ \psi = f \sqrt\mu $ is $ |\psi|^2 = |f|^2 \cdot \mu $.
\end{proof}

If $ \cM $ is shift-invariant then the continuous product constructed
above is homogeneous, and leads to an Arveson system (recall
\ref{3c9}). If, in addition, \almost{\cM} all $ C $ are perfect then
the constructed Arveson system has no other units, that is, is of type
$ II_0 $.

\begin{corollary}
(a)
Every shift-invariant factorizing measure type $ \cM $ on $ \Comp(\R)
$ is equal to $ \cM^u $ for some unit $ u $ of some Arveson system.

(b) If, in addition, $ \cM $ is concentrated on (the set of all)
perfect subsets of $ \R $ then $ \cM = \cM^\U $ for some Arveson system
of type $ II_0 $.
\end{corollary}

(See also \cite[Th.~3]{Li} for a stronger result.)
Thus, Arveson systems of type $ II_0 $ are at least as  diverse as
shift-invariant factorizing measure types on the space of all perfect
subsets of $ \R $.

The set of zeros of a Brownian motion is an example of such measure
type. More exactly, we may consider the random set $ \{ s \in (0,t) :
a+B_s = 0 \} $ for given $ t,a \in (0,\infty) $; the distribution of
the random set depends on $ a $, but its measure type does not. The
corresponding shift-invariant factorizing measure type on $ \Comp(\R)
$ exists and is unique. The random set is perfect, of Hausdorff
dimension $ \frac12 $ (unless empty).

Similarly, an example of a random set of any desired Hausdorff
dimension between $ 0 $ and $ 1 $ is given by zeros of a Bessel
process \cite[Sect.~3]{Ts00}. See also \cite[Sect.~4.4]{Li}.

However, random sets are much more diverse than Bessel processes. For
every perfect set of Hausdorff dimension less than $ \frac12 $ there
exists a random set (I mean, a shift-invariant factorizing measure
type) obtained from the given (nonrandom) set by a random perturbation
preserving almost all the microstructure of the given set. It may be
called the \emph{barcode construction,}\index{barcode}
see \cite[Sect.~6]{Ts02}. In
contrast to Bessel zeros, random sets obtained from the barcode
construction are in general not invariant under time reversal ($ t
\mapsto -t $) and time rescaling ($ t \mapsto ct $ for $ c \in
(0,\infty) $, $ c \ne 1 $); they gives us a continuum of mutually
non-isomorphic asymmetric Arveson systems of type $ II_0 $
\cite[Th.~7.3]{Ts02}.

\mysubsection{From off-white noises to Hilbert spaces}
\label{10b}

Random compact subsets of $ \R $ are one out of many sources of
continuous products of measure classes. One may use random compact (or
closed) subsets of $ \R \times L $ for some locally compact space $ L $
\cite[Sect.~4.1 and Sect.~11 (question 3)]{Li}, random measures
\cite[Note 6.8 and Sect.~8.3]{Li}, etc. But first of all we should try
Gaussian processes, for several reasons: they occupy a prominent place
among random processes; relations between sub-\sif s reduce to
relations between subspaces of a Hilbert space; the white noise is a
Gaussian process. Indeed, random sets used before generalize the
Poisson process, while off-white noises used below generalize the
white noise, and appear to lead to type $ III $ Arveson systems (see
\ref{10b6}).

Of course, the white noise cannot be treated as a random function on $
\R $. Gaussian random variables correspond to test functions rather
than points, which is harmless; we need only sub-\sif s $ \F_{s,t} $
that correspond to intervals $ (s,t) $, not points. The same
holds for other Gaussian processes considered here. Being stationary,
such process is described by its \emph{spectral measure,}%
\index{spectral measure!of a stationary Gaussian process}
a positive $
\si $-finite measure $ \nu $ on $ [0,\infty) $ such that the Gaussian
random variable corresponding to a test function $ f $ has mean $ 0 $
and variance $ \int | \hat f (\la) |^2 \, \nu(\D\la) $; here $ \hat f
(\la) = (2\pi)^{-1/2} \int f(t) \E^{-\I\la t} \, \D t $ is the Fourier
transform of $ f $. We restrict ourselves to measures $ \nu $
majorized by Lebesgue measure (that is, $ \nu(\D\la) \le C \D\la $ for
some $ C $; see \cite{Ts02off} and \cite[Sect.~9]{Ts02} for the
general case). Thus, every $ f \in L_2(\R) $ is an admissible test
function. The space $ G $ of Gaussian random variables may be
identified with the Hilbert space $ L_2(\nu) $. Each interval $ (s,t)
\subset \R $ leads to a subspace $ G_{s,t} \subset G $%
\index{zzg@$ G_{s,t} $, Gaussian subspace}
defined as the
closure of $ \{ \hat f : f \in L_2(s,t) \} $, and the corresponding
sub-\sif\ $ \F_{s,t} \subset \F $.

For the white noise the `past' and `future' spaces $ G_{-\infty,0} $,
$ G_{0,\infty} $ are orthogonal; the `past' and `future' sub-\sif\ $
\F_{-\infty,0} $, $ \F_{0,\infty} $ are independent, and we have a
continuous product of \emph{probability spaces.} More generally, in
order to give us a continuous product of \emph{measure classes,} these
sub-\sif s should be independent for some equivalent measure. A
necessary and sufficient condition is well-known (see
\cite[Th.~3.2]{Ts02off}, \cite[Th.~9.7]{Ts02}): $ \nu(\D\la) =
\E^{\phi(\la)} \D\la $ for some $ \phi : [0,\infty) \to \R $ such that
\begin{equation}\label{10b1}
\int_0^\infty \int_0^\infty \frac{ | \phi (\la_1) - \phi (\la_2) |^2
}{ | \la_1 - \la_2 |^2 } \, \D\la_1 \, \D\la_2 < \infty \, .
\end{equation}
A sufficient condition for \eqref{10b1} is available (see
\cite[Prop.~3.6(b)]{Ts02off}, \cite[(9.11)]{Ts02}):
\begin{equation}\label{10b2}
\text{$ \phi $ is continuously differentiable, and} \quad
\int_0^\infty \bigg| \frac{\D}{\D\la} \phi (\la) \bigg|^2 \la \,
 \D\la < \infty \, .
\end{equation}
In particular, the sufficient condition is satisfied by any strictly
positive smooth function $ \la \mapsto \E^{\phi(\la)} = \nu(\D\la) /
\D\la $ such that for $ \la $ large enough, one of the following
equalities holds:
\begin{gather}
\frac{\nu(\D\la)}{\D\la} = \ln^{-\al} \la \, , \quad 0 \le \al <
\infty \, ; \label{10b3} \\
\frac{\nu(\D\la)}{\D\la} = \exp ( - \ln^\be \la ) \, , \quad 0 < \be <
 \frac12 \, ; \label{10b4}
\end{gather}
see \cite[Examples 3.11, 3.12]{Ts02off},
\cite[(9.12)--(9.13)]{Ts02}. Every such $ \nu $ leads to a continuous
product of measure classes $ (\Om,\F,\cM), (\F_{s,t})_{s<t} $. Namely,
$ \cM $ is the equivalence class containing the Gaussian measure $ \ga
$ whose spectral measure is $ \nu $. The group $ (T_h)_{h\in\R} $ of
time shifts leaves invariant the equivalence class $ \cM $ and
moreover, the measure $ \ga $. The corresponding Hilbert spaces $
H_{s,t} = L_2 (\Om,\F_{s,t},\cM) $, being a homogeneous continuous
product of Hilbert spaces (over the time set $ [-\infty,\infty] $),
lead to an Arveson system.

The white noise, contained in \eqref{10b3} as the case $ \al = 0 $,
leads to a classical continuous product of probability spaces $
(\Om,\F,\ga), (\F_{s,t})_{s<t} $ and classical (type $ I $) Arveson
system. Generally, $ (\Om,\F,\ga), (\F_{s,t})_{s<t} $ is not a
continuous product of probability spaces, since the past and the
future are not independent on $ (\Om,\F,\ga) $. However, a
decomposable vector $ \psi = f \sqrt\mu $, $ \|\psi\|=1 $ (if any)
gives us a probability measure $ |\psi|^2 = |f|^2 \cdot \mu $,
decomposable in the sense that $ (\Om,\F,|\psi|^2), (\F_{s,t})_{s<t} $
is a continuous product of probability spaces. Especially, the measure
$ |\psi|^2 $ makes independent the pair of random variables
corresponding to such `comb' test functions $ f_n, g_n $ (for any
given $ n $):
\[
\begin{gathered}\includegraphics{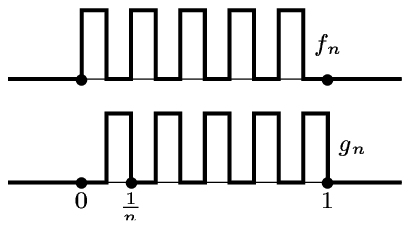}\end{gathered}
\qquad
\begin{gathered}
f_n(t) + g_n(t) = 1 \, , \\
f_n(t) - g_n(t) = \operatorname{sgn} \sin \pi n t \\
\qquad \text{for } t \in (0,1) \, .
\end{gathered}
\]
If $ \nu $ satisfies the condition
\begin{equation}\label{10b5}
\frac{ \nu(\D\la) }{ \D\la } \to 0 \quad \text{as } \la \to \infty
\end{equation}
(which excludes the white noise), then $ \| \hat f_n - \hat g_n
\|_{L_2(\nu)} \to 0 $ (see \cite[10.2]{Ts02}), therefore the
\emph{independent} random variables on $ (\Om,\F,|\psi|^2) $
corresponding to $ f_n $ and $ g_n $ converge (as $ n \to \infty $, in
probability) to \emph{the same} random variable $ Z $ corresponding to
the test function $ \frac12 \cdot \One_{(0,1)} $. We see that $ Z $ is
constant on $ (\Om,\F,|\psi|^2) $, therefore, has an atom on $
(\Om,\F,\ga) $ (since the measure $ |\psi|^2 $ is absolutely
continuous w.r.t.\ $ \ga $). However, the \emph{normal} distribution
of $ Z $ on $ (\Om,\F,\ga) $ is evidently nonatomic! The conclusion
follows.

\begin{proposition}\label{10b6}
\cite[10.3]{Ts02}
If $ \nu $ satisfies \eqref{10b5} then the corresponding continuous
product of Hilbert spaces has no decomposable vectors, and the
corresponding Arveson system is of type $ III $.
\end{proposition}

Every $ \al > 0 $ in \eqref{10b3} (as well as every $ \be $ in
\eqref{10b4}) gives us a nonclassical Arveson system. One may guess
that, the larger the parameter $ \al $, the more nonclassical the
system. Striving to an invariant able to confirm the guess, we
introduce `spaced comb' test functions $ f_{n,\eps} $,
\[
\begin{gathered}\includegraphics{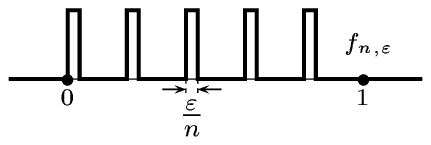}\end{gathered}
\qquad
\begin{gathered}
f_{n,\eps}(t) = 1 \text{ if } (nt \bmod 1) \in (0,\eps) ,
\end{gathered}
\]
and the corresponding Gaussian random variables $ Z_{n,\eps} $ on $
(\Om,\F,\ga) $. It is instructive to consider the correlation
coefficient
\[
\rho_{n,\eps} = \frac{ \int Z_{n,\eps} Z \, \D\ga }{ ( \int
Z^2_{n,\eps} \, \D\ga )^{1/2} ( \int Z^2 \, \D\ga )^{1/2} }
\]
between $ Z_{n,\eps} $ and the random variable $ Z $ corresponding to
the test function $ \One_{(0,1)} $. For the white noise, $
\rho_{n,\eps} $ does not depend on $ n $ (in fact, $ \rho_{n,\eps} =
\sqrt\eps $). However, \eqref{10b5} implies $ \rho_{n,\eps} \to 1 $ as
$ n \to \infty $ (for every $ \eps \in (0,1) $). On the other hand, $
\rho_{n,\eps} \to 0 $ as $ \eps \to 0 $ (for every $ n $). It is more
interesting to take $ \lim_n \rho_{n,\eps_n} $ when $ \eps_n \to 0
$. Especially, for $ \nu $ of the form \eqref{10b3} with $ \al > 0 $,
\begin{gather*}
\text{if } \eps_n \ln^\al n \to \infty \text{ then } \rho_{n,\eps_n}
 \to 1 \, , \\
\text{if } \eps_n \ln^\al n \to 0 \text{ then } \rho_{n,\eps_n} \to 0
 \, ,
\end{gather*}
which gives us a clue to a useful invariant. We should consider
`spaced comb' sets $ E_n \subset (0,1) $ (namely, $ E_n = \{ t \in
(0,1) : f_{n,\eps_n}(t) = 1 \} $), the decompositions $ (0,1) = E_n
\cup E_n^\co $ (where $ E_n^\co = (0,1) \setminus E_n $) of the
interval $ (0,1) $, and the corresponding decompositions (see also
\eqref{11a1})
\begin{gather*}
\cM_{0,1} = \cM_{E_n} \times \cM_{E_n^\co} \, , \\
G_{0,1} = G_{E_n} \oplus G_{E_n^\co} \, , \\
H_{0,1} = H_{E_n} \otimes H_{E_n^\co} \, , \\
\A_{0,1} = \A_{E_n} \otimes \A_{E_n^\co}
\end{gather*}
of the measure class $ \cM_{0,1} $, the Gaussian space $ G_{0,1} $ (a
linear subspace of $ G $), the Hilbert space $ H_{0,1} = L_2
(\cM_{0,1}) $ and the algebra $ \A_{0,1} $ of operators on $ H_{0,1}
$. Their asymptotic behavior (as $ n \to \infty $) should be sensitive
to the asymptotic behavior (as $ \la \to \infty $) of $ \nu(\D\la) /
\D\la $.

The rest of the story, sketched below, belongs to functional analysis
rather than probability. (In fact, the whole story is translated into
the language of analysis by Bhat and Srinivasan \cite{BS}.) The norm
on the Hilbert space $ H_{0,1} $ is singled out by the Gaussian
measure $ \ga $. However, the equivalence class $ \cM_{0,1} $ contains
many Gaussian measures; $ \ga $ is just one of them. Accordingly, $ G
$ should not be treated as a Hilbert space. Its natural structure is
given by an equivalence class of norms (rather than a single norm),
but the equivalence is much stronger than topological, it may be
called FHS\nobreakdash-\hspace{0pt}equivalence, and $ G $ may be
called an FHS\nobreakdash-\hspace{0pt}space%
\index{FHS space}
\cite[8.5]{Ts02}. The decomposition $ G_{0,1} = G_{E_n} \oplus
G_{E_n^\co} $ is orthogonal in the FHS sense, that is, orthogonal in
\emph{some} (depending on $ n $) norm of the given class. However, the
decomposition $ \A_{0,1} = \A_{E_n} \otimes \A_{E_n^\co} $ is treated
as usual; $ \A_{E_n^\co} $ is the commutant of $ \A_{E_n} $. In fact,
\emph{every} FHS space $ G $ leads to a Hilbert space $ H = \Exp(G) $,
and \emph{every} orthogonal decomposition of the FHS space $ G $ leads
to a decomposition of the operator algebra of $ H $ into tensor
product.

The desired invariant (of an Arveson system) is the set of all
sequences $ \eps_n \to 0 $ such that
\begin{equation}\label{10b7}
\limsup_{n\to\infty} \A_{E_n} \quad \text{is trivial} \, ;
\end{equation}
the latter means that for all $ A_1 \in \A_{E_1}, A_2 \in \A_{E_2},
\dots $ all limit points (in the weak operator topology) of the
sequence $ A_1, A_2, \dots $ are scalar operators. Condition
\eqref{10b7}, taken from \cite[Def.~26 and Th.~30]{BS}, is equivalent
to the condition \cite[2.2]{Ts02}: for every trace-class operator $ R
: H_{0,1} \to H_{0,1} $ satisfying $ \trace(R) = 0 $,
\begin{equation}\label{10b8}
\sup_{A\in\A_{E_n},\|A\|\le1} | \trace(AR) | \to 0 \quad \text{as } n
 \to \infty \, .
\end{equation}
(Note that the sequence $ (E_n)_n $ is not decreasing, in contrast to
\cite[Prop.~10.1 and Cor.~10.2]{Li}.) Fortunately, the condition can
be reformulated in terms of $ G_{E_n} $ and $ G_{E_n^\co} $.

\begin{proposition}\label{10b9}
Condition \eqref{10b7} holds if and only if
\begin{equation}\label{10b10}
\limsup_{n\to\infty} G_{E_n} = \{ 0 \} \quad \text{and} \quad
\liminf_{n\to\infty} G_{E_n^\co} = G \, .
\end{equation}
\end{proposition}

\begin{sloppypar}
The relation $ \limsup G_{E_n} = \{ 0 \} $ means that for all $ g_1
\in G_{E_1}, g_2 \in G_{E_2}, \dots $ the only possible limit point of
the sequence $ g_1, g_2, \dots $ is $ 0 $. The relation $ \liminf
G_{E_n^\co} = G $ means that every $ g \in G $ is the limit of some
sequence $ g_1, g_2, \dots $ such that $ g_1 \in G_{E_1^\co}, g_2 \in
G_{E_2^\co}, \dots $
\end{sloppypar}

Proposition \ref{10b9} appeared first in \cite[11.3]{Ts02} with a
long, complicated proof (occupying Sections 11 and 12 of
\cite{Ts02}). For a substantially simpler proof see
\cite[Th.~30]{BS}.

\begin{proposition}
(See \cite[13.10]{Ts02}.)
Let spectral measures $ \nu_1, \nu_2 $ satisfy \eqref{10b3} with
parameters $ \al_1, \al_2 $ respectively, $ 0 < \al_1 < \al_2 $. Then
there exists a sequence $ (\eps_n)_n $ such that $ \eps_n \to 0 $ and
the corresponding `spaced comb' sets $ E_n $ satisfy \eqref{10b7} for
the Arveson system corresponding to $ \nu_1 $ but not $ \nu_2 $.
\end{proposition}

In fact, the sequence $ \eps_n = \ln^{-c} n $ fits for $ \al_1 < c \le
\al_2 $. Of course, \eqref{10b10} is checked instead of
\eqref{10b7}. Still, the proof uses tedious calculations
\cite[Sect.~13]{Ts02}. Here is the conclusion.

\begin{theorem}\index{theorem}
(Tsirelson \cite[13.11]{Ts02})
There is a continuum of mutually non-isomorphic Arveson systems of
type $ III $.
\end{theorem}

\mysection{Beyond the one-dimensional time}
\label{sec:11}
\mysubsection{Boolean base}
\label{11a}

The definitions of continuous products (\ref{2b6}, \ref{3a1},
\ref{6d6}, \ref{10a3}) center round the relations
\[
\Om_{r,t} = \Om_{r,s} \times \Om_{s,t} \, , \quad H_{r,t} =
H_{r,s} \otimes H_{s,t} \quad \text{etc.}
\]
for $ r < s < t $. These are special cases (for $ A = (r,s) $ and $ B
= (s,t) $) of more general relations
\begin{equation}\label{11a1}
\Om_{A \uplus B} = \Om_A \times \Om_B \, , \quad H_{A \uplus B} =
H_A \otimes H_B  \quad \text{etc.}
\end{equation}
for \emph{elementary sets} $ A, B $ satisfying $ A \cap B = \emptyset
$. By an elementary set I mean a union of finitely many intervals,
treated modulo finite sets. (For example, $ (-5,1) \cup \{ 2 \} \cup
[9,\infty) $ is an elementary set, and $ [-5,1) \cup (9,\infty) $ is
the same elementary set.) The disjoint union $ A \uplus B $ is just $
A \cup B $ provided that $ A \cap B = \emptyset $. Of course,
\[
\Om_{(r,s)\cup(t,u)} = \Om_{r,s} \times \Om_{t,u} \, , \quad
H_{(r,s)\cup(t,u)} = H_{r,s} \otimes H_{t,u} \quad \text{etc.}
\]
for $ r < s < t < u $; the same for any finite number of intervals.

Elementary sets are a Boolean algebra. More generally, we may consider
an arbitrary Boolean algebra $ \A $ (instead of the time set $ T $)
and define continuous products (of probability spaces, Hilbert spaces
etc.) over $ \A $ by requiring \eqref{11a1} for all disjoint $ A,B \in
\A $. This `boolean base' approach is used in \cite{AW}, \cite{Fe},
\cite[Sect.~1]{TV}, \cite{Ts99}. Early works \cite{AW}, \cite{Fe}
concentrate on \emph{complete} Boolean algebras (which means that
every subset of the Boolean algebra has a supremum in the algebra),
striving to prove that all continuous products (satisfying appropriate
continuity conditions) are classical. More recent works
\cite{TV}, \cite[Sect.~2]{Ts99} prefer incomplete Boolean algebras and
nonclassical continuous products.

The following result answers a question of Feldman \cite[1.9]{Fe}.

\begin{theorem}\label{11a2}\index{theorem}
(Tsirelson \cite[6c7]{Ts03}, see also \cite[3.2]{Ts99})
A continuous product of probability spaces, satisfying the upward
continuity condition, is classical if and only if the map $ E \mapsto
\F_E $ can be extended from the algebra of elementary sets to the
Borel \sif, satisfying $ \F_{A \uplus B} = \F_A \otimes \F_B $ and the
upward continuity ($ A_n \uparrow A $ implies $ \F_{A_n} \uparrow \F_A
$) for Borel sets $ A,B,A_n $.
\end{theorem}

\mysubsection{Two-dimensional base}
\label{11b}

The two black noises considered in Sections \ref{7f}, \ref{7j} are
scaling limits of discrete models driven by two-dimensional arrays of
independent random variables. Their one-dimensional time is just one
of the two dimensions. The sub-\sif\ $ \F_{s,t} $ corresponds to the
strip $ (s,t) \times \R \subset \R^2 $. It should be possible to
define sub-\sif s $ \F_A $ for more general sets $ A \subset \R^2
$. The whole Borel \sif\ of $ \R^2 $ is too big (recall Theorem
\ref{11a2}); the Boolean algebra generated by rectangles $ (s,t)
\times (a,b) $ is a modest choice. No such theory is available for
now.

The class of appropriate sets $ A \subset \R^2 $ should depend on the
model. The same may be said about the one-dimensional time, if we do
not restrict ourselves to intervals (as in Sections
\ref{sec:4}--\ref{sec:10}) or elementary sets (as in
Sect.~\ref{11a}). The Hausdorff dimension $ \dim(\pd A) $ of the
boundary of $ A $ could be relevant. See also \cite[end of
Sect.~2d]{TV}. It could be related to the Hausdorff dimension of
spectral sets (recall \ref{9b}, \ref{9c}).

The model of Sect.~\ref{7f}, being a kind of oriented percolation, is
much simpler than the true percolation.

\begin{question}\index{question}
(See also \cite[8a1]{Ts03}.)
For the (conformally invariant) scaling limit of the critical site
percolation on the triangular lattice, invent an appropriate
conformally invariant Boolean algebra of sets on the plane and define
the corresponding sub-\sif s $ \F_A $ satisfying $ \F_{A \uplus B} =
\F_A \otimes \F_B $. Is it possible?
\end{question}

Hopefully, the answer is affirmative, that is, the two-dimensional
noise of percolation\index{percolation (noise of?)}
will be defined. Then it should appear to be a
(two-dimensional) black noise, see \cite[8a2]{Ts03}.

It would be the most important example of a black noise!

\bigskip\textbf{Acknowledgment.}
I thank the anonymous referee whose detailed comments have lead to better
readability.

\printindex

\end{document}